%% file: ArxivGNPIV.tex
\documentclass[hidelinks, 11 pt, a4paper]{article}
\usepackage{amsfonts}
\usepackage{amsmath}
\usepackage{amssymb}
\usepackage{amsthm}
\usepackage{subcaption}
\usepackage[pdftex]{graphicx}
\usepackage{natbib}
\usepackage[hypertexnames=false]{hyperref}
\usepackage{setspace}
\usepackage{pdfsync}
\usepackage{lscape}
\usepackage{titling}
\usepackage[hmargin=1.25in,vmargin=1.25in]{geometry}
\usepackage{bibunits}
\usepackage{mathrsfs,dsfont}
\usepackage[dvipsnames,svgnames, x11names]{xcolor}
\usepackage{pgfplots}
\usepackage{caption}
\usepackage{footmisc}
\usepackage{rotating}
\usepackage{algpseudocode}
\usepackage{algorithm}

\newcommand\DerP{\mathbb D_{P}}
\newcommand\PSigma{\Sigma_P}
\newcommand\Up{U_P}
\newcommand\UpS{U_P^\star}
\newcommand\IDset{\Theta_0}
\newcommand\IDsetRsieve{\Theta_{0n}^{\rm r}}
\newcommand\IDsetUsieve{\Theta_{0n}^{\rm u}}
\newcommand\IDpoint{\theta_0}
\newcommand\Gemp{\mathbb G_n}
\newcommand\WP{\mathbb W_P}
\newcommand\Bemp{\hat{\mathbb W}_n}
\newcommand\WPT{{\mathbb W}_P^\star}
\newcommand\Iso{\mathbb G_P}
\newcommand\IsoW{{\mathbb G}_P^\star}
\newcommand\Wboot{\hat{\mathbb G}_n}

\usetikzlibrary{patterns}
\hypersetup{colorlinks, linkcolor = {teal}, citecolor = {blue}, urlcolor ={blue!80!black}}
\DeclareGraphicsExtensions{.jpg}
\defaultbibliographystyle{ims}
\defaultbibliography{andres}
\makeatletter

\renewcommand\paragraph{\@startsection{paragraph}{4}{\z@}%
            {-2.5ex\@plus -1ex \@minus -.25ex}%
            {1.25ex \@plus .25ex}%
            {\normalfont\normalsize\bfseries}}
\makeatother
\def\qed{\rule{2mm}{2mm}}

\newtheorem{theorem}{Theorem}[section]
\newtheorem{lemma}{Lemma}[section]
\newtheorem{definition}{Definition}[section]

\newtheorem{corollary}{Corollary}[section]
\newtheorem{remark}{Remark}[section]
\newtheorem{assumption}{Assumption}[section]

\newenvironment{packed_enum}{\begin{enumerate} \setlength{\itemsep}{1pt}\setlength{\parskip}{0pt}\setlength{\parsep}{0pt}}{\end{enumerate}}

\pgfplotsset{compat=1.17} 

\begin{document}

\begin{bibunit}

\title{Constrained Conditional Moment Restriction Models\thanks{We thank Riccardo D'amato for excellent research assistance. We are also indebted to three anonymous referees and numerous seminar participants for their valuable comments.} }

{\author{Victor Chernozhukov \\ 
M.I.T. \\ vchern@mit.edu
\and
Whitney K. Newey\thanks{Research supported by NSF Grant 1757140.} \\ 
M.I.T. \\ wnewey@mit.edu
\and
Andres Santos\thanks{Research supported by NSF Grant SES-1426882.}\\ 
U.C.L.A.\\ andres@econ.ucla.edu}}

\date{First Draft: September, 2015 \\ This Draft: April 2022}
\maketitle

\begin{abstract}
Shape restrictions have played a central role in economics as both testable implications of theory and sufficient conditions for obtaining informative counterfactual predictions.
In this paper we provide a general procedure for inference under shape restrictions in identified and partially identified models defined by conditional moment restrictions.
Our test statistics and proposed inference methods are based on the minimum of the generalized method of moments (GMM) objective function with and without shape restrictions.
Uniformly valid critical values are obtained through a bootstrap procedure that  approximates a subset of the true local parameter space.
In an empirical analysis of the effect of childbearing on female labor supply, we show that employing shape restrictions in linear instrumental variables (IV) models can lead to shorter confidence regions for both local and average treatment effects.
Other applications we discuss include inference for the variability of quantile IV treatment effects and for bounds on average equivalent variation in a demand model with general heterogeneity. 
We find in Monte Carlo examples that the critical values are conservatively accurate and that tests about objects of interest have good power relative to unrestricted GMM.
\end{abstract}

\begin{center}
\end{center}

\thispagestyle{empty}

\newpage

\pagenumbering{arabic}

\section{Introduction}

Shape restrictions have played a central role in economics as both testable implications of classical theory and sufficient conditions for obtaining informative counterfactual predictions \citep{topkis1998supermodularity}.
A long tradition in applied and theoretical econometrics has as a result studied shape restrictions, their ability to aid in identification, estimation, and inference, and the possibility of testing for their validity \citep{matzkin:1994, chetverikov2018econometrics}.
A canonical example of this interplay between theory and practice is consumer demand analysis, where theoretical predictions such as Slutsky conditions have been extensively tested for and employed in estimation \citep{hausman:newey:1995, hausman2016individual, blundell:horowitz:parey:2012, dette2016testing}.
The empirical analysis of shape restrictions, however, goes well beyond this important application with recent examples including studies into the monotonicity of the state price density \citep{jackwerth2000recovering, ait2003nonparametric}, the presence of ramp-up and start-up costs \citep{wolak2007quantifying, reguant2014complementary}, and the existence of complementarities in demand \citep{gentzkow:2007} and organizational design \citep{athey1998empirical, kretschmer2012competitive}.

Shape restrictions are often equivalent to inequality restrictions on parameters of interest and on certain unknown functions.
For example, Slutsky negative semi-definiteness and monotonicity require that certain functions satisfy inequality restrictions. 
Inference with inequality restrictions is difficult. 
Such restrictions lead to discontinuities in (pointwise) limiting distributions where the inequality restrictions are ``close" to binding, which makes inference challenging due to non-pivotal and potentially unreliable pointwise asymptotic approximations \citep{andrews2000inconsistency, andrews2001testing}. 
Limit discontinuities further make it difficult to construct confidence intervals with uniform coverage.

We address these challenges by obtaining critical values through a bootstrap procedure that uniformly approximates a subset of the local parameter space.
The proposed critical values simultaneously deliver uniformly valid inference and pointwise limiting rejection probabilities that equal the nominal level of the test in many applications.
Our results apply to a class of conditional moment restriction models \citep{ai:chen:2007,ai2012semiparametric} that encompasses parametric \citep{hansen:1982}, semiparametric \citep{ai:chen:2003}, and nonparametric \citep{newey:powell:2003} instrumental variable (IV) models, as well as panel data applications \citep{chamberlain1992comment}, and the study of plug-in functionals.
For parametric IV our results deliver novel uniformly valid tests of inequality and equality restrictions as well as confidence intervals for parameters of interest in the presence of inequality restrictions in both identified and partially identified models.

Our test statistics and proposed inference methods are based on the difference of the minimum of a generalized method of moments (GMM) objective function with and without inequality restrictions.
The value of the test statistic increases when more binding constraints are imposed. 
To ensure uniform validity, critical values are obtained through a bootstrap procedure that acknowledges that some inequalities that do not bind in the sample could have bound under a different draw of the sample. 
Intuitively, in the bootstrap, we impose the inequalities that are within a region of the boundary that shrinks slightly slower than the convergence rate of the shape restricted estimator. 
The bootstrap procedure can further be set to ignore inequalities that are outside this shrinking region, leading to pointwise rejection probabilities that equal the nominal level in many applications.
As always, uniformity is essential for confidence intervals to be asymptotically valid over a set of unknown parameter values. The resulting inference is powerful in exploiting the large amount of information that inequality restrictions can provide in many cases relevant for applications. 

Our tests and confidence intervals remain valid under partial identification. 
In this setting, the tests and confidence intervals give an accurate and computationally feasible method of doing inference for a subvector of parameters under partial identification. 
Indeed, these methods have already been used by \cite{torgovitsky2019nonparametric} to construct informative confidence intervals for various partially identified state dependence parameters in the presence of unobserved heterogeneity.
Also, \cite{kline2021reasonable} used these methods to test shape constraints implied by a model of callback probabilities for employment applications.
By incorporating nuisance parameters into the definition of the parameter space, our results can further be applied to partially identified semi(non)-parametric models defined by conditional moment inequalities.

We demonstrate the usefulness of this approach in an empirical application.
Specifically, we conduct inference on the causal effect of childbearing on female labor force participation by relying on the instrumental variables approach of \cite{angrist1998children}. 
We find that monotonicity of the local average treatment effect (LATE) in education is not rejected by the data and neither is monotonicity and negativity -- these restrictions were discussed, but not formally tested, by \cite{angrist1998children}.
We further find that imposing these shape restrictions yields narrower confidence intervals for the LATE at different schooling levels.
Finally, we obtain similar results for the partially identified average treatment effect (ATE), though the data is less informative about the ATE because of the low proportion of compliers.

The inequalities associated with nonparametric shape restrictions necessitate consideration of parameter spaces that are sufficiently general yet endowed with enough structure to ensure a fruitful asymptotic analysis. 
An important theoretical insight of this paper is that this simultaneous flexibility and structure is possessed by sets defined by inequality restrictions on Abstract M (AM) spaces; i.e. Banach lattices whose norm obeys a condition discussed in Section \ref{sec:gentheory}.
We also introduce potentially regularized approximations to the local parameter spaces in order to account for the curvature present in nonlinear constraints. 
While aspects of our analysis are specific to models defined by conditional moment restrictions, the role of the local parameter space is solely dictated by the shape restrictions. 
As such, we expect the insights of the set up here to be applicable to the study of shape restrictions in alternative models as well. 
The critical values are shown to be uniformly asymptotically valid by developing strong approximations to both the test and bootstrap statistics.
Sufficient conditions are provided by adapting the coupling of \cite{koltchinskii:1994}.
Our coupling arguments and the use of AM spaces are key features of the theory that enable us to show that inference is uniformly valid and that partial identification is permitted.

We illustrate the general applicability of our analysis by obtaining novel uniformly valid inference results in a variety of problems.
Specifically, we:
(i) Conduct inference about partially identified sets of average equivalent variation and other objects of interest in demand estimation with general heterogeneity and smooth demand functions;
(ii) Test and impose shape restrictions on structural functions identified through quantile conditional moment restrictions; and
(iii) Impose the Slutsky restrictions to conduct inference in a linear conditional moment restriction model.
Additionally, while we do not pursue further examples in detail for conciseness, we note our results may be applied to conduct tests of homogeneity, supermodularity, and economies of scale or scope.

In a small Monte Carlo study, we examine instrumental variables estimation of a nonlinear structural function and consider the power of imposing monotonicity and/or convexity on the structural function.
We find rejection frequencies for our test that are conservatively accurate when testing a point null hypothesis about the value or derivative of the structural function.
In addition, we find that imposing shape restrictions leads to large increases in power relative to employing an unrestricted estimator, in moderately large samples.
Our Monte Carlo analysis further examines the performance of our test in a partially identified parametric IV model with discrete data.
In that context, we find that shape restrictions have substantial identifying power and that our test provides valid inference on the value of a function at a point.
A similar partially identified IV setting was previously studied by \cite{freyberger:horowitz:2012}, who also provide an inference procedure.
However, their procedure is based on limiting distributions that are discontinuous in true parameters leading to nonuniform inference.

Our paper contributes to an extensive literature studying semiparametric and nonparametric models under partial identification \citep{manski:2003, molinari2020microeconometrics}.
When specialized to finite dimensional models, our results enable us to conduct inference on functionals of the identified set in models defined by moment (in)equalities \citep{canay2017practical, ho2017partial}.
In that context, our results are complementary to those of \cite{bugni2017inference} and \cite{kaido2019confidence}, who provide uniformly valid procedures for subvector inference.
Their analysis is focused on convex models and can thus be invalid or conservative when conducting inference on nonlinear functionals or imposing non-convex restrictions -- we emphasize, however, that their analysis is also motivated by a different set of models than the ones we consider.
Our analysis is further related to \cite{hong2011inference}, \cite{santos2012inference}, \cite{tao:2014}, and \cite{chen2011sensitivity} who study inference on functionals of potentially partially identified structural functions, but do not allow for shape constraints as we do.

Following the original version of this paper, \cite{zhu2019inference} and \cite{fang2019general} have proposed inference methods for convex restrictions which, while applicable to an important class of problems, rule out inference on nonlinear functionals or tests of certain shape restrictions.
Also related is \cite{freyberger2018inference} who have more recently developed uniform inference for functionals under shape restrictions while imposing point identification.
Our paper is of course related to a large literature on shape restrictions; see \cite{samworth2018special} and \cite{chetverikov2018econometrics} for recent reviews.
We highlight here an important literature on linear Gaussian models focused on adaptivity (which we do not establish), but not applicable to many of the models that motivate us \citep{dumbgen2001multiscale, cai2013adaptive, armstrong:2015}.

The results here are also highly complementary to \cite{chetverikov2017nonparametric} in providing inference for nonparametric IV under shape restrictions while they showed that imposing monotonicity can greatly improve the convergence rate of the estimator -- an observation that additionally motivates our use of test statistics based on shape constrained (instead of unconstrained) estimators.
Finally, we note that our results do not lend themselves computationally for the construction of uniform confidence bands for shape restricted functions -- a problem that has been addressed in different contexts by \cite{rearrange} and \cite{horowitz2017nonparametric}.

The remainder of the paper is organized as follows.
In Section \ref{sec:testex} we show how to implement our tests in a linear instrumental variables model with inequality restrictions under both point and partial identification.
Section \ref{sec:testex} further illustrates our results by revisiting the analysis of \cite{angrist1998children}.
Section \ref{sec:gentheory} contains our main theoretical results, while Section \ref{sec:hetero} applies them to conduct inference in the heterogenous demand model of \cite{hausman2016individual}.
Finally, Section \ref{sec:mc} contains a brief simulation study.
All mathematical derivations are included in a series of appendices; see in particular Appendix \ref{sec:examples} for applications of our general results and Appendix \ref{sec:unifcoupling} for coupling results based on \cite{koltchinskii:1994}.

\section{Application for Linear Instrumental Variables} \label{sec:testex}

To fix ideas, we first describe our test in a linear instrumental variables model and illustrate its implementation by revisiting the analysis of \cite{angrist1998children}.
We reserve until later the full mathematical framework and focus here on implementation.

\subsection{Linear Instrumental Variables}\label{sec:testexiv}

As perhaps the simplest possible example, we first consider a linear instrumental variable model in which $\IDpoint\in \Theta\subseteq \mathbf R^{d_\theta}$ is identified through the moment conditions
\begin{equation*}
E_P[(Y - W^\prime \IDpoint)Z] = 0,
\end{equation*}
where $Y$ is a scalar, $W$ and $Z$ are vectors, and $P$ denotes the distribution of $V \equiv (Y,W,Z)$.
We are interested in testing whether $\IDpoint$ belongs to a set $R$ characterized by
\begin{equation}\label{eq:liniv1}
R=\{\theta \in \mathbf R^{d_\theta} : F \theta = f, ~ G \theta \leq g\},
\end{equation}
for known matrices $F$ and $G$ and known vectors $f$ and $g$.

We consider tests based on minimizing the norm of the weighted sample moments as in \cite{sargan:1958} and \cite{hansen:1982}.
To this end, we define the criterion
\begin{equation}\label{eq:liniv0p5}
Q_n(\theta) \equiv \|\hat \Sigma_n \{\frac{1}{n}\sum_{i=1}^n (Y_i - W_i^\prime \theta )Z_i\}\|_2,
\end{equation}
where $\|\cdot\|_2$ is the standard Euclidean norm and $\hat \Sigma_n$ is consistent for $(E[ZZ^\prime U^2])^{-1/2}$ for $U \equiv Y - W^\prime \IDpoint$.
Our analysis then enables us to employ tests based on the statistics
\begin{equation}\label{eq:liniv1p5}
I_n(R) \equiv \min_{\theta \in \Theta \cap  R} \sqrt n Q_n(\theta) \hspace{0.5 in} I_n(\Theta) \equiv \min_{\theta \in \Theta} \sqrt n Q_n(\theta);    
\end{equation}
e.g., we may consider a test that rejects for large values of $I_n(R) - I_n(\Theta)$.
In what follows it will also be helpful to let $\hat \theta_n$ and $\hat \theta_n^{\rm u}$ denote the minimizers of $Q_n$ over $\Theta \cap R$ and $\Theta$ respectively -- i.e.\ $\hat\theta_n$ and $\hat \theta_n^{\rm u}$ are the constrained and unconstrained estimators.

We construct critical values by relying on the multiplier bootstrap \citep{ledoux:talagrand:1988}.
Specifically, let $b \in \{1,\ldots, B\}$ index a bootstrap draw, $\{\omega_i^b\}_{i=1}^n$ be i.i.d.\ independent of the data with $\omega_i^b\sim N(0,1)$, and for any $\theta \in \mathbf R^{d_\theta}$ define
$$\hat {\mathbb W}^b_n(\theta) \equiv \frac{1}{\sqrt n}\sum_{i=1}^n \omega_i^b \{(Y_i - W_i^\prime \theta)Z_i - \frac{1}{n}\sum_{j=1}^n (Y_j - W_j^\prime \theta)Z_j\}, $$
which is a simulated draw of the true (centered) moment functions.\footnote{We follow previous work \citep{LEWBEL1995379, hansen1996inference} in considering Gaussian weights $\{\omega_i\}_{i=1}^n$ because it simplifies the proofs of our main results in Section \ref{sec:gentheory}. We expect our analysis extends to alternative specifications for the distribution of $\{\omega_i\}_{i=1}^n$ -- e.g., for $\omega_i$ following an exponential distribution, which results in a version of the Bayesian bootstrap advocated by \cite{chamberlain2003nonparametric}.} 
We also require an estimator of the derivative of the moment conditions, and to this end we set
$$\hat {\mathbb D}_n[h] \equiv -\frac{1}{n}\sum_{i=1}^n Z_iW_i^\prime h.$$
Here, we can think of $h$ as a local parameter, representing the possible values that the random variable $\sqrt n\{\hat \theta_n-\IDpoint\}$ may take (recall $\hat \theta_n$ is the minimizer of $Q_n$ over $\Theta \cap R$).

Finally, we need to enforce the inequality constraints in the bootstrap in a way that delivers a uniformly valid critical value. 
To this end, we account for the variation in $G_j \hat \theta_n  - g_j$ for each $j$, where $G_j$ is the $j^{th}$ row of $G$ and $g_j$ the $j^{th}$ coordinate of $g$. 
That is, we account for the likelihood that a constraint will bind at the restricted estimator $\hat \theta_n$ when computing $I_n(R) = \sqrt n Q_n(\hat \theta_n)$. 
For this purpose we introduce the set 
\begin{equation}\label{eq:ivVdef}
\hat V_n(\hat \theta_n,R) \equiv \{h \in \mathbf R^{d_\theta} : Fh = 0,~ G_j h \leq \sqrt n \max\{0, -(r_n + G_j \hat \theta_n  - g_j)\} \text{ for all } j\},
\end{equation}
where $r_n>0$ is a slackness parameter whose choice we discuss shortly.
The set $\hat V_n(\hat \theta_n,R)$ can be thought of as a local version of $R$, approximating the set of values $h$ that could equal $\sqrt n\{\hat \theta_n - \IDpoint\}$.
Our bootstrap approximations to $I_n(R)$ and $I_n(\Theta)$ are then 
\begin{align}
\hat U_n^b(R) & \equiv \min_{h \in \hat V_n(\hat \theta_n,R)} \|\hat \Sigma_n\{\hat {\mathbb W}_n^b(\hat \theta_n) + \hat {\mathbb D}_n[h]\}\|_2  \label{eq:linvi2} \\
\hat U_n^b(\Theta) & \equiv \min_{h \in \mathbf R^{d_\theta}} \|\hat \Sigma_n\{\hat {\mathbb W}_n^b(\hat \theta_n^{\rm u}) + \hat {\mathbb D}_n[h]\}\|_2.  \label{eq:linvi3}
\end{align}
Thus, we may obtain a level $\alpha$ test by rejecting whenever the test statistic $I_n(R)-I_n(\Theta)$ exceeds the $1-\alpha$ quantile of $\hat U_n^b(R) - \hat U_n^b(\Theta)$ across the $B$ bootstrap draws.
The main assumption required for the test to be asymptotically valid is that $\IDpoint$ be strongly identified -- i.e.\ $\IDpoint$ can be consistently estimated uniformly in $P$.

The critical value depends on the choice of $r_{n}$. 
When applied to linear instrumental variables, our asymptotic theory requires that $r_{n}$ tend to zero slower than the convergence rate of the restricted estimator, which is $1/\sqrt n$.
Heuristically, when $r_{n}$ tends to zero any constraint that is not binding at $\IDpoint$ will also not be binding in the bootstrap with probability approaching one (under pointwise in $P$ asymptotics). 
Consequently inference is not asymptotically conservative for a fixed data generating process. 
Setting $r_n \to 0$ while satisfying $r_n \sqrt n \to \infty$ leads to uniformly valid inference with constraints only being conservatively enforced when they are within order $1/\sqrt n$ of binding at $\IDpoint$.
Setting $r_{n}=+\infty$ is always theoretically valid, but it may be conservative and  result in a loss of power. 
Other, smaller choices of $r_{n}$ can lead to smaller, valid critical values and so may result in more powerful tests and tighter confidence intervals than $r_{n}=+\infty$. 

Intuitively, $r_{n}$ is meant to quantify the sampling uncertainty in $G\{\hat \theta_n - \IDpoint\}$.
Since the distribution of $\hat \theta_n$ cannot be uniformly consistently estimated, we suggest linking $r_n$ to the degree of sampling uncertainty in $G\{\hat \theta_n^{\rm u} - \IDpoint\}$ instead.
Specifically, for $\hat \theta_n^{\rm u\star}$ a ``bootstrap" analogue of $\hat \theta_n^{\rm u}$ and some $\gamma_n \to 0$, we recommend setting $r_n$ to satisfy
\begin{equation}\label{eq:liniv4}
 P(\max_{j} G_{j}\{\hat{\theta}^{\rm u}_n-\hat{\theta}^{\rm u\star}_n \}\leq r_{n}|\text{Data})=1-\gamma_n.
 \end{equation}
This approach changes the problem of selecting $r_{n}$ into the problem of selecting $\gamma_n$.
However, $\gamma_n$ is more interpretable: If we employed $\hat V_n(\hat \theta_n^{\rm u},R)$ in place of $\hat V_n(\hat \theta_n,R)$ in \eqref{eq:linvi2}, then a Bonferroni bound implies that the test that rejects whenever $I_n(R)-I_n(\Theta)$ exceeds the $1-\alpha$ quantile of $\hat U_n^b(R)-\hat U_n^b(\Theta)$ has asymptotic size at most $\alpha + \gamma_n$ even if $\gamma_n$ is fixed with $n$.\footnote{While we may replace $\hat V_n(\hat \theta_n,R)$ with $\hat V_n(\hat \theta_n^{\rm u},R)$ in identified models, in partially identified models we employ $\hat V_n(\hat \theta_n,R)$ due to the identified set potentially not being a subset of $R$ under the null hypothesis.} 
In particular, if we employed the $1-\alpha+\gamma_n$ quantile of $\hat U_n^b(R)-\hat U_n^b(\Theta)$ as a critical value instead, then the resulting test would have asymptotic size at most $\alpha$ (even if $\gamma_n$ is fixed). 
In simulations, however, we find the described bound to be pessimistic in that, when setting $r_n$ according to \eqref{eq:liniv4}, our test has a rejection probability under the null hypothesis of at most $\alpha$ for a wide range of choices of $\gamma_n$.

\begin{remark}\label{rm:linearCI} \rm
Our results may be employed to obtain confidence regions for a coordinate of $\IDpoint$ while imposing restrictions of the form $G\IDpoint \leq g$ on $\IDpoint$ (e.g., sign or monotonicity restrictions on  $w\mapsto w^\prime \IDpoint$).
For example, for $\theta_k$ the $k^{th}$ coordinate of $\theta \in \mathbf R^{d_\theta}$ let
$$R_\lambda = \{\theta \in \mathbf R^{d_\theta} : \theta_k = \lambda, ~ G\theta \leq g\},$$
which is a special case of \eqref{eq:liniv1}.
We may then obtain a confidence region for the $k^{th}$ coordinate of $\IDpoint$ by conducting test inversion in $\lambda$ employing the test based on $I_n(R_\lambda)-I_n(\Theta)$; see also Remark \ref{rm:combine} for alternative constructions based on our analysis. \qed
\end{remark}

\begin{remark}\label{rm:stud} \rm
In certain applications it may be desirable to studentize the constraints in our bootstrap approximation -- i.e.\ replace $G_j$ and $g_j$ by $G_j/\hat \sigma_j$ and $g_j/\hat \sigma_j$ everywhere in \eqref{eq:ivVdef} (and in \eqref{eq:liniv4} if employed).
In the empirical analysis below we proceed in this manner by setting $\hat \sigma_j^2$ to be an estimate of the asymptotic variance of $\sqrt n G_j\{\hat \theta_n^{\rm u} - \theta_0\}$. \qed
\end{remark}

\subsubsection{Fertility and Labor Supply: LATE}\label{sec:angristevans}

We illustrate the preceding discussion by revisiting the study by \cite{angrist1998children} on the causal effect of childbearing on female labor force participation.
Like \cite{angrist1998children}, we employ the 1980 Census Public Use Micro Sample restricted to mothers aged 21-35 with at least two children, and set: (i) $D\in \{0,1\}$ to indicate whether a mother has more than two children (the treatment); (ii) $Y\in \{0,1\}$ to indicate whether a mother is employed (the outcome of interest); and (iii) $Z\in \{0,1\}$ to indicate whether the first two children are of the same sex (the instrument).
We further adopt the heterogeneous treatment effects model of \cite{imbens1994identification} and let $Y_d$ denote the potential outcome under treatment status $d\in\{0,1\}$ and employ ``$\texttt{C}$," ``\texttt{NT}," and ``\texttt{AT}" to denote compliers, never takers, and always takers.

\begin{figure}[t!]
\begin{subfigure}{1 \textwidth}{\begin{center}
    \begin{tikzpicture}[scale = 0.7]
    \pgfplotsset{ymin = -0.6, ymax = 0.5}
    \begin{axis}[legend columns = 2, /pgf/number format/1000 sep={}, xtick = data,  xticklabels = {$<9$,9,10,11,12,13,14,15,16, $>16$}, width = \textwidth, height = 0.5\textwidth]
    \addplot [mark= square, violet, mark options = {solid}] table [x=XPos, y=2SLS]{Data4Graphs/MonConfidenceIntervalLATE.txt};
    \addlegendentry{2SLS Estimate};
    \addplot [mark = triangle, red, mark options = {solid}] table [x=XPos, y=R2SLS]{Data4Graphs/MonConfidenceIntervalLATE.txt};
    \addlegendentry{Restr. Estimate};
    \end{axis}
\end{tikzpicture}\end{center}}
\end{subfigure}\\
\begin{subfigure}{1 \textwidth}{\begin{center}
    \begin{tikzpicture}[scale = 0.7]
    \pgfplotsset{ymin = -0.6, ymax = 0.5}
    \begin{axis}[legend columns = 2, /pgf/number format/1000 sep={}, xtick = data,  xticklabels = {$<9$,9,10,11,12,13,14,15,16, $>16$}, width = \textwidth, height = 0.5\textwidth]
    \addplot [mark= square, violet, dashed, mark options = {solid}] table [x=XPos, y=LATEUB]{Data4Graphs/MonConfidenceIntervalLATE.txt};
    \addplot [mark = triangle, red, dashed, mark options = {solid}] table [x=XPos, y=RESLB]{Data4Graphs/MonConfidenceIntervalLATE.txt};
    \addplot [mark= square, violet, dashed, mark options = {solid}] table [x=XPos, y=LATELB]{Data4Graphs/MonConfidenceIntervalLATE.txt};
    \addplot [mark = triangle, red, dashed, mark options = {solid}] table [x=XPos, y=RESUB]{Data4Graphs/MonConfidenceIntervalLATE.txt};
    \addlegendentry{2SLS CI};
    \addlegendentry{Mon. Restr. CI};
    \end{axis}
\end{tikzpicture}\end{center}}
\end{subfigure}\\
\begin{subfigure}{1 \textwidth}{\begin{center}
    \begin{tikzpicture}[scale = 0.7]
    \pgfplotsset{ymin = -0.6, ymax = 0.5}
    \begin{axis}[legend columns = 2, /pgf/number format/1000 sep={}, xtick = data,  xticklabels = {$<9$,9,10,11,12,13,14,15,16, $>16$}, width = \textwidth, height = 0.5\textwidth]
    \addplot [mark= square, violet, dashed, mark options = {solid}] table [x=XPos, y=LATEUB]{Data4Graphs/PosConfidenceIntervalLATE.txt};
    \addplot [mark = triangle, red, dashed, mark options = {solid}] table [x=XPos, y=RESLB]{Data4Graphs/PosConfidenceIntervalLATE.txt};
    \addplot [mark= square, violet, dashed, mark options = {solid}] table [x=XPos, y=LATELB]{Data4Graphs/PosConfidenceIntervalLATE.txt};
    \addplot [mark = triangle, red, dashed, mark options = {solid}] table [x=XPos, y=RESUB]{Data4Graphs/PosConfidenceIntervalLATE.txt};
    \addlegendentry{2SLS CI};
    \addlegendentry{Mon.+Neg. Restr. CI};
    \end{axis}
\end{tikzpicture}\end{center}}
\end{subfigure}
\caption{First Panel: Unconstrained and shape restricted LATE estimates (imposing monotonicty or monotonicity and negativity yield the same estimates). Second and Third Panels: 95$\%$ Confidence intervals for LATE at different education levels. }\label{fig:LATEci}
\end{figure}

\cite{angrist1998children} document that the impact of childbearing on labor force participation depends on observable characteristics.
In particular, their two stage least squares (2SLS) estimates suggest a negative impact of childbearing on labor force participation across different levels of schooling, but that the magnitude of the impact decreases with schooling -- a phenomenon that may reflect that more educated mothers have a stronger attachment to the labor force.
To formally examine this claim, we introduce dummy variables $S$ for each year of schooling between 9 and 16 and for the categories ``less than 9" and ``more than 16."
Defining the local average treatment effects
\begin{equation*}
\text{LATE}(S) \equiv E[Y_1 - Y_0|S,\texttt{C}]
\end{equation*}
we then test whether: (i) $\text{LATE}(\cdot)$ is increasing in schooling, and (ii) $\text{LATE}(\cdot)$ is increasing in schooling and nonpositive.
Both hypotheses fall within the framework of the preceding section because $\text{LATE}(\cdot)$ is identified through linear moment restrictions and the hypothesized restrictions are linear in $\text{LATE}(\cdot)$.
Employing five thousand bootstrap replications and setting $r_n = +\infty$ or $r_n$ as suggested in \eqref{eq:liniv4} with $\gamma_n = 0.05$ yields in this case equal $p$-values that fail to reject either null hypothesis. The $p$-values for $\text{LATE}(\cdot)$ being nondecreasing is 0.21 and for it being nondecreasing and nonpositive is 0.394.

In Figure \ref{fig:LATEci} we study the values of $\text{LATE}(S)$ at different schooling levels $S$.
The first panel displays the unconstrained 2SLS estimates and their monotonicity restricted counterparts -- the latter are negative and hence additionally demanding nonpositivity does not change the estimates.
Unfortunately, two sided confidence regions based on the (pointwise in $P$) asymptotic distribution of the shape-restricted 2SLS estimator can asymptotically undercover the true parameter.
In the second panel of Figure \ref{fig:LATEci} we instead proceed as in Remark \ref{rm:linearCI} to obtain $95\%$ confidence intervals while imposing monotonicity and again selecting $r_n$ by setting $\gamma_n = 0.05$ in \eqref{eq:liniv4}.
Employing the monotonicity restriction in this manner yields confidence intervals that are sometimes substantially shorter than their 2SLS counterparts.
Notably, we observe lower upper ends for the restricted confidence intervals at the lower education levels and higher lower ends at higher education levels.
As shown in the third panel of Figure \ref{fig:LATEci}, additionally imposing that $\text{LATE}(\cdot)$ be nonpositive mostly reduces the upper bound of our confidence intervals at higher education levels.

\subsection{Partial Identification} \label{sec:testexpi}

We next illustrate the implementation of our results in a partially identified setting.
With an eye towards extending the preceding empirical analysis to study average treatment effects (ATEs), we maintain that the parameter of interest $\IDpoint\in \Theta\subseteq \mathbf R^{d_\theta}$ satisfies
\begin{equation}\label{eq:liniv5}
E_P[(Y - W^\prime \IDpoint)Z] = 0,    
\end{equation}
but no longer assume $\IDpoint$ is identified by \eqref{eq:liniv5}.
Instead, we define the identified set
\begin{equation}\label{eq:liniv6}
\IDset \equiv \{\theta \in \Theta : E_P[(Y - W^\prime \theta)Z] = 0\}    
\end{equation}
and consider the problem of testing whether the intersection of $\IDset$ and $R$ is nonempty (i.e.\ $\IDset \cap R\neq \emptyset$).
Such hypotheses can be employed, for instance, to build confidence regions for functionals of the identified set; see Remark \ref{rm:IMCI} below.
We also now set
\begin{equation}\label{eq:liniv7}
R=\{\theta \in \mathbf R^{d_\theta} : \Upsilon_F(\theta) = 0, ~ G \theta \leq g\},
\end{equation}
for $\Upsilon_F$ a known possibly nonlinear function -- e.g., $\Upsilon_F(\theta) = F\theta - f$ recovers \eqref{eq:liniv1}.

We continue to rely on the statistics $I_n(R)$ and $I_n(\Theta)$ (as in \eqref{eq:liniv1p5}) for inference.
However, since in many settings in which $\IDpoint$ fails to be identified by \eqref{eq:liniv5} we will have that the dimension of $Z$ is smaller than that of $W$, in what follows we assume for ease of exposition that $I_n(\Theta) = 0$ (almost surely); see Section \ref{subsec:critrec} for a general discussion.
Another distinction relative to Section \ref{sec:testexiv} is that the choice of $\hat \Sigma_n$ (as in \eqref{eq:liniv0p5}) may need to be modified in settings in which $U \equiv Y - W^\prime \IDpoint$ cannot be consistently estimated due to $\IDpoint$ being partially identified.
In such instances we may, for example, set
$$\hat \Sigma_n \equiv (\frac{1}{n}\sum_{i=1}^n Z_iZ_i^\prime(Y_i - W_i^\prime \hat \theta_n^{\rm u})^2)^{-1/2},$$
where we now interpret $\hat \theta_n^{\rm u}$ as the minimum norm minimizer of $Q_n$ over $\Theta$. 
While the choice of $\hat \Sigma_n$ has an impact on how local power is directed, we note that the test has correct asymptotic size provided $\hat \Sigma_n$ converges in probability to a non-stochastic limit.

Our bootstrap procedure requires two modifications relative to our preceding discussion.
First, because in \eqref{eq:liniv7} we consider nonlinear equality constraints, we now set
$$\hat V_n(\theta,R) \equiv \{h \in \mathbf R^{d_\theta} : \Upsilon_F(\theta + \frac{h}{\sqrt n}) = 0, ~ G_jh\leq \sqrt n\max\{0,-(r_n + G_j \theta - g_j)\} \text{ for all } j\}$$
(notice that if $\Upsilon_F(\theta) = F\theta - f$, then we recover \eqref{eq:ivVdef}).
A distinction with Section \ref{sec:testexiv} is that if one aims to employ \eqref{eq:liniv4} to select  $r_n$, then an alternative to an unrestricted estimator $\hat \theta_n^{\rm u}$ may be necessary; see Section \ref{sec:angristevanspi} for an example.
Second, our bootstrap approximation employs an estimator $\hat \Theta_n^{\rm r}$ for $\IDset \cap R$.
To this end, we set
\begin{equation*}
    \hat \Theta_n^{\rm r} \equiv \{\theta \in \Theta \cap R: Q_n(\theta) \leq \inf_{\theta \in \Theta \cap R} Q_n(\theta) + \tau_n\} 
\end{equation*}
where $\tau_n \geq 0$ is a bandwidth whose choice we discuss shortly -- i.e.\ $\hat \Theta_n^{\rm r}$ is the set of ``near" minimizers of $Q_n$ over $\Theta \cap R$.
Our bootstrap approximation to $I_n(R)$ then equals
$$\hat U_n^b(R) \equiv \min_{\theta \in \hat \Theta_n^{\rm r}} \min_{h \in \hat V_n(\theta,R)} \|\hat \Sigma_n\{\hat {\mathbb W}_n^b(\theta) + \hat {\mathbb D}_n[h]\}\|_2.$$
Thus, to obtain a level $\alpha$ test we reject the null hypothesis whenever $I_n(R)$ exceeds the $1-\alpha$ quantile of $\hat U_n^b(R)$ across bootstrap draws. 
Paralleling Section \ref{sec:testexiv}, a principal assumption for the test to be asymptotically valid is that $\IDset$ be strongly identified. 

When specialized to the current setting, our asymptotic theory requires that $\tau_n$ tend to zero.
It is theoretically valid to set $\tau_n = 0$, which simplifies the computation of our bootstrap statistic -- e.g., let $\hat \Theta_n^{\rm r} = \{\hat \theta_n\}$ for any $\hat \theta_n$ minimizing $Q_n$ over $\Theta\cap R$ to recover \eqref{eq:linvi2}.
However, setting $\tau_n = 0$ can result in lower power in applications for which the corresponding $\hat \Theta_n^{\rm r}$ is not consistent for $\IDset \cap R$ (in the Hausdorff metric) --  
to ensure consistency, $\tau_n$ must in addition satisfy $\tau_n\sqrt n \to \infty$.
For applications in which it is desirable to set $\tau_n > 0$, we propose a procedure inspired by \cite{romano:shaikh:2010}.
Specifically, for any set $K\subseteq \Theta \cap R$ we define the quantile $\hat q_n(K)$ according to
$$P(\sup_{\theta \in K} \|\hat \Sigma_n \Bemp(\theta)\|_2 \leq \hat q_n(K)|\text{Data}) = 1-\gamma_n$$
where $\gamma_n \in (0,1)$.
Letting $S_1 \equiv \Theta \cap R$, we then inductively define $S_{j+1} \equiv \{\theta \in \Theta\cap R: \sqrt n Q_n(\theta) \leq \hat q_n(S_j)\}$ noting that by construction $S_{j+1}\subseteq S_j$.
To select $\tau_n$, we proceed inductively until we find $S_j = \emptyset$, in which case we set $\tau_n = 0$, or $S_{j+1} = S_j \neq \emptyset$, in which case we set $\tau_n = \hat q_n(S_j)$.
Heuristically, under such a choice of $\tau_n$, the set $\hat \Theta_n^{\rm r}$ may be interpreted as a $1-\gamma_n$ confidence region for $\IDset \cap R$.
While power considerations suggest setting $\gamma_n$ to tend to zero, for practical considerations we suggest simply setting $1-\gamma_n$ to be a high quantile fixed with $n$ (e.g., $1-\gamma_n = 0.8$).

\begin{remark}\label{rm:IMCI} \rm
The introduced test can be employed to obtain confidence regions for functionals of the identified set satisfying the coverage requirement advocated by \cite{imbens:manski:2004}. 
Specifically, given a functional $\Upsilon_F:\Theta \to \mathbf R$ we may set
$$R_\lambda = \{\theta \in \mathbf R^{d_\theta} : \Upsilon_F(\theta) = \lambda, G\theta \leq g\}$$
and obtain the desired confidence region by conducting test inversion in $\lambda$ of the null hypothesis that the set $\IDset \cap R_\lambda$ is not empty. \qed
\end{remark}

\subsubsection{Fertility and Labor Supply: ATE}\label{sec:angristevanspi}

Returning to our analysis of the causal impact of fertility on female labor force participation, we next turn to estimating the average treatment effect at different education levels $S$ (denoted $\text{ATE}(S)$).
Following the literature, we decompose $\text{ATE}(S)$ into 
\begin{equation}\label{eq:atede}
\text{LATE}(S)P(\texttt{C}|S) + E[Y_1-Y_0|S,\texttt{AT}]P(\texttt{AT}|S) + E[Y_1-Y_0|\texttt{NT},S]P(\texttt{NT}|S),
\end{equation}
where recall \texttt{C}, \texttt{AT}, and \texttt{NT} denote ``compliers," ``always takers," and ``never takers."
With the exception of $E[Y_0|\texttt{AT},S]$ and $E[Y_1|\texttt{NT},S]$, all terms in \eqref{eq:atede} can be identified through linear moment restrictions.\footnote{Technically, the moment equations have the structure $E_P[(Y_\jmath - W_{\jmath}^\prime \IDpoint)Z_\jmath]=0$ with the instruments $Z_\jmath$ not being common across all $1\leq \jmath \leq \mathcal J$ equations. The bootstrap implementation in this case, formally studied in Section \ref{sec:gentheory}, is identical with only $\Bemp$ and $\hat {\mathbb D}_n$ being modified in the natural way.}
Because $S$ has ten support points, we obtain sixty moments and eighty parameters so that $I_n(\Theta) = 0$ almost surely.

\begin{figure}[t!]
\begin{center}
    \begin{tikzpicture}[scale = 0.9]
    \pgfplotsset{ymin = -0.7, ymax = 0.7}
    \begin{axis}[legend columns = 3, /pgf/number format/1000 sep={}, xtick = data,  xticklabels = {$<9$,9,10,11,12,13,14,15,16, $>16$}, ytick = {-0.6, -0.4, -0.2, 0, 0.2, 0.4, 0.6}, width = \textwidth, height = 0.5\textwidth]
    \addplot [mark= square, violet, dashed, mark options = {solid}] table [x=XPos, y=UNLBCI]{Data4Graphs/MonConfidenceIntervalATE.txt};
    \addplot [mark=triangle, red, dashed, mark options = {solid}] table [x=XPos, y=RESLBCI]{Data4Graphs/MonConfidenceIntervalATE.txt};
    \addplot [mark=x, blue, dashed, mark options = {solid}] table [x=XPos, y=RESLBCI]{Data4Graphs/PosConfidenceIntervalATE.txt};
    \addlegendentry{Unr.};
    \addlegendentry{Mon. Restr.};
    \addlegendentry{Mon.+Neg. Restr.};
    \addplot [mark= square, violet, dashed, mark options = {solid}] table [x=XPos, y=UNUBCI]{Data4Graphs/MonConfidenceIntervalATE.txt};
    \addplot [mark=triangle, red, dashed, mark options = {solid}] table [x=XPos, y=RESUBCI]{Data4Graphs/MonConfidenceIntervalATE.txt};
    \addplot [mark=x, blue, dashed, mark options = {solid}] table [x=XPos, y=RESUBCI]{Data4Graphs/PosConfidenceIntervalATE.txt};
    \end{axis}
\end{tikzpicture}\end{center}\vspace{-0.1 in}\caption{95$\%$ Confidence intervals for $\text{ATE}$ at different education levels. ``Unr." uses bounds implied by $Y_d\in\{0,1\}$; ``Mon. Restr." adds that average treatment effects be increasing in education for all types; ``Mon.+Neg. Restr." also requires they be  negative.}\label{fig:ATEci}
\end{figure}

Following our analysis of $\text{LATE}(S)$ we conduct inference on $\text{ATE}(S)$ under three increasingly stringent set of (linear) restrictions: (i) The logical bounds implied by $Y_d\in \{0,1\}$; (ii) Adding to (i) that the average treatment effect be increasing in schooling among all types (i.e.\ \texttt{C}, \texttt{NT}, and \texttt{AT}); (iii) Adding to (ii) that average treatment effects be nonpositive for all levels of education and types.
Figure \ref{fig:ATEci} reports the resulting $95\%$ confidence regions obtained through the approach described in Remark \ref{rm:IMCI} -- here, the restriction $G\theta \leq g$ imposes the described shape constraints while the nonlinear restriction $\Upsilon_F(\theta) = 0$ corresponds to imposing a hypothesized value for $\text{ATE}(S)$ through \eqref{eq:atede}.
In our bootstrap approximation, we set $\tau_n = 0$ and selected $r_n$ according to \eqref{eq:liniv4} with $\gamma_n = 0.05$ and where, when necessary, we used the distribution of estimators of identified parameters for their partially identified counterparts.\footnote{E.g., for the constraint $E[Y_1|\texttt{NT},S]\leq 1$ we substituted the corresponding $G_j\{\hat \theta_n^{\rm u} - \hat \theta_n^{\rm u^\star}\}$ term in \eqref{eq:liniv4} with a mean zero normal distribution with the variance of the estimator for $E[Y_0|\texttt{NT},S]$.}
We do not report estimates of the identified sets for $\text{ATE}(S)$ as they are very close to the obtained confidence intervals: On average the bounds of the confidence intervals exceed the bounds of estimates of the identified set by 0.011.
Nonetheless, the unrestricted confidence intervals are large as the estimates for the identified set are themselves large -- a result driven by the low proportion of compliers ($5\%$ on average across schooling levels).
Imposing monotonicity across types carries identifying information on the upper end of the identified set at low levels of education and on the lower end of the identified set at high levels of education.
Additionally imposing nonpositivity sharpens the upper bound of the identified set at all schooling levels. 
The resulting confidence regions sign $\text{ATE}(S)$ at all education levels (weakly) smaller than 12 as strictly negative, though very close to zero.

Finally, as a preview of our general analysis in Section \ref{sec:gentheory}, in Table \ref{table:aggate} we employ the same shape restrictions to report estimates and $95\%$ confidence intervals for the identified sets of the average treatment effects for: High School Dropouts ($\text{edu}\in [9,12)$), College Dropouts ($\text{edu}\in [13,15)$), College Graduates ($\text{edu} \geq 16$) and the overall average treatment effect.
These confidence regions are obtained through test inversion after noting that a hypothesized value for the average treatment effect of a subgroup can be written as a nonlinear moment restriction in $\IDpoint$ through \eqref{eq:atede} -- nonlinear moment restrictions fall within our general framework but outside the scope of Section \ref{sec:testexpi}. Overall the impact of imposing shape restrictions parallels the results in Figure \ref{fig:ATEci}.

\begin{table}[t!] \centerline {\small
\begin{tabular}{c cc cc cc}\hline \hline
                    & \multicolumn{2}{c}{Unrestricted}      & \multicolumn{2}{c}{Mon. Restr.}     & \multicolumn{2}{c}{Mon.+Neg Restr.} \\ \cline{2-3} \cline{4-5} \cline{6-7}
Subgroup            & Estimate          & $95\%$ CI         &  Estimate       & $95\%$ CI         & Estimate        & $95\%$ CI\\ \hline
HS   Drop           & [-0.520,0.426]    & [-0.526,0.432]    &  [-0.489,0.346] & [-0.500,0.356]    & [-0.489,-0.008] & [-0.501,-0.003]   \\
Coll. Drop          & [-0.561,0.380]    & [-0.566,0.385]    &  [-0.447,0.325] & [-0.460,0.337]    & [-0.447,-0.004] & [-0.462,0.000]    \\
Coll. Grad          & [-0.579,0.375]    & [-0.586,0.382]    &  [-0.446,0.328] & [-0.462,0.339]    & [-0.446,-0.002] & [-0.464,0.000]    \\
All                 & [-0.545,0.395]    & [-0.547,0.398]    &  [-0.467,0.328] & [-0.477,0.338]    & [-0.467,-0.008] & [-0.478,-0.003]   \\
\hline \hline
\end{tabular}}\caption{Point Estimates and $95\%$ confidence intervals for the average treatment effect at different groups defined by schooling levels under different shape restrictions.}\label{table:aggate}
\end{table}

\section{General Analysis}\label{sec:gentheory}

We next develop a general inferential framework that encompasses the tests discussed in Section \ref{sec:testex}.
The class of models we consider are those in which the parameter of interest $\IDpoint \in \Theta$ satisfies a finite number $\mathcal J$ of conditional moment restrictions
\begin{equation*}
E_P[\rho_\jmath(X,\IDpoint)|Z_\jmath] = 0 \text{ for } 1 \leq \jmath \leq \mathcal J
\end{equation*}
with $\rho_\jmath : \mathbf X \times \Theta \to \mathbf R$, $X\in \mathbf X$, and $Z_\jmath \in \mathbf Z_\jmath$.
For notational simplicity, we also let $Z \equiv (Z_1,\ldots, Z_{\mathcal J})$ and $V \equiv (X, Z)$ with $V\sim P\in \mathbf P$.
In some of the applications that motivate us, the parameter $\IDpoint$ is not identified.
We therefore define the identified set
\begin{equation*}
\IDset \equiv \{\theta \in \Theta : E_P[\rho_{\jmath}(X,\theta)|Z_{\jmath}] = 0 \text{ for } 1\leq \jmath \leq \mathcal J \}
\end{equation*}
and employ it as the basis of our statistical analysis -- we emphasize that $\IDset$ depends on $P$, but leave such dependence implicit to simplify notation.
For a set $R$ of parameters satisfying a conjectured restriction, we develop a test for the hypothesis
\begin{equation}\label{def:hyp}
H_0: \IDset \cap R \neq \emptyset \hspace{0.5 in} H_1: \IDset \cap R = \emptyset ;
\end{equation}
i.e.\ we devise a test of whether at least one element of the identified set satisfies the posited constraint.
In what follows, we denote the set of distributions $P\in \mathbf P$ satisfying the null hypothesis in \eqref{def:hyp} by $\mathbf P_0$.
We also note that in an identified model, a test of \eqref{def:hyp} is equivalent to a test of whether $\IDpoint$ itself satisfies the hypothesized constraint.

The defining elements determining the type of applications encompassed by \eqref{def:hyp} are the choices of $\Theta$ and $R$.
In imposing restrictions on $\Theta$ and $R$ we therefore aim to allow for a general framework while simultaneously ensuring enough structure for a fruitful asymptotic analysis.
To this end, we require $\Theta$ to be a subset of a complete vector space $\mathbf B$ with norm $\|\cdot\|_{\mathbf B}$ (i.e.\ $(\mathbf B,\|\cdot\|_{\mathbf B}$) is a Banach space) and consider sets $R$ satisfying
\begin{equation}\label{def:Rset}
R = \{\theta \in \mathbf B : \Upsilon_F(\theta) = 0 \text{ and } \Upsilon_G(\theta) \leq 0\},
\end{equation}
where $\Upsilon_F:\mathbf B \rightarrow \mathbf F$ and $\Upsilon_G:\mathbf B \rightarrow \mathbf G$ are known maps.
Our first assumption formalizes the basic structure of the hypothesis testing problem we study.

\begin{assumption}\label{ass:param}
(i) $\{V_i\}_{i=1}^n$ is i.i.d.\ with $V\sim P\in \mathbf P$;
(ii) $\Theta \subseteq \mathbf B$, where $(\mathbf B,\|\cdot\|_{\mathbf B})$ is a Banach space; 
(iii) $\Upsilon_F:\mathbf B \rightarrow \mathbf F$ and $\Upsilon_G:\mathbf B \rightarrow \mathbf G$, where $(\mathbf F,\|\cdot\|_{\mathbf F})$ is a Banach space and $(\mathbf G,\|\cdot\|_{\mathbf G})$ is an AM space with order unit $\mathbf {1_G}$.
\end{assumption}

Through Assumption \ref{ass:param}(i) we focus on the i.i.d.\ setting, though extensions to other sampling frameworks are feasible.
Assumption \ref{ass:param}(ii) allows us to address parametric, semiparametric, and nonparametric models, while Assumption \ref{ass:param}(iii) allows $\Upsilon_F$ to impose both finite dimensional or infinite dimensional equality restrictions.
Assumption \ref{ass:param}(iii) further requires that $\Upsilon_G$ take values in an AM space $\mathbf G$ -- we provide an overview of AM spaces in the supplemental appendix.
Heuristically, the key properties of $\mathbf G$ are: (i) $\mathbf G$ is a vector space equipped with a partial order ``$\leq$"; (ii) The partial order and the vector space operations interact in the same manner they do on $\mathbf R$ (e.g.\ if $\theta_1 \leq \theta_2$, then $\theta_1 + \theta_3 \leq \theta_2 + \theta_3$); and (iii) The order unit $\mathbf {1_G} \in \mathbf G$ is an element such that for any $\theta \in \mathbf G$ there exists a scalar $\lambda > 0$ satisfying $|\theta| \leq \lambda \mathbf {1_G}$ (e.g.\ when $\mathbf G = \mathbf R^d$ we may set $\mathbf {1_G} \equiv (1,\ldots, 1)^\prime \in \mathbf R^d$).
These properties of an AM space will prove instrumental in our analysis.
In particular, the order unit ${\bf 1_G}$ will provide a crucial link between the partial order ``$\leq$", the norm $\|\cdot\|_{\mathbf G}$, and (through smoothness of $\Upsilon_G$) allow us to leverage a rate of convergence in $\mathbf B$ to build a suitable sample analogue to the local parameter space.

\subsection{Main Results}\label{sec:strongapp}

Our analysis centers around a statistic $I_n(R)$ that constitutes a ``building block" for different tests of \eqref{def:hyp} -- e.g., it may be employed to implement a generalization of the $J$-test of \cite{sargan:1958} and \cite{hansen:1982} or the incremental $J$-test of \cite{eichenbaum1988time}.
In this section we first introduce $I_n(R)$, obtain an approximation to its finite sample distribution, and devise a bootstrap procedure for estimating its quantiles. 
Together, these results allow us to establish the asymptotic validity of different tests.

\subsubsection{The Building Block}

We first introduce the statistic $I_n(R)$ that we employ to build different tests.
To this end, for each instrument $Z_{\jmath}$ we consider transformations $\{q_{k,\jmath}\}_{k=1}^{k_{n,\jmath}}$ and let $q_{\jmath}^{k_{n,\jmath}}(z_{\jmath}) \equiv (q_{1,\jmath}(z_\jmath),\ldots,q_{k_{n,\jmath},\jmath}(z_\jmath))^\prime$.
Recalling that $Z \equiv (Z_{1},\ldots, Z_{\mathcal J})$, we further set $k_n \equiv \sum_{\jmath =1}^{\mathcal J} k_{n,\jmath}$, $q^{k_n}(z) \equiv (q_{1}^{k_{n,1}}(z_1)^\prime,\ldots, q_{\mathcal J}^{k_{n,\mathcal J}}(z_{\mathcal J})^\prime)^\prime$, $\rho(x,\theta) \equiv (\rho_1(x,\theta),\ldots,\rho_{\mathcal J}(x,\theta))^\prime$, and let
\begin{equation*}
\rho(X_i,\theta)*q^{k_n}(Z_i) \equiv \left(\begin{array}{c}  \rho_{1}(X_i,\theta) q_{1}^{k_{n,1}}(Z_{i,1}) \\ \vdots \\  \rho_{\mathcal J}(X_i,\theta) q_{\mathcal J}^{k_{n,\mathcal J}}(Z_{i,\mathcal J})\end{array}\right);
\end{equation*}
i.e.\ for each $\theta$ we take the product of each ``residual" $\rho_\jmath(X,\theta)$ with the transformations of its respective instrument $Z_{\jmath}$.
For a $k_n\times k_n$ matrix $\hat \Sigma_n$, we then define
\begin{equation*}
Q_n(\theta) \equiv \|\frac{1}{n} \sum_{i=1}^n \rho(X_i,\theta)*q^{k_n}(Z_i)\|_{\hat \Sigma_n,p} ,
\end{equation*}
where $\|a\|_{\hat \Sigma_n,p} \equiv \|\hat \Sigma_n a\|_p$ and $\|\cdot\|_p$ is the $p$-norm on $\mathbf R^{k_n}$ for any $p\geq 2$ -- i.e.\ $\|a\|_p \equiv  (\sum_{i =1}^d |a^{(i)}|^p)^{1/p}$ for any $a \equiv (a^{(1)},\ldots,a^{(d)})^\prime \in \mathbf R^d$.
Letting $\Theta_n \cap R$ be a finite dimensional subset of $\Theta\cap R$ that grows dense in $\Theta \cap R$ \citep{chen:2006}, we then define $I_n(R)$ to equal
\begin{equation*}
I_n(R) \equiv \inf_{\theta \in \Theta_n \cap R} \sqrt n Q_n(\theta) .
\end{equation*}
We note that setting $p =2$ is often computationally attractive.
However, we allow for $p>2$ because higher values of $p$ enable us to establish distributional approximations under weaker conditions on the number of unconditional moments $k_n$.

Heuristically, $\sqrt n Q_n$ should diverge to infinity when evaluated at any $\theta \notin \IDset$ and remain ``stable" when evaluated at a $\theta \in \IDset$.
Thus, examining the minimum of $\sqrt nQ_n$ over $R$ should reveal whether there is a $\theta$ that simultaneously makes $\sqrt nQ_n(\theta)$ ``stable" ($\theta \in \IDset$) and satisfies the conjectured restriction ($\theta \in R$).
This intuition suggests $I_n(R)$ may be employed as a test statistic that is similar in spirit to the $J$-statistic of \cite{hansen:1982}.
Alternatively, we may build a test by considering the recentered test statistic
$$I_n(R) - I_n(\Theta),$$
which aims power in a different direction than $I_n(R)$ \citep{chen2018overidentification}.
Conceptually, it is important to note that $I_n(\Theta)$ is a special case of $I_n(R)$ (i.e.\ set $R = \Theta$).
We refer to $I_n(R)$ as a ``building block" in the sense that, together with closely related variants like $I_n(\Theta)$, it may be employed to obtain a variety of different tests.

\subsubsection{Strong Approximation}

We first obtain a strong approximation to statistics of the form $I_n(R)$. 
Before proceeding, we introduce some additional notation.
First, we define the class
\begin{equation}\label{not5}
\mathcal F_n \equiv \{\rho_{\jmath}(\cdot,\theta) : \theta \in \Theta_n \cap R \text{ and } 1 \leq \jmath \leq \mathcal J\} .
\end{equation}
The ``size" of $\mathcal F_n$ plays a crucial role, and we control it through the bracketing integral
\begin{equation*}
J_{[\hspace{0.02 in}]}(\delta,\mathcal F_n,\|\cdot\|_{P,2}) \equiv \int_0^\delta \sqrt{1 + \log N_{[\hspace{0.02 in}]}(\epsilon,\mathcal F_n,\|\cdot\|_{P,2})}d\epsilon ,
\end{equation*}
where $\|f\|_{P,2}^2 \equiv E_P[f^2(V)]$ and $N_{[\hspace{0.02 in}]}(\epsilon,\mathcal F_n,\|\cdot\|_{P,2})$ is the smallest number of $\epsilon$-brackets (under $\|\cdot\|_{P,2}$) required to cover $\mathcal F_n$.
Finally, we denote the empirical process by
$$\Gemp(\theta) \equiv \frac{1}{\sqrt n}\sum_{i=1}^n \{\rho(X_i,\theta)*q^{k_n}(Z_i) - E_P[\rho(X,\theta)*q^{k_n}(Z)]\}.$$
Our next assumptions imposes requirements on $\Theta_n \cap R$ and the transformation $q^{k_n}(Z)$.

\begin{assumption}\label{ass:startreg}
(i) $\max_{1\leq \jmath \leq \mathcal J} \max_{1\leq k \leq k_{n,\jmath}} \|q_{k,\jmath}\|_{\infty} \leq B_n$ with $B_n \geq 1$;
(ii) The eigenvalues of $E_P[q^{k_{n,\jmath}}_{\jmath}(Z_{\jmath})q^{k_{n,\jmath}}_{\jmath}(Z_{\jmath})^\prime]$ are bounded uniformly in $k_{n,\jmath}$ and $P\in\mathbf P$;
(iii) $\mathcal F_n$ has envelope $F_n$, $\sup_{P\in \mathbf P} \|F_n\|_{P,2} <\infty$, and $\sup_{P\in \mathbf P}J_{[\hspace{0.02 in}]}(\|F_n\|_{P,2},\mathcal F_n,\|\cdot\|_{P,2}) \leq J_n$ with $J_n <\infty$.
\end{assumption}

\begin{assumption}\label{ass:coupling}
(i) $\sup_{\theta \in \Theta_n \cap R}\|\Gemp(\theta) - \WP(\theta)\|_p = o_P(a_n)$ uniformly in $P\in \mathbf P$ for some $a_n = o(1)$ and Gaussian $\WP$ satisfying $E[\WP(\theta)] = 0$ and ${\rm Cov}\{\WP(\theta),\WP(\theta^\prime)\} = {\rm Cov}_P\{\Gemp(\theta),\Gemp(\theta^\prime)\}$;
(ii) There is a norm $\|\cdot\|_{\mathbf E}$, $\kappa_\rho > 0$, and $K_\rho< \infty$ such that $E_P[\|\rho(X,\theta_1)-\rho(X,\theta_2)\|_2^2] \leq K_\rho^2\|\theta_1 - \theta_2\|_{\mathbf E}^{2\kappa_\rho}$ for all $\theta_1,\theta_2 \in \Theta_n \cap R$ and $P\in \mathbf P$.
\end{assumption}

Assumptions \ref{ass:startreg}(i)(ii) impose standard requirements on the transformations $q^{k_n}$ -- e.g., Assumption \ref{ass:startreg}(i) holds with $B_n = 1$ for trigonometric series and $B_n \asymp \sqrt{k_n}$ for normalized $B$-splines.
Assumption \ref{ass:startreg}(iii) controls the ``size" of $\mathcal F_n$.
We allow $J_n$ to depend on $n$ to accommodate non-compact parameter spaces \citep{chen:pouzo:2012, chen:pouzo:2015}.
Assumption \ref{ass:coupling}(i) requires that the empirical process be approximately Gaussian.
The sequence $\{a_n\}_{n=1}^\infty$ denotes a bound on the rate of coupling, which in turn characterizes the rate of convergence of our strong approximation.
In the appendix, we verify Assumption \ref{ass:coupling}(i) by relying on existing results \citep{yurinskii:1977, zhai2018high} or a novel extension of \cite{koltchinskii:1994}.
Assumption \ref{ass:coupling}(ii) imposes a mild restriction on the moment functions that ensures $\WP$ is equicontinuous with respect to $\|\cdot\|_{\mathbf E}$.

In establishing our strong approximation to $I_n(R)$, it is helpful to derive the rate of convergence of the minimizer of $Q_n$ over $\Theta_n \cap R$.
To this end, we follow the literature on set estimation \citep{chernozhukov:hong:tamer:2007, beresteanu:molinari:2008, santos:2011, Kaido_Santos2013} and for any sets $A$ and $B$ we define
$$\overrightarrow{d}_{H}(A,B, \|\cdot\|_{\mathbf E})  \equiv \sup_{a\in A}\inf_{b\in B} \|a - b\|_{\mathbf E}, $$
which is known as the directed Hausdorff distance.
For each $\theta \in \Theta \cap R$, we further let $\Pi_n \theta$ denote its approximation on $\Theta_n \cap R$ and denote the approximation to $\IDset \cap R$ by
\begin{equation}\label{def:thetaset}
\IDsetRsieve \equiv \{\Pi_n\theta : \theta \in \IDset\cap R\}.
\end{equation}
Our next assumption enables us to obtain a rate of convergence (under $\|\cdot\|_{\mathbf E}$) to $\IDsetRsieve$.

\begin{assumption}\label{ass:keycons}
There are $\mathcal V_{n}(P)\subseteq \Theta_n \cap R$ and a sequence constants $\{\nu_{n}\}$ with $0< \nu^{-1}_{n} = O(1)$ such that (i) For any $\theta \in \mathcal V_n(P)$ it holds that
$$\nu_{n}^{-1}  \overrightarrow d_H(\theta, \IDsetRsieve,\|\cdot\|_{\mathbf E})
\leq \sup_{\tilde \theta \in \IDsetRsieve} \|E_P[(\rho(X,\theta)-\rho(X,\tilde \theta))*q^{k_n}(Z)]\|_{\PSigma,p};$$
(ii) There is a $\hat \theta_n \in \mathcal V_n(P)$ satisfying $Q_n(\hat \theta_n) \leq \inf_{\theta \in \Theta_n \cap R} Q_n(\theta) + o(a_n/\sqrt n)$ with probability tending to one uniformly in $P\in \mathbf P_0$.
\end{assumption}

Assumption \ref{ass:keycons}(ii) requires that an approximate minimum of $Q_n$ over $\Theta_n \cap R$ be attained at a point $\hat \theta_n$ in a set $\mathcal V_n(P)$ with high probability.
Typically, $\mathcal V_n(P)$ may be taken to equal the entire sieve in convex models, or it may be taken to equal a local neighborhood of $\IDsetRsieve$ after establishing the consistency of $\hat \theta_n$ through standard arguments; see, e.g., Lemma \ref{lm:setcons} in the appendix.
Assumption \ref{ass:keycons}(i) introduces a local identification condition on $\mathcal V_n(P)$ by requiring that the moments ``change" at a rate $\nu_n^{-1}$ as $\theta$ moves away from $\IDsetRsieve$.
The parameter $\nu_n^{-1}$, which implicitly depends on $k_n$ and the choice of sieve $\Theta_n\cap R$, is conceptually related to sieve measure of ill-posedness \citep{blundell:chen:kristensen:2007}.

By employing Assumption \ref{ass:keycons}, we are able to show that with arbitrarily high probability, $\hat \theta_n$ is contained in a $\|\cdot\|_{\mathbf E}$-neighborhood of $\IDsetRsieve$ that shrinks at a rate
\begin{equation}\label{def:Rn}
\mathcal R_n \equiv \nu_n\{\frac{k_n^{1/p}\sqrt{\log(1 + k_n)} J_n B_n}{\sqrt n}\},
\end{equation}
where recall $B_n$ and $J_n$ where introduced in Assumption \ref{ass:startreg}.
Under assumptions on the (Hausdorff) distance between $\IDsetRsieve$ and $\IDset \cap R$, the triangle inequality can yield a rate of convergence of $\hat \theta_n$ to $\IDset \cap R$.
Heuristically, we focus on convergence to $\IDsetRsieve$ (instead of $\IDset\cap R$) because our strong approximation will rely on undersmoothing.

In our final assumptions, we follow the literature and accommodate non-differentiable moment functions by requiring that their conditional expectations be differentiable \citep{chen:pouzo:2009,chen:pouzo:2012}.
Specifically, for each $1\leq \jmath\leq \mathcal J$ and $\theta\in \Theta$ we set
$$m_{P,\jmath}(\theta)(Z_\jmath) \equiv E_P[\rho_\jmath(X,\theta)|Z_\jmath];$$
i.e.\ $m_{P,\jmath}$ maps each $\theta \in \Theta$ to a square integrable function of $Z_\jmath$.
Letting $\mathbf B_n$ denote the vector subspace generated by $\Theta_n\cap R$, we then impose the following:

\begin{assumption}\label{ass:driftlin}
There is a norm $\|\cdot\|_{\mathbf L}$ on $\mathbf B_n$, linear maps $\nabla m_{P,\jmath}(\theta) : \mathbf B \rightarrow L^2_P$, and constants $\epsilon > 0$ and $K_m, M < \infty$ such that for all $P\in \mathbf P$, $h\in \mathbf B_n$, and elements $\theta_1,\theta_2 \in \{\theta \in \Theta_n \cap R : \overrightarrow{d}_H(\theta,\IDsetRsieve, \|\cdot\|_{\mathbf E}) \leq \epsilon\}$ we have: (i) $\|m_{P,\jmath}(\theta_1) - m_{P,\jmath}(\theta_2) - \nabla m_{P,\jmath}(\theta_2)[\theta_1-\theta_2]\|_{P,2} \leq K_m\|\theta_1-\theta_2\|_{\mathbf L}\|\theta_1-\theta_2\|_{\mathbf E}$; (ii) $\|\nabla m_{P,\jmath}(\theta_1)[h] - \nabla m_{P,\jmath}(\theta_2)[h]\|_{P,2} \leq K_m\|\theta_1 - \theta_2\|_{\mathbf L}\|h\|_{\mathbf E}$;  
(iii) $\|\nabla m_{P,\jmath}(\theta_2)[h]\|_{P,2} \leq M \|h\|_{\mathbf E}$.
\end{assumption}

\begin{assumption}\label{ass:locrates}
(i) $k_n^{1/p}\sqrt{\log(1+k_n)}B_n \sup_{P\in \mathbf P} J_{[\hspace{0.02 in}]}({\mathcal R}_n^{\kappa_\rho},\mathcal F_n,\|\cdot\|_{P,2}) =  o(a_n)$;
(ii) $\sup_{P\in \mathbf P_0} \sup_{\theta \in \IDsetRsieve} \sqrt n\|E_P[\rho(X,\theta)*q^{k_n}(Z)]\|_{\PSigma,p} = o(a_n)$.
\end{assumption}

\begin{assumption}\label{ass:weights}
(i) For each $P\in \mathbf P$ there is a $k_n\times k_n$ matrix $\PSigma > 0$ such that $\|\hat \Sigma_n - \PSigma\|_{o,p} = o_P(1 \wedge a_n \{k_n^{1/p}\sqrt{\log(1+k_n)}B_nJ_n\}^{-1})$ uniformly in $P\in \mathbf P$;
(ii) $\|\PSigma\|_{o,p}$ and $\|\PSigma^{-1}\|_{o,p}$ are uniformly bounded in $k_n$ and $P\in \mathbf P$.
\end{assumption}

Assumption \ref{ass:driftlin}(i) ensures $m_{P,\jmath}$ is approximated by linear maps $\nabla m_{P,\jmath}$ with an approximation error that is controlled by $\|\cdot\|_{\mathbf E}$ and a potentially stronger norm $\|\cdot\|_{\mathbf L}$.
In turn, Assumptions \ref{ass:driftlin}(ii)(iii) impose continuity conditions on $\nabla m_{P,\jmath}$ -- these assumptions are not used in this section, but will be needed for our bootstrap results.
Assumption \ref{ass:locrates} contains our key rate restrictions. 
Assumption \ref{ass:locrates}(i) ensures the rate of convergence $\mathcal R_n$ (as in \eqref{def:Rn}) is sufficiently fast to overcome an asymptotic loss of equicontinuity -- a requirement that can hold even when $\mathcal R_n$ is slower than the traditional $o(n^{-1/4})$ rate employed to linearize nonlinear models.
Assumption \ref{ass:locrates}(ii) is an undersmoothing assumption, which ensures that $I_n(R)$ is properly centered under the null hypothesis.
Finally, Assumption \ref{ass:weights} requires $\hat \Sigma_n$ to converge to an invertible matrix $\PSigma$ at a suitable rate -- here, $\|\cdot\|_{o,p}$ denotes the operator norm when $\mathbf R^{k_n}$ is endowed with $\|\cdot\|_p$.

The introduced assumptions suffice for obtaining a strong approximation through a local reparametrization.
Formally, we denote the local deviations from $\theta\in \Theta_n \cap R$ by
\begin{equation*}
 V_n(\theta,R|\ell) \equiv \{h\in \mathbf B_n : \theta + \frac{h}{\sqrt n} \in \Theta_n \cap R \text{ and } \|\frac{h}{\sqrt n}\|_{\mathbf E} \leq \ell\}.
\end{equation*}
Recall $\mathbf B_n$ denotes the vector subspace generated by $\Theta_n \cap R$ and for any $h \in \mathbf B_n$ set
$$\DerP(\theta)[h] \equiv E_P[\nabla m_P(\theta)[h](Z)*q^{k_n}(Z)],$$
where $\nabla m_P(\theta)[h](Z) \equiv (\nabla m_{P,1}(\theta)[h](Z_1),\ldots, \nabla m_{P,\mathcal J}(\theta)[h](Z_{\mathcal J}))^\prime$.
For any given sequence $\ell_n$, we then define a sequence of random variables $\Up(R|\ell_n)$ to be given by
\begin{equation}\label{eq:hatudef}
\Up(R|\ell_n) \equiv \inf_{\theta \in \IDsetRsieve} \inf_{h \in V_{n}(\theta,R|\ell_n)} \|\WP(\theta) + \DerP(\theta)[h]\|_{\PSigma,p}.
\end{equation}
As a final piece of notation, for any two norms $\|\cdot\|_{\mathbf A_1}$ and $\|\cdot\|_{\mathbf A_2}$ defined on $\mathbf B_n$, we set
\begin{equation*}
\mathcal S_n(\mathbf A_1,\mathbf A_2) \equiv \sup_{b\in \mathbf B_n} \frac{\|b\|_{\mathbf A_1}}{\|b\|_{\mathbf A_2}},
\end{equation*}
which we note depends on the sample size $n$ only through the choice of sieve $\Theta_n \cap R$.

The next result establishes the relation between $\Up(R|\ell_n)$ and $I_n(R)$.
It is helpful to recall here that the norm $\|\cdot\|_{\mathbf L}$ and constants $K_m$, introduced in Assumption \ref{ass:driftlin}, control the linearization of the moments and that $K_m = 0$ for linear models.

\begin{theorem}\label{th:localdrift}
Let Assumptions \ref{ass:param}(i), \ref{ass:startreg}, \ref{ass:coupling}, \ref{ass:keycons}, \ref{ass:driftlin}(i), \ref{ass:locrates}, and \ref{ass:weights} hold. Then:
(i) For any $\ell_n \downarrow 0$ satisfying $k_n^{1/p}\sqrt{\log(1+k_n)}B_n\times \sup_{P\in \mathbf P}J_{[\hspace{0.02 in}]}(\ell_n^{\kappa_\rho},\mathcal F_n,\|\cdot\|_{P,2}) =  o(a_n)$ and $K_m\ell_n^2\times  \mathcal S_n(\mathbf L,\mathbf E) = o(a_nn^{-1/2})$ it follows uniformly in $P\in \mathbf P_0$ that:
$$I_{n}(R) \leq \Up(R|\ell_n) + o_P(a_n) . $$
(ii) If in addition $K_m\mathcal R_n^2\times\mathcal S_n(\mathbf L,\mathbf E) =  o(a_nn^{-1/2})$, then $\ell_n$ may be additionally chosen to satisfy $\mathcal R_n = o(\ell_n)$, in which case it follows uniformly in $P\in \mathbf P_0$ that:
$$I_{n}(R) = \Up(R|\ell_n) +  o_P(a_n) . $$
\end{theorem}

Theorem \ref{th:localdrift} is perhaps best understood as establishing the validity of a family (indexed by $\{\ell_n\}$) of strong approximations that differ on the size of the local neighborhoods of $\IDsetRsieve$ that they employ.
Its proof crucially relies on the linearization
\begin{equation}\label{eq:local1}
\DerP(\theta)[h] \approx \sqrt n \{E_P[\rho(X,\theta+\frac{h}{\sqrt n})*q^{k_n}(Z)] -  E_P[\rho(X,\theta)*q^{k_n}(Z)]\},
\end{equation}
which holds for nonlinear moments ($K_m\neq 0$) when $h/\sqrt n$ is sufficiently small.
In particular, if the infimum defining $I_n(R)$ is attained at a point $\hat\theta_n$ that converges to $\IDsetRsieve$ sufficiently fast, then we may apply \eqref{eq:local1} to establish Theorem \ref{th:localdrift}(ii).
Regrettably, in certain models the rate of convergence of $\hat \theta_n$ may be too slow to apply the approximation in \eqref{eq:local1} to $\hat \theta_n$.
In such instances, we may instead rely on the inequality
\begin{equation}\label{eq:local2}
I_n(R) = \inf_{\theta \in \Theta_n \cap R} \sqrt nQ_n(\theta) \leq \inf_{(\theta,h) \in (\IDsetRsieve,V_n(\theta,R|\ell_n))} \sqrt n Q_n(\theta + \frac{h}{\sqrt n})
\end{equation}
and successfully couple the right hand side of \eqref{eq:local2} by restricting attention to sequences $\ell_n$ for which \eqref{eq:local1} is accurate.
Thus, by regularizing the local parameter space through a norm bound, we obtain in Theorem \ref{th:localdrift}(i) a distributional approximation that, while potentially conservative, holds under weaker requirements on the rate of convergence.

\subsubsection{Bootstrap Approximation} \label{sec:boot}

Theorem \ref{th:localdrift} shows that the distribution of $U_P(R|\ell_n)$ is a suitable approximation for the distribution of $I_n(R)$. 
We next develop a bootstrap procedure for estimating the distribution of $U_P(R|\ell_n)$ with the goal of obtaining valid critical values.

We estimate the distribution of $\Up(R|\ell_n)$ by replacing population parameters with suitable sample analogues.
The key ingredients are: (i) A random variable $\Bemp$ whose distribution conditional on the data is consistent for the distribution of $\WP$; (ii) An estimator $\hat {\mathbb D}_n(\theta)$ for $\DerP(\theta)$; (iii) An estimator $\hat \Theta_n^{\text{r}}$ for $\IDsetRsieve$ (as in \eqref{def:thetaset}); and (iv) A sample analogue $\hat V_n(\theta,R| \ell_n)$ for the local parameter space $V_{n}(\theta,R|\ell_n)$.
We then approximate the distribution of $\Up(R|\ell_n)$ by the distribution (conditional on the data) of 
\begin{equation*}
\hat U_n(R|\ell_n) \equiv \inf_{\theta \in \hat \Theta_n^{\rm r}} \inf_{h \in \hat V_n(\theta,R|\ell_n)} \|\Bemp(\theta) + \hat {\mathbb D}_n(\theta)[h]\|_{\hat \Sigma_n,p}.
\end{equation*}

For concreteness, we employ the following sample analogues in our construction.

\noindent {\bf Gaussian Distribution:} We estimate the distribution of $\WP$ with the multiplier bootstrap.
Specifically, for i.i.d.\ $\{\omega_i\}_{i=1}^n$ with $\omega_i \sim N(0,1)$ independent of $\{V_i\}_{i=1}^n$ we let
\begin{equation*}
\Bemp(\theta) \equiv \frac{1}{\sqrt n} \sum_{i=1}^n \omega_i\{\rho(X_i,\theta)*q^{k_n}(Z_i) - \frac{1}{n}\sum_{j=1}^n \rho(X_j,\theta)*q^{k_n}(Z_j)\}.
\end{equation*}
We focus on the multiplier bootstrap due to its theoretical tractability, though we note that alternative bootstrap approaches can also be valid. \qed

\noindent {\bf The Derivative:} We estimate $\DerP(\theta)$ by employing a construction that is applicable to non-differentiable moments.
Specifically, for any $\theta \in \Theta_n \cap R$ and $h\in \mathbf B_n$ we set
\begin{equation*}
\hat{\mathbb D}_n(\theta)[h] \equiv \frac{1}{\sqrt n} \sum_{i=1}^n (\rho(X_i,\theta + \frac{h}{\sqrt n}) - \rho(X_i,\theta))*q^{k_n}(Z_i).
\end{equation*}
We employ $\hat {\mathbb D}_n(\theta)$ due to its general applicability, though alternative approaches may be preferable in some applications. 
In particular, if moments are differentiabile, then using
\begin{equation*}
 \frac{1}{n}\sum_{i=1}^n  \nabla_\theta \rho(X_i, \theta)[h]* q^{k_n}(Z_i) 
\end{equation*}
as an estimator for $\DerP(\theta)[h]$ leads to a computationally simpler bootstrap statistic. \qed

\noindent {\bf The Identified Set:} We estimate the identified set by employing the set of (approximate) minimizers of $Q_n$ on $\Theta_n \cap R$. 
Formally, for a sequence $\tau_n \downarrow 0$, we let
\begin{equation}\label{def:setest}
\hat \Theta_n^{\text{r}} \equiv \{\theta \in \Theta_n \cap R: Q_n(\theta) \leq \inf_{\theta \in \Theta_n \cap R} Q_n(\theta) + \tau_n\}.
\end{equation}
We may set $\tau_n = 0$ in identified models, in which case $\hat \Theta_n^{\rm r}$ reduces to the minimizer of $Q_n$. In partially identified models, $\hat \Theta_n^{\rm r}$ can be shown to asymptotically lie in a shrinking neighborhood of $\IDsetRsieve$ provided $\tau_n \to 0$.
In order for $\hat \Theta_n^{\rm r}$ to additionally be Hausdorff consistent for $\IDsetRsieve$, however, $\tau_n$ must not tend to zero too fast; see Lemma \ref{lm:setcons}. \qed

\noindent {\bf Local Parameter Space:} We account for the role inequality constraints play in determining the local parameter space by estimating ``binding" sets in analogy to approaches pursued in the moment inequalities literature \citep{chernozhukov:hong:tamer:2007, andrews:soares:2010}.
Specifically, for a sequence $r_n$ and any $\theta \in \Theta_n \cap R$ we define
\begin{equation*}
G_n(\theta) \equiv  \{h \in \mathbf B_n : \Upsilon_G(\theta + \frac{h}{\sqrt n}) \leq (\Upsilon_G(\theta) - K_gr_n\|\frac{h}{\sqrt n}\|_{\mathbf B}\mathbf {1_G}) \vee (-r_n\mathbf {1_G}) \} ,
\end{equation*}
where recall $\mathbf {1_G}$ is the order unit in $\mathbf G$ and $g_1 \vee g_2$ represents the supremum of any $g_1,g_2\in \mathbf G$.
The constant $K_g$, formally introduced in Assumption \ref{ass:locineq} below, is related to the curvature of $\Upsilon_G$ and equals zero for linear $\Upsilon_G$.
For any $\ell_n$ we then define
\begin{equation}\label{eq:Vndef}
\hat V_n(\theta,R| \ell_n) \equiv\{h \in \mathbf B_n : h \in G_n(\theta), ~ \Upsilon_F(\theta + \frac{h}{\sqrt n}) = 0 \text{ and } \|\frac{h}{\sqrt n}\|_{\mathbf B} \leq \ell_n\},
\end{equation}
i.e.\ in comparison to $V_n(\theta,R|\ell_n)$ we: (i) Replace $\Upsilon_G(\theta + h/\sqrt n) \leq 0$ by $h \in G_n(\theta)$; (ii) Retain $\Upsilon_F(\theta + h/\sqrt n) = 0$; and (iii) Substitute $\|h/\sqrt n \|_{\mathbf E} \leq \ell_n$ with $\|h/\sqrt n\|_{\mathbf B} \leq \ell_n$. \qed

Before establishing the asymptotic validity of the proposed bootstrap procedure, we require some additional notation. 
For any set $A \subseteq \mathbf B_n$, we denote its $\epsilon$-neighborhood by
$$(A)^\epsilon \equiv \{\theta \in \mathbf B_n : \inf_{a \in A} \|a - \theta\|_{\mathbf B} \leq \epsilon\}.$$
We further denote the closure of the linear span of $\Upsilon_F(\mathbf B_n)$ by $\mathbf F_n$, and for any linear map $\Gamma$ on $\mathbf B$ we let $\mathcal N(\Gamma) \equiv \{h \in \mathbf B: \Gamma(h) = 0\}$ denote its null space.
In the assumptions that follow, it is helpful to recall that $\Theta_{0n}^{\rm r}$ is implicitly a function of $P$.

\begin{assumption}\label{ass:locineq}
For some $K_g,M < \infty$, $\epsilon > 0$ and all $n$, $P\in \mathbf P_0$, $\theta_1,\theta_2 \in (\IDsetRsieve)^\epsilon$ (i) $\Upsilon_G$ is Fr\'echet differentiable with $\|\Upsilon_G(\theta_1) - \Upsilon_G(\theta_2) - \nabla\Upsilon_G(\theta_1)[\theta_1-\theta_2]\|_{\mathbf G} \leq K_g \|\theta_1 - \theta_2\|_{\mathbf B}^2$; (ii) $\|\nabla \Upsilon_G(\theta_1) - \nabla \Upsilon_G(\theta_2)\|_o \leq K_g\|\theta_1 - \theta_2\|_{\mathbf B}$; (iii) $\|\nabla \Upsilon_G(\theta_1)\|_o \leq M$.
\end{assumption}

\begin{assumption}\label{ass:loceq}
For some $K_f,M < \infty$, $\epsilon > 0$ and all $n$, $P\in \mathbf P_0$, $\theta_1,\theta_2 \in  (\IDsetRsieve)^\epsilon$ (i) $\Upsilon_F$ is Fr\'echet differentiable with $\|\Upsilon_F(\theta_1) - \Upsilon_F(\theta_2) - \nabla\Upsilon_F(\theta_1)[\theta_1-\theta_2]\|_{\mathbf F} \leq K_f \|\theta_1 - \theta_2\|_{\mathbf B}^2$; (ii) $\|\nabla \Upsilon_F(\theta_1) - \nabla \Upsilon_F(\theta_2)\|_o \leq K_f\|\theta_1 - \theta_2\|_{\mathbf B}$; (iii) $\|\nabla \Upsilon_F(\theta_1)\|_o \leq M$; (iv) $\nabla \Upsilon_F(\theta_1):\mathbf B_n \rightarrow \mathbf F_n$ admits a right inverse $\nabla\Upsilon_F(\theta_1)^{-}$ with $K_f\|\nabla \Upsilon_F(\theta_1)^{-}\|_o \leq M$.
\end{assumption}

\begin{assumption}\label{ass:ineqlindep}
Either (i) $\Upsilon_F : \mathbf B \rightarrow \mathbf F$ is affine, or (ii) There are constants $\epsilon > 0$, $M < \infty$ such that for every $P\in \mathbf P_0$, $n$, and $\theta \in \IDsetRsieve$ there exists a $h \in \mathbf B_n \cap \mathcal N(\nabla \Upsilon_F(\theta))$ satisfying $\Upsilon_G(\theta) + \nabla \Upsilon_G(\theta)[h] \leq -\epsilon \mathbf {1_G}$ and $\|h\|_{\mathbf B} \leq M$.
\end{assumption}

Assumption \ref{ass:locineq} imposes that $\Upsilon_G$ be Fr\'echet differentiable. 
The constant $K_g$, employed in the construction of $\hat V_n(\theta,R|\ell_n)$, may be interpreted as a bound on the second derivative of $\Upsilon_G$ and equals zero when $\Upsilon_G$ is linear.
Assumptions \ref{ass:loceq} and \ref{ass:ineqlindep} mark an important difference between hypotheses in which $\Upsilon_F$ is linear and those in which $\Upsilon_F$ is nonlinear -- note linear $\Upsilon_F$ automatically satisfy Assumptions \ref{ass:loceq} and \ref{ass:ineqlindep}.
This distinction reflects that when $\Upsilon_F$ is linear its impact on the local parameter space is known and need not be estimated.\footnote{For linear $\Upsilon_F$, the requirement $\Upsilon_F(\theta + h/\sqrt n) = 0$ is equivalent to $\Upsilon_F(h) = 0$ for any $\theta \in R$.} 
Thus, while Assumptions \ref{ass:loceq}(i)-(iii) impose conditions analogous to those required of $\Upsilon_G$, Assumption \ref{ass:loceq}(iv) additionally demands that $\nabla \Upsilon_F(\theta)$ posses a norm bounded right inverse on ($\IDsetRsieve)^\epsilon$ -- the existence of a right inverse is equivalent to a classical rank condition.\footnote{Recall for a linear map $\Gamma : \mathbf B_n \rightarrow \mathbf F_n$, its right inverse is a map $\Gamma^{-} : \mathbf F_n \rightarrow \mathbf B_n$ such that $\Gamma \Gamma^{-}(h) = h$ for any $h \in \mathbf B_n$. The right inverse $\Gamma^{-}$ need not be unique if $\Gamma$ is not bijective, in which case Assumption \ref{ass:loceq}(iv) is satisfied as long as it holds for some right inverse of $\nabla \Upsilon_F(\theta):\mathbf B_n \rightarrow \mathbf F_n$.}
Finally, for nonlinear $\Upsilon_F$, Assumption \ref{ass:ineqlindep}(ii) requires the existence of a local perturbation to any $\theta\in \IDsetRsieve$ that relaxes  ``active" inequality constraints without a first order effect on the equality restrictions.


We impose a final set of assumptions in order to couple our bootstrap statistic.

\begin{assumption}\label{ass:bootcoupling}
$\sup_{\theta \in \Theta_n \cap R} \|\Bemp(\theta) - \WPT(\theta)\|_p =  o_P(a_n)$ uniformly in $\Phi\times P$ with $P\in \mathbf P$ for $\Phi$ the standard normal distribution, $a_n = o(1)$, and $\WPT$ independent of $\{V_i\}_{i=1}^n$ and having the same distribution as $\WP$.
\end{assumption}

\begin{assumption}\label{ass:extra}
(i) For some $M < \infty$, $\|h\|_{\mathbf E}\leq M\|h\|_{\mathbf B}$ for all $h\in \mathbf B_n$;
(ii) There is an $\epsilon > 0$ such that $P((\hat \Theta_n^{\text{\rm r}})^{\epsilon} \subseteq \Theta_n)$ tends to one uniformly in $P\in \mathbf P_0$;
(iii) For $\mathcal V_n(P)$ as in Assumption \ref{ass:keycons}, $P(\hat \Theta_n^{\rm r} \subseteq \mathcal V_n(P))$ tends to one uniformly in $P\in \mathbf P_0$.
\end{assumption}

\begin{assumption}\label{ass:bootrates}
(i) Either $\Upsilon_F$ and $\Upsilon_G$ are affine or $(\mathcal R_n + \nu_n\tau_n) \times \mathcal S_n(\mathbf B,\mathbf E) = o(1)$;
(ii) The sequences $\ell_n,\tau_n$ satisfy $k_n^{1/p}\sqrt{\log(1+k_n)}B_n\times \sup_{P\in \mathbf P}J_{[\hspace{0.03 in}]}(\ell_n^{\kappa_\rho}\vee(\nu_n\tau_n)^{\kappa_\rho},\mathcal F_n,\|\cdot\|_{P,2}) =  o(a_n)$, $K_m \ell_n( \ell_n + \mathcal R_n+ \nu_n\tau_n)\times \mathcal S_n(\mathbf L, \mathbf E) =  o(a_nn^{-1/2})$, and $\ell_n( \ell_n + \{\mathcal R_n+\nu_n\tau_n\} \times \mathcal S_n(\mathbf B,\mathbf E))1\{K_f > 0\} =  o(a_nn^{-1/2})$;
(iii) The sequence $r_n$ satisfies $\limsup_{n\rightarrow \infty}  1\{K_g>0\}\ell_n/r_n < 1/2$ and $(\mathcal R_n +\nu_n\tau_n) \times \mathcal S_n(\mathbf B, \mathbf E) = o(r_n)$.
\end{assumption}

Assumption \ref{ass:bootcoupling} demands that $\Bemp$ be coupled with a Gaussian $\WPT$ independent of $\{V_i\}_{i=1}^n$.
This condition implies the multiplier bootstrap is valid in our potentially non-Donsker setting; see Appendix \ref{sec:bootcoup} for sufficient conditions.
More generally, we note that our analysis remains valid if the multiplier bootstrap is replaced with any other re-sampling scheme (e.g., nonparametric bootstrap) satisfying a coupling requirement like Assumption \ref{ass:bootcoupling}.
Assumption \ref{ass:extra}(i) ensures that $\|\cdot\|_{\mathbf B}$ is (weakly) stronger than $\|\cdot\|_{\mathbf E}$.
Assumption \ref{ass:extra}(ii) demands that $\hat \Theta_n^{\rm r}$ be asymptotically contained in the interior of $\Theta_n$.
This requirement does not rule out that parameter space restrictions be binding at $\IDsetRsieve$ -- instead, Assumption \ref{ass:extra}(ii) requires that all such restrictions be stated through $R$.
Together with Assumption \ref{ass:keycons}(i), Assumption \ref{ass:extra}(iii) enables us to obtain a rate of convergence for $\hat \Theta_n^{\rm r}$ and may be verified in the same manner as Assumption \ref{ass:keycons}(ii).

Assumption \ref{ass:bootrates} contains our main rate requirements.
In particular, Assumption \ref{ass:bootrates}(i) ensures the one sided Hausdorff convergence of $\hat \Theta_n^{\rm r}$ to $\IDsetRsieve$ under $\|\cdot\|_{\mathbf B}$ whenever $\Upsilon_F$ or $\Upsilon_G$ are nonlinear. 
The main conditions on $\ell_n$, employed in constructing $\hat V_n(\theta,R|\ell_n)$, are contained in Assumption \ref{ass:bootrates}(ii).
These conditions ensure the consistency of $\hat {\mathbb D}_n(\theta)[h]$, the applicability of Theorem \ref{th:localdrift}, and that $\hat V_n(\theta,R|\ell_n)$ be well approximated by the true local parameter space. 
Heuristically, whenever the rate of convergence $\mathcal R_n$ is too slow, regularizing the local parameter space by selecting a small $\ell_n$ can ensure the asymptotic validity of the test.
As in Section \ref{sec:testex}, however, we note that whenever the rate of convergence $\mathcal R_n$ is sufficiently fast such regularization is unnecessary and it is possible to set $\ell_n = +\infty$ -- in such applications, setting $\ell_n$ to be too small can lead to a loss of power.
In turn, Assumption \ref{ass:bootrates}(iii) requires that $r_n$ not decrease to zero faster than the $\|\cdot\|_{\mathbf B}$-rate of convergence of $\hat \Theta_n^{\rm r}$. 
Assumption \ref{ass:bootrates}(iii) is always satisfied if $r_n = +\infty$, though setting $r_n \to 0$ can improve power against certain alternatives.
Similarly, we note that the requirements on $\tau_n$ imposed by Assumption \ref{ass:bootrates} can always be satisfied by setting $\tau_n = 0$, but such a choice can lead to a loss of power in certain partially identified models (recall the discussion in Section \ref{sec:testexpi}).

Our next result provides a coupling result for our bootstrap statistic.
In its statement, $\UpS(R|\ell_n)$ is defined identically to $\Up(R|\ell_n)$ but with $\WPT$ in place of $\WP$.

\begin{theorem}\label{th:coupsmooth}
If Assumptions \ref{ass:param}, \ref{ass:startreg}, \ref{ass:coupling}, \ref{ass:keycons}(i),  \ref{ass:driftlin}, \ref{ass:locrates}(ii), \ref{ass:weights}, \ref{ass:locineq}, \ref{ass:loceq}, \ref{ass:ineqlindep}, \ref{ass:bootcoupling}, \ref{ass:extra}, \ref{ass:bootrates} hold, then there is $\tilde \ell_n \asymp \ell_n$ so that uniformly in $P\in \mathbf P_0$
$$\hat U_n(R|\ell_n) \geq \UpS(R| \tilde \ell_n) +  o_P(a_n) .$$
\end{theorem}

Theorem \ref{th:coupsmooth} shows that with unconditional probability tending to one uniformly on $P\in \mathbf P_0$ our bootstrap statistic is bounded from below by a random variable that is independent of the data.
The significance of this result lies in that the lower bound is equal in distribution to the coupling to $I_n(R)$ obtained in Theorem \ref{th:localdrift}.
Thus, Theorems \ref{th:localdrift} and \ref{th:coupsmooth} provide the basis for constructing tests that employ increasing functions of $I_n(R)$ as a test statistic and the analogous bootstrap quantiles of $\hat U_n(R|\ell_n)$ as critical values.
The resulting tests may be conservative, however, whenever the inequalities in Theorems \ref{th:localdrift} and \ref{th:coupsmooth} are not ``sharp."
In particular, in order for the pointwise (in $P$) rejection probability to equal the nominal level of the test under the null hypothesis we require: (i) The rate of convergence $\mathcal R_n$ must be sufficiently fast for Theorem \ref{th:localdrift}(ii) to apply (in which case setting $\ell_n = +\infty$ is often valid); (ii) We should select $r_n$ to tend to zero with the sample size; and (iii) In partially identified settings, $\tau_n$ must tend to zero sufficiently slowly so that $\hat \Theta_n^{\rm r}$ is Hausdorff consistent for $\IDsetRsieve$.

\subsection{The Tests}

We next employ the theoretical results of Section \ref{sec:strongapp} to establish the asymptotic validity of different tests of the null hypothesis defined in \eqref{def:hyp}.
In what follows, for any statistic $\hat T_n$ that is a function of $\{V_i\}_{i=1}^n$ and the bootstrap weights $\{\omega_i\}_{i=1}^n$, we let $\hat q_{\tau}(\hat T_n)$ denote its conditional $\tau$ quantile given $\{V_i\}_{i=1}^n$.
For example, we have that
\begin{equation*}
\hat q_{1-\alpha}(\hat U_n(R|\ell_n)) = \inf\{ u : P(\hat U_n(R|\ell_n) \leq u |\{V_i\}_{i=1}^n) \geq 1- \alpha\}.
\end{equation*}

\subsubsection{Tests Based on \texorpdfstring{$I_n(R)$}{TEXT} }\label{subsec:critunc}

We first examine a test that employs $I_n(R)$ as a test statistic and a bootstrap quantile of $\hat U_n(R|\ell_n)$ as a critical value.
As has been shown in the literature, uniform consistent estimation of approximating distributions is not sufficient for characterizing the asymptotic size of a test \citep{romano:shaikh:2012}.
Heuristically, to establish the asymptotic validity of a test the approximating distributions must additionally be suitably uniformly continuous.
Our next assumption suffices for verifying this final requirement.

\begin{assumption}\label{ass:equic}
There is $\eta \geq 0$ and $\varrho_n = o(a_n^{-1})$ such that for $\hat c_n = \hat q_{1-\alpha}(\hat U_n(R|\ell_n))$ and any $\tilde \ell_n \asymp \ell_n$: (i) $P(I_n(R) > \hat c_n) = P( I_n(R) > \hat c_n\vee \eta) +o(1)$ uniformly in $P\in \mathbf P_0$, and (ii) $\sup_{P\in \mathbf P_0}\sup_{t\in(\eta - a_n,+\infty)} P(|\Up(R|\tilde \ell_n) - t|\leq \epsilon) \leq \varrho_n(\epsilon \wedge 1) + o(1)$.
\end{assumption}

Assumption \ref{ass:equic}(i) trivially holds with $\eta = 0$ since both $I_n(R)$ and $\hat U_n(R|\ell_n)$ are (weakly) positive almost surely.
However, in some applications it is possible to verify Assumption \ref{ass:equic}(i) in fact holds with $\eta > 0$ by arguing that the bootstrap quantiles of $\hat U_n(R|\ell_n)$ are suitably bounded away from zero when $I_n(R)$ is strictly positive.
Establishing Assumption \ref{ass:equic}(i) holds with $\eta > 0$ eases the verification of Assumption \ref{ass:equic}(ii), which intuitively requires that $\Up(R|\tilde \ell_n)$ be continuously distributed on $(\eta-a_n,+\infty)$ with a density bounded by a, possibly diverging, $\varrho_n$.
Because $U_P(R|\tilde \ell_n)$ is a functional of the Gaussian measure $\WP$, Assumption \ref{ass:equic}(ii) can in some applications be verified using available results in the literature \citep{davydov:lifshits:smorodina:1998}.
For instance, when $\Up(R|\tilde \ell_n)$ is a convex function of $\WP$, as in the application of Section \ref{sec:angristevans}, the distribution of $\Up(R|\tilde \ell_n)$ can readily be shown to be continuous in $(0,+\infty)$.
We refer the reader to \cite{chernozhukov2014comparison} for further discussion and motivation of conditions such as Assumption \ref{ass:equic}(ii), called \textit{anti-concentration} conditions.

The next result establishes the asymptotic validity of a test based on $I_n(R)$.

\begin{corollary}\label{th:critval}
Let Assumption \ref{ass:equic} hold and the conditions of Theorem \ref{th:localdrift}(i) and Theorem \ref{th:coupsmooth} be satisfied.
If $\hat c_n = \hat q_{1-\alpha}(\hat U_n(R|\ell_n))$, then it follows that:
$$\limsup_{n\rightarrow \infty} \sup_{P \in \mathbf P_0} P(I_n(R) > \hat c_n) \leq \alpha .$$
\end{corollary}

In Algorithm \ref{alg:cap} below we describe how to obtain p-values for the test described in Corollary \ref{th:critval} when the moments are differentiable.
We note that if there are no inequality constraints, then it is possible to show that the test in Corollary \ref{th:critval} is similar and its asymptotic size equals the nominal level $\alpha$ whenever the conditions of Theorem \ref{th:localdrift}(ii) are satisfied.
The consistency of the test against any $P\in \mathbf P\setminus \mathbf P_0$ for which $\max_{\jmath} \|E_P[\rho_\jmath(X,\theta)|Z_\jmath]\|_{P,2}$ is bounded away from zero (in $\theta \in \Theta\cap R$) is also straightforward to establish under suitable conditions. 
Finally, we also note that if we instead employ the critical value $\hat c_n = \hat q_{1-\alpha + \delta}(\hat U_n(R|\ell_n)) +\delta$ for any $\delta > 0$, then the conclusion of Corollary \ref{th:critval} holds without needing to impose Assumption \ref{ass:equic}; see Corollary \ref{cor:infactor}.
This modification to the critical value was originally proposed in a different context by \cite{andrews:shi:2013}, who suggest setting $\delta = 10^{-6}$.

\begin{algorithm}[hbt!]
\caption{Computing p-values for test based on $I_n(R)$}\label{alg:cap}
\begin{algorithmic}[1]
\Require $\Theta_n$, $\Upsilon_F$, $\Upsilon_G$, $\{\rho(X_i,\theta)*q^{k_n}(Z_i)\}_{i=1}^n$, $\hat \Sigma_n$, $r_n$, $\tau_n$, $\ell_n$

\item[]

\Statex \(\triangleright\) Compute the Test Statistic
\State $Q_n(\theta) \gets \| \hat \Sigma_n  \{\frac{1}{n}\sum_{i=1}^n \rho(X_i,\theta)*q^{k_n}(Z_i)\}\|_{p}$ \Comment{Criterion function}

\State $R\gets \{\theta : \Upsilon_F(\theta) = 0, \Upsilon_G(\theta) \leq 0\}$ \Comment{Constraint Set}

\State $I_n(R) \gets \min_{\theta \in \Theta_n} \sqrt nQ_n(\theta) \text{ s.t. } \theta \in R$ \Comment{ Test Statistic}

\item[]

\Statex \(\triangleright\) Prepare variables for bootstrap problem

\State $\hat{\mathbb D}_n(\theta)[h] \gets \frac{1}{n}\sum_{i=1}^n \nabla_\theta \rho(X_i,\theta)[h] *q^{k_n}(Z_i)$ \Comment{Moments Derivative}

\State $\hat \Theta_n^{\rm r}\gets \{\theta \in \Theta_n \cap R:  Q_n(\theta) \leq I_n(R)/\sqrt n + \tau_n \}$ \Comment{Boot Constraint $\theta$}

\State $G_n(\theta) \gets \{h : \Upsilon_G(\theta + h/\sqrt n) \leq (\Upsilon_G(\theta) - K_g r_n\|h/\sqrt n\|_{\mathbf B} \mathbf{1_G}) \vee (-r_n \mathbf {1_G})\}$

\State $\hat V_n(\theta,R|\ell_n) \gets \{h \in G_n(\theta) : \Upsilon_F(\theta + h/\sqrt n) = 0, \|h\|_{\mathbf B} \leq \ell_n \sqrt n\}$ \Comment{Boot Constraint $h$}

\item[]

\Statex \(\triangleright\) Compute $B$ bootstrap statistics and obtain p-value

\For{$b = 1 \text{ to } B$}

\State $\{\omega_i^b\}_{i=1}^n \gets $ \text{Generate i.i.d.\ sample of $N(0,1)$ variables}

\State $\hat {\mathbb W}_n^b(\theta) \gets \frac{1}{\sqrt n}\sum_{i=1}^n \omega_i^{b} \{\rho(X_i,\theta)*q^{k_n}(Z_i) - \frac{1}{n}\sum_{j=1}^n \rho(X_j,\theta)*q^{k_n}(Z_j)\}$

\State $F^b_n(\theta,h) \gets \|\hat \Sigma_n\{\hat {\mathbb W}_n^b(\theta) + \hat{\mathbb D}_n(\theta)[h]\}\|_p$ \Comment{Boot Criterion}

\State {\rm Boot}[b] $\gets \min_{\theta,h} F^b(\theta,h) \text{ s.t. } \theta \in \hat \Theta_n^{\rm r}, h\in \hat V_n(\theta,R|\ell_n)$ \Comment{Boot Statistic}

\EndFor

\State pval $\gets \frac{1}{B}\sum_{b=1}^B 1\{I_n(R) \leq {\rm Boot}[b]\}$ \Comment{Compute p-value}

\end{algorithmic}
\end{algorithm}

\begin{remark} \rm \label{rm:combine}
Suppose $\theta_0$ is identified, we aim to test whether $\Upsilon_F(\IDpoint) = 0$, and we are confident $\IDpoint$ satisfies $\Upsilon_G(\IDpoint) \leq 0$. 
We could then set $R$ to equal $R_1$ or $R_2$, where
\begin{align*}
R_1 & = \{\theta \in \mathbf B :\Upsilon_G(\theta) \leq 0 \text{ and } \Upsilon_F(\theta) = 0\}\\
R_2 & = \{\theta \in \mathbf B : \Upsilon_F(\theta) = 0\}.
\end{align*}
The power functions of the corresponding tests are not necessarily ranked.
As a result, it can be desirable to combine both tests by, for instance, using the test statistic
$T_n \equiv \max\{F_1(I_n(R_1)),F_2(I_n(R_2))\}$ for $F_1,F_2$ increasing functions, and the quantiles of $\max\{F_1(\hat U_n(R_1|\ell_n)),F_2(\hat U_n(R_2|\ell_n))\}$ as critical values -- e.g., $F_j$ may be c.d.f.\ of $\hat U_n(R_j|\ell_n)$ conditional on the data.
The asymptotic validity of such a test follows from Theorems \ref{th:localdrift} and \ref{th:coupsmooth} under a suitable modification of Assumption \ref{ass:equic}. \qed
\end{remark}

\subsubsection{Tests Based on \texorpdfstring{$I_n(R)-I_n(\Theta)$}{TEXT}}\label{subsec:critrec}

We next establish the asymptotic validity of a test based on $I_n(R) - I_n(\Theta)$ by also relying on Theorems \ref{th:localdrift} and \ref{th:coupsmooth}.
In what follows, we signify parameters associated with setting $R = \Theta$ by a ``${\rm u}$" superscript -- e.g.\ $\mathcal F_n^{\rm u}$ is understood to be as in \eqref{not5} but with $R = \Theta$.

In order to obtain a distributional approximation to the recentered statistic, we may simply apply Theorem \ref{th:localdrift}(i) to $I_n(R)$ and Theorem \ref{th:localdrift}(ii) to $I_n(\Theta)$ to conclude that
\begin{equation}\label{eq:critrec1}
I_n(R) - I_n(\Theta) \leq \Up(R|\ell_n) - \Up(\Theta|\ell_n^{\rm u}) + o_P(a_n).
\end{equation}
Moreover, by Theorem \ref{th:coupsmooth} we may approximate the distribution of $\Up(R|\ell_n)$ by using $\hat U_n(R|\ell_n)$.
Similarly, to obtain a bootstrap approximation to $\Up(\Theta|+\infty)$, we define 
$$\hat \Theta_n^{\rm u}    \equiv \{\theta \in \Theta_n : Q_n(\theta) \leq \inf_{\theta \in \Theta_n} Q_n(\theta) + \tau_n^{\text{u}}\};  $$
i.e.\ $\hat \Theta_n^{\rm u}$ is simply the set estimator in \eqref{def:setest} applied with $\Theta = R$.
For $\mathbf B_n^{\rm u}$ the closed linear span of $\Theta_n$, we then approximate the law of $\Up(\Theta|\ell_n^{\rm u})$ by employing
$$\hat U_{n}(\Theta|+\infty)  \equiv \inf_{\theta \in \hat \Theta_n^{\text{u}}} \inf_{h \in \mathbf B_n^{\rm u}} \| \Bemp(\theta) + \hat {\mathbb D}_n(\theta)[h]\|_{\hat \Sigma_n,p};$$
i.e.\ the bootstrap approximation equals that of Theorem \ref{th:coupsmooth}, with the local parameter space being unconstrained due to the absence of equality or inequality restrictions.

The preceding discussion suggests that the quantiles of $\hat U_n(R|\ell_n) -\hat U_n(\Theta|+\infty)$ conditional on the data provide valid critical values for the recentered statistic.
Our next result formally establishes that the resulting test is indeed asymptotically valid.

\begin{corollary}\label{th:incJcritval}
Let the conditions of Theorems \ref{th:localdrift}(i) and \ref{th:coupsmooth} hold with $R$ as in \eqref{def:Rset}, the conditions of Theorems \ref{th:localdrift}(ii) and \ref{th:coupsmooth} hold with $R = \Theta$, and Assumption \ref{ass:equic} hold with $I_n(R) - I_n(\Theta)$, $\hat U_n(R|\ell_n)-\hat U_n(\Theta|+\infty)$, and $\Up(R|\tilde \ell_n) - \Up(\Theta|\tilde \ell_n^{\rm u})$ in place of $I_n(R)$, $\hat U_n(R|\ell_n)$, and $\Up(R|\tilde \ell_n)$ with $\tilde \ell_n^{\rm u}$ satisfying $\mathcal R_n^{\rm u} = o(\tilde \ell_n^{\rm u})$ and Assumption \ref{ass:bootrates}(ii) with $R = \Theta$.
If $\tau_n^{\rm u}\downarrow 0$ satisfies $J_n^{\text{\rm u}}B_n k_n^{1/p}\sqrt{\log(1+k_n)/n} = o(\tau_n^{\text{\rm u}})$ and $\nu_n^{\rm u}\tau_n^{\rm u} \times \mathcal S_n^{\rm u}(\mathbf B, \mathbf E) = o(1)$, then for $\hat c_n \equiv \hat q_{1-\alpha}(\hat U_n(R|\ell_n)-\hat U_n(\Theta|+\infty))$ it follows that
$$\limsup_{n\rightarrow \infty} \sup_{P \in \mathbf P_0} P(I_n(R) - I_n(\Theta) > \hat c_n) \leq \alpha .$$
\end{corollary}

It is worth emphasizing that in coupling $I_n(\Theta)$ we must rely on Theorem \ref{th:localdrift}(ii) instead of Theorem \ref{th:localdrift}(i) in order to ensure that \eqref{eq:critrec1} holds.
As a result, whenever moments are nonlinear, Corollary \ref{th:incJcritval} requires the rate of convergence of the unconstrained estimator to be sufficiently fast for Theorem \ref{th:localdrift}(ii) to apply.
Similarly, in coupling $\hat U_n(\Theta|+\infty)$ it is important that $\hat \Theta_n^{\rm u}$ be consistent in the Hausdorff metric.
Thus, while we may set $\tau_n^{\rm u} = 0$ in identified models, in partially identified models we require that $\tau_n^{\rm u}$ not tend to zero too fast; see Theorem \ref{app:th:setrates}.
Finally, we note that in identified models, it is possible to employ either $\Bemp(\hat \theta_n)$ or $\Bemp(\hat \theta_n^{\rm u})$ in constructing both $\hat U_n(R|\ell_n)$ and $\hat U_n(\Theta |+\infty)$ -- a change that results in an asymptotically equivalent coupling but ensures that the bootstrap statistic $\hat U_n(R|\ell_n)-\hat U_n(\Theta|+\infty)$ is (weakly) positive.


\section{Heterogeneity and Demand Analysis}\label{sec:hetero}

For our final example, we illustrate how to conduct inference in the heterogeneous demand model of  \cite{hausman2016individual} -- alternative models of demand under conditional moment restrictions include the analysis in \cite{hausman:newey:1995}, \cite{blundell:horowitz:parey:2012}, and \cite{chen2018optimal}.
Specifically, for $Y\in  [0,1]$ equal to the expenditure share on a commodity, $W\in {\mathbf{W}}$ a vector of prices, income, and covariates, and $\eta$ representing unobserved individual heterogeneity we suppose
\begin{equation} \label{ex:hetero1}
Y=g(W,\eta )
\end{equation}
where $g$ is a known function of $(W,\eta )$.
The unobserved heterogeneity $\eta$ can potentially be infinite dimensional.
For instance, \cite{hausman2016individual} set $\eta =\{\beta _{j}\}_{j=1}^{\infty }$ to be a random variable in the sequence space $\ell ^{2}\equiv \{\{a_{j}\}_{j=1}^{\infty }:\sum_{j}a_{j}^{2}<\infty \}$, and let
\begin{equation} \label{ex:hetero2}
g(W,\eta) = \sum_{j=1}^{\infty }\psi_{j}(W)\beta_{j},
\end{equation}
where $\{\psi_{j}\}_{j=1}^{\infty }$ is a known basis satisfying $\sum_{j=1}^{\infty }\psi _{j}^{2}(W) < \infty$ almost surely (in $W$).

If the covariates $W$ are independent of $\eta$, then for any $c\in {\mathbf{R}}$ it follows that
\begin{equation} \label{ex:hetero3}
P(Y\leq c|W) = P(g(W,\eta)\leq c|W)=\int 1\{g(W,\eta )\leq c\}\mu_0(d\eta )
\end{equation}
where $\mu_0$ denotes the unknown distribution of $\eta $.
Result \eqref{ex:hetero3} restricts the possible values of $\mu_0$ and hence the identified set for functionals of $\mu_0$, such as average exact consumer surplus or average share.
Specifically, for $\Psi (g,\eta)$ an object of interest for preferences denoted by $\eta $, such as equivalent variation, \cite{hausman2016individual} study functionals
\begin{equation}\label{ex:hetero4}
\int \Psi (g,\eta )\mu_0(d\eta ),
\end{equation}
which is the average across individuals.
By evaluating the set of values of \eqref{ex:hetero4} which can be generated by a distribution $\mu_0$ satisfying \eqref{ex:hetero3} at a grid $\{c_{\jmath }\}_{\jmath =1}^{{\mathcal{J}}}$,
\cite{hausman2016individual} provide estimates of the identified set for the functional of interest.
We further note bounds on the distribution of $\Psi(g,\eta)$ under $\mu_0$ can be obtained by replacing $\Psi (g,\eta)$  in \eqref{ex:hetero4} with an indicator that $\Psi (g,\eta )$ be less than or equal to some number.

In what follows, we apply our results to conduct inference on functionals as in \eqref{ex:hetero4}.
To this end, we let $F_{P}(c_{\jmath}|W)\equiv P(Y\leq c_{\jmath }|W)$ for a given grid $\{c_{\jmath }\}_{\jmath=1}^{{\mathcal{J}}}$. 
To define ${\mathbf{B}}$, we suppose $\eta \in \Omega $ for some known Hausdorff space $\Omega$, set ${\mathcal{B}}$ to be the Borel $\sigma $-algebra on $\Omega $, let $\mathcal M$ be the space of regular signed Borel measures on $\Omega$, and let $\|\cdot\|_{TV}$ denote the total variation norm. 
Assuming $F_{P}(c_\jmath|\cdot)\in C_{B}({\mathbf{W}})$ for $C_{B}({\mathbf{W}})$ the space of continuous and bounded functions on ${\mathbf{W}}$, we set $\mathbf{B} = (\bigotimes_{\jmath =1}^{{\mathcal{J}}}C_{B}({\mathbf{W}}))\times {\mathcal M }$, for any $(\{F(c_\jmath|\cdot)\}_{\jmath =1}^{{\mathcal{J}}},\mu)=\theta \in {\mathbf{B}}$ let $\Vert \theta \Vert _{{\mathbf{B}}}=\sum_{\jmath =1}^{{\mathcal{J}}}\Vert F(c_\jmath|\cdot)\Vert _{\infty }+\Vert \mu\Vert _{TV}$, and set
\begin{equation}\label{ex:hetero6}
\Theta = \{(\{F(c_\jmath|\cdot)\}_{\jmath =1}^{{\mathcal{J}}},\mu )=\theta \in {\mathbf{B}}:\max_{1\leq \jmath \leq {\mathcal{J}}}\Vert F(c_\jmath|\cdot)\Vert _{\infty }\leq 2\},
\end{equation}
where the ``2" norm bound is simply selected to ensure $\IDset$ is in the interior of $\Theta$.

Letting $X=(Y,W)$ and setting $Z_\jmath = W$ for every $1\leq \jmath \leq {\mathcal{J}}$ we then define
\begin{equation}\label{ex:hetero7}
\rho _{\jmath }(X,\theta )=1\{Y\leq c_{\jmath }\}-F(c_{\jmath }|W),
\end{equation}
which yields conditional moment restrictions that identify $F_{P}(c_{\jmath }|W)$ -- note, however, that $\mu_0$ is potentially partially identified.
For a grid $\{w_l\}_{l=1}^{\mathcal L}\subseteq \mathbf W$ we test whether a hypothesized value $\lambda$ belongs to the identified set for the functional in \eqref{ex:hetero4} by setting
\begin{multline}\label{ex:hetero7paux}
    R = \Big\{(\{F(c_\jmath|\cdot)\}_{\jmath =1}^{{\mathcal{J}}},\mu ): \mu(\Omega) = 1, ~ \mu(B) \geq 0 \text{ for all } B\in \mathcal B, ~ \int \Psi(g,\eta) \mu(d\eta) = \lambda, \\
    \text{ and } F(c_\jmath|w_l) = \int 1\{g(w_l,\eta)\leq c_\jmath\}\mu(d\eta) \text{ for all } 1\leq \jmath\leq \mathcal J, 1 \leq l \leq \mathcal L\Big\}.
\end{multline}
Thus, the null hypothesis that $\Theta_0 \cap R$ be nonempty corresponds to requiring that there exist a distribution $\mu$ for $\eta$ satisfying the restrictions in \eqref{ex:hetero3} at the points $(c_\jmath,w_l)$ and yielding a value for the functional in \eqref{ex:hetero4} of $\lambda$.
By conducting test inversion in $\lambda$ we can obtain a confidence region for the desired functional.
To map $R$ into the framework of Section \ref{sec:gentheory}, we set $\mathbf G = \ell^\infty(\mathcal B)$ for $\ell^\infty(\mathcal B)$ the set of bounded functions on $\mathcal B$ and for any $(\{F(c_{\jmath }|\cdot )\}_{\jmath =1}^{{\mathcal{J}}},\mu )=\theta \in {\mathbf{B}}$ let $\Upsilon _{G}:{\mathbf{B}}\rightarrow \ell ^{\infty }({\mathcal{B}})$ be given by
\begin{equation}\label{ex:hetero9p0}
\Upsilon _{G}(\theta )(B)=-\mu (B).
\end{equation}
Finally, we  set $\Upsilon _{F}:{\mathbf{B}}\rightarrow {\mathbf{R}}^{{\mathcal{J}}{\mathcal{L}}+2}$ to equal $\Upsilon _{F}(\theta )=(\Upsilon _{F}^{({\text{e}})}(\theta ),\Upsilon_{F}^{(\mu )}(\theta ),\Upsilon _{F}^{({\text{s}})}(\theta ))$, where
\begin{align}
\Upsilon _{F}^{({\text{e}})}(\theta )& =\{F(c_{\jmath }|w_{l})-\int 1\{g(w_{l},\eta )\leq c_{\jmath }\}\mu (d\eta )\}_{1\leq \jmath \leq {\mathcal{J}},1\leq l\leq {\mathcal{L}}}  \notag \\ 
\Upsilon _{F}^{(\mu )}(\theta )& =\mu (\Omega )-1 \notag \\
\Upsilon _{F}^{({\text{s}})}(\theta )& =\int \Psi (g,\eta )\mu (d\eta)-\lambda. \label{ex:hetero9p3}
\end{align}
Given these definitions, we may then map $R$ (as introduced in \eqref{ex:hetero7paux}) into the framework of Section \ref{sec:gentheory} by noting that $R = \{\theta \in \mathbf B : \Upsilon_F(\theta) = 0 \text{ and } \Upsilon_G(\theta) \leq 0 \}$.

As in \cite{hausman2016individual}, we can impose utility maximization by requiring that the support $\Omega$ consist only of $\eta$ such that $g(\cdot ,\eta )$ satisfies the Slutsky conditions.
One may sample from $\Omega$ by drawing randomly from sets of $\eta$ that satisfy Slutsky symmetry and only keeping those where the compensated price effects matrix is negative semidefinite on a grid.
This is the procedure followed in \cite{hausman2016individual} for two goods.
Importantly, we emphasize that because the utility maximization restrictions are imposed through $\Omega$, they do not affect the  basic structure of $\Upsilon_F$ and $\Upsilon_G$ -- i.e., $\Upsilon_F$ and $\Upsilon_G$ remain linear maps  satisfying Assumptions \ref{ass:locineq}-\ref{ass:ineqlindep}.
In this sense, as long as they are imposed through the support $\Omega$ of $\eta$, our procedure allows us to accommodate a wide array of shape restrictions on individual demand $g(\cdot,\eta)$.

Given a collection of orthogonal probability measures $\{\delta_{s}\}_{s=1}^{s_{n}}\subseteq {\mathcal M}$ we employ
\begin{equation*}
{\mathcal{M}}_{n}= \{\mu \in {\mathcal{M}}:\mu =\sum_{s=1}^{s_{n}}\alpha_{s}\delta_{s}{\text{ for some }}\{\alpha_{s}\}_{s=1}^{s_{n}}\in {\mathbf{R}}^{s_{n}}\}
\end{equation*}%
as a sieve for ${\mathcal{M}}$.
Employing orthogonal measures, such as distinct Dirac measures, is computationally attractive as it simplifies imposing the nonnegativity constraint on any $\mu \in \mathcal M_n$.
As a sieve for $\{F_{P}(c_{\jmath }|\cdot )\}_{\jmath=1}^{{\mathcal{J}}}$, we employ approximating functions $\{p_{j}\}_{j=1}^{j_{n}}$.
In particular, setting $p^{j_{n}}(w) = (p_{1}(w),\ldots ,p_{j_{n}}(w))^{\prime }$, we set as our sieve
\begin{equation*}
\Theta _{n}= \{(\{p^{j_{n}\prime }\beta _{\jmath }\}_{\jmath =1}^{{%
\mathcal{J}}},\mu ):\mu \in {\mathcal{M}}_{n}{\text{ and }}%
\max_{1\leq \jmath \leq {\mathcal{J}}}\Vert p^{j_{n}\prime }\beta
_{\jmath }\Vert _{\infty }\leq 2\}.  \label{ex:hetero10}
\end{equation*}
Similarly, for a sequence $\{q_{k}\}_{k=1}^{k_{n}}$ and $k_n \times k_n$ positive definite matrices $\{\hat \Sigma _{\jmath ,n}\}_{\jmath=1}^{{\mathcal{J}}}$, we set  $q^{k_{n}}(w)= (q_{1}(w),\ldots ,q_{k_{n}}(w))^{\prime}$ and for any $(\{F(c_{\jmath }|\cdot )\}_{\jmath =1}^{{\mathcal{J}}},\mu )=\theta$ define
\begin{equation}\label{ex:hetero11}
Q_{n}(\theta )=\{\sum_{\jmath =1}^{{\mathcal{J}}}\Vert \frac{1}{n}\sum_{i=1}^{n}(1\{Y_{i}\leq c_{\jmath }\}-F(c_{\jmath}|W_{i}))q^{k_{n}}(W_{i})\Vert _{\hat{\Sigma}_{\jmath ,n},2}^{2}\}^{1/2}.
\end{equation}
The statistics $I_n(R)$ and $I_n(\Theta)$ then equal the minimums of $\sqrt nQ_n$ over $\Theta_n \cap R$ and $\Theta_n$.

Our next set of assumptions enable us to couple $I_n(R)$ and $I_n(R)-I_n(\Theta )$.

\begin{assumption}\label{ass:heterosieve}
(i) $\{Y_i,W_i\}_{i=1}^n$ is i.i.d.\ with $(Y,W)\sim P\in {\mathbf{P}}$;
(ii) $\sup _{w}\|p^{j_n}(w)\|_2\lesssim \sqrt{j_n}$;
(iii) $E_P[p^{j_n}(W)p^{j_n}(W)^{\prime }]$ has eigenvalues bounded away from zero and infinity uniformly in $P\in {\mathbf{P}}$ and $j_n$;
(iv) For each $P\in {\mathbf{P}}_0$ and $\theta \in \IDset\cap R$, there exists a $\Pi_n\theta =(\{F_n(c_{\jmath }|\cdot )\}_{\jmath=1}^{{\mathcal{J}}},\mu _n)\in \Theta _n\cap R$ such that $\sum _{\jmath=1}^{{\mathcal{J}}}\|E_P[(F_n(c_{\jmath }|W)-F_P(c_{\jmath}|W))q^{k_n}(W)]\|_2=O((n\log (n))^{-1/2})$ uniformly in $P\in {\mathbf{P}}_0$ and $\theta \in \IDset\cap R$.
\end{assumption}

\begin{assumption}\label{ass:heteromoments}
(i) $\max _{1\leq k\leq k_n}\|q_{k}\|_{\infty}\lesssim \sqrt{k_n}$;
(ii) $E_P[q^{k_n}(W)q^{k_n}(W)^{\prime }]$ has eigenvalues bounded uniformly in $P\in {\mathbf{P}}$ and $k_n$;
(iii) $E_P[q^{k_n}(W)p^{j_n}(W)^{\prime }]$ has singular values bounded away from zero uniformly in $P\in {\mathbf{P}}$ and $(k_n,j_n)$; (iv) $k_n^2j_n\log ^3(n)=o(n^{1/2})$.
\end{assumption}

\begin{assumption}\label{ass:heterosigma}
For all $1\leq \jmath \leq {\mathcal{J}}$: (i) $\|\hat {\Sigma }_{\jmath ,n}-\Sigma _{\jmath,P }\|_{o,2}=o_P(1/k_n\sqrt{j_n}\log ^2(n))$ uniformly in $P\in {\mathbf{P}}$;
(ii) The $k_n \times k_n$ matrices $\Sigma _{\jmath,P}$ are invertible and $\|\Sigma _{\jmath,P}\|_{o,2}$ and $\|\Sigma _{\jmath,P}^{-1}\|_{o,2}$ are bounded uniformly in $P\in {\mathbf{P}}$ and $k_n$.
\end{assumption}

Assumptions \ref{ass:heterosieve}(ii)-(iv) state the conditions on $\Theta_n$, with Assumptions \ref{ass:heterosieve}(ii)(iii) being satisfied by standard choices such as B-Splines or wavelets.
Assumption \ref{ass:heterosieve}(iv) is an asymptotic unbiasedness requirement -- a condition that is eased by noting no requirements are imposed on the approximating space for $\mu_0$.
The requirements on $\{q_{k}\}_{k=1}^{k_n}$ are imposed in Assumption \ref{ass:heteromoments}(i)(iii) and are again satisfied by standard choices. 
Assumption \ref{ass:heteromoments}(iv) states a rate condition that suffices for verifying the coupling requirements of Theorem \ref{th:localdrift}.
Assumption \ref{ass:heterosigma} imposes the requirements on the weighting matrices.

Our next result employs Theorem \ref{th:localdrift}(ii) to obtain strong approximations.

\begin{theorem}\label{th:heteroapprox}
 Let Assumptions \ref{ass:heterosieve}, \ref{ass:heteromoments}, \ref{ass:heterosigma} hold, $a_n=(\log (n))^{-1/2}$, and for any $\theta = (\{F(c_\jmath|\cdot)\}_{\jmath =1}^{{\mathcal{J}}},\mu )\in \mathbf B$ let $\Vert \theta\Vert _{{\mathbf{E}}}=\sum_{\jmath = 1}^{\mathcal J} \sup_{P\in \mathbf P} \|F(c_\jmath|\cdot)\|_{P,2}$.
If $\ell _n,\ell _n^{{\text {\rm u}}}\downarrow 0$ satisfy $k_n\sqrt{j_n}\log ^2(n)(\ell _n\vee \ell _n^{{\text {\rm u}}})=o(1)$, $k_n\sqrt{j_n}\log (n)/\sqrt{n}=o(\ell _n\wedge \ell _n^{{\text {\rm u}}})$, then uniformly in $P\in {\mathbf{P}}_0$:
\begin{align*}
I_n(R)&=\Up(R|\ell _n)+o_P(a_n) \\
I_n(R)-I_n(\Theta )&=\Up(R|\ell _n)-\Up(\Theta |\ell _n^{{\text {\rm u}}})+o_P(a_n).
\end{align*}
\end{theorem}

In order to conduct inference, we next aim to estimate the distributions of $\Up(R|\ell_n)$ and $\Up(\Theta|\ell_n^{\rm u})$.
To this end, we note that $\IDsetRsieve$ (as in \eqref{def:thetaset}) is potentially non-singleton and we therefore employ a set estimator $\hat \Theta_n^{{\rm r}}$ (as in \eqref{def:setest}) to estimate the distribution of $\Up(R|\ell_n)$.
In contrast, since $\Up(\Theta|\ell_n^{\rm u})$ only depends on the identified component $\{F_P(c_{\jmath }|\cdot )\}_{\jmath =1}^{{\mathcal{J}}}$, for the unconstrained problem we employ
any minimizer $\hat \theta_n^{\rm u}$ of $Q_n$ over $\Theta_n$.
With regards to the local parameter space, we note that in this application
\begin{equation}\label{ex:hetero20}
G_{n}(\theta )=\{(\{p^{j_n\prime }\beta _{\jmath,h}\}_{\jmath =1}^{{\mathcal{J}} },\mu _h) :\mu _{h}(B)\geq \sqrt n \min \{r_{n}-\mu (B),0\}{\text{ for all }}B\in {\mathcal{B}}\}
\end{equation}
for any $\theta = (\{F(c_{\jmath }|\cdot )\}_{\jmath =1}^{{\mathcal{J}}},\mu )$. 
Computationally, since any $\mu ,\mu _{h}\in {\mathcal{M}}_{n}$ has the structure $\mu = \sum_{s=1}^{s_n} \alpha_s\delta_{s}$ and $\mu_h = \sum_{s=1}^{s_n} \alpha_{sh}\delta_{s}$ it follows that the constraints in \eqref{ex:hetero20} reduce to $\alpha_{sh}\geq \min \{r_{n}-\alpha_{s},0\}$ for all $1\leq s \leq s_n$ whenever $\{\delta_{s}\}_{s=1}^{s_{n}}$ are orthogonal.
Furthermore, since moments and restrictions are linear, we may let $\ell_n = +\infty$ and set
\begin{equation}\label{ex:hetero21}
\hat{V}_{n}(\theta ,R|+\infty )= \{(\{p^{j_{n}\prime }\beta _{\jmath,h}\}_{\jmath =1}^{{\mathcal{J}}},\mu _{h}): h\in G_{n}(\theta ),~\Upsilon_F(h) = 0 \}.
\end{equation}
For each $\theta \in \Theta_n$, we denote the bootstrap process for the $\jmath^{th}$ conditional moment by
$$\hat{\mathbb W}_{\jmath,n}(\theta)= \frac{1}{\sqrt n}\sum_{i=1}^n \omega_i\{\rho_\jmath(X_i,\theta)q^{k_n}(W_i) - \frac{1}{n}\sum_{j=1}^n \rho_\jmath(X_j,\theta)q^{k_n}(W_j)\}.$$
Similarly, we set  $\hat {\mathbb D}_{\jmath,n}[h] = -\sum_{i=1}^n q^{k_n}(W_i)p^{j_n}(W_i)^\prime \beta_{\jmath,h}/n$ for any $h = (\{p^{j_n \prime} \beta_{\jmath, h}\}_{\jmath = 1}^{\mathcal J}, \mu_h)$. 
Thus, the estimators of the strong approximations obtained in Theorem \ref{th:heteroapprox} equal
\begin{align*}
\hat{U}_{n}(R|+\infty )& = \inf_{\theta \in \hat{\Theta}_{n}^{{\text{\textrm{r}}}}}\inf_{h\in \hat{V}_{n}(\theta,R|+\infty )}\{\sum_{\jmath =1}^{{\mathcal{J}}}\Vert \hat{{\mathbb{W}}}_{\jmath,n}(\theta) +\hat{{\mathbb{D}}}_{\jmath ,n}[h]\Vert _{\hat{\Sigma}_{\jmath ,n},2}\}^{1/2} \\
\hat{U}_{n}(\Theta |+\infty )& = \inf_{h }\{\sum_{\jmath =1}^{{\mathcal{J}}}\Vert \hat{{\mathbb{W}}}_{\jmath,n}(\hat{\theta}_{n}^{{\text{\textrm{u}}}})+\hat{{\mathbb{D%
}}}_{\jmath ,n}[h]\Vert _{\hat{\Sigma}_{\jmath ,n},2}\}^{1/2}.
\end{align*}

Before stating our final assumption, we need an auxiliary result. 
To this end, define
\begin{equation}  \label{eq:hetero23}
\Gamma _n(\theta )\equiv \{\tilde {\mu }\in {\mathcal{M}}_n: \tilde {\theta }=(\{F(c_{\jmath }|\cdot )\}_{\jmath =1}^{{\mathcal{J}}},%
\tilde {\mu }){\text { satisfies }}\Upsilon _F(\tilde {\theta })=0,~\Upsilon
_G(\tilde {\theta })\leq 0\}
\end{equation}
for any $\theta = (\{F(c_{\jmath }|\cdot )\}_{\jmath =1}^{{\mathcal{J}}},\mu )$ -- i.e.\ $\Gamma _n(\theta )$ is the set of distributions of $\eta$ that agree with the restrictions implied by $\{F(c_{\jmath }|\cdot )\}_{\jmath =1}^{{\mathcal{J}}}$.
Our next result bounds the $\|\cdot \|_{TV}$-Hausdorff distance between $\Gamma_n(\theta _1)$ and $\Gamma _n(\theta_2)$, which we denote by $d_H(\Gamma_n(\theta_1),\Gamma_n(\theta_2),\|\cdot\|_{TV})$.

\begin{lemma}\label{lm:heterolipschitz}
If the probability measures $\{\delta_{s}\}_{s=1}^{s_n}$ are orthogonal, then for every $n$ there exists a constant $\zeta _n<\infty $ independent of $\mathbf P$ such that 
\begin{equation*}
d_H(\Gamma _n(\theta _1),\Gamma _n(\theta _2),\|\cdot \|_{TV})\leq \zeta _n\sum _{\jmath =1}^{{\mathcal{J}}}\|F_1(c_{\jmath }|\cdot )-F_2(c_{\jmath }|\cdot )\|_{\infty }
\end{equation*}
for any $(\{F_1(c_{\jmath }|\cdot )\}_{\jmath =1}^{{\mathcal{J}}},\mu _1)=\theta _1\in \Theta _n\cap R$ and $(\{F_2(c_{\jmath }|\cdot )\}_{\jmath =1}^{{\mathcal{J}}},\mu _2)=\theta _2\in \Theta _n\cap R$.

\end{lemma}

We introduce our final assumption to show the validity of our bootstrap procedure.

\begin{assumption}\label{ass:heteroboot}
(i) $\Psi (g,\cdot )$ is bounded on $\Omega $;
(ii) The probability measures $\{\delta_{s}\}_{s=1}^{s_{n}}$ are orthogonal;
(iii) $k_{n}^{4}j_{n}^{5}\log ^{5}(n)/n=o(1)$;
(iv) $\Pi _{n}\theta=(\{F_{n}(c_{\jmath }|\cdot )\}_{\jmath =1}^{{\mathcal{J}}},\mu_n )$ satisfies $\Vert F_{n}(c_{\jmath }|\cdot)-F_{P}(c_{\jmath }|\cdot )\Vert _{\infty }=o(1)$ uniformly in $\theta\in\IDset \cap R$ and $P\in {\mathbf{P}}_{0}$;
(v) $k_{n}\sqrt{j_{n}}\log ^{2}(n)\tau _{n}=o(1)$, and $\zeta _{n}(k_{n}j_{n}\log (n)/\sqrt{n}+\sqrt{j_{n}}\tau _{n})=o(r_{n})$.
\end{assumption}

The boundedness of $\Psi (g,\cdot )$ on $\Omega $ ensures $\Upsilon _{F}^{({\text{\textrm{s}}})}$ (as in \eqref{ex:hetero9p3}) is continuous, while Assumption \ref{ass:heteroboot}(ii) allows us to apply Lemma \ref{lm:heterolipschitz}.
Assumption \ref{ass:heteroboot}(iii) is a low level sufficient condition for verifying the bootstrap coupling requirement of Assumption \ref{ass:bootcoupling}.
These rate requirements could be improved under smoothness conditions on $F_{P}(c_{\jmath }|\cdot )$.
Finally, Assumption \ref{ass:heteroboot}(iv) imposes a mild requirement on the sieve, while Assumption \ref{ass:heteroboot}(v) states conditions on $\tau _{n}$ and $r_{n}$ -- note $\tau _{n}=0$ and $r_{n}=+\infty $ are always valid, though such choices can lead to lower local power against certain alternatives.

Our final result obtains a coupling for our bootstrap approximations.

\begin{theorem}\label{th:heteroboot}
Let the conditions of Theorem \ref{th:heteroapprox} hold and Assumption \ref{ass:heteroboot} be satisfied.
Then: there are sequences $\ell _n,\ell _n^{{\rm u}}\downarrow 0$ satisfying $k_n\sqrt {j_n}\log (n)/\sqrt{n}=o(\ell _n\wedge \ell _n^{\rm u})$ and $k_n\sqrt{j_n}\log ^2(n)(\ell _n\vee \ell
_n^{{\rm u}})=o(1)$ such that uniformly in $P\in {\mathbf{P}}_0$
\begin{align*}
\hat {U}_n(R|+\infty )&\geq \UpS(R|\ell _n)+o_P(a_n) \\
\hat {U}_n(R|+\infty )-\hat {U}_n(\Theta |+\infty )&\geq \UpS(R|\ell _n)-\UpS(\Theta |\ell _n^{{{\rm u}}})+o_P(a_n).
\end{align*}
\end{theorem}

In particular, since the conditions on $\ell_n$ and $\ell_n^{\rm u}$ imposed in Theorems \ref{th:heteroapprox} and \ref{th:heteroboot} are the same, it follows that we may employ the quantiles of $\hat {U}_n(R|+\infty )$ and $\hat {U}_n(R|+\infty )-\hat {U}_n(\Theta |+\infty )$ conditional on the data as critical values for $I_n(R)$ and $I_n(R)-I_n(\Theta )$.


\section{Simulation Evidence} \label{sec:mc}

To conclude, we study the finite sample performance of our inference procedure by revisiting the simulation design in \cite{chetverikov2017nonparametric}.

\subsection{Identified Model}\label{sec:mcid}

We first consider a nonparametric instrumental variable model in which, for some unknown function $\IDpoint$, the distribution of $(Y,W,Z) \in \mathbf R^3$ satisfies the restriction
\begin{equation}\label{mc:eq1}
Y = \IDpoint(W) + \varepsilon \hspace{0.3 in} E[\varepsilon|Z] = 0;
\end{equation}
see Appendix \ref{sec:examples} for a formal study of this model.
Following \cite{chetverikov2017nonparametric}, we set $\IDpoint(w) \equiv 0.2w+w^2$ and for $(\epsilon,\zeta,\nu)$ independent standard normal random variables we let $Z = \Phi(\zeta)$, $W = \Phi(0.3 \zeta + \sqrt{1-(0.3)^2} \epsilon)$, and $\varepsilon = (0.3 \epsilon + \sqrt{1-(0.3)^2} \nu)/2$ for $\Phi$ the cumulative distribution function of a standard normal.
All reported results are based on five thousand replications employing five hundred bootstrap draws each.

In what follows, we utilize the restriction $\Upsilon_F(\IDpoint) = 0$ to impose a hypothesized value on the the level or the derivative of $\IDpoint$ at the point $w_0 = 0.5$ and use $\Upsilon_G(\IDpoint)\leq 0$ to impose that $\theta_0$ be either monotonically increasing or monotonically increasing and convex.
We employ the test statistic $I_n(R)-I_n(\Theta)$ with $p = 2$ and $\hat \Sigma_n$ an estimate of the optimal weighting matrix based on a first stage unconstrained estimator.
The implementation of the test is similar to that of the linear model of Section \ref{sec:testexiv}, with the difference that we must select the sieve $\Theta_n = \{p^{j_n\prime}\beta : \beta \in \mathbf R^{j_n}\}$ and $q^{k_n}$.
We follow \cite{chetverikov2017nonparametric} in employing continuously differentiable piecewise quadratic splines with equally spaced knots for both $p^{j_n}$ and $q^{k_n}$.

In computing critical values we set $\ell_n = +\infty$ since the model is linear and $\tau_n = 0$ since the model is identified.
We select $r_n$ by proceeding as in Section \ref{sec:testexiv}.
Specifically, the choice of sieve implies that, for any $\theta = p^{j_n\prime}\beta$, the restriction $\Upsilon_G(\theta )\leq 0$ is equivalent to $G\beta \leq 0$ for a known matrix $G$.
For $p^{j_n\prime}\hat \beta_n^{\rm u}$ the minimizer of $I_n(\Theta)$ and $p^{j_n\prime}\hat \beta_n^{\rm u\star}$ its score bootstrap analogue \citep{kline2012score}, we therefore set $r_n$ to satisfy
\begin{equation}\label{mc:eq2}
P(\max_{j} G_j\{\hat \beta_n^{\rm u\star \prime} - \hat \beta_n\} \leq r_n |\{V_i\}_{i=1}^n) = 1-\gamma_n
\end{equation}
where $\gamma_n \in (0,1)$ and the vectors $G_j\in \mathbf R^{j_n}$ depend on the shape restriction being imposed.
We emphasize that the sequence $\gamma_n$ must tend to zero in order for $r_n$ to satisfy our assumptions.
Finally, we employ the minimizer of $I_n(R)$ in obtaining bootstrap draws for both $\hat U_n(R|+\infty)$ and $\hat U_n(\Theta|+\infty)$; see discussion following Corollary \ref{th:incJcritval}.

\begin{table}[t!] \centerline{ \small
\begin{tabular}{cc c cccc c cccc} \hline \hline
         &                           & & \multicolumn{4}{c}{Imposed: Mon.}                      & & \multicolumn{4}{c}{Imposed: Mon.+ Conv.} \\ \cline{4-7} \cline{9-12}
         &                           & & \multicolumn{2}{c}{Level}  &  \multicolumn{2}{c}{Derivative}   & & \multicolumn{2}{c}{Level}& \multicolumn{2}{c}{Derivative} \\ \cline{2-12}
         & $r_n / (j_n,k_n)$         & & (4,4)  & (4,6)             & (4,4) & (4,6)                     & & (4,4) & (4,6)            & (4,4) & (4,6)    \\   \hline
         &   $\infty$                & & 1.90   & 1.72	            & 1.88	& 2.02	                    & & 1.44  & 1.52             & 2.74  & 2.84     \\
         &   95$\%$                  & & 1.74	& 1.68	            & 1.90	& 2.08	                    & & 1.46  & 1.54             & 2.68  & 2.84     \\
$n=500$  &   50$\%$                  & & 1.74	& 1.70	            & 1.90	& 2.10	                    & & 1.46  & 1.54             & 2.68  & 2.84     \\
         &   $5\%$                   & & 2.18	& 2.90	            & 2.20	& 2.96	                    & & 1.52  & 1.82             & 2.74  & 2.98     \\
         &   0                       & & 5.30	& 5.10	            & 4.62	& 4.48	                    & & 5.42  & 5.36             & 5.08  & 4.84     \\
\hline
         &   $\infty$                & & 1.56   & 1.82	            & 1.68	& 1.94	                    & & 1.40  & 1.54             & 2.26  & 2.32     \\
         &   95$\%$                  & & 1.52	& 1.84	            & 1.64	& 1.86	                    & & 1.36  & 1.44             & 2.04  & 2.26     \\
$n=1000$ &   50$\%$                  & & 1.52	& 1.86	            & 1.64	& 1.86	                    & & 1.36  & 1.44             & 2.04  & 2.26     \\
         &   $5\%$                   & & 2.02	& 2.84	            & 2.06	& 3.06	                    & & 1.44  & 1.86             & 2.14  & 2.38     \\
         &   0                       & & 4.54	& 4.56	            & 4.58	& 4.68	                    & & 4.62  & 4.78             & 4.38  & 4.20     \\
\hline
         &   $\infty$                & & 1.34   & 1.58	            & 1.26	& 1.52	                    & & 1.04  & 1.36             & 1.36  & 1.58     \\
         &   95$\%$                  & & 1.40	& 1.50	            & 1.32	& 1.62	                    & & 1.06  & 1.42             & 1.36  & 1.62     \\
$n=5000$ &   50$\%$                  & & 1.42	& 1.52	            & 1.32	& 1.62	                    & & 1.06  & 1.42             & 1.36  & 1.62     \\
         &   $5\%$                   & & 2.20	& 3.62	            & 2.36	& 3.36	                    & & 1.42  & 2.38             & 1.46  & 1.86     \\
         &   0                       & & 3.98	& 4.56	            & 4.68	& 4.50	                    & & 4.10  & 4.74             & 3.98  & 4.06     \\ \hline \hline
\end{tabular}}\caption{Empirical rejection probabilities for $5\%$-level tests based on $I_n(R)-I_n(\Theta)$.
Value of $r_n$ set to a percentile corresponds to choice of $1-\gamma_n$ in \eqref{mc:eq2}.}\label{table:mcsize1}
\end{table}

Table \ref{table:mcsize1} reports empirical rejection probabilities under the null hypothesis for $5\%$-level tests on the derivative and level of $\IDpoint$ at $w_0 = 0.5$ under different shape restrictions.
With regards to $r_n$, we examine the extreme possible values ($0$ and $\infty$) and choices corresponding to \eqref{mc:eq2} for different $\gamma_n$.
In accord to theory, which requires $\gamma_n \downarrow 0$, we find that the rejection probability is no larger than the nominal level except for very small values of $1-\gamma_n$.
Overall, we find the general lack of sensitivity to different choices of bandwidths to be reassuring for empirical practice.

\begin{figure}[t!]
\begin{subfigure}[t]{0.5 \textwidth}{
\begin{tikzpicture}[scale = 0.5]
            \pgfplotsset{every axis legend/.append style={scale only axis, at={(0.5,-0.15)},anchor=north}, ymin = 0, ymax = 1, xmin = -0.2, xmax = 1.1}
            \def\axisdefaultwidth{5in}
            \def\axisdefaultheight{2in}
            \begin{axis}[legend columns = 3, /pgf/number format/1000 sep={}, title = {\Large Power Curve for Level: $(j_n,k_n) = (4,4)$, $n = 1000$}]
                \addplot [mark = triangle, red] table [x=Grid, y=prn95]{Data4Graphs/PowerLev44MonN1000.txt};
                \addlegendentry{Mon};
                \addplot [mark=x, blue] table [x=Grid, y=prn95]{Data4Graphs/PowerLev44ConvN1000.txt};
                \addlegendentry{Mon+Conv};
                \addplot [mark= square, violet] table [x=Grid, y=GMM]{Data4Graphs/PowerLev44MonN1000.txt};
                \addlegendentry{Unres.};
            \end{axis}
\end{tikzpicture}}
\end{subfigure}
\begin{subfigure}[t]{0.5 \textwidth}{
\begin{tikzpicture}[scale = 0.5]
            \pgfplotsset{every axis legend/.append style={scale only axis, at={(0.5,-0.15)},anchor=north}, ymin = 0, ymax = 1, xmin = -0.2, xmax = 1.1}
            \def\axisdefaultwidth{5in}
            \def\axisdefaultheight{2in}
            \begin{axis}[legend columns = 3, /pgf/number format/1000 sep={}, title = {\Large  Power Curve for Level: $(j_n,k_n) = (4,6)$, $n = 1000$}]
                \addplot [mark = triangle, red] table [x=Grid, y=prn95]{Data4Graphs/PowerLev46MonN1000.txt};
                \addlegendentry{Mon};
                \addplot [mark=x, blue] table [x=Grid, y=prn95]{Data4Graphs/PowerLev46ConvN1000.txt};
                \addlegendentry{Mon+Conv};
                \addplot [mark= square, violet] table [x=Grid, y=GMM]{Data4Graphs/PowerLev46MonN1000.txt};
                \addlegendentry{Unres.};
            \end{axis}
\end{tikzpicture}}
\end{subfigure} \vspace{-0.05 in} \\
\begin{subfigure}{0.5 \textwidth}{
\begin{tikzpicture}[scale = 0.5]
            \pgfplotsset{every axis legend/.append style={scale only axis, at={(0.5,-0.15)},anchor=north}, ymin = 0, ymax = 1, xmin = -0.2, xmax = 1.1}
            \def\axisdefaultwidth{5in}
            \def\axisdefaultheight{2in}
            \begin{axis}[legend columns = 3, /pgf/number format/1000 sep={}, title = {\Large  Power Curve for Level: $(j_n,k_n) = (4,4)$, $n = 5000$}]
                \addplot [mark = triangle, red] table [x=Grid, y=prn95]{Data4Graphs/PowerLev44MonN5000.txt};
                \addlegendentry{Mon};
                \addplot [mark=x, blue] table [x=Grid, y=prn95]{Data4Graphs/PowerLev44ConvN5000.txt};
                \addlegendentry{Mon+Conv};
                \addplot [mark= square, violet] table [x=Grid, y=GMM]{Data4Graphs/PowerLev44MonN5000.txt};
                \addlegendentry{Unres.};
            \end{axis}
\end{tikzpicture}}
\end{subfigure}
\begin{subfigure}{0.5 \textwidth}{
\begin{tikzpicture}[scale = 0.5]
            \pgfplotsset{every axis legend/.append style={scale only axis, at={(0.5,-0.15)},anchor=north}, ymin = 0, ymax = 1, xmin = -0.2, xmax = 1.1}
            \def\axisdefaultwidth{5in}
            \def\axisdefaultheight{2in}
            \begin{axis}[legend columns = 3, /pgf/number format/1000 sep={}, title = {\Large  Power Curve for Level: $(j_n,k_n) = (4,6)$, $n = 5000$}]
                \addplot [mark = triangle, red] table [x=Grid, y=prn95]{Data4Graphs/PowerLev46MonN5000.txt};
                \addlegendentry{Mon};
                \addplot [mark=x, blue] table [x=Grid, y=prn95]{Data4Graphs/PowerLev46ConvN5000.txt};
                \addlegendentry{Mon+Conv};
                \addplot [mark= square, violet] table [x=Grid, y=GMM]{Data4Graphs/PowerLev46MonN5000.txt};
                \addlegendentry{Unres.};
            \end{axis}
\end{tikzpicture}}
\end{subfigure}\\ \vspace{0.1 in} \\
\begin{subfigure}{0.5 \textwidth}{
\begin{tikzpicture}[scale = 0.5]
            \pgfplotsset{every axis legend/.append style={scale only axis, at={(0.5,-0.15)},anchor=north}, ymin = 0, ymax = 1, xmin = 0, xmax = 4}
            \def\axisdefaultwidth{5in}
            \def\axisdefaultheight{2in}
            \begin{axis}[legend columns = 3, /pgf/number format/1000 sep={}, title = {\Large Power Curve for Derivative: $(j_n,k_n) = (4,4)$, $n = 1000$}]
                \addplot [mark = triangle, red] table [x=Grid, y=prn95]{Data4Graphs/PowerDer44MonN1000.txt};
                \addlegendentry{Mon};
                \addplot [mark=x, blue] table [x=Grid, y=prn95]{Data4Graphs/PowerDer44ConvN1000.txt};
                \addlegendentry{Mon+Conv};
                \addplot [mark= square, violet] table [x=Grid, y=GMM]{Data4Graphs/PowerDer44MonN1000.txt};
                \addlegendentry{Unres.};
            \end{axis}
\end{tikzpicture}}
\end{subfigure}
\begin{subfigure}{0.5 \textwidth}{
\begin{tikzpicture}[scale = 0.5]
            \pgfplotsset{every axis legend/.append style={scale only axis, at={(0.5,-0.15)},anchor=north}, ymin = 0, ymax = 1, xmin = 0, xmax = 4}
            \def\axisdefaultwidth{5in}
            \def\axisdefaultheight{2in}
            \begin{axis}[legend columns = 3, /pgf/number format/1000 sep={}, title = {\Large Power Curve for Derivative: $(j_n,k_n) = (4,6)$, $n = 1000$}]
                \addplot [mark = triangle, red] table [x=Grid, y=prn95]{Data4Graphs/PowerDer46MonN1000.txt};
                \addlegendentry{Mon};
                \addplot [mark=x, blue] table [x=Grid, y=prn95]{Data4Graphs/PowerDer46ConvN1000.txt};
                \addlegendentry{Mon+Conv};
                \addplot [mark= square, violet] table [x=Grid, y=GMM]{Data4Graphs/PowerDer46MonN1000.txt};
                \addlegendentry{Unres.};
            \end{axis}
\end{tikzpicture}}
\end{subfigure}\\ \vspace{0.1 in} \\
\begin{subfigure}{0.5 \textwidth}{
\begin{tikzpicture}[scale = 0.5]
            \pgfplotsset{every axis legend/.append style={scale only axis, at={(0.5,-0.15)},anchor=north}, ymin = 0, ymax = 1, xmin = 0, xmax = 4}
            \def\axisdefaultwidth{5in}
            \def\axisdefaultheight{2in}
            \begin{axis}[legend columns = 3, /pgf/number format/1000 sep={}, title = {\Large Power Curve for Derivative: $(j_n,k_n) = (4,4)$, $n = 5000$}]
                \addplot [mark = triangle, red] table [x=Grid, y=prn95]{Data4Graphs/PowerDer44MonN5000.txt};
                \addlegendentry{Mon};
                \addplot [mark=x, blue] table [x=Grid, y=prn95]{Data4Graphs/PowerDer44ConvN5000.txt};
                \addlegendentry{Mon+Conv};
                \addplot [mark= square, violet] table [x=Grid, y=GMM]{Data4Graphs/PowerDer44MonN5000.txt};
                \addlegendentry{Unres.};
            \end{axis}
\end{tikzpicture}}
\end{subfigure}
\begin{subfigure}{0.5 \textwidth}{
\begin{tikzpicture}[scale = 0.5]
            \pgfplotsset{every axis legend/.append style={scale only axis, at={(0.5,-0.15)},anchor=north}, ymin = 0, ymax = 1, xmin = 0, xmax = 4}
            \def\axisdefaultwidth{5in}
            \def\axisdefaultheight{2in}
            \begin{axis}[legend columns = 3, /pgf/number format/1000 sep={}, title = {\Large Power Curve for Derivative: $(j_n,k_n) = (4,6)$, $n = 5000$}]
                \addplot [mark = triangle, red] table [x=Grid, y=prn95]{Data4Graphs/PowerDer46MonN5000.txt};
                \addlegendentry{Mon};
                \addplot [mark=x, blue] table [x=Grid, y=prn95]{Data4Graphs/PowerDer46ConvN5000.txt};
                \addlegendentry{Mon+Conv};
                \addplot [mark= square, violet] table [x=Grid, y=GMM]{Data4Graphs/PowerDer46MonN5000.txt};
                \addlegendentry{Unres.};
            \end{axis}
\end{tikzpicture}}
\end{subfigure}
\caption{Rejection probabilities for $5\%$-level tests on conjectured value of $\IDpoint(0.5)$ (true value 0.35) and $\IDpoint^\prime(0.5)$ (true value 1.2). Tests implemented with $1-\gamma_n = 0.05$ in \eqref{mc:eq2}. }\label{fig:powall}
\end{figure}

In Figure \ref{fig:powall} we report power curves for different $5\%$-level tests concerning the value of $\IDpoint$ and its derivative at $w_0 = 0.5$.
For conciseness, we focus on the sample sizes $n\in \{1000,5000\}$ and $r_n$ chosen as in \eqref{mc:eq2} with $1-\gamma_n = 0.95$.
The curves labeled ``Mon" and ``Mon+Conv" correspond to tests based on $I_n(R)- I_n(\Theta)$ with $R$ imposing monotonicity and monotonicity and convexity while changing the conjectured value of $\IDpoint$ and its derivative at $w_0 = 0.5$.
The curve labeled ``Unres." corresponds to a Wald test based on the unrestricted estimator.
For all designs we find that imposing shape restrictions can improve power.
The effect of imposing shape restrictions, however, depend on both the sampling uncertainty and how ``close" the shape restrictions are to binding \citep{chetverikov2018econometrics}.
Since our design is fixed with $n$ and $\IDpoint$ is strictly increasing and convex, in our simulations we see the advantages of imposing shape restrictions decrease with $n$ as sample uncertainty decreases.
Similarly, since estimating the derivative is a harder than estimating the level, we observe larger power gains when imposing shape restrictions in the former problem. 

\subsection{Partially Identified Model}\label{sec:mcpid}

We next examine the performance of our test in a partially identified setting by discretizing the simulation design in \cite{chetverikov2017nonparametric}.
Concretely, we generate $(W,Z,\epsilon)\in [0,1]^2\times \mathbf R$ as in Section \ref{sec:mcid}, divide $[0,1]$ into $S_w$ and $S_z$ equally spaced segments, and generate dummy variables $D_w$ and $D_z$ for the segment to which $W$ and $Z$ belong -- e.g.\ if $(S_w,S_z) = (3,2)$, then $D_w(W) \equiv (1\{W\in [0,1/3]\}, 1\{W\in (1/3,2/3]\}, 1\{W\in (2/3,1]\})^\prime$ and $D_z(Z) \equiv (1\{Z\in [0,1/2]\}, 1\{Z\in (1/2,1]\})^\prime$.
The outcome $Y$ is generated according to \eqref{mc:eq1} but employing $D_w$ in place of $W$.

The discretized design is characterized by $S_z$ linear unconditional moment restrictions in $S_w$ unknowns.
For conciseness, we focus on imposing that $\IDpoint$ be monotonically increasing and convex while conducting inference on the value of $\IDpoint$ at the point $d_0 \equiv D_w(0.5)$ -- e.g, if $S_w = 3$, then $d_0 = (0,1,0)^\prime$.
The parameter $\IDpoint(d_0)$ is generically not identified whenever $S_w > S_z$ but, as we report in Table \ref{table:idset}, imposing a shape restriction on $\IDpoint$ partially identifies $\IDpoint(d_0)$.
A similar setting was previously studied by \cite{freyberger:horowitz:2012} who develop confidence regions for parameters such as $\IDpoint(d_0)$. Their leading procedure is computationally simpler than ours, but can suffer from size distortions, for example, when the identified set for $\IDpoint(d_0)$ is ``small."

\begin{table}[h!] \centerline {\small
\begin{tabular}{c c ccc}\hline \hline
                    &   & \multicolumn{3}{c}{$(S_w,S_z)$} \\ \cline{3-5}
Restriction on $\IDpoint$                    &   & $(3,2)$   & $(4,2)$   & $(3,2)$ \\ \hline
Mon.+Convex    &   & [0.059, 0.252]      & [0.100, 0.412]     & [0.310, 0.388]  \\
No Restriction      &   & $(-\infty,+\infty)$   & $(-\infty,+\infty)$   & $(-\infty,+\infty)$ \\  \hline \hline
\end{tabular}}\caption{Identified sets for $\IDpoint(d_0)$ with and without shape restrictions.}\label{table:idset}
\end{table}

We test whether a value $\lambda$ belongs to the identified set for $\IDpoint(d_0)$ by setting $\Upsilon_F(\theta) = \theta(d_0) - \lambda$ and employ the constraint $\Upsilon_G(\theta) \leq 0$ to impose that $\theta$ be monotonically increasing and convex.
We base inference on $I_n(R)$ with $p= 2$, $\hat \Sigma_n$ the sample analogue to $E[D_zD_z^\prime]$, all moment restrictions ($k_n = S_z$), and a saturated model for $\IDpoint$ ($j_n = S_w$).
To compute critical values we set $\ell_n = +\infty$ and $\tau_n = 0$ -- though note $\hat \Theta_n^{\rm r}$ need not be a singleton when $\tau_n = 0$ because $j_n > k_n$.
We select $r_n$ by modifying the approach employed in Section \ref{sec:mcid}.
Specifically, we note that the constraint $\Upsilon_G(\theta)\leq 0$ may be written as $G\theta \leq 0$ for some matrix $G$, and for $\hat \theta^{L}_n$ and $\hat \theta^U_n$ the minimizer and maximizers of $\theta(d_0)$ over the set of $\theta$ that are monotonically increasing, convex, and minimize $\| \sum_{i=1}^n (Y_i -\theta(D_{w,i}))D_{z,i}/n\|_\infty$, we set $r_n$ according to
\begin{equation}\label{mc:eq4}
P(\max_{ j} \max\{G_j(\hat \theta^{L\star}_n - \hat \theta^L_n), G_j(\hat \theta^{U\star}_n - \hat \theta^U_n)\} \leq r_n |\{V_i\}_{i=1}^n) = 1-\gamma_n,
\end{equation}
where $\hat \theta^{L\star}_n$ and $\hat\theta^{U\star}_n$ are again computed employing the score bootstrap.
As in our previous analysis, $\gamma_n$ must tend to zero with $n$ in order for $r_n$ to satisfy our assumptions.

\begin{table}[t!] \centerline{ \small
\begin{tabular}{cc c ccc c ccc c ccc} \hline \hline
            &               & & \multicolumn{3}{c}{Lower Endpoint} & & \multicolumn{3}{c}{Midpoint}    & & \multicolumn{3}{c}{Upper Endpoint} \\ \cline{4-6} \cline{8-10} \cline{12-14}
            &               & & \multicolumn{3}{c}{$(S_w,S_z)$} & & \multicolumn{3}{c}{$(S_w,S_z)$} & & \multicolumn{3}{c}{$(S_w,S_z)$} \\ \cline{2-14}
            &   $r_n$       & & (3,2)   & (4,2) & (4,3)         & & (3,2)   & (4,2) & (4,3)         & & (3,2)   & (4,2) & (4,3) \\ \hline
            &   $\infty$    & & 1.96    & 3.34  & 1.48          & &	0.10    & 0.02  & 1.48          & & 1.88    & 3.10  & 2.00  \\
            &   95$\%$      & & 3.64    & 4.70  & 1.46          & & 0.10    & 0.02  & 1.46          & & 2.26    & 3.12  & 1.98  \\
$n=500$     &   50$\%$      & & 5.34    & 5.24  & 1.46          & & 0.50    & 0.06  & 1.50          & & 5.22    & 5.02  & 2.04  \\
            &   $5\%$       & & 5.36    & 5.24  & 3.56          & & 0.50    & 0.06  & 3.44          & & 5.24    & 5.02  & 3.54  \\
            &   0           & & 5.34    & 5.26  & 4.64          & & 0.50    & 0.06  & 4.48          & & 5.24    & 5.16  & 4.60  \\
\hline
            &   $\infty$    & & 1.84    & 3.06  & 1.12          & & 0.00    & 0.00  & 1.10          & & 1.96    & 2.90	& 1.34  \\
            &   95$\%$      & & 4.98    & 4.84  & 1.12          & & 0.02    & 0.00  & 1.08          & & 2.98    & 2.90  & 1.34  \\
$n=1000$    &   50$\%$      & & 5.10    & 4.88  & 1.20          & & 0.12    & 0.00  & 1.14          & & 5.00    & 4.86  & 1.44  \\
            &   $5\%$       & & 5.10    & 4.88  & 3.48          & & 0.12    & 0.00  & 3.12          & & 5.00    & 4.86  & 2.78  \\
            &   0           & & 5.28    & 4.88  & 4.42          & & 0.08    & 0.00  & 4.14          & & 5.10    & 4.86  & 3.82  \\
\hline
            &   $\infty$    & & 1.98    & 4.40  & 1.34          & & 0.00    & 0.00  & 1.22          & & 1.98    & 2.80	& 1.36  \\
            &   95$\%$      & & 5.08    & 6.76  & 1.34          & & 0.00    & 0.00  & 1.26          & & 4.56    & 4.86  & 1.34  \\
$n=5000$    &   50$\%$      & & 5.08    & 8.30  & 1.48          & & 0.00    & 0.00  & 1.44          & & 4.58    & 4.84  & 1.52  \\
            &   $5\%$       & & 5.08    & 9.00  & 4.28          & & 0.00    & 0.00  & 4.14          & & 4.58    & 4.84  & 3.58  \\
            &   0           & & 4.96    & 8.84  & 4.70          & & 0.00    & 0.00  & 4.38          & & 4.64    & 5.02  & 4.46  \\ \hline \hline
\end{tabular}}\caption{Empirical rejection probabilities for $5\%$-level tests based on $I_n(R)$ for different points in the null hypothesis.
Lower and upper endpoints correspond to Table \ref{table:idset}.}\label{table:mcsize2}
\end{table}

Table \ref{table:mcsize2} reports empirical rejection rates for testing whether $\lambda$ belongs to the identified set, with the lower and upper endpoint columns corresponding to setting $\lambda$ to equal the lower and upper endpoints in Table \ref{table:idset}.
All tests are conducted at a $5\%$ nominal level.
Across designs, we find that setting $r_n = +\infty$ always delivers tests with rejection probabilities below their nominal level.
Setting $r_n$ according to \eqref{mc:eq4} with $1-\gamma_n = 0.95$ also delivers adequate size control, with the exception of $n = 5000$ and $(S_w,S_z) = (4,2)$ where we see a modest over-rejection at the lower endpoint of the identified set.
Overall, the degree of sensitivity to the choice of $r_n$ is similar to that found in Section \ref{sec:mcid}.


\phantomsection
\addcontentsline{toc}{section}{References}

{\small
\singlespace
\putbib
}

\end{bibunit}

\newpage


\begin{bibunit}

\setcounter{footnote}{0}

\renewcommand{\thesection}{A.\arabic{section}}
\renewcommand{\theequation}{A.\arabic{equation}}
\renewcommand{\thelemma}{A.\arabic{section}.\arabic{lemma}}
\renewcommand{\thecorollary}{A.\arabic{section}.\arabic{corollary}}
\renewcommand{\thetheorem}{A.\arabic{section}.\arabic{theorem}}
\renewcommand{\theassumption}{A.\arabic{section}.\arabic{assumption}}
\setcounter{lemma}{0}
\setcounter{theorem}{0}
\setcounter{corollary}{0}
\setcounter{equation}{0}
\setcounter{remark}{0}
\setcounter{section}{0}
\setcounter{assumption}{0}

\emptythanks

\title{Supplemental Appendix I}

{\author{Victor Chernozhukov \\ Department of Economics \\ M.I.T. \\ vchern@mit.edu
\and
Whitney K. Newey\thanks{Research supported by NSF Grant 1757140.} \\ Department of Economics \\ M.I.T. \\ wnewey@mit.edu
\and
Andres Santos\thanks{Research supported by NSF Grant SES-1426882.}
\\ Department of Economics\\ U.C.L.A.\\ andres@econ.ucla.edu}}

\date{April, 2022}

\maketitle

This Supplemental Appendix to ``Constrained Conditional Moment Restriction Models" is organized as follows. 
Sections \ref{sec:appAM} provides a review of AM spaces.
Section \ref{sec:examples} specializes the general results of Section \ref{sec:gentheory} to three additional examples: (i) GMM, (ii) Quantile Treatment Effects, and (iii) The Slutsky restriction in a partially linear model.
The proofs for all results are included in Supplmental Appendix II.

\thispagestyle{empty}

\newpage



\input{Appendix/AppAM}


\input{Appendix/ExamplesMain/AppExamples}


\phantomsection
\addcontentsline{toc}{section}{References}

{\small
\singlespace
\putbib }

\end{bibunit}

\newpage


\begin{bibunit}

\setcounter{footnote}{0}

\renewcommand{\thesection}{S.\arabic{section}}
\renewcommand{\theequation}{S.\arabic{equation}}
\renewcommand{\thelemma}{S.\arabic{section}.\arabic{lemma}}
\renewcommand{\thecorollary}{S.\arabic{section}.\arabic{corollary}}
\renewcommand{\thetheorem}{S.\arabic{section}.\arabic{theorem}}
\renewcommand{\theassumption}{S.\arabic{section}.\arabic{assumption}}
\setcounter{lemma}{0}
\setcounter{theorem}{0}
\setcounter{corollary}{0}
\setcounter{equation}{0}
\setcounter{remark}{0}
\setcounter{section}{0}
\setcounter{assumption}{0}

\emptythanks

\title{Supplemental Appendix II\\ Not Intended For Publication}

{\author{Victor Chernozhukov \\ Department of Economics \\ M.I.T. \\ vchern@mit.edu
\and
Whitney K. Newey\thanks{Research supported by NSF Grant 1757140.} \\ Department of Economics \\ M.I.T. \\ wnewey@mit.edu
\and
Andres Santos\thanks{Research supported by NSF Grant SES-1426882.}
\\ Department of Economics\\ U.C.L.A.\\ andres@econ.ucla.edu}}

\date{April, 2022}

\maketitle

This Supplemental Appendix to ``Constrained Conditional Moment Restriction Models" contains the proofs for all results.
Section \ref{sec:apprate} derives rate of convergence results that are employed in our strong and bootstrap approximations.
In Section \ref{sec:appstrong} we establish Theorem \ref{th:localdrift}, while the proofs for all remaining results concerning our bootstrap approximation and test are contained in Section \ref{sec:appboot}.
Section \ref{sec:exproofs} includes the proofs of the results stated in Section \ref{sec:hetero} and the examples discussed in Supplemental Appendix I.
Finally, Sections \ref{sec:localparam}, \ref{sec:unifcoupling}, and \ref{sec:bootcoup} develop results that may be of independent interest, and include the analysis of the local parameter space, empirical process coupling results based on \cite{koltchinskii:1994}, and bootstrap coupling results.

\thispagestyle{empty}

\newpage




\input{Appendix/AppRate}        

\input{Appendix/AppStrong}      

\input{Appendix/AppBoot}        


\input{Appendix/ExamplesProofs/AppExamplesProofs}


\input{Appendix/AppLocParam}     

\input{Appendix/AppKoltCoup}     

\input{Appendix/AppBootCoup}     


\phantomsection
\addcontentsline{toc}{section}{References}

{\small
\singlespace
\putbib }

\end{bibunit}

\end{document}

%% file: Appendix/AppAM.tex

\section{AM Spaces} \label{sec:appAM}

We provide a brief introduction to AM spaces and refer the reader to  Chapters 8 and 9 of \cite{aliprantis:border:2006} for a more detailed exposition.
Before proceeding, we first recall the definitions of a partially ordered set and a lattice.

\begin{definition}\label{def:ordered} \rm
A \emph{partially ordered set} $(\mathbf G, \geq)$ is a set $\mathbf G$ with a partial order relationship $\geq$ defined on it -- i.e.\ $\geq$ is a transitive $(x\geq y \text{ and } y \geq z \text{ implies } x\geq z)$, reflexive $(x\geq x)$, and antisymmetric $(x\geq y \text{ implies the negation of } y\geq x)$ relation. \qed
\end{definition}

\begin{definition}\label{def:lattice} \rm
A \emph{lattice} is a partially ordered set $(\mathbf G, \geq)$ such that any pair $x,y \in \mathbf G$ has a least upper bound (denoted $x\vee y$) and a greatest lower bound (denoted $x \wedge y$). \qed
\end{definition}

Whenever $\mathbf G$ is both a vector space and a lattice, it is possible to define objects that depend on both the vector space and lattice operations. In particular, for $x\in \mathbf G$ we define the positive part $x^+ \equiv x\vee 0$, the negative part $x^- \equiv (-x)\vee 0$, and the absolute value $|x| \equiv x\vee (-x)$.
It is also natural to demand that the order relation $\geq$ interact with the algebraic operations in a manner analogous to that of $\mathbf R$ -- i.e.\  to have
\begin{align}
x\geq y & \text{ implies } x+ z \geq y+ z \text{ for each } z\in \mathbf G \label{appA:eq2} \\
x\geq y & \text{ implies } \alpha x \geq \alpha y \text{ for each }  0\leq \alpha \in \mathbf R ~. \label{appA:eq3}
\end{align}
A complete normed vector space that shares these familiar properties of $\mathbf R$ under a given order relation $\geq$ is referred to as a \emph{Banach lattice}. Formally, we define:

\begin{definition}\label{def:riesz} \rm
A Banach space $\mathbf G$ with norm $\|\cdot\|_{\mathbf G}$ is a \emph{Banach lattice} if (i) $\mathbf G$ is a lattice under $\geq$, (ii) $\| x \|_{\mathbf G} \leq \|y\|_{\mathbf G}$ when $|x|\leq |y|$, (iii) \eqref{appA:eq2} and \eqref{appA:eq3} hold. \qed
\end{definition}

An AM space is a Banach lattice in which the maximum of the norms of any two positive elements is equal to the norm of the maximums of the two elements.

\begin{definition}\label{def:am} \rm
A Banach lattice $\mathbf G$ is called an AM space if for any elements $0\leq x, y \in \mathbf G$ it follows that $\|x\vee y\|_{\mathbf G} = \max\{\|x\|_{\mathbf G},\|y\|_{\mathbf G}\}$. \qed
\end{definition}

In certain Banach lattices there may exist an element $\mathbf {1_G} > 0$ called an \emph{order unit} such that for any $x \in \mathbf G$ there exists a $0 < \lambda \in \mathbf R$ for which $|x| \leq \lambda  \mathbf {1_G}$ -- for example, in $\mathbf R^d$ the vector $(1,\ldots, 1)^\prime$ is an order unit. The order unit $\mathbf {1_G}$ can be used to define
\begin{equation}\label{appA:eq4}
\|x\|_\infty \equiv\inf \{\lambda > 0 : |x| \leq \lambda \mathbf {1_G}\},
\end{equation}
which is a norm on $\mathbf G$.
In principle, $\|\cdot\|_\infty$ need not be related to the original norm $\|\cdot\|_{\mathbf G}$. 
However, if $\mathbf G$ is an AM space, then $\|\cdot\|_{\mathbf G}$ and $\|\cdot\|_\infty$ are equivalent in that they generate the same topology. 
Hence,  we refer to $\mathbf G$ as an \emph{AM space with unit} $\mathbf {1_G}$ if: (i) $\mathbf G$ is an AM space, (ii) $\mathbf {1_G}$ is an order unit in $\mathbf G$, and (iii) The norm of $\mathbf G$ equals $\|\cdot\|_\infty$.

%% file: Appendix/ExamplesMain/AppExamples.tex
\section{Illustrative Examples} \label{sec:examples}

In this Section, we examine special cases of our general analysis and illustrate both how to implement our procedure and verify the assumptions in the main text.

\input{Appendix/ExamplesMain/ExGMM}      

\input{Appendix/ExamplesMain/ExNPIV}     

\input{Appendix/ExamplesMain/ExQuant}    

%% file: Appendix/ExamplesMain/ExGMM.tex

\subsection{Generalized Method of Moments}\label{sec:gmmex}

Our first example concerns the generalized method of moments (GMM) model of \cite{hansen:1982}.
We assume the parameter of interest $\IDpoint$ is identified as the unique solution to 
\begin{equation}\label{eq:appgmm0}
E_P[\rho(X,\IDpoint)] = 0,
\end{equation}
where $X\in \mathbf X$ is distributed according to $P\in \mathbf P$ and $\rho: \mathbf X \times \Theta \to \mathbf R^{\mathcal J}$.
This model maps into our general framework by letting $Z_\jmath = 1$ for all $1\leq \jmath \leq \mathcal J$.
Moreover, since we have assumed $\IDpoint$ is identified, the hypothesis testing problem simplifies to
$$H_0:\IDpoint\in R \hspace{0.5 in} H_1:\IDpoint \notin R.$$

The set $R$ is, as in the main text, defined by equality and inequality restrictions.
In particular, for known functions $\Upsilon_F : \mathbf R^{d_\theta} \to \mathbf R^{d_F}$ and $\Upsilon_G : \mathbf R^{d_\theta} \to \mathbf R^{d_G}$ we set
\begin{equation}\label{eq:appgmm0p5}
R \equiv \{\theta \in \mathbf R^{d_\theta} : \Upsilon_F(\theta) = 0 \text{ and } \Upsilon_G(\theta) \leq 0\}.
\end{equation}
To verify Assumptions \ref{ass:param}(ii)(iii), note $\mathbf R^d$ is a Banach space under any norm $\|\cdot\|_p$ with $1\leq p \leq \infty$, so for concreteness we set $\mathbf B = \mathbf R^{d_\theta}$, $\mathbf F = \mathbf R^{d_F}$, and $\|\cdot\|_{\mathbf B} = \|\cdot\|_{\mathbf F}  = \|\cdot\|_2$.
The space $\mathbf R^{d}$ is in addition a lattice under the standard pointwise partial order
\begin{equation}\label{eq:appgmm1}
a \leq b \text{ if and only if } a_i \leq b_i \text{ for all } 1\leq i \leq d
\end{equation}
for any $(a_1,\ldots, a_d)^\prime = a$ and $(b_1,\ldots, b_d)^\prime = b$ in $\mathbf R^d$, while the least upper bound equals
$$a \vee b  = (\max\{a_1,b_1\},\ldots, \max\{a_d,b_d\})^\prime. $$
The vector $(1,\ldots, 1)^\prime$ is an order unit in $\mathbf R^d$ under the partial order in \eqref{eq:appgmm1}. 
As discussed in Section \ref{sec:appAM} of this Supplemental Appendix, the order unit induces the norm
$$\{\inf \lambda > 0 : |a|\leq \lambda (1,\ldots, 1)^\prime\} = \max_{1\leq i \leq d} |a_i|,$$
which corresponds to the usual $\|\cdot\|_\infty$ norm on $\mathbf R^d$.
Hence, by setting $\mathbf G = \mathbf R^{d_G}$, $\|\cdot\|_{\mathbf G} = \|\cdot\|_\infty$, and $\mathbf {1_G} = (1,\ldots, 1)^\prime$ we verify the requirements of Assumption \ref{ass:param}(ii)(iii).

Since the parameter space $\Theta$ is finite dimensional and all moment restrictions are unconditional, we may set $\Theta_n = \Theta$ and $k_n = \mathcal J$ for all $n$.
We base our test statistic on quadratic forms in the moments ($p =2$), which implies $Q_n(\theta)$ is given by
$$Q_n(\theta) \equiv \| \hat \Sigma_n \{\frac{1}{n}\sum_{i=1}^n \rho(X_i,\theta)\}\|_{2}.$$
In what follows, we consider tests based on both the un-centered statistic $I_n(R)$ and the re-centered statistic $I_n(R) - I_n(\Theta)$.
To this end, we impose the following:

\begin{assumption}\label{ass:gmmbasic}
(i) $\{X_i\}_{i=1}^n$ is i.i.d.\ with $X_i \sim P\in \mathbf P$;
(ii) For each $P\in \mathbf P_0$ there exists a unique $\IDpoint \in \Theta$ solving \eqref{eq:appgmm0}; (iii) $\Theta$ is convex and compact.
\end{assumption}

\begin{assumption}\label{ass:gmmrho}
(i) The function $\rho(x,\cdot) : \Theta \to \mathbf R^{\mathcal J}$ is twice differentiable for all $x$;
(ii) $E_P[\sup_{\theta \in \Theta} \|\rho(X,\theta)\|_2^3]$, $E_P[\sup_{\theta \in \Theta} \|\nabla_\theta \rho(X,\theta)\|_{o,2}^{2}]$, $E_P[\sup_{\theta \in \Theta}\|\nabla^2_\theta \rho_\jmath(X,\theta)\|^{1+\delta}_{o,2}]$ are finite and bounded uniformly in $P\in \mathbf P$ for some $\delta > 0$.
\end{assumption}

\begin{assumption}\label{ass:gmmid}
(i) $\inf_{P\in \mathbf P_0} \inf_{\theta \in \Theta : \|\theta- \IDpoint\|_2 \geq \epsilon} \|E_P[\rho(X,\theta)]\|_2 > 0$ for all $\epsilon > 0$;
(ii) The singular values of $E_P[\nabla_\theta \rho(X,\IDpoint)]$ are bounded away from zero in $P\in \mathbf P_0$.
\end{assumption}

\begin{assumption}\label{ass:gmmsigma}
(i) $\|\hat \Sigma_n - \PSigma\|_{o,2} = O_P(n^{-1/2})$ uniformly in $P\in \mathbf P$;
(ii) $\PSigma$ is invertible and $\|\PSigma\|_{o,2}$ and $\|\PSigma^{-1}\|_{o,2}$ are bounded uniformly in $P\in \mathbf P$.
\end{assumption}

In Assumption \ref{ass:gmmrho} we focus on differentiable moments for simplicity.
Assumption \ref{ass:gmmid} essentially imposes strong identification of $\IDpoint$ and hence guarantees that $\IDpoint$ can be consistently estimated uniformly in $P\in \mathbf P_0$ -- recall that $\IDpoint$ depends on $P$ through \eqref{eq:appgmm0}, though the dependence is left implicit in the notation.
Finally, Assumption \ref{ass:gmmsigma} states the requirements on the $\mathcal J \times 
\mathcal J$ weighting matrix $\hat \Sigma_n$. 

In what follows, we set the local parameter spaces $V_n(\theta,R|\ell)$ and $V_n(\theta,\Theta|\ell)$ to equal
\begin{align*}
V_n(\theta,R|\ell) & = \{h \in \mathbf R^{d_\theta} : \theta+ h/\sqrt n \in \Theta \cap R \text{ and } \|h/\sqrt n\|_2 \leq \ell\} \\
V_n(\theta,\Theta|\ell) & = \{h \in \mathbf R^{d_\theta} : \theta+ h/\sqrt n \in \Theta \text{ and } \|h/\sqrt n\|_2 \leq \ell\} .
\end{align*}
Setting $\DerP(\IDpoint)[h]\equiv E_P[\nabla_\theta \rho(X,\IDpoint)]h$ and letting $\WP(\IDpoint) \sim N(0, \text{Var}_P\{\rho(X,\IDpoint)\})$ we then denote the variables to which $I_n(R)$ and $I_n(\Theta)$ will be coupled to by
\begin{align*}
\Up(R|\ell_n) & \equiv \inf_{h \in V_n(\IDpoint,R|\ell_n)} \|\mathbb W_{P}(\IDpoint) + \DerP(\IDpoint)[h] \|_{\PSigma,2} \\ 
\Up(\Theta|\ell_n) & \equiv \inf_{h \in V_n(\IDpoint,\Theta|\ell_n)} \|\mathbb W_{P}(\IDpoint) + \DerP(\IDpoint)[h] \|_{\PSigma,2} .  
\end{align*}

Our distributional approximations follow immediately from Theorem \ref{th:localdrift}(ii).

\begin{theorem}\label{th:gmmapprox}
Let Assumptions \ref{ass:gmmbasic}, \ref{ass:gmmrho}, \ref{ass:gmmid}, and \ref{ass:gmmsigma} hold, $\Upsilon_F$ and $\Upsilon_G$ be continuous, and set $a_n = \sqrt{\log(n)}/n^{\frac{1}{10 + 5d_\theta}}$.
Then: For any $\ell_n,\ell_n^{\text{\rm u}} \downarrow 0$ satisfying $(\ell_n\vee \ell_n^{\text{\rm u}}) \sqrt{\log(1/\ell_n\vee\ell_n^{\text{\rm u}})} = o(a_n)$ and $n^{-1/2} = o(\ell_n\vee\ell_n^{\text{\rm u}})$ we have uniformly in $P\in \mathbf P_0$
\begin{align*}
I_n(R) & = \Up (R|\ell_n) + o_P(a_n) \notag \\
I_n(R) - I_n(\Theta) & = \Up(R|\ell_n) - \Up(\Theta|\ell_n^{\text{\rm u}}) + o_P(a_n).
\end{align*}
\end{theorem}


The rate of coupling $a_n = \sqrt{\log(n)}/n^{\frac{1}{10+5d_\theta}}$ obtained in Theorem \ref{th:gmmapprox} suffices for both the empirical process and bootstrap coupling; see Lemmas \ref{lm:gmmcoupling} and \ref{lm:bootgmmcoupling} in Supplemental Appendix II.
While the rate is adequate for our purposes, it can be improved under additional moment restrictions.
Here, we rely in \cite{yurinskii:1977} both to illustrate the diversity of coupling arguments that can be employed to verify Assumption \ref{ass:coupling}(i) and to impose only the weak third moment restriction of Assumption \ref{ass:gmmrho}(ii).

Our next goal is to obtain bootstrap approximations to the distributions of $\Up(R|\ell_n)$ and $\Up(\Theta|\ell_n^{\rm u})$. 
To this end, we write $\Upsilon_F(\theta) = (\Upsilon_{F,1}(\theta),\ldots, \Upsilon_{F,d_F}(\theta))^\prime$ and $\Upsilon_G(\theta) = (\Upsilon_{G,1}(\theta),\ldots, \Upsilon_{G,d_G}(\theta))^\prime$, for any $\epsilon > 0$ we define $B^\epsilon \equiv \bigcup_{P\in\mathbf P_0} \{\theta : \|\theta - \IDpoint\|_2 \leq \epsilon\}$ (where recall $\IDpoint$ implicitly depends on $P$ through \eqref{eq:appgmm0}), and impose:

\begin{assumption}\label{ass:gmmbootp1}
For some $\epsilon > 0$:
(i)  $B^\epsilon \subseteq \Theta$;
(ii) $\Upsilon_{F}$ and $\Upsilon_{G}$ are twice differentiable on $B^\epsilon$;
(iii) $ \|\nabla \Upsilon_F(\theta)\|_{o,2}$ and $\|\nabla \Upsilon_G(\theta)\|_{o,2}$ are bounded on $B^\epsilon$;
(iv) $\|\nabla^2 \Upsilon_{F,j}(\theta)\|_{o,2}$ is bounded on $B^\epsilon$ for $1\leq j \leq d_F$;
(v)  $\|\nabla^2 \Upsilon_{G,j}(\theta)\|_{o,2}$ is bounded on $B^\epsilon$ for $1\leq j \leq d_G$;
(vi) $\nabla \Upsilon_F(\theta)$ has full row-rank on $B^\epsilon$.
\end{assumption}

\begin{assumption}\label{ass:gmmbootp2}
Either (i) $\Upsilon_F : \mathbf R^{d_\theta} \to \mathbf R^{d_F}$ is affine, or (ii) There is an $\epsilon > 0$ and $M < \infty$ such that the singular values of $\nabla \Upsilon_F(\theta)^\prime$ are bounded away from zero uniformly in $\theta \in B^\epsilon$, and for every $P\in \mathbf P_0$ there is an $h\in \mathcal N(\nabla \Upsilon_F(\IDpoint))$ with $\|h\|_2 \leq M$ satisfying $\Upsilon_{G,j}(\IDpoint) + \nabla \Upsilon_{G,j}(\IDpoint)[h] \leq -\epsilon$ for all $1\leq j \leq d_G$.
\end{assumption}

In order to describe our bootstrap procedure in this application, we let $\hat \theta_n$ and $\hat \theta_n^{\rm u}$ denote the minimizers of $Q_n$ over $\Theta \cap R$ and $\Theta$ respectively.
Employing $\hat \theta_n$ and $\hat \theta_n^{\text{u}}$ we obtain estimators for the distribution of  $\mathbb W_{P}(\IDpoint)$ and for $\DerP(\IDpoint)$ by evaluating 
\begin{align}
\hat{\mathbb W}_n (\theta) & \equiv \frac{1}{\sqrt n} \sum_{i=1}^n \omega_i\{\rho(X_i,\theta) - \frac{1}{n}\sum_{j=1}^n \rho(X_j,\theta)\}  \label{eq:appgmm2}\\
\hat{\mathbb D}_n(\theta ) & \equiv \frac{1}{n}\sum_{i=1}^n \nabla_\theta \rho(X_i, \theta), \label{eq:appgmm3}
\end{align}
at $\theta = \hat \theta_n$ and $\theta = \hat \theta_n^{\text{u}}$, where recall $\{\omega_i\}_{i=1}^n$ is an i.i.d.\ sample independent of $\{X_i\}_{i=1}^n$ with $\omega_i \sim N(0,1)$.
We note that because moments are differentiable, we employ an analytical derivative in \eqref{eq:appgmm3} instead of the numerical derivative studied in Section \ref{sec:gentheory}.

With regards to the local parameter space, we note that the construction of $\hat V_n(\theta,R|\ell)$ requires the bound $K_g$ on the second derivative of $\Upsilon_G$ (as specified in Assumption \ref{ass:locineq}).
In particular, Assumption \ref{ass:gmmbootp1}(v) implies Assumption \ref{ass:locineq} is satisfied with
$$K_g \equiv \max_{1\leq j \leq d_G} \sup_{\theta \in B^\epsilon} \|\nabla^2_\theta \Upsilon_{G,j}(\theta)\|_{o,2}$$
(see Lemma \ref{lm:gmm4ver}).
If an a-priory bound on the second derivative is not available, then it is also possible to simply substitute $K_g$ with the data driven choice
$$\hat K_g \equiv \max_{1\leq j \leq d_G} \sup_{\theta \in \Theta : \|\theta - \hat \theta_n\|_2 \leq r_n} \|\nabla^2_\theta \Upsilon_{G,j}(\theta)\|_{o,2},$$
where we discuss the choice of $r_n$ below.
Given $K_g$ (or $\hat K_g$), we set $G_n(\theta)$ to equal
$$G_n(\theta) = \{h \in \mathbf R^{d_\theta} : \Upsilon_{j,G}(\theta + \frac{h}{\sqrt n}) \leq \max\{\Upsilon_{j,G}(\theta) - K_gr_n \|\frac{h}{\sqrt n}\|_2, - r_n\} \text{ for all } j\}$$
In this application we may additionally specify $\ell_n$ to be infinite, and hence we set
$$\hat V_n(\theta,R|+\infty) = \{h \in \mathbf R^{d_\theta} : h \in G_n(\theta) \text{ and } \Upsilon_F(\theta + \frac{h}{\sqrt n}) = 0\}.$$

The approximations to the distributions of $I_n(R)$ and $I_n(\Theta)$ are then given by the laws of $\hat U_n(R|+\infty)$ and $\hat U_n(\Theta|+\infty)$ conditional on the data, where
\begin{align*}
\hat U_n(R|+\infty) & \equiv \inf_{h\in \hat V_n(\hat\theta_n,R|+\infty)} \|\hat{\mathbb W}_n(\hat \theta_n) + \hat {\mathbb D}_n(\hat \theta_n) [h]\|_{\hat \Sigma_n,2} \\
\hat U_n(\Theta|+\infty) & \equiv \inf_{h\in \mathbf R^{d_\theta}} \|\hat{\mathbb W}_n (\hat \theta_n^{\text{u}}) + \hat {\mathbb D}_n(\hat \theta_n^{\text{u}})  [h]\|_{\hat \Sigma_n,2}.
\end{align*}
The validity of these distributional approximations follows from Theorem \ref{th:coupsmooth}.

\begin{theorem}\label{th:gmmboot}
Let Assumptions \ref{ass:gmmbasic}, \ref{ass:gmmrho}, \ref{ass:gmmid}, \ref{ass:gmmsigma}, \ref{ass:gmmbootp1}, and \ref{ass:gmmbootp2} hold, set $a_n = \sqrt{\log(n)}/n^{\frac{1}{10 + 5d_\theta}}$, and let $n^{-1/2} = o(r_n)$.
Then: there are sequences $\ell_n,\ell_n^{\text{\rm u}} \downarrow 0$ satisfying $(\ell_n\vee\ell_n^{\text{\rm u}})^2 \sqrt{\log(1/(\ell_n\vee \ell_n^{\text{\rm u}}))} = o(a_n n^{-\frac{1}{2}})$, $\ell_n  = o(r_n)$, and $n^{-\frac{1}{2}} = o(\ell_n\wedge \ell_n^{\text{\rm u}})$ for which it follows uniformly in $P\in \mathbf P_0$ that
\begin{align*}
\hat U_n(R|+\infty)  & \geq \Up^\star(R|\ell_n) + o_P(a_n) \\
\hat U_n(R|+\infty) - \hat U_n(\Theta|+\infty) & \geq \UpS (R|\ell_n) - \UpS (\Theta|\ell_n^{\text{\rm u}}) + o_P(a_n).
\end{align*}
\end{theorem}


Crucially, note that any sequences $\ell_n$ and $\ell_n^{\text{\rm u}}$ satisfying the conditions of Theorem \ref{th:gmmboot} also satisfy the conditions of Theorem \ref{th:gmmapprox}.
Therefore, Theorems \ref{th:gmmboot} and \ref{th:gmmapprox} together establish the validity of employing the laws of $\hat U_n(R|+\infty)$ and $\hat U_n(\Theta|+\infty)$ conditional on the data to approximate the laws of $I_n(R)$ and $I_n(\Theta)$.
In particular, for a level $\alpha$ test we may compare the test statistic $I_n(R)$ to the critical value
$$\hat q_{1-\alpha}(\hat U_n(R|+\infty)) \equiv \inf \{c : P(\hat U_n(R|+\infty) \leq c |\{X_i\}_{i=1}^n) \geq 1-\alpha\}.$$
Similarly, for the re-centered statistic $I_n(R) - I_n(\Theta)$, valid critical values are given by:
\begin{multline*}
\hat q_{1-\alpha}(\hat U_n(R|+\infty) - \hat U_n(\Theta|+\infty)) \\ \equiv \inf \{c : P(\hat U_n(R|+\infty) - \hat U_n(\Theta|+\infty) \leq c |\{X_i\}_{i=1}^n) \geq 1-\alpha\}.
\end{multline*}

These approximations are valid under the requirement that $r_n$ satisfy $r_n \sqrt n \to \infty$.
Intuitively, the bandwidth $r_n$ is meant to reflect a bound on the distance between $\hat \theta_n$ and $\IDpoint$.
For a data driven choice of $r_n$ we may therefore employ a bootstrap estimate of an upper quantile of the distribution of the \emph{unconstrained} estimator.
Specifically, for $\hat \theta_n^{\rm u \star}$ the bootstrapped version of $\hat \theta_n^{\text{u}}$, we may set $\hat r_n$ to be given by
$$\hat r_n \equiv\inf \{c: P(\|\hat \theta_n^{\rm u \star} - \hat \theta_n^{\text{u}} \|_2 \leq c|\{X_i\}_{i=1}^n)\} \geq 1-\gamma_n$$
for $\gamma_n \to 0$ as the sample size $n$ tends to infinity, and employ $\hat r_n$ in place of $r_n$.

%% file: Appendix/ExamplesMain/ExNPIV.tex

\subsection{Consumer Demand}\label{sec:exconsumer}

We base our next example on a long-standing literature aiming to replace parametric assumptions with shape restrictions implied by economic theory \citep{matzkin:1994}.
Specifically, suppose that quantity demanded by individual $i$, denoted $Q_i$, satisfies
$$Q_i = g_0(S_i, Y_i) + W_i^\prime \gamma_0 + U_i,$$
where $S_i \in \mathbf R_+$ denotes price, $Y_i \in \mathbf R_+$ denotes income, and $W_i \in \mathbf R^{d_w}$ is a set of covariates.
In addition, we assume there is an instrument $Z_i$ yielding the restriction
\begin{equation}\label{eq:appslutsky0}
E_P[Q - g_0(S,Y) - W^\prime \gamma_0 |Z] = 0.
\end{equation}
For instance, under exogeneity of prices we may let $Z = (S,Y,W^\prime)^\prime$ as in \cite{blundell:horowitz:parey:2012}.
Alternatively, if there is a concern that prices are endogenous, then we may set $Z = (I,Y,W^\prime)^\prime$ for $I$ an instrument for $S$, as in \cite{blundell2017nonparametric}.

Our goal is to conduct inference on the level of demand at particular price income pair $(s_0,y_0)$ while imposing that the function $g_0$ satisfies the Slutsky restriction
\begin{equation}\label{eq:appslutsky1}
\frac{\partial}{\partial s} g_0(s,y) + g_0(s,y)\frac{\partial}{\partial y} g_0(s,y) \leq 0 .
\end{equation}
To map this problem into our framework, we assume that for some set $\Omega$, $(S,Y)\in \Omega \subseteq \mathbf R_+^2$ with probability one for all $P\in \mathbf P$ and impose that $g_0 \in C^1_B(\Omega)$, where
$$C^{m}_B(\Omega) \equiv \{g: \Omega \to \mathbf R \text{ s.t. } \|g\|_{m,\infty} < \infty\} \hspace{0.3 in} \|g\|_{m,\infty} \equiv\sup_{0\leq \alpha \leq m} \sup_{(s,y)\in \Omega} |\nabla^\alpha g(s,y)|.$$
Since $\theta_0 \equiv (g_0,\gamma_0)$ with $\gamma_0 \in \mathbf R^{d_w}$, we set $\mathbf B = C^1_B(\Omega)\times \mathbf R^{d_w}$ and for any $(g,\gamma) = \theta \in \mathbf B$ let $\|\theta\|_{\mathbf B} = \max\{\|g\|_{1,\infty},\|\gamma\|_2\}$.
We also note that $X = (Q,S,Y,W)$ and
\begin{equation}\label{eq:appslutsky1p5}
\rho(X,\theta) = Q - g(S,Y) - W^\prime \gamma.
\end{equation}
We will assume $\theta_0 \equiv (g_0,\gamma_0)$ is identified by \eqref{eq:appslutsky0}. 
Hence, we can think of $\theta_0$ as a function of $P$ through \eqref{eq:appslutsky0}, though we leave such dependence implicit in the notation.

In order to impose the Slutsky restriction in \eqref{eq:appslutsky1} we let $\mathbf G = C^0_B(\Omega)$ and $\|\cdot\|_{\mathbf G} = \|\cdot\|_\infty$, where with some abuse of notation we write $\|\cdot\|_\infty$ in place of $\|\cdot\|_{0,\infty}$.
The space $C^0_B(\Omega)$ is a Banach lattice under the standard pointwise ordering given by
\begin{equation}\label{eq:appslutsky2}
a \leq b \text{ if and only if } a(s,y) \leq b(s,y) \text{ for all } (s,y) \in \Omega
\end{equation}
for any $a,b\in C^0_B(\Omega)$. 
The constant function ${ \text{\bf c}}\in C_B^0(\Omega)$ satisfying $\text{\bf c}(s,y) = 1$ for all $(s,y)\in \Omega$ is an order unit under the partial ordering in \eqref{eq:appslutsky2}.
Its induced norm is
$$\{\inf \lambda > 0 : |a| \leq \lambda \text{\bf c}\} = \sup_{(s,y)\in \Omega} |a(s,y)|,$$
which coincides with the norm $\|\cdot\|_\infty$ on $C_B^0(\Omega)$, and we therefore set $\mathbf {1_G} = \text{\bf c}$.
To encode the Slutsky restriction in \eqref{eq:appslutsky1} we then let the map $\Upsilon_G : \mathbf B \to \mathbf G$ equal
\begin{equation}\label{eq:appslutsky3}
\Upsilon_G(\theta)(s,y) = \frac{\partial}{\partial s} g(s,y) + g(s,y)\frac{\partial}{\partial y} g(s,y)
\end{equation}
for any $\theta = (g,\gamma)\in \mathbf B$.
Finally, to test whether the level of demand at a prescribed price $s_0$ and income $y_0$ equals a hypothesized value $c_0$, we set $\mathbf F = \mathbf R$, $\|\cdot\|_{\mathbf F} = |\cdot|$, and
\begin{equation}\label{eq:appslutsky4}
\Upsilon_F(\theta) = g(s_0,y_0) - c_0
\end{equation}
for any $\theta = (g,\gamma) \in \mathbf B$.
By setting $R = \{\theta \in \mathbf B : \Upsilon_G(\theta) \leq 0 \text{ and } \Upsilon_F(\theta) = 0\}$ and conducting test inversion (over different values of $c_0$) of the null hypothesis
$$H_0: \IDpoint \in R \hspace{0.5 in } H_1: \IDpoint \notin R$$
we may obtain a confidence region for the level of demand at price $s_0$ and income $y_0$.

We set the parameter space to be a ball in $\mathbf B$ under $\|\cdot\|_{\mathbf B}$ by letting $\Theta$ be equal to
\begin{equation}\label{eq:slutskyTheta}
\Theta \equiv \{(g,\gamma) \in C^1_B(\Omega) \times \mathbf R^{d_w}: \|g\|_{1,\infty} \leq C_0 \text{ and } \|\gamma\|_2 \leq C_0\}
\end{equation}
for some $C_0 < \infty$.
Given a sequence of approximating functions $\{p_{j}\}_{j=1}^{j_n}$, we then let $p^{j_n}(s,y) \equiv (p_{1}(s,y),\ldots, p_{j_n}(s,y))^\prime$ and set the sieve $\Theta_n$ to equal
$$\Theta_n \equiv \{(p^{j_n\prime}\beta, \gamma) : \|p^{j_n\prime}\beta\|_{1,\infty} \leq C_0 \text{ and } \|\gamma\|_2 \leq C_0\}.$$
Similarly, for a sequence $\{q_{k}\}_{k=1}^{k_n}$ of transformations of the conditioning variable $Z$, we let $q^{k_n}(z) \equiv (q_{1}(z),\ldots, q_{k_n}(z))^\prime$.
We base our test statistic on the quadratic forms
$$Q_n(\theta) \equiv \|\hat  \Sigma_n \{\frac{1}{n}\sum_{i=1}^n (Q_i - g(S_i,Y_i)-W_i^\prime \gamma)q^{k_n}(Z_i)\}\|_2$$
for some $k_n\times k_n$ weighting matrix $\hat \Sigma_n$ and every $(g,\gamma) = \theta \in \Theta$.
The statistics $I_n(R)$ and  $I_n(\Theta)$ simply equal the minimums of $\sqrt n Q_n(\theta)$ over $\Theta_n \cap R$ and $\Theta_n$ respectively.

The next assumptions suffice for obtaining a strong approximation.
In their statement, the notation $\underline{\text{sing}}\{A\}$ denotes the smallest singular value of a matrix $A$.

\begin{assumption}\label{ass:slutskyproc}
(i) $\{X_i,Z_i\}_{i=1}^n$ is i.i.d. with $(X,Z)$ distributed according to $P\in \mathbf P$;
(ii) For $\Theta$ as in \eqref{eq:slutskyTheta} and each $P\in \mathbf P_0$ there exists a unique $\IDpoint \in \Theta$ satisfying $E_P[\rho(X,\IDpoint)|Z] = 0$;
(iii) The support of $(Q,W)$ is bounded uniformly in $P\in \mathbf P$.
\end{assumption}

\begin{assumption}\label{ass:slutskysieve}
(i) $\sup_{(s,y)} \|p^{j_n}(s,y)\|_2 \lesssim \sqrt {j_n}$;
(ii) $\sup_{(s,y)} \|\partial_a p^{j_n}(s,y)\|_2 \lesssim j_n^{3/2}$ for $a\in\{s,y\}$;
(iii) The eigenvalues of $E_P[p^{j_n}(S,Y)p^{j_n}(S,Y)^\prime]$ are bounded away from zero and infinity uniformly in $P\in \mathbf P$ and $j_n$;
(iv) For each $P\in \mathbf P_0$ there is a $\Pi_n\IDpoint = (g_n,\gamma_0) \in \Theta_n \cap R$ with $\sup_{P\in \mathbf P_0}\|E_P[(g_0(S,Y) - g_n(S,Y))q^{k_n}(Z)\|_{2} = o((n\log(n))^{-1/2})$.
\end{assumption}

\begin{assumption}\label{ass:slutskymoments}
(i) $\max_{1\leq k \leq k_n} \|q_{k}\|_\infty \lesssim \sqrt{k_n}$;
(ii) $E_P[q^{k_n}(Z)q^{k_n}(Z)^\prime]$ has eigenvalues bounded uniformly in $P\in \mathbf P$, $k_n$;
(iii) $\text{\rm s}_n \equiv  \inf_{P\in \mathbf P} \text{\rm \underline{sing}}\{E_P[q^{k_n}(Z)(p^{j_n}(S,Y)^\prime ~ W^\prime)]\}$ satisfies $0 < \text{\rm s}_n = O(1)$;
(iv) $j_n^2 k_n^3 \log^3(n) = o(n)$ and $k_n^2 j_n\log^{3/2}(1+k_n)/(\text{\rm s}_n\sqrt{n})(1\vee \sqrt{\log(\text{\rm s}_n \sqrt n/k_n)}) = o((\log(n))^{-1/2})$.
\end{assumption}

\begin{assumption}\label{ass:slutskysigma}
(i) $\|\hat \Sigma_n - \PSigma\|_{o,2} = o_P((k_n\sqrt{j_n}\log^{3/2}(n))^{-1})$ uniformly in $P\in \mathbf P$;
(ii) $\PSigma$ is invertible and $\|\PSigma\|_{o,2}$ and $\|\PSigma^{-1}\|_{o,2}$ are bounded in $P\in \mathbf P$ and $k_n$.
\end{assumption}

Assumption \ref{ass:slutskyproc}(iii) requires $(Q,W)$ to be bounded, which enables us to apply the recent coupling results by \cite{zhai2018high}.
Alternatively, Assumption \ref{ass:slutskyproc}(iii) can be relaxed under appropriate tail conditions.
Assumptions \ref{ass:slutskysieve}(i)-(iii) are standard requirements on $\Theta_n$ that can be satisfied by, e.g., tensor product wavelets or B-splines \citep{newey:1997, chen:2006, belloni2015some, chen2018optimal}.
Assumption \ref{ass:slutskysieve}(iv) pertains the approximating requirements on the sieve; see Remarks \ref{rm:slutskybiaswell} and \ref{rm:slutskybiasill} below.
In turn, Assumption \ref{ass:slutskymoments}(i)(ii) imposes standard requirements on $\{q_{k}\}_{k=1}^{k_n}$.
Assumption \ref{ass:slutskymoments}(iii)(iv) contains the required rate conditions, which are governed by $\text{\rm s}_n$ -- a parameter that is proportional to $\nu_n^{-1}$ (as in Assumption \ref{ass:keycons}) and is closely linked the degree of ill-posedness; see Remark \ref{rm:slutskybiasill} below.
Finally, Assumption \ref{ass:slutskysigma} states the conditions on the weighting matrix $\hat \Sigma_n$.

In this application, we may set $\|\theta\|_{\mathbf E} = \sup_{P\in \mathbf P} \|g\|_{P,2} + \|\gamma\|_2$ for any $(g,\gamma) \in \Theta$.
Since in addition any $\theta = (g,\gamma) \in \Theta_n\cap R$ has the structure $g = p^{j_n\prime}\beta$, we have
\begin{align}
V_n(\theta,R|\ell) = \Big\{(p^{j_n\prime}\beta_h,\gamma_h) : ~
& \|g + \frac{p^{j_n\prime}\beta_h}{\sqrt n}\|_{1,\infty} \leq C_0 \text{ and } \|\gamma + \frac{\gamma_h}{\sqrt n}\|_2 \leq C_0 \label{Vndef1} \\
& p^{j_n}(s_0,y_0)^\prime\beta_h = 0 \label{Vndef2} \\
& \frac{\partial}{\partial s}(g + \frac{p^{j_n\prime}\beta_h}{\sqrt n}) + (g + \frac{p^{j_n\prime}\beta_h}{\sqrt n}) \frac{\partial}{\partial y} (g + \frac{p^{j_n\prime}\beta_h}{\sqrt n}) \leq 0 \label{Vndef3}\\
& \sup_{P\in \mathbf P} \|p^{j_n\prime}\beta_h\|_{P,2} + \|\gamma_h\|_2 \leq \ell \sqrt n \Big\}, \label{Vndef4}
\end{align}
where constraint \eqref{Vndef1} corresponds to $(\theta + h/\sqrt n) \in \Theta_n$, constraints \eqref{Vndef2} and \eqref{Vndef3} impose $\theta + h/\sqrt n \in R$, and constraint \eqref{Vndef4} imposes $\|h/\sqrt n\|_{\mathbf E} \leq \ell$. Similarly,
\begin{align}
V_n(\theta,\Theta|\ell) = \Big\{(p^{j_n\prime}\beta_h,\gamma_h) : ~
& \|g + \frac{p^{j_n\prime}\beta_h}{\sqrt n}\|_{1,\infty} \leq C_0 \text{ and } \|\gamma + \frac{\gamma_h}{\sqrt n}\|_2 \leq C_0 \label{VnUdef1} \\
& \sup_{P\in \mathbf P} \|p^{j_n\prime}\beta_h\|_{P,2} + \|\gamma_h\|_2 \leq \ell \sqrt n \Big\}. \label{VnUdef2}
\end{align}
Finally, recall that $\WP(\theta) \sim N(0,\text{Var}_P\{\rho(X,\theta)q^{k_n}(Z)\})$ and define $\DerP$ to be given by
$$\DerP[h] \equiv -E_P[q^{k_n}(Z)(p^{j_n}(S,Y)^\prime \beta_h + W^\prime \gamma_h)]$$
for any $h = (p^{j_n\prime}\beta_h,\gamma_h)$.
Given these definitions, note that for any $\ell_n$ we have that
\begin{align*}
\Up (R|\ell_n) & \equiv \inf_{h \in V_n(\Pi_n \IDpoint,R|\ell_n)} \|\WP(\Pi_n \IDpoint) + \DerP[h] \|_{\PSigma,2} \notag\\
\Up (\Theta|\ell_n) & \equiv \inf_{h \in V_n(\Pi_n\IDpoint,\Theta|\ell_n)} \|\WP(\Pi_n \IDpoint) + \DerP[h] \|_{\PSigma,2}.
\end{align*}

Theorem \ref{th:localdrift}(ii) immediately yields the following distributional approximations.

\begin{theorem}\label{th:slutskyapprox}
Let Assumptions \ref{ass:slutskyproc}-\ref{ass:slutskysigma} hold, and $a_n = (\log(n))^{-1/2}$.
Then: for any $\ell_n,\ell_n^{\text{\rm u}} \downarrow 0$ satisfying $k_n\sqrt{j_n\log(1+k_n)}(\ell_n\vee \ell_n^{\text{\rm u}})\sqrt{\log(\sqrt{j_n}/(\ell_n\vee \ell_n^{\text{\rm u}}))} = o(a_n)$ and $k_n\sqrt{j_n}\log(1+k_n)/\text{\rm s}_n \sqrt{n} = o(\ell_n\wedge \ell_n^{\text{\rm u}})$ it follows uniformly in $P\in \mathbf P_0$ that
\begin{align*}
I_n(R) & = \Up (R|\ell_n) + o_P(a_n) \\
I_n(R) - I_n(\Theta) & = \Up (R|\ell_n) - \Up (\Theta|\ell_n^{\text{\rm u}}) + o_P(a_n).
\end{align*}
\end{theorem}


To obtain bootstrap estimates of the distributional approximations in Theorem \ref{th:slutskyapprox} we let $\hat \theta_n$ and $\hat \theta_n^{\rm u}$ denote the minimizers of $Q_n$ over $\Theta_n \cap R$ and $\Theta_n$ respectively.
For $\rho(\cdot,\theta)$ as in \eqref{eq:appslutsky1p5}, we approximate the law of $\WP(\Pi_n \IDpoint)$ by evaluating
$$\Bemp(\theta) \equiv \frac{1}{\sqrt n}\sum_{i=1}^n \omega_i\{q^{k_n}(Z_i)\rho(X_i, \theta) - \frac{1}{n}\sum_{j=1}^n q^{k_n}(Z_j) \rho(X_j, \theta)\},$$
at $\theta = \hat \theta_n$ and $\theta = \hat \theta_n^{\text{u}}$, where $\{\omega_i\}_{i=1}^n$ is an i.i.d.\ sample independent of the data satisfying $\omega_i \sim N(0,1)$.
As our estimator for $\DerP[h]$, for any $h = (p^{j_n\prime}\beta_h,\gamma_h)$, we let
$$\hat{\mathbb D}_n[h] = -\frac{1}{n}\sum_{i=1}^n q^{k_n}(Z_i)(W_i^\prime \gamma_h + p^{j_n}(S_i,Y_i)^\prime \beta_h).$$

With regards to the local parameter space, we note that in this application Assumptions \ref{ass:locineq}(i)(ii) are satisfied with $K_g = 2$ (see Lemma \ref{lm:slutsky4ver}).
Therefore, we have
\begin{multline}
G_n(\hat \theta_n) = \Big\{h : \frac{\partial}{\partial s} p^{j_n}(s,y)^\prime(\hat \beta_n + \frac{\beta_h}{\sqrt n}) + p^{j_n}(s,y)^\prime(\hat \beta_n+\frac{\beta_h}{\sqrt n})\frac{\partial}{\partial y} p^{j_n}(s,y)^\prime(\hat \beta_n + \frac{\beta_h}{\sqrt n}) \\
\leq \max\{\frac{\partial}{\partial s} p^{j_n}(s,y)^\prime\hat \beta_n + p^{j_n}(s,y)^\prime\hat \beta_n\frac{\partial}{\partial y} p^{j_n}(s,y)^\prime \hat \beta_n - 2r_n\|\frac{p^{j_n\prime}\beta_h}{\sqrt n}\|_{1,\infty}, -r_n\}\Big\}.
\end{multline}
Moreover, because $\rho(X,\cdot)$ and $\Upsilon_F$ are linear, we may set $\ell_n = +\infty$ and obtain that
$$\hat V_n(\hat \theta_n,R|+\infty) = \{h = (p^{j_n\prime} \beta_h, \gamma_h) : h \in G_n(\hat \theta_n) \text{ and } p^{j_n}(s_0,y_0)^\prime\beta_h = 0\}.$$
Given the introduced notation, we define the statistics $\hat U_n(R|+\infty)$ and $\hat U_n(\Theta|+\infty)$ by
\begin{align*}
\hat U_n(R|+\infty) & \equiv \inf_{h \in \hat V_n(\hat \theta_n,R|+\infty)} \|\Bemp(\hat \theta_n) + \hat {\mathbb D}_n[h]\|_{\hat \Sigma_n,2} \\
\hat U_n(\Theta|+\infty) & \equiv \inf_{h = (p^{j_n\prime}\beta_h,\gamma_h) } \|\Bemp(\hat \theta_n^{\text{u}}) + \hat {\mathbb D}_n[h]\|_{\hat \Sigma_n,2}.
\end{align*}

We impose one final assumption to establish the validity of the bootstrap.
In the requirements below, it is helpful to recall $\theta_0$ is implicitly a function of $P$ through \eqref{eq:appslutsky0}.

\begin{assumption}\label{ass:slutskyboot}
(i) There is an $\epsilon > 0$ such that $\|g_0\|_{1,\infty}\vee\|\gamma_0\|_2 \leq C_0 - \epsilon$ for all $P\in \mathbf P_0$;
(ii) $\Pi_n\IDpoint = (g_n,\gamma_0) \in \Theta_n \cap R$ satisfies $\|g_n - g_0\|_{1,\infty} = o(1)$ uniformly in $P\in \mathbf P_0$;
(iii) The sequence $r_n \downarrow 0$ satisfies $k_n j_n^2\sqrt{\log(1+k_n)}/\text{\rm s}_n\sqrt{n} = o(r_n/\sqrt{\log(n)})$;
(iv) $k_n j_n^{3/4}(\mathcal E_n\vee \sqrt{\log(k_n)})\log^{1/4}(1+k_n) = o(n^{1/4}/\sqrt{\log(n)})$, where $\mathcal E_n \equiv \int_0^\infty\sqrt{\log(\epsilon,\mathcal C_n,\|\cdot\|_2)}d\epsilon$ and $\mathcal C_n \equiv \{\beta  : \|p^{j_n\prime}\beta\|_{1,\infty} \leq C_0\}$.
\end{assumption}

Assumptions \ref{ass:slutskyboot}(i)(ii) suffice for verifying Assumption \ref{ass:extra}(ii).
These requirements may be dropped at the expense of modifying $\hat V_n(\hat \theta_n,R|+\infty)$ to reflect the possible impact of $\Pi_n\IDpoint$ being ``near" the boundary of $\Theta_n$.
Assumption \ref{ass:slutskyboot}(iii) imposes the rate conditions on $r_n$.
Finally, Assumption \ref{ass:slutskyboot}(iv) controls the ``size" of the set of coefficients $\beta$ corresponding to elements $p^{j_n\prime }\beta \in \Theta_n$ and suffices for verifying the bootstrap coupling requirement of Assumption \ref{ass:bootcoupling}.
For instance, $\mathcal E_n \asymp j_n^{1/4}$ for tensor product B-splines (see Lemma \ref{lm:auxsplines}), which implies a sufficient condition for Assumption \ref{ass:slutskyboot}(iv) is that $k_n^4 j_n^4 \log^4(k_n) = o(n/\log^2(n))$.
The rate requirements for a bootstrap coupling can be weakened if the test statistic is based on the $\|\cdot\|_\infty$-norm (see Lemma \ref{lm:slutskzybootcoup}) or under additional smoothness assumptions (see Theorem \ref{th:mainbootcoup}(ii)).

Our next result characterizes the properties of the proposed bootstrap statistics.

\begin{theorem}\label{th:slutskyboot}
Let Assumptions \ref{ass:slutskyproc}, \ref{ass:slutskysieve}, \ref{ass:slutskymoments}, \ref{ass:slutskysigma}, \ref{ass:slutskyboot} hold, and $a_n = (\log(n))^{-1/2}$.
Then: there are sequences $\ell_n,\ell_n^{\text{\rm u}} \downarrow 0$ satisfying $k_nj_n^2\log(1+k_n)/\text{\rm s}_n\sqrt n = o(\ell_n\wedge \ell_n^{\text{\rm u}})$, $\ell_n = o(r_n)$, and $k_n\sqrt{j_n\log(1+k_n)} (\ell_n\vee \ell_n^{\text{\rm u}}) \sqrt{\log(\sqrt{j_n}/(\ell_n\vee\ell_n^{\text{\rm u}}))} = o(a_n)$ for which it follows that uniformly in $P\in \mathbf P_0$ we have
\begin{align*}
\hat U_n(R|+\infty) & \geq \UpS (R|\ell_n) + o_P(a_n) \\
\hat U_n(R|+\infty) - \hat U_n(\Theta|+\infty) & \geq \UpS (R|\ell_n) - \UpS (\Theta|\ell_n^{\text{\rm u}}) + o_P(a_n).
\end{align*}
\end{theorem}


Importantly, any sequences $\ell_n$ and $\ell_n^{\text{\rm u}}$ satisfying the requirements of Theorem \ref{th:slutskyboot} also satisfy the requirements of Theorem \ref{th:slutskyapprox}.
Hence, we may employ
$$\hat q_{1-\alpha}(\hat U_n(R|+\infty))  \equiv \inf \{c : P(\hat U_n(R|+\infty) \leq c | \{V_i\}_{i=1}^n) \geq 1-\alpha \}$$
as a critical value for $I_n(R)$.
Similarly, for the statistic $I_n(R) - I_n(\Theta)$ we may employ
\begin{multline*}
\hat q_{1-\alpha}(\hat U_n(R|+\infty) - \hat U_n(\Theta|+\infty)) \\ \equiv \inf \{c : P(\hat U_n(R|+\infty) - \hat U_n(\Theta|+\infty) \leq c | \{V_i\}_{i=1}^n) \geq 1-\alpha \}.
\end{multline*}

\begin{remark}\label{rm:slutskybiaswell} \rm
Suppose for notational simplicity that there are no covariates $W$ and let the marginal distribution of $(S,Y,Z)$ be constant in $P\in \mathbf P$. If $Z = (S,Y)$ (i.e.\ $(S,Y)$ is exogenous), we may set $q^{k_n}(Z) = p^{k_n}(S,Y)^\prime$ for some $k_n \geq j_n$.
The singular value $\text{s}_n$ can then be assumed to be bounded away from zero, and a sufficient condition for Assumption \ref{ass:slutskymoments}(iv) is that $k_n^4 j_n^2 \log^5(n) = o(n)$.
In order to appreciate the content of Assumption \ref{ass:slutskysieve}(iv), suppose $\{p_{j}\}_{j=1}^\infty$ is an orthonormal basis such that
$$g_0 = \sum_{j=1}^{\infty} \beta_j p_j \text{ with } |\beta_j| = O(j^{-\gamma_\beta}).$$
Setting $\Pi_n^{\text{u}} g_0 = \sum_{j=1}^{j_n} p_j \beta_j$, we obtain from a standard integral bound for a sum that
\begin{equation}\label{rm:slutskybiaswell1}
\|E_P[(g_0(S,Y) - \Pi_n^{\text{u}} g_0(S,Y))q^{k_n}(Z)]\|_2^2 \lesssim \sum_{j=j_n+1}^{k_n} \frac{1}{j^{2\gamma_\beta}} \lesssim \frac{1}{j_n^{2\gamma_\beta -1}} - \frac{1}{k_n^{2\gamma_\beta -1}}.
\end{equation}
For instance, if $k_n - j_n = O(1)$, then the bound in \eqref{rm:slutskybiaswell1} is of order $1/j_n^{2\gamma_\beta}$.
Hence, provided the approximation error by $\Pi_n^{\text{u}}g_0$ and $g_n$ (as in Assumption \ref{ass:slutskysieve}(iv)) are of the same order when $g_0 \in R$, we obtain that Assumption \ref{ass:slutskysieve}(iv) is equivalent to $\sqrt{n\log(n)}/j_n^{\gamma_\beta} = o(1)$ when $k_n - j_n = O(1)$. This approximation requirement is compatible with the condition $k_n^4 j_n^2 \log^5(n) = o(n)$ provided $\gamma_\beta > 3$. \qed
\end{remark}

\begin{remark}\label{rm:slutskybiasill} \rm
Building on Remark \ref{rm:slutskybiaswell}, suppose again there are no covariates $W$ and the marginal distribution of $(S,Y,Z)$ is constant in $P\in \mathbf P$, but now let $(S,Y)$ be endogenous.
A standard benchmark for nonparametric models with endogeneity is to assume the operator $g\mapsto E_P[g(S,Y)|Z]$ is compact, in which case there are orthonormal sequences of functions $\{\phi_j\}_{j=1}^\infty$ of $(S,Y)$ and $\{\psi_j\}_{j=1}^\infty$ of $Z$ satisfying
\begin{equation*}
E_P[\phi_j(S,Y)|Z] = \lambda_j \psi_j(Z) \hspace{0.5 in} E_P[\psi_j(Z)|S,Y] = \lambda_j \phi_j(S,Y)
\end{equation*}
where $\lambda_j >0$ tends to zero. In addition suppose $g_0$ admits for an expansion satisfying
\begin{equation*}
g_0 = \sum_{j=1}^\infty \beta_j \phi_j \text{ with } |\beta_j| = O(j^{-\gamma_\beta}),
\end{equation*}
and let $p^{j_n} = (\phi_1,\ldots, \phi_{j_n})^\prime$, $q^{k_n} = (\psi_1,\ldots, \psi_{k_n})^\prime$ with $k_n \geq j_n$ and $k_n - j_n = O(1)$, and set $\Pi_n^{\text{u}} g_0 = \sum_{j=1}^{j_n} \phi_j \beta_j$.
Provided the approximation error of $\Pi_n^{\text{u}} g_0$ and $g_n$ (as in Assumption \ref{ass:slutskysieve}(iv)) are of the same order when $g_0 \in R$, we then obtain
$$\|E_P[(g_0(S,Y) - g_n(S,Y))q^{k_n}(Z)]\|_2 \lesssim \frac{\lambda_{j_n}}{j_n^{\gamma_\beta}}.$$
Moreover, direct calculation shows $\text{s}_n$, which is proportional to $\nu_n^{-1}$ as in Assumption \ref{ass:keycons}, satisfies $\text{s}_n = \lambda_{j_n}$ and hence equals the reciprocal of the sieve measure of ill-posedness \citep{blundell:chen:kristensen:2007}.
It follows that if $\lambda_{j} \asymp j^{-\gamma_\lambda}$, and $\gamma_\beta > 3$, then Assumptions \ref{ass:slutskysieve}(iv) and \ref{ass:slutskymoments}(iv) can be satisfied by setting $j_n \asymp n^\kappa$ with $(\gamma_\lambda + \gamma_\beta)^{-1} < 2\kappa < (3 + \gamma_\lambda)^{-1}$ and $k_n - j_n = O(1)$.
Alternatively, if $\lambda_j = \exp\{-\gamma_\lambda j\}$, then Assumption \ref{ass:slutskysieve}(iv) and \ref{ass:slutskymoments}(iv) can be satisfied when $\gamma_\beta > 4$ by setting, for example, $j_n = (\log(n) - \kappa\log(\log(n)))/2\gamma_\lambda$ with $7 < \kappa < 2\gamma_\beta -1$ and $k_n - j_n = O(1)$. \qed
\end{remark}

%% file: Appendix/ExamplesMain/ExQuant.tex

\subsection{Quantile Treatment Effects} \label{sec:exquant}
 
For our next example, we study a nonparametric quantile treatment effect (QTE) model.
Specifically, for an outcome $Y\in \mathbf R$, treatment $D\in [0,1]$, instrument $Z \in \mathbf R$, and quantile $\tau \in (0,1)$, we assume the parameter of interest $\IDpoint$ satisfies
\begin{equation}\label{ex:quant1}
P(Y\leq \IDpoint(D)|Z) = \tau.
\end{equation}
If $D$ is randomly assigned, then we may set $D = Z$ and interpret $\nabla \IDpoint$ as the $\tau^{th}$ quantile treatment effect (QTE).
Alternatively, if $D\neq Z$, then we obtain the QTE model of \cite{chernozhukov:hansen:2005}.
To map \eqref{ex:quant1} into our framework, we set 
\begin{equation}\label{ex:quant6}
\rho(X,\theta) = 1\{Y \leq \theta(D)\} - \tau,
\end{equation}
where $X = (Y,D)\in \mathbf X \equiv \mathbf R\times [0,1]$.
In order to illustrate our conditions in a number of different settings, we focus on conducting inference on a nonlinear function of $\IDpoint$.
Specifically, we conduct inference on the variance of the quantile treatment effects: 
\begin{equation*}\label{ex:quant2}
\int_0^1 (\nabla \IDpoint(u))^2 du - (\int_0^1 \nabla \IDpoint(u) du)^2
\end{equation*}
while imposing that the QTE be increasing in treatment intensity (i.e.\ $d\mapsto \nabla \IDpoint(d)$ is increasing).
To map this problem into our framework we define
\begin{equation*}\label{ex:quant3}
C^m_B([0,1])\equiv  \{\theta :[0,1] \to \mathbf R \text{ s.t. } \|\theta\|_{m,\infty} < \infty\} \hspace{0.3 in} \|\theta\|_{m,\infty} \equiv \sup_{0\leq \alpha \leq m} \sup_{d\in [0,1]} |\nabla^\alpha \theta(d)|,
\end{equation*}
and set $\mathbf B = C^2_B([0,1])$ and $\|\cdot\|_{\mathbf B} = \|\cdot\|_{2,\infty}$.
We impose the restriction that the quantile treatment effect be increasing in the intensity of treatment by letting $\mathbf G = C^0_B([0,1])$, $\|\cdot\|_{\mathbf G} = \|\cdot\|_{\infty}$ (where we write $\|\cdot\|_\infty$ in place of $\|\cdot\|_{0,\infty}$), and defining
\begin{equation}\label{ex:quant4}
\Upsilon_G(\theta) \equiv -\nabla^2 \theta.
\end{equation}
As shown in Section \ref{sec:exconsumer}, $\mathbf G$ is a lattice with order unit $\mathbf {1_G} = \mathbf c$ for $\mathbf c$ the constant function $\mathbf c(d) = 1$ for all $d\in [0,1]$.
Setting $\mathbf F = \mathbf R$ with $\|\cdot\|_{\mathbf F} = |\cdot|$, we test whether the variance of the quantile treatment effects equals a hypothesized value $\lambda \neq 0$ by setting
\begin{equation}\label{ex:quant5}
\Upsilon_F(\theta) = \int_0^1 (\nabla \theta(u))^2 du - (\int_0^1 \nabla \theta(u)du)^2 - \lambda.
\end{equation}

For the parameter space for $\IDpoint$ we employ a ball in $\mathbf B$ and we thus set $\Theta$ to equal
\begin{equation}\label{eq:quantTheta}
\Theta \equiv \{ \theta \in C^2_B([0,1]) \text{ s.t. } \|\theta\|_{2,\infty} \leq C_0\}
\end{equation}
for some $C_0 < \infty$.
For a sequence of approximating functions $\{p_{j}\}_{j=1}^{j_n}$ defined on $[0,1]$ we then let $p^{j_n}(d) \equiv (p_{1}(d),\ldots, p_{j_n}(d))^\prime$ and define $\Theta_n$ to equal
\begin{equation}\label{eq:quantThetaN}
\Theta_n \equiv \{ p^{j_n\prime} \beta \in C^2_B([0,1]): \|p^{j_n\prime} \beta\|_{2,\infty} \leq C_0\}.
\end{equation}
Similarly for a sequence $\{q_{k}\}_{k=1}^{k_n}$, we set $q^{k_n}(z) \equiv (q_{1}(z),\ldots, q_{k_n}(z))^\prime$ and define
$$Q_n(\theta) \equiv \|\hat \Sigma_n \{\frac{1}{n}\sum_{i=1}^n (1\{Y_i \leq \theta(D_i)\} - \tau)q^{k_n}(Z_i)\}\|_p$$
for some $2\leq p \leq \infty$ and weighting matrix $\hat \Sigma_n$.
The statistics $I_n(R)$ and $I_n(\Theta)$ then equal the minimums of $\sqrt nQ_n$ over $\Theta_n \cap R$ and $\Theta_n$ respectively. 

In what follows, we will assume for simplicity that $\IDpoint$ is identified.
As a result, we may think of $\IDpoint$ as a function of $P$ through \eqref{ex:quant1}, though we leave such dependence implicit in the notation.
We next impose the following assumptions:

\begin{assumption}\label{ass:quantproc}
(i) $\{Y_i,D_i,Z_i\}_{i=1}^n$ is i.i.d. with $(Y,D,Z) \in \mathbf R\times [0,1]\times \mathbf R$ distributed according to $P\in \mathbf P$;
(ii) For $\Theta$ as in \eqref{eq:quantTheta} and each $P\in \mathbf P_0$ there exists a unique $\IDpoint \in \Theta$ satisfying \eqref{ex:quant1};
(iii) The distribution of $Y$ conditional on $(D,Z)$ is absolutely continuous with density $f_{Y|DZ,P}(\cdot|D,Z)$ that is bounded and Lipschitz uniformly in $(D,Z)$ and $P\in \mathbf P$;
(iv) Assumptions \ref{ass:coupreg} and \ref{ass:supportreg} hold.
\end{assumption}

\begin{assumption}\label{ass:quantsieve}
(i) $\sup_{d} \|p^{j_n}(d)\|_2 \lesssim \sqrt {j_n}$;
(ii) $E_P[p^{j_n}(D)p^{j_n}(D)^\prime]$ has eigenvalues bounded away from zero and infinity uniformly in $P\in \mathbf P$ and $j_n$;
(iii) For each $P\in \mathbf P_0$ there is a $\Pi_n \IDpoint \in \Theta_n \cap R$ satisfying $\sup_{P\in \mathbf P_0}\|E_P[(1\{Y\leq \Pi_n \IDpoint(D)\} -1\{Y \leq \IDpoint(D)\})q^{k_n}(Z)]\|_{p} = O((n\log(n))^{-1/2})$ and $\sup_{P\in \mathbf P_0}\|\IDpoint - \Pi_n\IDpoint\|_{1,\infty} = o(1)$.
\end{assumption}

\begin{assumption}\label{ass:quantid}
(i) $\inf_{P\in \mathbf P_0}\inf_{\theta\in \Theta : \|\theta-\IDpoint\|_{1,\infty} \geq \epsilon} E_P[(P(Y\leq \theta(D)|Z)-\tau)^2] > 0$  for every $\epsilon > 0$;
(ii) There are $\epsilon$ and $\text{\rm s}_n > 0$ satisfying for all $P\in \mathbf P_0$ and $\|\theta-\Pi_n \IDpoint\|_{1,\infty} \leq \epsilon$, $\text{\rm s}_n \leq \underline{\text{\rm sing}}\{E_P[f_{Y|D,Z}(\theta(D)|D,Z)q^{k_n}(Z)p^{j_n}(D)^\prime]\}$ and $\text{\rm s}_n = O(1)$.
\end{assumption}

\begin{assumption}\label{ass:quantmoments}
(i) $\max_{1\leq k \leq k_n} \|q_{k}\|_\infty = O(1)$;
(ii) $\max_{1\leq k \leq k_n} \|q_{k}\|_{1,\infty} = O(k_n)$;
(iii) $E_P[q^{k_n}(Z)q^{k_n}(Z)^\prime]$ has eigenvalues bounded away from zero and infinity uniformly in $P\in \mathbf P$ and $k_n$;
(iv) For each $\theta \in \Theta$ there is a $\pi_n(\theta) \in \mathbf R^{k_n}$ with $E_P[(E_P[\rho(X,\theta)|Z] - q^{k_n}(Z)^\prime \pi_n(\theta))^2] = o(1)$ uniformly in $P\in \mathbf P$ and $\theta \in \Theta$;
(v) $k_n^{1/p}\sqrt{j_n}\log^{3/2}(n)(n^{1/6}\vee k_n)/n^{1/3} = o(1)$ and $j_n\log^{3/2}(1+k_n) k_n^{2/p + 1/2}/\text{\rm s}_n\sqrt n = o((\log(n))^{-2})$.
\end{assumption}

\begin{assumption}\label{ass:quantsigma}
(i) $\|\hat \Sigma_n - \PSigma\|_{o,p} = o_P((k_n^{1/p}\log(n))^{-1})$ uniformly in $P\in \mathbf P$;
(ii) $\PSigma$ is invertible and $\|\PSigma\|_{o,p}$ and $\|\PSigma^{-1}\|_{o,p}$ are bounded in $P\in \mathbf P$ and $k_n$.
\end{assumption}

Assumption \ref{ass:quantproc} imposes regularity conditions on the distribution $P$ that enable us to apply the empirical process coupling results of Appendix \ref{sec:unifcoupling}.
Assumption \ref{ass:quantsieve} states the requirements on $\Theta_n$, including demanding an asymptotically negligible bias in Assumption \ref{ass:quantsieve}(iii).
Assumption \ref{ass:quantid}(i) holds pointwise in $P\in \mathbf P_0$ due to $\Theta$ being compact under $\|\cdot\|_{1,\infty}$, and hence the uniformity in $P\in \mathbf P_0$ demanded by Assumption \ref{ass:quantid}(i) corresponds to imposing strong identification.
Assumption \ref{ass:quantid}(ii) enables us to obtain a uniform rate of convergence under $\|\cdot\|_{\mathbf E} = \sup_{P\in \mathbf P} \|\cdot\|_{P,2}$.
As in Section \ref{sec:exconsumer}, $\text{s}_n$ can be shown to be related to the degree of ill-posedness.
Assumptions \ref{ass:quantmoments}(i)-(iv) impose conditions on $\{q_{k}\}_{k=1}^{k_n}$ including that they be bounded -- this requirement can be relaxed at the cost of more stringent rate restrictions to ensure a coupling of the empirical process (see Lemma \ref{lm:quantcoup}).
Finally, Assumption \ref{ass:quantmoments}(v) states our rate restrictions, which we note are easier to satisfy for higher values of $p$.

For any $\theta = p^{j_n\prime}\beta \in \Theta_n \cap R$, in this application the local parameter space equals
\begin{align}\label{ex:quantrestdef}
V_n(\theta,R|\ell) = \Big\{ h = p^{j_n\prime}\beta_h :&  \|\theta + \frac{h}{\sqrt n}\|_{2,\infty}\leq C_0,~ \sup_{P\in \mathbf P} \|h\|_{P,2}\leq \ell \sqrt n, \notag \\
& \int_0^1(\nabla \theta(u) + \frac{\nabla h(u)}{\sqrt n})^2du - (\int_0^1 \{\nabla \theta(u) +  \frac{\nabla h(u)}{\sqrt n}\}du)^2 = \lambda, \notag \\
& - \nabla^2 \theta(d) - \frac{\nabla^2 h(d)}{\sqrt n} \leq 0  \text{ for all } d\in [0,1]\Big\},
\end{align}
where the first two constraints impose that $\theta + h/\sqrt n \in \Theta_n$ and $\|h/\sqrt n \|_{\mathbf E} \leq \ell$, while the final two constraints require that $\theta+ h/\sqrt n \in R$.
Similarly, here
$$V_n(\theta,\Theta|\ell) = \Big\{ h = p^{j_n\prime}\beta_h : \|\theta + \frac{h}{\sqrt n}\|_{2,\infty}\leq C_0 \text{ and } \sup_{P\in \mathbf P} \|h\|_{P,2}\leq \ell \sqrt n \Big\}.$$
Also recall that $\WP(\theta) \sim N(0,\text{Var}_P\{\rho(X,\theta)q^{k_n}(Z)\})$ and for any $h = p^{j_n\prime}\beta_h $ define 
\begin{equation}\label{ex:quatderdef}
\DerP(\theta)[h] \equiv E_P[q^{k_n}(Z)f_{Y|DZ,P}(\theta(D)|D,Z)p^{j_n}(D)^\prime \beta_h].
\end{equation}
The random variables to which $I_n(R)$ and $I_n(\Theta)$ will be coupled are then given by
\begin{align*}
\Up (R|\ell_n) & \equiv \inf_{h \in V_n(\Pi_n \IDpoint,R|\ell_n)} \|\WP(\Pi_n \IDpoint) + \DerP(\Pi_n\IDpoint)[h] \|_{\PSigma,2} \notag\\
\Up (\Theta|\ell_n) & \equiv \inf_{h \in V_n(\Pi_n\IDpoint,\Theta|\ell_n)} \|\WP(\Pi_n \IDpoint) + \DerP(\Pi_n\IDpoint)[h] \|_{\PSigma,2}.
\end{align*}

Our next result obtains distributional approximations by applying Theorem \ref{th:localdrift}.

\begin{theorem}\label{th:quantapprox}
Let Assumptions \ref{ass:quantproc}, \ref{ass:quantsieve}, \ref{ass:quantid}, \ref{ass:quantmoments}, and \ref{ass:quantsigma} hold, $a_n = (\log(n))^{-1/2}$, and $\ell_n \downarrow 0$ satisfy $k_n^{1/p}\sqrt{j_n\ell_n\log(1+k_n)\log(1/\ell_n)} = o((\log(n))^{-1/2})$ and $\ell_n^2\sqrt{nj_n\log(n)} = o(1)$.
Then: (i) Uniformly in $P\in \mathbf P_0$ it follows that
$$I_n(R) \leq \Up (R|\ell_n) + o_P(a_n).$$
(ii) If in addition $k_n\log(1+k_n)\sqrt{j_n \log(n)}/\text{\rm s}_n^2 \sqrt n = o(1)$, then for any $\ell_n^{\text{\rm u}}\downarrow 0$ satisfying $k_n^{1/p}\sqrt{j_n\ell_n^{\text{\rm u}}\log(1+k_n)\log(1/\ell_n^{\text{\rm u}})} = o((\log(n))^{-1/2})$, $(\ell_n^{\text{\rm u}})^2\sqrt{nj_n\log(n)} = o(1)$, and $\sqrt{k_n\log(1+k_n)}/\text{\rm s}_n \sqrt n = o(\ell_n^{\text{\rm u}})$, it follows uniformly in $P\in \mathbf P_0$ that
$$I_n(R) - I_n(\Theta) \leq \Up (R|\ell_n) - \Up (\Theta|\ell_n^{\text{\rm u}}) + o_P(a_n).$$
\end{theorem}


Theorem \ref{th:quantapprox}(i) obtains an upper bound for $I_n(R)$ by relying on Theorem \ref{th:localdrift}(i).
In order to approximate the re-centered statistic $I_n(R) - I_n(\Theta)$, we cannot rely on an upper bound for $I_n(\Theta)$ as the resulting approximation could fail to control size.
Therefore, Theorem \ref{th:quantapprox}(ii) instead relies on Theorem \ref{th:localdrift}(ii).
Applying Theorem \ref{th:localdrift}(ii), however, requires an additional rate condition in order to establish the linearization of the moment conditions is asymptotically valid.
We also note that the conclusion of Theorem \ref{th:quantapprox}(ii) in fact holds with equality if $\ell_n$ satisfies the same rate restrictions as $\ell_n^{\text{\rm u}}$.

In order to provide bootstrap estimates for these distributional approximations, we let $\hat \theta_n$ and $\hat \theta_n^{\rm u}$ denote minimizers of $Q_n$ over $\Theta_n \cap R$ and $\Theta_n$ respectively.
Our bootstrap approximation estimates the law of $\WP(\IDpoint)$ and the derivative $\DerP(\IDpoint)$ by evaluating
$$\Bemp(\theta) \equiv \frac{1}{\sqrt n}\sum_{i=1}^n \omega_i \{q^{k_n}(Z_i)(1\{Y_i \leq \theta(D_i)\} - \tau) - \frac{1}{n}\sum_{j=1}^n q^{k_n}(Z_j)(1\{Y_j \leq \theta(D_j)\} - \tau)\}$$
$$\hat {\mathbb D}_n(\theta)[h] \equiv \frac{1}{\sqrt n}\sum_{i=1}^n q^{k_n}(Z_i)(1\{Y_i \leq \theta(D_i) + \frac{h(D_i)}{\sqrt n}\} - 1\{Y_i \leq \theta(D_i)\})$$
at $\hat \theta_n$ and $\hat \theta_n^{\rm u}$.
An unappealing feature of $\hat {\mathbb D}_n(\theta)$ is that it is not linear in $h$, which complicates computation.
Alternatively, a plug-in estimator based on \eqref{ex:quatderdef} could be used, though at the expense of having to estimate the density $f_{Y|DZ,P}$.

With regards to the local parameter space, we note that in this application 
$$G_n(\hat \theta_n) \equiv \{h = p^{j_n\prime}\beta_h : -\nabla^2 \hat \theta_n(d) - \frac{\nabla^2 h(d)}{\sqrt n} \leq  \max\{- \nabla^2 \hat \theta_n(d) \vee -r_n\} \text{ for all } d\in [0,1]\}.$$
Employing that $\|\cdot\|_{\mathbf B} = \|\cdot\|_{2,\infty}$ and the expression for $\Upsilon_F$ in \eqref{ex:quant5}, we  obtain that
\begin{align*}
\hat V_n(\hat \theta_n,R|\ell_n) = \Big\{ h & =  p^{j_n\prime}\beta_h  : h \in G_n(\hat \theta_n), ~ \|\frac{h}{\sqrt n}\|_{2,\infty} \leq \ell_n \notag \\
& \int_0^1 (\nabla \hat \theta_n(u) + \frac{\nabla h(u)}{\sqrt n})^2 du - (\int_0^1  (\nabla\hat \theta_n(u) + \frac{\nabla h(u)}{\sqrt n})du)^2 = \lambda \Big\},
\end{align*}
where $\ell_n$ is chosen to satisfy conditions stated below.
The bootstrap statistics $\hat U_n(R|\ell_n)$ and $\hat U_n(\Theta|+\infty)$ for approximating the distributions in Theorem \ref{th:quantapprox} are then
\begin{align*}
\hat U_n(R|\ell_n) & \equiv \inf_{h \in \hat V_n(\hat \theta_n,R|\ell_n)} \|\Bemp(\hat \theta_n) + \hat{\mathbb D}_n(\hat \theta_n)[h]\|_{\hat \Sigma_n,p} \\
\hat U_{n}(\Theta|+\infty) & \equiv \inf_{ h = p^{j_n\prime}\beta_h} \|\Bemp(\hat \theta_n^{\text{\rm u}}) + \hat{\mathbb D}_n(\hat \theta_n^{\text{\rm u}})[h]\|_{\hat \Sigma_n,p}.
\end{align*}

The following final assumption will enable us to establish bootstrap validity.
In the requirements below, it is helpful to recall $\IDpoint$ is implicitly a function of $P$ through \eqref{ex:quant1}.

\begin{assumption}\label{ass:4bootquant}
(i) The functions $\theta(d) = 1$, $\theta(d) = d^2$ are in $\mathbf B_n$;
(ii) $\|\IDpoint - \Pi_n \IDpoint\|_{2,\infty} = o(1)$ uniformly in $P\in \mathbf P_0$ and $\sup_{P\in \mathbf P_0} \|\IDpoint\|_{2,\infty} < C_0$;
(iii) $k_n$ satisfies $k_n^{1/p + 12/26} = o(n^{1/26}/\log(n))$;
(iv) $\sup_d \|\nabla^2 p^{j_n}(d)\|_2 \vee \|\nabla p^{j_n}(d)\|_2 \lesssim j_n^{5/2}$;
(v) $r_n,\ell_n$ satisfy $k_n^{1/p} \sqrt{j_n \ell_n\log(1+k_n)\log(1/\ell_n) } = o((\log(n))^{-1/2})$, $j_n^{5/2}\sqrt{k_n\log(1+k_n)}/\text{\rm s}_n \sqrt n = o(1 \wedge r_n)$, and
$ \ell_n (\sqrt {j_n n} \ell_n + j_n^{5/2}\sqrt{k_n\log(1+k_n)}/\text{\rm s}_n ) = o((\log(n))^{-1/2})$.
\end{assumption}

Assumption \ref{ass:4bootquant}(i) requires that the quadratic functions belong to $\mathbf B_n$ -- a condition that holds if quadratic functions belong to the span of $\{p_j\}_{j=1}^{j_n}$.
Assumption \ref{ass:4bootquant}(ii) implies that $\IDpoint$ and its approximation $\Pi_n\IDpoint$ belong to the interior of $\Theta$.
Assumption \ref{ass:4bootquant}(iii) enables us to verify the bootstrap coupling requirement of Assumption \ref{ass:bootcoupling} by applying the results in Appendix \ref{sec:bootcoup} to a Haar basis expansion.
While condition \ref{ass:4bootquant}(iii) suffices for verifying Assumption \ref{ass:bootcoupling} in both the endogenous ($Z \neq D$) and exogenous ($Z =D$) settings, we note that in both cases better rate conditions can be obtained.\footnote{For instance under endogeneity, a better rate could be obtained by conducting a basis expansion using the tensor product of a Haar Basis for $(Y,D)$ and the functions $\{q_{k}\}_{k=1}^{k_n}$.}
Finally, Assumption \ref{ass:4bootquant}(iv) ensures $\mathcal S_n(\mathbf B,\mathbf E) \asymp j_n^{5/2}$, while Assumption \ref{ass:4bootquant}(v) imposes the requirements on $\ell_n$ and $r_n$.

The next theorem establishes the validity of the bootstrap procedure.

\begin{theorem}\label{th:quantboot}
Let Assumptions \ref{ass:quantproc}, \ref{ass:quantsieve}, \ref{ass:quantid}, \ref{ass:quantmoments}, \ref{ass:quantsigma}, and \ref{ass:4bootquant} hold and $a_n = (\log(n))^{-1/2}$.
Then, there is a sequence $\tilde \ell_n \asymp \ell_n$ satisfying
$$\hat U_n(R|\ell_n)  \geq \UpS (R|\tilde \ell_n) + o_P(a_n)$$
uniformly in $P\in \mathbf P_0$.
(ii)  If in addition $k_n\log(1+k_n)\sqrt{j_n \log(n)}/\text{\rm s}_n^2 \sqrt n = o(1)$, then for any $\tilde \ell_n^{\rm u}$ satisfying the conditions of Theorem \ref{th:quantapprox}(ii) we have uniformly in $P\in \mathbf P_0$
$$\hat U_n(R|\ell_n) - \hat U_n(\Theta|+\infty)  \geq \UpS (R|\tilde \ell_n) - \UpS (\Theta|\tilde \ell_n^{\text{\rm u}}) + o_P(a_n).$$
\end{theorem}


Theorems \ref{th:quantapprox}(i) and \ref{th:quantboot}(i) imply that as critical value for $I_n(R)$ we may employ
$$\hat q_{1-\alpha}(\hat U_n(R|\ell_n)) \equiv \inf\{c: P(\hat U_n(R|\ell_n) \leq c |\{V_i\}_{i=1}^n) \geq 1-\alpha \}.$$
If in addition $k_n\log(1+k_n)\sqrt{j_n \log(n)}/\text{\rm s}_n^2 \sqrt n = o(1)$, then Theorems \ref{th:quantapprox}(ii) and \ref{th:quantboot}(ii) imply a valid test can be obtained by rejecting whenever $I_n(R) - I_n(\Theta)$ exceeds 
$$\hat q_{1-\alpha}(\hat U_n(R|\ell_n) - \hat U_n(\Theta|+\infty)) \equiv \inf\{c: P(\hat U_n(R|\ell_n) - \hat U_n(\Theta|+\infty)\leq c |\{V_i\}_{i=1}^n) \geq 1-\alpha \}.$$

Our critical values depend on the choices of $r_n$ and $\ell_n$.
The slackness parameter $r_n$ again measures sampling uncertainty in whether constraints ``bind." 
Following the discussion in Section \ref{sec:testexiv}, for $\hat \theta_n^{\rm u\star}$ a ``bootstrap" analogue to $\hat \theta_n^{\rm u}$, we may thus set
\begin{equation}\label{eq:refrl1}
P(\max_{d\in [0,1]} \nabla^2 \hat \theta_n^{\rm u}(d) - \nabla^2 \hat \theta_n^{\rm u\star}(d) \leq r_n | \{V_i\}_{i=1}^n) = 1-\gamma_n
\end{equation}
with $\gamma_n \to 0$.
With regards to $\ell_n$, we note that its main role in this application is to ensure that $\hat V_n(\hat \theta_n,R|\ell_n)$ is well approximated by the true local parameter space despite the nonlinearity of $\Upsilon_F$.
To this end, the requirements on $\ell_n$ imposed in Assumption \ref{th:quantboot} ensure $\sqrt n \ell_n \|\hat \theta_n - \Pi_n\IDpoint\|_{\mathbf B} = o_P(a_n)$ uniformly in $P\in \mathbf P_0$.
Since $\|\cdot\|_{\mathbf B} = \|\cdot\|_{2,\infty}$ in this application, we may select $\ell_n$ in a data driven way by setting it to satisfy
\begin{equation}\label{eq:refrl2}
P(\max_{d\in [0,1]} |\nabla^2 \hat \theta_n^{\rm u}(d) - \nabla^2 \hat \theta_n^{\rm u\star}(d)| \leq \frac{1}{\sqrt n \ell_n} |\{V_i\}_{i=1}^n) = 1-\gamma_n
\end{equation}
for some $\gamma_n\to 0$. 
While we set $\gamma_n$ in \eqref{eq:refrl1} and \eqref{eq:refrl2} to be the same, it is worth noting they could be different.
In fact, $r_n$ and $\ell_n$ do not ``interact" in the requirements of Assumption \ref{ass:4bootquant}(v) and, in this sense, can be set independently. 
We also note that in settings in which the rate of convergence is sufficiently fast, \eqref{eq:refrl2} should deliver a ``large" $\ell_n$ in the sense that $\hat U_n(R|\ell_n)$ and $\hat U_n(R|+\infty)$ are asymptotically equivalent. 
Moreover, in applications in which we expect the rate of convergence of $\hat \theta_n$ to be sufficiently fast, we may directly set $\ell_n = +\infty$; see Lemma \ref{lm:lnotbind}.

\begin{remark} \rm \label{ex:quantprim}
To illustrate the role of $\ell_n$, it is helpful to conduct a pointwise (in $P$) analysis, set $p=2$, and connect our assumptions to the literature on estimation of conditional moment restriction models \citep{chen:pouzo:2012}.
We follow the literature in imposing a local curvature assumption, which in our application, corresponds to
\begin{multline}\label{ex:quantprim1}
\|E_P[(P(Y \leq h(D)|Z) - \tau)q^{k_n}(Z)]\|_2 \\ \asymp \|E_P[f_{Y|DZ,P}(\bar \theta(D)|D,Z)(\IDpoint(D) - h(D))q^{k_n}(Z)]\|_2
\end{multline}
for all $h\in \Theta_n$ and $\bar \theta \in \Theta$ that are in a neighborhood of $\IDpoint$. 
We further suppose the linear operator $h\mapsto E_P[f_{Y|DZ,P}(\IDpoint(D)|D,Z)h(D)|Z]$ is compact, in which case there exist orthonormal bases $\{\psi_j\}$ and $\{\phi_k\}$ and a sequence $\lambda_j \downarrow 0$ satisfying
\begin{equation}\label{ex:quantprim2}
E_P[f_{Y|DZ,P}(\IDpoint(D)|D,Z)\phi_j(D)|Z] = \lambda_j \psi_j(Z).
\end{equation}
Setting $k_n \geq j_n$ with $k_n - j_n = O(1)$, $p^{j_n} = (\phi_1,\ldots, \phi_{j_n})^\prime$, $q^{k_n} = (\psi_1,\ldots, \psi_{k_n})^\prime$, and $\Pi_n^{\rm u} \IDpoint = \sum_{j=1}^{j_n} \phi_j \beta_j$, we also suppose $\IDpoint$ admits an expansion
\begin{equation}\label{ex:quantprim3}
\IDpoint = \sum_{j=1}^\infty \beta_j\phi_j \text{ with } |\beta_j| = O(j^{-\gamma_\beta}).
\end{equation}
Provided that the approximation error of $\Pi_n\IDpoint$ (as in Assumption \ref{ass:quantsieve}(iii)) and $\Pi_n^{\rm u}\IDpoint$ are of the same order, it then follows from \eqref{ex:quantprim1} and \eqref{ex:quantprim2} that
\begin{equation}\label{ex:quantprim4}
\|E_P[(1\{Y \leq \Pi_n \IDpoint(D)\} - 1\{Y \leq \IDpoint(D)\})q^{k_n}(Z)]\|_2 \lesssim \frac{\lambda_{j_n}}{j_n^{\gamma_\beta}}
\end{equation}
and ${\rm s}_n \asymp \lambda_{j_n}$ -- i.e.\ ${\rm s}_n$ is proportional to the reciprocal of the sieve measure of ill-posedness \citep{chen:pouzo:2012}.
As a result, if $\lambda_j \asymp j^{-\gamma_\lambda}$ and $\gamma_\beta > \max\{5/2,3-\gamma_\lambda\}$, then Theorem \ref{th:quantapprox} may be applied to couple $I_n(R)$ by setting $j_n \asymp n^{\kappa}$ with $(2(\gamma_\lambda+\gamma_\beta))^{-1} < \kappa < \min\{(5+2\gamma_\lambda)^{-1},1/6\}$, while coupling $I_n(R) - I_n(\Theta)$ additionally requires $\gamma_\beta > 3/2 + \gamma_\lambda$ and $\kappa < (3+4\gamma_\lambda)^{-1}$.
In contrast, in the severely ill-posed case in which $\lambda_j \asymp \exp\{- \gamma_\lambda j\}$, the conditions of Theorem \ref{th:quantapprox} for coupling $I_n(R) - I_n(\Theta)$ are not satisfied.
However, the conditions for coupling $I_n(R)$ can still be met provided $\gamma_\beta > 4$ by setting $j_n = (\log(n) - \kappa(\log(\log(n))))/2\gamma_\lambda$ with $7 < \kappa < 2\gamma_\beta -1$.
Thus, while in the ill-posed case the rate of convergence is too slow to apply Theorem \ref{th:quantapprox}(ii), Theorem \ref{th:quantapprox}(i) is still able to deliver a coupling upper bound for suitable  $\ell_n$. \qed
\end{remark}

%% file: Appendix/AppRate.tex

\section{Rate of Convergence} \label{sec:apprate}

This section contains consistency and rate of convergence results for $\hat \Theta_n^{\rm r}$.
The assumptions in the main text, which are designed to deliver a strong approximation, are stronger than needed for deriving the results in this section.
We therefore next introduce a weaker set of assumptions that suffice for obtaining a rate of convergence.
To this end, we set
\begin{equation}\label{sec:apprate:eq1}
Q_P(\theta) \equiv \|E_P[\rho(X,\theta)*q^{k_n}(Z)]\|_{\PSigma,p};
\end{equation}
i.e.\ $Q_P$ is the population analogue to the criterion function $Q_n$.
In addition, we define
\begin{align*}
\overrightarrow{d}_{H}(A,B, \|\cdot\|_{\mathbf E}) & \equiv \sup_{a\in A}\inf_{b\in B} \|a - b\|_{\mathbf E} \\ 
d_{H}(A,B,\|\cdot\|_{\mathbf E}) & \equiv \max\{\overrightarrow{d}_{H}(A,B,\|\cdot\|_{\mathbf E}),\overrightarrow{d}_{H}(B,A,\|\cdot\|_{\mathbf E})\} , 
\end{align*}
which constitute the directed Hausdorff and the Hausdorff distance (under $\|\cdot\|_{\mathbf E}$) between two sets $A$ and $B$. Given these definition, we introduce the following requirements:

\begin{assumption}\label{app:ass:weights}
(i) There are $k_n\times k_n$ matrices $\PSigma > 0$ with $\|\hat \Sigma_n - \PSigma\|_{o,p} = o_P(1)$ uniformly in $P\in \mathbf P$;
(ii) $\|\PSigma\|_{o,p}\vee\|\PSigma^{-1}\|_{o,p}$ is uniformly bounded in $k_n$ and $P\in \mathbf P$.
\end{assumption}

\begin{assumption}\label{app:ass:rates}
Define the sequence $\eta_n \equiv J_nB_nk_n^{1/p}\sqrt{\log(1+k_n)/n}$. Then: (i) $\sup_{\theta\in \IDsetRsieve} Q_P(\theta)\times \|\hat\Sigma_n - \PSigma\|_{o,p} = O_P(\eta_n)$ uniformly in $P\in \mathbf P_0$; (ii) $\sup_{\theta \in \IDsetRsieve} Q_P(\theta) = \inf_{\theta \in \Theta_n \cap R} Q_P(\theta) + O(\eta_n)$ uniformly in $P\in \mathbf P_0$.
\end{assumption}

\begin{assumption}\label{app:ass:keycons}
There are sets $\mathcal V_{n}(P)\subseteq \Theta_n \cap R$ and a sequence $\{\nu_{n}\}_{n=1}^\infty$ with $\nu^{-1}_{n} = O(1)$, such that $\hat \Theta_n^{\text{\rm r}} \subseteq \mathcal V_n(P)$ with probability tending to one uniformly in $P\in \mathbf P_0$ and for any $\theta \in \mathcal V_n(P)$ and $\eta_n \equiv J_nB_nk_n^{1/p} \sqrt{\log(1+k_n)/n}$ it follows that
$$\nu_{n}^{-1}  \overrightarrow d_H(\{\theta\}, \IDsetRsieve,\|\cdot\|_{\mathbf E})
\leq \{Q_{P}(\theta) - \inf_{\tilde \theta \in \Theta_n \cap R} Q_{P}(\tilde \theta) \}+ O(\eta_n).$$
\end{assumption}

In particular, note Assumption \ref{app:ass:weights} is implied by Assumption \ref{ass:weights}.
Similarly, Assumption \ref{app:ass:rates} follows from Assumptions \ref{ass:weights}(i) and \ref{ass:locrates}(ii), while Assumption \ref{app:ass:keycons} will be verified by relying on Assumptions \ref{ass:keycons}(i), \ref{ass:keycons}(ii) or \ref{ass:extra}(iii) (depending on the choice of $\tau_n)$, and \ref{ass:locrates}(ii).
Given these assumptions, we next establish a consistency (Lemma \ref{lm:setcons}) and rate of convergence results (Theorem \ref{app:th:setrates}) for $\hat\Theta_n^{\rm r}$.

\begin{lemma}\label{lm:setcons}
Let Assumptions \ref{ass:param}(i), \ref{ass:startreg}(i)(iii), \ref{app:ass:weights}, \ref{app:ass:rates}(i), $\|\cdot\|_{\mathbf A}$ be a norm on $\mathbf B_n$ and for $\epsilon > 0$ let $\mathcal V_n(P) \equiv \{\theta \in \Theta_n \cap R : \overrightarrow{d}_H(\{\theta\},\IDsetRsieve,\|\cdot\|_{\mathbf A}) \leq \epsilon\}$, and define
$$S_n(\epsilon) \equiv \inf_{P\in \mathbf P_0} \{\inf_{\theta \in (\Theta_n \cap R)\setminus \mathcal V_n(P)} Q_P(\theta) - \inf_{\theta \in \Theta_n \cap R} Q_P(\theta)\}.$$
(i) If $\eta_n \vee \tau_n = o(S_n(\epsilon))$ for  $\eta_n \equiv k_n^{1/p}\sqrt{\log(1+k_n)}J_nB_n/\sqrt n$, then $\hat \Theta_n^{\text{\rm r}} \subseteq \mathcal V_n(P)$ with probability tending to one uniformly in $P\in \mathbf P_0$.
(ii) If Assumption \ref{app:ass:rates}(ii) holds and $\eta_n = o(\tau_n)$, then $\IDsetRsieve \subseteq \hat \Theta_n^{\text{\rm r}}$ with probability tending to one uniformly in $P\in \mathbf P_0$.
\end{lemma}

\noindent {\sc Proof:} For a given $\epsilon > 0$ first notice that by definition of $\hat \Theta_n^{\text{r}}$ and $\mathcal V_n(P)$ we have
\begin{equation}\label{lm:setcons1}
P (\overrightarrow d_H(\hat \Theta_n^{\text{r}}, \IDsetRsieve,\|\cdot\|_{\mathbf A}) > \epsilon)\\
\leq  P(\inf_{\theta \in (\Theta_n\cap R)\setminus \mathcal V_n(P)}  Q_{n}(\theta) \leq \inf_{\theta \in \Theta_n \cap R}Q_n(\theta) + \tau_n)
\end{equation}
for all $P\in \mathbf P_0$.
Setting $\hat Q_P(\theta) \equiv \| E_P[\rho(X,\theta)*q^{k_n}(Z)]\|_{\hat \Sigma_n,p}$ then note that Lemma \ref{aux:suppQ} and $\|\hat \Sigma_n\|_{o,p} = O_P(1)$ uniformly in $P\in \mathbf P_0$ by Lemma \ref{aux:weights} allow us to conclude
\begin{equation}\label{lm:setcons2}
\inf_{\theta \in (\Theta_n\cap R)\setminus \mathcal V_n(P)} \hat Q_P(\theta) \leq \inf_{\theta \in (\Theta_n\cap R)\setminus \mathcal V_n(P)} Q_{n}(\theta) + O_P(\eta_n)
\end{equation}
uniformly in $P\in \mathbf P_0$.
In addition, by similar arguments we obtain uniformly in $P\in \mathbf P_0$
\begin{equation}\label{lm:setcons3}
\inf_{\theta \in \Theta_n \cap R}  Q_n(\theta)  \leq \inf_{\theta \in \Theta_n \cap R} \hat Q_P(\theta) + O_P(\eta_n).
\end{equation}
Next note that for any $a\in \mathbf R^{k_n}$ we have $ \|\PSigma a\|_p \leq \|\PSigma \hat \Sigma_n^{-1}\|_{o,p} \|\hat \Sigma_n a\|_p$, and therefore 
\begin{align}\label{lm:setcons4}
\inf_{\theta \in (\Theta_n\cap R)\setminus \mathcal V_n(P)} \hat Q_P(\theta) & \geq \|\PSigma \hat \Sigma_n^{-1}\|_{o,p}^{-1} \inf_{\theta \in (\Theta_n\cap R)\setminus \mathcal V_n(P)} Q_P(\theta) \notag \\
& \geq \|\PSigma \hat \Sigma_n^{-1}\|_{o,p}^{-1}\{S_n(\epsilon) + \inf_{\theta \in \Theta_n \cap R} Q_P(\theta)\}
\end{align}
by definition of $S_n(\epsilon)$.
Similarly, employing that $\|\hat \Sigma_n a\|_p \leq \|\hat \Sigma_n \PSigma^{-1}\|_{o,p} \|\PSigma a\|_p$ yields
\begin{multline}\label{lm:setcons5}
\inf_{\theta \in \Theta_n \cap R} \hat Q_P(\theta) - \|\PSigma \hat \Sigma_n^{-1}\|_{o,p}^{-1} \inf_{\theta \in \Theta_n \cap R} Q_P(\theta) \\
\leq \{\|\hat \Sigma_n \PSigma^{-1}\|_{o,p} - \|\PSigma \hat \Sigma_n^{-1}\|_{o,p}^{-1}\}\inf_{\theta \in \Theta_n \cap R} Q_P(\theta).
\end{multline}
For $I_{k_n}$ the $k_n \times k_n$ identity matrix, then note that $\|I_{k_n}\|_{o,p} = 1$ implies the bound
\begin{multline}\label{lm:setcons6}
|\|\PSigma\hat \Sigma_n^{-1}\|_{o,p} - 1| = |\|\PSigma\hat \Sigma_n^{-1}\|_{o,p} - \|I_{k_n}\|_{o,p}| \leq \|(\PSigma - \hat \Sigma_n)\hat \Sigma_n^{-1}\|_{o,p}\\
\leq \|\hat \Sigma_n^{-1}\|_{o,p} \|\PSigma - \hat \Sigma_n\|_{o,p} = O_P(\|\PSigma - \hat \Sigma_n\|_{o,p}),
\end{multline}
where the final equality holds uniformly in $P\in \mathbf P_0$ by Lemma \ref{aux:weights}.
By identical arguments it follows that $|\|\hat \Sigma_n \PSigma^{-1}\|_{o,p} -1| =O_P(\|\hat \Sigma_n - \PSigma\|_{o,p})$ uniformly in $P\in \mathbf P_0$, and therefore \eqref{lm:setcons5}, $\IDsetRsieve\subseteq \Theta_n \cap R$, and Assumption \ref{app:ass:rates}(i) imply that
\begin{equation}\label{lm:setcons7}
\inf_{\theta \in \Theta_n \cap R} \hat Q_P(\theta) - \|\PSigma \hat \Sigma_n^{-1}\|_{o,p}^{-1} \inf_{\theta \in \Theta_n \cap R} Q_P(\theta) \\
\leq O_P(\eta_n)
\end{equation}
uniformly in $P\in \mathbf P_0$.
Therefore, \eqref{lm:setcons1}, \eqref{lm:setcons2}, \eqref{lm:setcons3}, \eqref{lm:setcons4}, and \eqref{lm:setcons7} yield that
\begin{multline*}
\limsup_{n\rightarrow \infty} \sup_{P\in \mathbf P_0} P (\overrightarrow d_H(\hat \Theta_n^{\text{r}}, \IDsetRsieve,\|\cdot\|_{\mathbf A}) > \epsilon)\\
\leq \limsup_{M \uparrow \infty} \limsup_{n\rightarrow \infty} \sup_{P\in \mathbf P_0} P(S_n(\epsilon) \leq \|\PSigma\|_{o,p}\|\hat \Sigma_n^{-1}\|_{o,p} M (\eta_n+\tau_n)) = 0,
\end{multline*}
where the equality follows from Lemma \ref{aux:weights}, Assumption \ref{app:ass:weights}(ii), and $\eta_n \vee \tau_n = o(S_n(\epsilon))$ by hypothesis.
Part (i) of the lemma then follows by definition of $\mathcal V_n(P)$.

In order to establish part (ii) of the lemma, note that the definition of $\hat \Theta_n^{\text{r}}$ implies
\begin{equation}\label{lm:setcons9}
P(\IDsetRsieve \subseteq \hat \Theta_n^{\text{r}}) \geq P(\sup_{\theta \in \IDsetRsieve} Q_n(\theta) \leq \inf_{\theta \in \Theta_n \cap R} Q_n(\theta) + \tau_n )
\end{equation}
for all $P\in \mathbf P_0$.
Moreover, applying Lemmas \ref{aux:suppQ} and \ref{aux:weights} together with $\|\hat \Sigma_n a\|_p \leq \|\hat \Sigma_n \PSigma^{-1}\|_{o,p}\|\PSigma a\|_p$ for any $a\in \mathbf R^{k_n}$ implies that uniformly in $P\in \mathbf P_0$
\begin{multline}\label{lm:setcons10}
\sup_{\theta \in \IDsetRsieve} Q_n(\theta) \leq \sup_{\theta \in \IDsetRsieve} \hat Q_P (\theta) + O_P(\eta_n) \\ \leq \|\hat \Sigma_n \PSigma^{-1}\|_{o,p} \sup_{\theta \in \IDsetRsieve} Q_P(\theta) + O_P(\eta_n) = \inf_{\theta \in \Theta_n \cap R} Q_P(\theta) + O_P(\eta_n),
\end{multline}
where the final equality follows from Assumption \ref{app:ass:rates}(ii), identical arguments to those in \eqref{lm:setcons6} implying $| \|\hat \Sigma_n \PSigma^{-1}\|_{o,p} -1| = O_P(\|\hat \Sigma_n - \PSigma\|_{o,p})$ uniformly in $P\in \mathbf P$, and Assumption \ref{app:ass:rates}(i).
Similarly, Lemmas \ref{aux:suppQ} and \ref{aux:weights}, $\|\PSigma a\|_p \leq \|\PSigma \hat \Sigma_n^{-1}\|_{o,p}\|\hat \Sigma_n a\|_p$ for any $a \in \mathbf R^{k_n}$, Assumption \ref{app:ass:rates}(i), and result \eqref{lm:setcons6} imply that uniformly in $P\in \mathbf P_0$
\begin{multline}\label{lm:setcons11}
\inf_{\theta \in \Theta_n \cap R} Q_n(\theta) \geq \inf_{\theta \in \Theta_n \cap R} \hat Q_{P}(\theta) - O_P(\eta_n) \\ \geq \|\PSigma\hat \Sigma_n^{-1}\|_{o,p}^{-1} \inf_{\theta \in \Theta_n\cap R} Q_{P}(\theta) - O_P(\eta_n)  = \inf_{\theta \in \Theta_n \cap R} Q_{P}(\theta) - O_P(\eta_n).
\end{multline}
Part (ii) of the lemma thus follows from \eqref{lm:setcons9}, \eqref{lm:setcons10}, \eqref{lm:setcons11}, and $\eta_n = o(\tau_n)$. \qed

\begin{theorem}\label{app:th:setrates}
Let Assumptions \ref{ass:param}(i), \ref{ass:startreg}(i)(iii), \ref{app:ass:weights}, \ref{app:ass:rates}, \ref{app:ass:keycons} hold, and 
\begin{equation}\label{th:setratesdef}
\mathcal R_n \equiv \nu_n\{\frac{k_n^{1/p}\sqrt{\log(1+k_n)}J_nB_n}{\sqrt n}\} .
\end{equation}
Then uniformly in $P\in \mathbf P_0$: (i) $\overrightarrow d_H(\hat \Theta_n^{\text{\rm r}}, \IDsetRsieve,\|\cdot\|_{\mathbf E}) = O_P(\mathcal R_n + \nu_n\tau_n)$; and
(ii) $d_H(\hat \Theta_n^{\text{\rm r}}, \IDsetRsieve,\|\cdot\|_{\mathbf E}) = O_P(\nu_n\tau_n)$ provided $J_nB_nk_n^{1/p}\sqrt{\log(1+k_n)/n} = o(\tau_n)$.
\end{theorem}

\noindent {\sc Proof:} Let $\eta_n \equiv k_n^{1/p}\sqrt{\log(1+k_n)} J_n B_n/\sqrt n$, $\delta_n^{-1} \equiv \nu_n(\eta_n + \tau_n)$, and $Q_P(\theta) \equiv \|E_P[\rho(X,\theta)*q^{k_n}(Z)]\|_{\PSigma,p}$.
In addition, we define $A_n \equiv A_{n1} \cap A_{n2} \cap A_{n3}$ where
\begin{align}\label{th:setrates1}
A_{n1} & \equiv \{\hat \Theta_n^{\text{r}} \subseteq \mathcal V_n(P)\} \nonumber\\
A_{n2} & \equiv \{\hat \Sigma_n^{-1} \text{ exists and } \|\hat \Sigma_n^{-1}\|_{o,p} \vee \|\hat \Sigma_n\|_{o,p} \vee \|\PSigma^{-1}\|_{o,p} \vee \|\PSigma\|_{o,p} < B\} \notag\\
A_{n3} & \equiv \{\sup_{\theta \in \IDsetRsieve} Q_P(\theta) \times \|\hat \Sigma_n - \PSigma\|_{o,p} \leq B \eta_n \text{ and } \|\hat \Sigma_n - \PSigma\|_{o,p} \leq \frac{1}{2B} \}.
\end{align}
Moreover, note that for any $\epsilon > 0$ and $B$ sufficiently large we can conclude that
\begin{equation}\label{th:setrates2}
\limsup_{n\rightarrow \infty} \sup_{P\in \mathbf P_0} P(A_n^c) < \epsilon
\end{equation}
due to Lemma \ref{aux:weights} and Assumptions \ref{app:ass:weights}(i), \ref{app:ass:rates}(i), and \ref{app:ass:keycons}.
Therefore, we obtain
\begin{multline}\label{th:setrates3}
\limsup_{n\rightarrow \infty} \sup_{P\in \mathbf P_0} P( \delta_n \overrightarrow d_H(\hat \Theta_n^{\text{r}}, \IDsetRsieve,\|\cdot\|_{\mathbf E}) > 2^M) \\
\leq \limsup_{n\rightarrow \infty} \sup_{P\in \mathbf P_0} P( \delta_n \overrightarrow d_H(\hat \Theta_n^{\text{r}}, \IDsetRsieve,\|\cdot\|_{\mathbf E}) > 2^M; ~ A_n) + \epsilon
\end{multline}
for any $M$.
For each $P\in \mathbf P_0$, next partition $\mathcal V_n(P)$ into subsets $S_{n,j}(P)$ defined by
\begin{equation*} 
S_{n,j}(P) \equiv \{\theta \in \mathcal V_n(P) : 2^{j-1} < \delta_n\overrightarrow d_H(\{\theta\},\IDsetRsieve,\|\cdot\|_{\mathbf E}) \leq 2^j\} .
\end{equation*}
Since $\hat \Theta_n^{\text{r}} \subseteq \mathcal V_n(P)$ under $A_n$, it follows from  the definition of $\hat \Theta_n^{\text{r}}$, and \eqref{th:setrates3} that
\begin{multline}\label{th:setrates5}
 \limsup_{n\rightarrow \infty} \sup_{P\in \mathbf P_0} P( \delta_n \overrightarrow d_H(\hat \Theta_n^{\text{r}}, \IDsetRsieve,\|\cdot\|_{\mathbf E}) > 2^M) \\
\leq \limsup_{n\rightarrow \infty} \sup_{P\in \mathbf P_0} \sum_{j \geq M}^\infty P( \inf_{\theta \in S_{n,j}(P)}Q_n(\theta) \leq \inf_{\theta \in \Theta_n \cap R} Q_n(\theta) + \tau_n; ~ A_n)+\epsilon .
\end{multline}
Letting $\hat Q_P(\theta) \equiv \|E_P[\rho(X,\theta) * q^{k_n}(Z)]\|_{\hat \Sigma_n,p}$, we then obtain from Lemma \ref{aux:suppQ} that
\begin{equation}\label{th:setrates6}
\inf_{\theta \in \Theta_n \cap R} Q_n(\theta) \leq \inf_{\theta \in \Theta_n \cap R} \hat Q_P(\theta) + \|\hat \Sigma_n\|_{o,p} \mathcal Z_{n,P} \leq \inf_{\theta \in \Theta_n \cap R} \hat Q_P(\theta) + B \mathcal Z_{n,P}
\end{equation}
where the final inequality holds under the event $A_n$ by \eqref{th:setrates1}.
Moreover, since for any $a\in \mathbf R^{k_n}$ we have $\|\hat \Sigma_n a\|_p \leq \|\hat \Sigma_n \PSigma^{-1}\|_{o,p} \|\PSigma a\|_p$, we obtain from $\IDsetRsieve \subseteq \Theta_n\cap R$ and the inequality $\|\hat \Sigma_n \PSigma^{-1}\|_{o,p} \leq \|\{\hat \Sigma_n - \PSigma\}\PSigma^{-1}\|_{o,p} + 1$ that under the event $A_n$ we have
\begin{multline}\label{th:setrates7}
\inf_{\theta \in \Theta_n \cap R} \hat Q_P(\theta)  \leq \|\hat \Sigma_n \PSigma^{-1}\|_{o,p} \inf_{\theta \in \Theta_n \cap R} Q_P(\theta)  \\
 \leq \{1 + \|\PSigma^{-1}\|_{o,p} \|\hat \Sigma_n - \PSigma\|_{o,p}\}\inf_{\theta \in \Theta_n \cap R}Q_P(\theta)  \leq \inf_{\theta \in \Theta_n \cap R} Q_P(\theta) + B^2\eta_n.
\end{multline}
In addition, note that by similar arguments we also obtain from Lemma \ref{aux:suppQ} and $\|\PSigma a\|_p \leq \|\PSigma \hat \Sigma_n^{-1}\|_{o,p} \|\hat\Sigma_n a\|_p$ that under the event $A_n$ we must have
\begin{multline}\label{th:setrates8}
\inf_{\theta \in S_{n,j}(P)} Q_n(\theta) \geq \inf_{\theta \in S_{n,j}(P)} \hat Q_P(\theta) - \|\hat \Sigma_n\|_{o,p} \mathcal Z_{n,P} \\
\geq \|\PSigma \hat \Sigma_n^{-1}\|_{o,p}^{-1} \inf_{\theta \in S_{n,j}(P)}  Q_P(\theta) - B \mathcal Z_{n,P} .
\end{multline}
Next, we note the triangle inequality, $\|(\PSigma - \hat \Sigma_n)\hat \Sigma_n^{-1}\|_{o,p} \leq \|\hat \Sigma_n^{-1}\|_{o,p} \|\hat \Sigma_n - \PSigma\|_{o,p}$, and $\|\hat \Sigma_n^{-1}\|_{o,p} \leq B$ under the event $A_n$ by \eqref{th:setrates1} yield the inequality
\begin{multline}\label{th:setrates9}
\|\PSigma \hat \Sigma_n^{-1}\|_{o,p}^{-1} - 1 \geq (\|(\PSigma - \hat \Sigma_n)\hat \Sigma_n^{-1}\|_{o,p} + 1)^{-1} - 1 \\
\geq - \|(\PSigma -\hat \Sigma_n)\hat \Sigma_n^{-1}\|_{o,p} \geq  - B \|\hat \Sigma_n - \PSigma\|_{o,p}.
\end{multline}
Therefore, combining results \eqref{th:setrates8} and \eqref{th:setrates9}, together with Assumption \ref{app:ass:keycons} and the definition of $S_{n,j}(P)$ we obtain for $B$ sufficiently large that under the event $A_n$ we have
\begin{align}\label{th:setrates10}
\inf_{\theta \in S_{n,j}(P)} Q_n(\theta) & \geq (1 - B\|\hat \Sigma_n - \PSigma\|_{o,p})\times \inf_{\theta \in S_{n,j}(P)} Q_P(\theta) - B\mathcal Z_{n,P} \notag \\
& \geq (1 - B\|\hat \Sigma_n - \PSigma\|_{o,p})(\inf_{\theta \in \Theta_n \cap R} Q_P(\theta) + \frac{2^{j-1}}{\nu_n \delta_n} - B\eta_n) - B\mathcal Z_{n,P} \notag \\
& \geq \inf_{\theta \in \Theta_n \cap R} Q_P(\theta) + \frac{2^{j-2}}{\nu_n \delta_n} - B(\mathcal Z_{n,P} + 2B\eta_n),
\end{align}
where the final inequality follows from $\IDsetRsieve\subseteq \Theta_n \cap R$ and the definition of the evnet $A_n$ in \eqref{th:setrates1}.
Hence, results \eqref{th:setrates5}, \eqref{th:setrates6}, \eqref{th:setrates7}, and \eqref{th:setrates10} yield
\begin{align}\label{th:setrates11}
\limsup_{M\uparrow \infty} & \limsup_{n\rightarrow \infty} \sup_{P\in \mathbf P_0} P( \delta_n \overrightarrow d_H(\hat \Theta_n^{\text{r}}, \IDsetRsieve,\|\cdot\|_{\mathbf E}) > 2^M) \nonumber \\
& \leq \limsup_{M\uparrow \infty} \limsup_{n\rightarrow \infty} \sup_{P\in \mathbf P_0} \sum_{j \geq M}^\infty P(\frac{2^{j-2}}{\nu_n\delta_n} \leq  3B(B\eta_n + \mathcal Z_{n,P}) + \tau_n; ~ A_n) + \epsilon\nonumber \\ & \leq \limsup_{M\uparrow \infty} \limsup_{n\rightarrow \infty} \sup_{P\in \mathbf P_0} \sum_{j \geq M}^\infty P(2^{(j-3)}(\eta_n + \tau_n) \leq 3B\mathcal Z_{n,P}) + \epsilon,
\end{align}
where in the final inequality we employed that we had defined $\delta_n^{-1} \equiv \nu_n(\eta_n + \tau_n)$.
Therefore, $\mathcal Z_{n,P}\in \mathbf R_+$, Lemma \ref{aux:suppQ}, $\tau_n \geq 0$, and Markov's inequality yield
\begin{multline}\label{th:setrates12}
\limsup_{M\uparrow \infty} \limsup_{n\rightarrow \infty} \sup_{P\in \mathbf P_0} \sum_{j\geq M}^\infty P(2^{(j-3)}(\eta_n + \tau_n) \leq 3B\mathcal Z_{n,P}) \\
\lesssim \limsup_{M\uparrow \infty} \limsup_{n\rightarrow \infty}  \sum_{j\geq M}2^{-j}\times \frac{1}{\eta_n}\frac{k_n^{1/p}\sqrt {\log(1+k_n)}J_nB_n}{\sqrt n} = 0 ,
\end{multline}
where in the final result we employed that $\eta_n \equiv k_n^{1/p}\sqrt{\log(1 + k_n)}J_nB_n/\sqrt n$.
The first claim of the theorem therefore follows from \eqref{th:setrates11}, \eqref{th:setrates12}, and $\epsilon$ being arbitrary.

To establish the second claim of the theorem, next define the event $A_{n4} \equiv \{\IDsetRsieve \subseteq \hat \Theta_n^{\text{r}}\}$.
Since $\overrightarrow d_H(\IDsetRsieve, \hat\Theta_n^{\text{r}},\|\cdot\|_{\mathbf E}) = 0$ whenever the event $A_{n4}$ occurs, we can conclude from Lemma \ref{lm:setcons}(ii) and part (i) of this theorem that
\begin{multline}\label{th:setrates13}
\limsup_{M\uparrow \infty}\limsup_{n\rightarrow \infty} \sup_{P\in \mathbf P_0} P( \delta_n d_H(\hat \Theta_n^{\text{r}}, \IDsetRsieve,\|\cdot\|_{\mathbf E}) > 2^M) \\
= \limsup_{M\uparrow \infty}\limsup_{n\rightarrow \infty} \sup_{P\in \mathbf P_0} P( \delta_n \overrightarrow d_H(\hat \Theta_n^{\text{r}}, \IDsetRsieve,\|\cdot\|_{\mathbf E}) > 2^M) = 0,
\end{multline}
and thus the theorem follows from $\delta_n^{-1} = \nu_n (\eta_n + \tau_n)$ and $\eta_n = o(\tau_n)$. \qed

\begin{corollary}\label{cor:setrates}
If Assumptions \ref{ass:param}(i), \ref{ass:startreg}(i)(iii), \ref{ass:coupling}(i), \ref{ass:keycons}, \ref{ass:locrates}(ii), and \ref{ass:weights} hold, then $\overrightarrow d_H(\hat \theta_n, \IDsetRsieve,\|\cdot\|_{\mathbf E}) = O_P(\mathcal R_n)$ uniformly in $P\in \mathbf P_0$.
\end{corollary}

\noindent {\sc Proof:} Follows from Theorem \ref{app:th:setrates}(i) applied with $\tau_n \equiv a_n/\sqrt n$ after noting that $a_n = o(1)$ (by Assumption \ref{ass:coupling}(i)) implies $\nu_n a_n/\sqrt n = o(\mathcal R_n)$ and: (i) Assumption \ref{app:ass:weights} holds by Assumption \ref{ass:weights}; (ii) Assumption \ref{app:ass:rates}(i) holds by Assumptions \ref{ass:locrates}(ii) and \ref{ass:weights}(i); (iii) Assumption \ref{app:ass:rates}(ii) holds by $Q_P(\theta) \geq 0$ and Assumption \ref{ass:locrates}(ii); and (iv) Assumption \ref{app:ass:keycons} holds with $\tau_n \equiv a_n/\sqrt n$ by Assumptions \ref{ass:keycons} and \ref{ass:locrates}(ii), the triangle inequality, and $\inf_{\theta \in \Theta_n \cap R} Q_P(\theta) \leq \sup_{\theta \in \IDsetRsieve} Q_P(\theta)$ due to $\IDsetRsieve\subseteq \Theta_n \cap R$. \qed

\begin{corollary}\label{cor:bootsetrates}
Let Assumptions \ref{ass:param}(i), \ref{ass:startreg}(i)(iii), \ref{ass:coupling}(i), \ref{ass:keycons}(i), \ref{ass:locrates}(ii), \ref{ass:weights}, and \ref{ass:extra}(iii) hold. Then uniformly in $P\in \mathbf P_0$: (i) $\overrightarrow d_H(\hat \Theta_n^{\rm r}, \IDsetRsieve,\|\cdot\|_{\mathbf E}) = O_P(\mathcal R_n + \nu_n \tau_n)$; and (ii) $d_H(\hat \Theta_n^{\rm r}, \IDsetRsieve,\|\cdot\|_{\mathbf E}) = O_P(\nu_n \tau_n)$ provided $J_nB_nk_n^{1/p}\sqrt{\log(1+k_n)/n} = o(\tau_n)$.
\end{corollary}

\noindent {\sc Proof:} Follows from Theorem \ref{app:th:setrates} after noting that $a_n = o(1)$ (by Assumption \ref{ass:coupling}(i)) implies: (i) Assumption \ref{app:ass:weights} holds by Assumption \ref{ass:weights}; (ii) Assumption \ref{app:ass:rates}(i) holds by Assumptions \ref{ass:locrates}(ii) and \ref{ass:weights}(i); (iii) Assumption \ref{app:ass:rates}(ii) holds by $Q_P(\theta) \geq 0$ and Assumption \ref{ass:locrates}(ii); and (iv) Assumption \ref{app:ass:keycons} holds by Assumptions \ref{ass:keycons}(i), \ref{ass:locrates}(ii), \ref{ass:extra}(iii), the triangle inequality, and $\IDsetRsieve\subseteq \Theta_n \cap R$. \qed

\begin{lemma}\label{aux:suppQ}
Let $\hat Q_{P}(\theta) \equiv \|E_P[\rho(X,\theta)*q^{k_n}(Z)]\|_{\hat \Sigma_n,p}$, and Assumptions \ref{ass:param}(i), \ref{ass:startreg}(i), and \ref{ass:startreg}(iii) hold. Then, for each $P\in \mathbf P$ there are random $\mathcal Z_{n,P} \in \mathbf R_+$ with
\begin{equation*} 
|Q_n(\theta) - \hat Q_{P}(\theta)| \leq \|\hat \Sigma_n\|_{o,p} \times \mathcal Z_{n,P} ,
\end{equation*}
for all $\theta \in \Theta_n \cap R$ and in addition $\sup_{P\in \mathbf P} E_P[\mathcal Z_{n,P}] = O(k_n^{1/p}\sqrt{\log(1+k_n)}J_nB_n/\sqrt n)$.
\end{lemma}

\noindent {\sc Proof:} Let $\mathcal G_n \equiv \{fq_{k,\jmath} : f\in\mathcal F_n, ~1\leq \jmath \leq \mathcal J \text{ and } 1\leq k \leq k_{n,\jmath}\}$.
Note that by Assumption \ref{ass:startreg}(i), $ \|q_{k,\jmath}\|_{\infty} \leq B_n$ for all $1\leq \jmath \leq \mathcal J$ and $1\leq k \leq k_{n,\jmath}$.
Hence, letting $F_n$ be the envelope for $\mathcal F_n$, as in Assumption \ref{ass:startreg}(iii), it follows that $G_n \equiv B_nF_n$ is an envelope for $\mathcal G_n$ satisfying $\sup_{P\in \mathbf P}E_P[G^2_n(V)] < \infty$.
Thus, we obtain
\begin{equation}\label{aux:suppQ1}
\sup_{P\in \mathbf P} E_P[\sup_{g \in \mathcal G_n} |\frac{1}{\sqrt n}\sum_{i=1}^n (g(V_i) - E_P[g(V)])|] \lesssim \sup_{P\in \mathbf P}J_{[\hspace{0.02 in}]}(\|G_n\|_{P,2},\mathcal G_n,\|\cdot\|_{P,2})
\end{equation}
by Theorem 2.14.2 in \cite{vandervaart:wellner:1996}.
Moreover, also notice that Lemma \ref{aux:simpbracket}, the change of variables $u = \epsilon/B_n$, and $B_n \geq 1$ imply
\begin{multline}\label{aux:suppQ2}
\sup_{P\in \mathbf P}J_{[\hspace{0.02 in}]}(\|G_n\|_{P,2},\mathcal G_n,\|\cdot\|_{P,2}) \leq \sup_{P\in \mathbf P} \int_0^{\|G_n\|_{P,2}} \sqrt{1 + \log (k_n N_{[\hspace{0.02 in}]}(\epsilon/B_n,\mathcal F_n,\|\cdot\|_{P,2}))}d\epsilon \\ \leq (1+\sqrt{\log (k_n)})B_n \times \sup_{P\in \mathbf P} J_{[\hspace{0.02 in}]}(\|F_n\|_{P,2},\mathcal F_n,\|\cdot\|_{P,2}) = O(\sqrt{\log(1+k_n)}B_nJ_n) ,
\end{multline}
where the final equality follows from Assumption \ref{ass:startreg}(iii).
Next define $\mathcal Z_{n,P} \in \mathbf R_+$ by
\begin{equation*} 
\mathcal Z_{n,P} \equiv \frac{k_n^{1/p}}{\sqrt n} \times \sup_{g\in \mathcal G_n} |\frac{1}{\sqrt n}\sum_{i=1}^n(g(V_i) - E_P[g(V)])|
\end{equation*}
and note \eqref{aux:suppQ1} and \eqref{aux:suppQ2} imply $\sup_{P\in \mathbf P} E_P[\mathcal Z_{n,P}] = O(k_n^{1/p}\sqrt {\log(1+k_n)}B_nJ_n/\sqrt n)$ as desired. 
Moreover, for any $\theta \in \Theta_n \cap R$, the definitions of $\Gemp(\theta)$, $\mathcal G_n$, and $\mathcal Z_{n,P}$ yield
\begin{multline*} 
|Q_n(\theta) - \hat Q_{P}(\theta)| \leq \frac{\|\hat \Sigma_n\|_{o,p}}{\sqrt n} \times \|\Gemp(\theta)\|_p \\ \leq \|\hat \Sigma_n\|_{o,p}\times \frac{k_n^{1/p}}{\sqrt n} \times  \sup_{g\in \mathcal G_n}|\frac{1}{\sqrt n}\sum_{i=1}^n(g(V_i) - E_P[g(V)])| \equiv \|\hat \Sigma_n\|_{o,p}\times \mathcal Z_{n,P}  ,
\end{multline*}
which establishes the claim of the lemma. \qed

\begin{lemma}\label{aux:simpbracket}
Let $\{g_j\}_{j=1}^{J}$ be functions satisfying $\max_{1\leq j \leq J} \|g_j\|_{\infty} \leq C < \infty$ and define $\mathcal G_n \equiv \{f g_j : f\in \mathcal F_n, ~ 1 \leq j \leq J\}$.
Then for any $P$ it follows that 
$$N_{[\hspace{0.02 in}]}(\epsilon, \mathcal G_n,\|\cdot\|_{P,2}) \leq J \times N_{[\hspace{0.02 in}]}(\epsilon/C, \mathcal F_n,\|\cdot\|_{P,2}).$$
\end{lemma}

\noindent {\sc Proof:} First define $g_{j}^+ \equiv g_j \vee 0$ and $g_{j}^- \equiv g_j \wedge 0$, where $\vee$ and $\wedge$ denote the pointwise maximums and minimums.
If $\{[f_{i,l}, f_{i,u}]\}_i$ is a collection of brackets for $\mathcal F_n$ satisfying
\begin{equation}\label{aux:simpbracket1}
\int (f_{i,u} - f_{i,l})^2dP \leq \epsilon^2
\end{equation}
for all $i$, then it follows that the following collection of brackets covers the class $\mathcal G_n$:
\begin{equation}\label{aux:simpbracket2}
\{[g_{j}^+f_{i,l} + g_{j}^-f_{i,u}, ~ g_{j}^-f_{i,l} + g_{j}^+f_{i,u}]\}_{i,j} .
\end{equation}
Moreover, since $|g_{j}| = g_{j}^+ - g_{j}^-$ by construction, we also obtain from result \eqref{aux:simpbracket1} that
\begin{equation}\label{aux:simpbracket3}
\int(g_{j}^+f_{i,u} + g_{j}^-f_{i,l} - g_{j}^+f_{i,l} - g_{j}^-f_{i,u})^2 dP 
 = \int (f_{i,u} - f_{i,l})^2|g_j|^2dP\leq  \epsilon^2C^2.
\end{equation}
Since there are $J \times N_{[\hspace{0.02 in}]}(\epsilon,\mathcal F_n,\|\cdot\|_{P,2})$ brackets in \eqref{aux:simpbracket2}, we conclude from \eqref{aux:simpbracket3} that $N_{[\hspace{0.02 in}]}(\epsilon,\mathcal G_n,\|\cdot\|_{P,2}) \leq J \times N_{[\hspace{0.02 in}]}(\epsilon/C,\mathcal F_n,\|\cdot\|_{P,2})$, which establishes the lemma. \qed

\begin{lemma}\label{aux:weights}
If Assumption \ref{app:ass:weights} holds, then there is a constant $B < \infty$ such that
\begin{equation*}
\liminf_{n\rightarrow \infty} \inf_{P\in \mathbf P} P( \hat \Sigma_n^{-1} \text{ exists and } \max\{\|\hat \Sigma_n\|_{o,p},\|\hat \Sigma_n^{-1}\|_{o,p}\} < B) = 1 .
\end{equation*}
\end{lemma}

\noindent {\sc Proof:} Recall that $\hat \Sigma_n$ and $\Sigma_P$ are $k_n\times k_n$ matrices, though the dependence on $k_n$ was suppressed from the notation. 
Then note that by Assumption \ref{app:ass:weights}(ii) there exists a constant $B < \infty$ such that for all $P\in \mathbf P$ and $k_n$ we have that
\begin{equation}\label{aux:weights1}
\max\{\|\PSigma\|_{o,p},\|\PSigma^{-1}\|_{o,p}\} < \frac{B}{2} .
\end{equation}
Next, let $I_{k_n}$ denote the $k_n \times k_n$ identity matrix and for each $P\in \mathbf P$ rewrite $\hat \Sigma_n$ as
\begin{equation}\label{aux:weights2}
\hat \Sigma_n = \PSigma\{I_{k_n} - \PSigma^{-1}(\PSigma - \hat \Sigma_n)\} .
\end{equation}
By Theorem 2.9 in \cite{kress:1999}, the matrix $\{I_{k_n} - \PSigma^{-1}(\PSigma - \hat \Sigma_n)\}$ is invertible and the operator norm of its inverse is bounded by two when $\|\PSigma^{-1}(\PSigma - \hat \Sigma_n)\|_{o,p} < 1/2$. 
Since $\PSigma$ is invertible by Assumption \ref{app:ass:weights}(i), result \eqref{aux:weights2} implies that $\hat \Sigma_n$ is invertible if and only if $\{I_{k_n} - \PSigma^{-1}(\PSigma - \hat \Sigma_n)\}$ is invertible, which yields that
\begin{multline}\label{aux:weights3}
P( \hat \Sigma_n^{-1} \text{ exists and } \|\{I_{k_n} - \PSigma^{-1}(\PSigma - \hat \Sigma_n)\}^{-1}\|_{o,p} < 2) \\ \geq  P( \|\PSigma^{-1}(\hat \Sigma_n - \PSigma)\|_{o,p} < \frac{1}{2}) \geq P(\|\hat \Sigma_n - \PSigma\|_{o,p} < \frac{1}{B}) ,
\end{multline}
where we employed $\|\PSigma^{-1}(\hat \Sigma_n - \PSigma)\|_{o,p} \leq \|\PSigma^{-1}\|_{o,p}\|\hat \Sigma_n - \PSigma\|_{o,p}$ and \eqref{aux:weights1}. 
Hence, since \eqref{aux:weights2} implies $\hat \Sigma_n^{-1} = \{I_{k_n} - \PSigma^{-1}(\PSigma - \hat \Sigma_n)\}^{-1} \PSigma^{-1}$ whenever $ \{I_{k_n} - \PSigma^{-1}(\PSigma - \hat \Sigma_n)\}^{-1}$ exists, the bound $\|\PSigma^{-1}\|_{o,p} < B/2$ and result \eqref{aux:weights3} allow us to conclude
\begin{equation}\label{aux:weights4}
P( \hat \Sigma_n^{-1} \text{ exists and } \|\hat \Sigma_n^{-1}\|_{o,p} < B) \geq P(\|\hat \Sigma_n - \PSigma\|_{o,p} < \frac{1}{B}) .
\end{equation}
Finally, since $\|\hat \Sigma_n\|_{o,p} \leq B/2 + \|\hat \Sigma_n - \PSigma\|_{o,p}$ by \eqref{aux:weights1}, result \eqref{aux:weights4} implies that
\begin{multline*}
\liminf_{n\rightarrow \infty} \inf_{P\in \mathbf P} P( \hat \Sigma_n^{-1} \text{ exists and } \max\{\|\hat \Sigma_n\|_{o,p},\|\hat \Sigma_n^{-1}\|_{o,p}\} < B) \\ \geq \liminf_{n\rightarrow \infty} \inf_{P\in \mathbf P} P(\|\hat \Sigma_n - \PSigma\|_{o,p} < \min\{\frac{B}{2},\frac{1}{B}\}) = 1 ,
\end{multline*}
where the equality, and hence the lemma, follows from Assumption \ref{app:ass:weights}(i). \qed

\begin{corollary}\label{cor:weights}
If Assumption \ref{ass:weights} holds, then for some $B < \infty$ it follows that:
$$\liminf_{n\rightarrow \infty} \inf_{P\in \mathbf P} P( \hat \Sigma_n^{-1} \text{ exists and } \max\{\|\hat \Sigma_n\|_{o,p},\|\hat \Sigma_n^{-1}\|_{o,p}\} < B) = 1 .$$
\end{corollary}

\noindent {\sc Proof:} Follows from Lemma \ref{aux:weights} and Assumption \ref{ass:weights} together with $a_n = o(1)$, which is imposed by Assumption \ref{ass:coupling}(i) (or  \ref{ass:bootcoupling}), implying Assumption \ref{app:ass:weights} holds. \qed

\begin{lemma}\label{lm:obviousineq}
If $a\in \mathbf R^d$, then $\|a\|_{\tilde p} \leq d^{(\frac{1}{\tilde p} - \frac{1}{p})_+}\|a\|_p$ for any $\tilde p, p\in [2,\infty]$.
\end{lemma}

\noindent {\sc Proof:} The case $p \leq \tilde p$ trivially follows from $\|a\|_{\tilde p} \leq \|a\|_p$ for all $a\in \mathbf R^d$.
For the case $p > \tilde p$, let $a = (a_1,\ldots, a_d)^\prime$ and note that by H\"{o}lder's inequality we obtain
\begin{multline}\label{lm:obviousineq1}
\|a\|_{\tilde p}^{\tilde p} = \sum_{i=1}^d \{|a_i|^{\tilde p} \times 1 \} \leq \{\sum_{i=1}^d (|a_i|^{\tilde p})^\frac{p}{\tilde p}\}^{\frac{\tilde p}{p}} \{\sum_{i=1}^d 1^{\frac{p}{p-\tilde p}}\}^{1 - \frac{\tilde p}{p}} = \{\sum_{i=1}^d |a_i|^p\}^{\frac{\tilde p}{p}} d^{1-\frac{\tilde p}{p}} .
\end{multline}
Thus, the claim of the lemma for $p > \tilde p$ follows from taking the $1/\tilde p$ power in \eqref{lm:obviousineq1}. \qed

%% file: Appendix/AppStrong.tex

\section{Strong Approximation} \label{sec:appstrong}

This Section contains the proof of Theorem \ref{th:localdrift} and supporting results.

\noindent{\sc Proof of Theorem \ref{th:localdrift}:}
First note that by Assumption \ref{ass:weights}(ii) there is a constant $C_0 < \infty$ such that $\|\PSigma\|_{o,p} \leq C_0$ for all $P\in \mathbf P_0$.
Hence, Assumption \ref{ass:locrates}(ii) and Lemma \ref{lm:obviousineq} imply that for all $P\in \mathbf P_0$, $\theta\in \IDsetRsieve$, and $h\in V_{n}(\theta,R|\ell_n)$ we have
\begin{multline}\label{th:localdrift2}
\| \sqrt n E_P[\rho(X,\theta + \frac{h}{\sqrt n}) *q^{k_n}(Z)] - \DerP(\theta)[h]\|_{\PSigma,p} \\
\leq C_0\|\sqrt n E_P[(\rho(X,\theta + \frac{h}{\sqrt n}) - \rho(X,\theta)) *q^{k_n}(Z)] - \DerP(\theta)[h]\|_2 + o(a_n) .
\end{multline}
Moreover, Lemma \ref{aux:bessel} and the maps $m_{P,\jmath}$ satisfying Assumption \ref{ass:driftlin}(i) imply that
\begin{align}\label{th:localdrift3}
\sum_{\jmath = 1}^{\mathcal J}\sum_{k=1}^{k_{n,\jmath}} \langle \sqrt n & \{m_{P,\jmath}(\theta + \frac{h}{\sqrt n}) - m_{P,\jmath}(\theta)\} - \nabla m_{P,\jmath}(\theta)[h], q_{k,\jmath}\rangle_{L^2_P}^2  \nonumber \\
& \leq  \sum_{\jmath = 1}^{\mathcal J} C_1\|\sqrt n\{m_{P,\jmath}(\theta + \frac{h}{\sqrt n}) - m_{P,\jmath} (\theta) - \nabla m_{P,\jmath}(\theta)[\frac{h}{\sqrt n}]\}\|_{P,2}^2 \nonumber\\
& \leq  \sum_{\jmath = 1}^{\mathcal J} C_1K_m^2\times n \times \|\frac{h}{\sqrt n}\|_{\mathbf L}^2 \times \|\frac{h}{\sqrt n}\|_{\mathbf E}^2
\end{align}
for some constant $C_1 < \infty$ and all $P\in \mathbf P_0$, $\theta \in \IDsetRsieve$, and $h\in V_{n}(\theta,R|\ell_n)$.
Therefore, by results \eqref{th:localdrift2} and \eqref{th:localdrift3}, the law of iterated expectations, the definition of $\mathcal S_n(\mathbf L,\mathbf E)$, and $K_m\ell^2_n\times \mathcal S_n(\mathbf L,\mathbf E) =  o(a_nn^{-1/2})$ by hypothesis, we obtain that
\begin{multline}\label{th:localdrift4}
\sup_{P\in \mathbf P_0}  \sup_{\theta \in \IDsetRsieve} \sup_{h \in V_{n}(\theta,R|\ell_n)} \|\sqrt n E_P[\rho(X,\theta + \frac{h}{\sqrt n}) * q^{k_n}(Z)] - \DerP(\theta)[h]\|_{\PSigma,p} \\
\lesssim K_m\times \sqrt n \ell_n^2\times \mathcal S_n(\mathbf L,\mathbf E) + o(a_n) =  o(a_n).
\end{multline}
Next, note that since $k_n^{1/p}\sqrt{\log(1+k_n)}B_n\times \sup_{P\in \mathbf P} J_{[\hspace{0.02 in}]}(\ell_n^{\kappa_\rho},\mathcal F_n,\|\cdot\|_{P,2}) = o(a_n)$, Assumption \ref{ass:locrates}(i) implies there is a sequence $\tilde \ell_n$ satisfying the conditions of Lemma \ref{lm:localprelim} and $\ell_n = o(\tilde \ell_n)$.
Therefore, applying Lemma \ref{lm:localprelim} we obtain uniformly in $P\in \mathbf P_0$
\begin{equation}\label{th:localdrift5}
I_{n}(R) = \inf_{\theta \in \IDsetRsieve} \inf_{h\in V_{n}(\theta,R|\tilde \ell_n)} \|\WP(\theta) + \sqrt n E_P[\rho(X,\theta + \frac{h}{\sqrt n})*q^{k_n}(Z)]\|_{\PSigma,p} +  o_P(a_n).
\end{equation}
Moreover, since $\ell_n = o(\tilde \ell_n)$ implies that $V_n(\theta,R|\tilde \ell_n) \subseteq V_n(\theta,R|\ell_n)$ for all $\theta \in \Theta_n \cap R$ for $n$ sufficiently large, we obtain uniformly in $P\in \mathbf P_0$ that
\begin{align}\label{th:localdrift6}
& \inf_{\theta \in \IDsetRsieve} \inf_{h \in V_{n}(\theta,R|\tilde \ell_n)} \|\WP(\theta) + \sqrt n E_P[\rho(X,\theta + \frac{h}{\sqrt n})*q^{k_n}(Z)]\|_{\PSigma,p} \nonumber \\
& \leq  \inf_{\theta \in \IDsetRsieve} \inf_{h \in V_{n}(\theta,R|\ell_n)} \|\WP(\theta) + \sqrt n E_P[\rho(X,\theta + \frac{h}{\sqrt n})*q^{k_n}(Z)]\|_{\PSigma,p} \nonumber \\
& = \inf_{\theta \in \IDsetRsieve} \inf_{h \in V_{n}(\theta,R|\ell_n)} \|\WP(\theta) + \DerP(\theta)[h]\|_{\PSigma,p} +  o(a_n),
\end{align}
where the final equality following from \eqref{th:localdrift4}, Assumption \ref{ass:weights}(ii) and Lemma \ref{lm:geninf}.
Thus, the first claim of the Theorem follows from \eqref{th:localdrift5} and \eqref{th:localdrift6}, while the second follows by noting that if $K_m\mathcal R_n^2\times \mathcal S_n(\mathbf L,\mathbf E) =  o(a_nn^{-1/2})$, then we may set $\ell_n$ to simultaneously satisfy the conditions of Lemma \ref{lm:localprelim} and $K_m\mathcal \ell_n^2\times\mathcal S_n(\mathbf L,\mathbf E) =  o(a_nn^{-1/2})$, which obviates the need to introduce $\tilde \ell_n$ in \eqref{th:localdrift5} and \eqref{th:localdrift6}. \qed

\begin{lemma}\label{lm:localprelim}
Let Assumptions \ref{ass:param}(i), \ref{ass:startreg}(i), \ref{ass:startreg}(iii), \ref{ass:coupling}, \ref{ass:keycons}, \ref{ass:locrates}, and \ref{ass:weights} hold.
Then, for any sequence $\{\ell_n\}$ satisfying $k_n^{1/p}\sqrt{\log(1+k_n)}B_n \sup_{P\in \mathbf P}J_{[\hspace{0.02 in}]}(\ell_n^{\kappa_\rho},\mathcal F_n,\|\cdot\|_{P,2}) = o(a_n)$ and $\mathcal R_n = o(\ell_n)$, we have uniformly in $P\in \mathbf P_0$ that:
$$I_{n}(R) = \inf_{\theta \in \IDsetRsieve} \inf_{h \in V_{n}(\theta,R|\ell_n)} \|\WP(\theta) + \sqrt n E_P[\rho(X,\theta + \frac{h}{\sqrt n})*q^{k_n}(Z)]\|_{\PSigma,p} + o_P(a_n).$$
\end{lemma}

\noindent {\sc Proof:} First note that the required sequence $\{\ell_n\}$ exists by Assumption \ref{ass:locrates}(i). 
Next, note that by Assumption \ref{ass:keycons}(ii) and Corollary \ref{cor:setrates} there is a $\hat \theta_n\in \Theta_n \cap R$ satisfying 
\begin{equation}\label{lm:localprelim1}
Q_n(\hat \theta_n) \leq \inf_{\theta \in \Theta_n \cap R} Q_n(\theta) +  o(a_n/\sqrt n)
\end{equation}
and $\overrightarrow d_H(\hat \theta_n,\IDsetRsieve,\|\cdot\|_{\mathbf E}) = O_P(\mathcal R_n)$ uniformly in $P\in \mathbf P_0$.
Hence, defining $(\IDsetRsieve)^{\ell_n}\equiv\{\theta \in \Theta_n \cap R : \overrightarrow d_H(\theta,\IDsetRsieve,\|\cdot\|_{\mathbf E})\leq \ell_n\}$, which implicitly depends on $P\in \mathbf P_0$, we obtain
\begin{equation}\label{lm:localprelim3}
I_n(R) = \inf_{\theta \in (\IDsetRsieve)^{\ell_n}} \sqrt n Q_n(\theta) + o_P(a_n)
\end{equation}
uniformly in $P\in \mathbf P_0$ due to $\mathcal R_n = o(\ell_n)$, $\overrightarrow d_H(\hat \theta_n,\IDsetRsieve,\|\cdot\|_{\mathbf E}) = O_P(\mathcal R_n)$, $(\IDsetRsieve)^{\ell_n} \subseteq \Theta_n \cap R$ by construction, result \eqref{lm:localprelim1}, and the definition of $I_n(R)$.
Next, note that by Assumption \ref{ass:coupling}(i), Corollary \ref{cor:weights}, and Lemma \ref{lm:geninf} it follows that
\begin{multline}\label{lm:localprelim5}
|\inf_{\theta \in (\IDsetRsieve)^{\ell_n}} \sqrt n Q_n(\theta) - \inf_{\theta \in (\IDsetRsieve)^{\ell_n}} \|\WP(\theta) + \sqrt n E_P[\rho(X,\theta)*q^{k_n}(Z)]\|_{\hat \Sigma_n,p}| \\ \leq \|\hat \Sigma_n\|_{o,p} \times \sup_{\theta \in\Theta_n \cap R} \|\Gemp(\theta) - \WP(\theta)\|_p =  o_P(a_n)
\end{multline}
uniformly in $P\in \mathbf P_0$.
Similarly, employing Corollary \ref{cor:weights}, Lemmas \ref{lm:locprocess}, \ref{lm:geninf}, and $\ell_n$ satisfying $k_n^{1/p}\sqrt{\log(1+k_n)}B_n\times \sup_{P\in \mathbf P}J_{[\hspace{0.02 in}]}(\ell_n^{\kappa_\rho},\mathcal F_n,\|\cdot\|_{P,2}) =  o(a_n)$  yields
\begin{multline*}
\inf_{\theta\in \IDsetRsieve}\inf_{h \in V_{n}(\theta,R|\ell_n)} \|\WP(\theta+\frac{h}{\sqrt n}) + \sqrt n E_P[\rho(X,\theta + \frac{h}{\sqrt n})*q^{k_n}(Z)]\|_{\hat \Sigma_n,p}\\
= \inf_{\theta\in \IDsetRsieve}\inf_{h \in V_{n}(\theta,R|\ell_n)} \|\WP(\theta) + \sqrt n E_P[\rho(X,\theta + \frac{h}{\sqrt n})*q^{k_n}(Z)]\|_{\hat \Sigma_n,p} +  o_P(a_n)
\end{multline*}
uniformly in $P\in \mathbf P_0$, which together with results \eqref{lm:localprelim3} and \eqref{lm:localprelim5}, and Lemma \ref{lm:locsigma} establish the claim of the lemma. \qed

\begin{lemma}\label{lm:locprocess}
Let Assumptions \ref{ass:startreg}(i) and \ref{ass:coupling}(ii) hold.
If $\{\delta_n\}$ is a sequence satisfying $k_n^{1/p}\sqrt{\log(1+k_n)}B_n \times \sup_{P\in \mathbf P}J_{[\hspace{0.02 in}]}(\delta_n^{\kappa_\rho},\mathcal F_n,\|\cdot\|_{P,2}) =  o(a_n)$, then uniformly in $P\in \mathbf P$:
$$\sup_{\theta\in \IDsetRsieve} \sup_{h \in V_{n}(\theta,R|\delta_n)} \|\WP(\theta + \frac{h}{\sqrt n}) - \WP(\theta)\|_{p} = o_P(a_n) .$$
\end{lemma}

\noindent {\sc Proof:} Since $\|q_{k,\jmath}\|_{\infty} \leq B_n$ for all $1\leq \jmath \leq \mathcal J$ and $1\leq k \leq k_{n,\jmath}$ by Assumption \ref{ass:startreg}(i), Assumption \ref{ass:coupling}(ii) yields for any $P\in \mathbf P$, $\theta \in \Theta_n \cap R$, and $h \in V_{n}(\theta,R|\delta_n)$ that
\begin{equation}\label{lm:locprocess1}
E_P[\|\rho(X,\theta + \frac{h}{\sqrt n}) - \rho(X,\theta)\|_2^2q_{k,\jmath}^2(Z)] \leq K_\rho^2B_n^2 \|\frac{h}{\sqrt n}\|_{\mathbf E}^{2\kappa_\rho} \leq K_\rho^2 B_n^2 \delta_n^{2\kappa_\rho} .
\end{equation}
Set $\mathcal G_n \equiv \{fq_{k,\jmath} \text{ for some } f \in \mathcal F_n, ~1\leq \jmath \leq \mathcal J, ~ 1\leq k \leq k_{n,\jmath}\}$ and let $\Iso$ be a Gaussian process on $\mathcal G_n$ satisfying $E[\Iso(g_1)] = 0$ and $E[\Iso (g_1)\Iso(g_2)] = \text{Cov}_P\{g_1(V),g_2(V)\}$ for any $g_1,g_2\in \mathcal G_n$.  
Since $\|a\|_p \leq  k_n^{1/p}\|a\|_\infty$ for any $a\in \mathbf R^{k_n}$, result \eqref{lm:locprocess1} yields 
\begin{multline}\label{lm:locprocess3}
E_P[\sup_{\theta\in \IDsetRsieve}  \sup_{h \in V_{n}(\theta,R|\delta_n)} \|\WP(\theta + \frac{h}{\sqrt n}) - \WP(\theta)\|_{p}] \\  \leq k_n^{1/p}\times E[\sup_{g_1,g_2 \in \mathcal G_n: \|g_1 - g_2\|_{P,2} \leq K_\rho B_n \delta_n^{\kappa_\rho}} |\Iso(g_1) - \Iso(g_2)| ].
\end{multline}
Moreover, Corollary 2.2.8 in \cite{vandervaart:wellner:1996} additionally implies that
\begin{align}\label{lm:locprocess4}
\sup_{P\in \mathbf P} E_P[& \sup_{g_1,g_2 \in \mathcal G_n : \|g_1 - g_2\|_{P,2} \leq K_\rho B_n \delta_n^{\kappa_\rho}} |\Iso(g_1)- \Iso(g_2)|] \notag\\
& \lesssim \sup_{P\in \mathbf P} \int_0^{K_\rho B_n\delta_n^{\kappa_\rho}} \sqrt{\log N_{[\hspace{0.02 in}]}(\epsilon,\mathcal G_n,\|\cdot\|_{P,2})}d\epsilon \notag \\
& \lesssim \sup_{P\in \mathbf P}\sqrt{\log(1+k_n)}B_n \int_0^{K_\rho \delta_n^{\kappa_\rho}} \sqrt{1+ \log N_{[\hspace{0.02 in}]}(u,\mathcal F_n,\|\cdot\|_{P,2})}du,
\end{align}
where the second inequality follows from Lemma \ref{aux:simpbracket} and the change of variables $u = \epsilon/B_n$.
However, note that since $N_{[\hspace{0.02 in}]}(u,\mathcal F_n,\|\cdot\|_{P,2})$ is decreasing in $u$, it follows that $J_{[\hspace{0.02 in}]}(K_\rho \delta_n^{\kappa_\rho},\mathcal F_n,\|\cdot\|_{P,2}) \leq K_\rho J_{[\hspace{0.02 in}]}(\delta_n^{\kappa_\rho},\mathcal F_n,\|\cdot\|_{P,2})$.
Therefore, the lemma follows from results \eqref{lm:locprocess3} and \eqref{lm:locprocess4}, the definition of $J_{[\hspace{0.02 in}]}(\epsilon,\mathcal F_n,\|\cdot\|_{P,2})$, and $k_n^{1/p}\sqrt{\log(1+k_n)}B_n\times \sup_{P\in \mathbf P}J_{[\hspace{0.02 in}]}(\delta_n^{\kappa_\rho},\mathcal F_n,\|\cdot\|_{P,2}) =  o(a_n)$ by hypothesis. \qed

\begin{lemma}\label{lm:locsigma}
Let Assumptions \ref{ass:startreg}(i), \ref{ass:startreg}(iii), \ref{ass:locrates}(ii), and \ref{ass:weights} hold with $a_n = o(1)$. For any positive sequence $\delta_n$ it then follows that uniformly in $P\in \mathbf P_0$ we have
\begin{multline*}
\inf_{\theta \in \IDsetRsieve} \inf_{h\in V_{n}(\theta,R|\delta_n)} \|\WP(\theta) + \sqrt n E_P[\rho(X,\theta + \frac{h}{\sqrt n})*q^{k_n}(Z)]\|_{\PSigma,p} \\ = \inf_{ \theta \in \IDsetRsieve} \inf_{h \in V_{n}(\theta,R|\delta_n)} \|\WP(\theta) + \sqrt n E_P[\rho(X,\theta + \frac{h}{\sqrt n})*q^{k_n}(Z)]\|_{\hat \Sigma_n,p} + o_P(a_n).
\end{multline*}
\end{lemma}

\noindent {\sc Proof:} First note that by Assumption \ref{ass:weights}(ii) there is a $C_0 < \infty$ such that $\|\PSigma\|_{o,p} \vee \|\PSigma^{-1}\|_{o,p} \leq C_0$ for all $P\in \mathbf P$.
Since $\|\hat \Sigma_n a\|_p \leq \|\hat \Sigma_n \PSigma^{-1}\|_{o,p} \|\PSigma a\|_p$ for any $a \in \mathbf R^{k_n}$, and $\|\hat\Sigma_n\PSigma^{-1}\|_{o,p} \leq \|\PSigma^{-1}\|_{o,p}\|\hat \Sigma_n - \PSigma\|_{o,p} + 1$ by the triangle inequality, we obtain
\begin{multline}\label{lm:locsigma1}
 \{C_0\|\hat \Sigma_n - \PSigma\|_{o,p} + 1\}\|\WP(\theta) + \sqrt n E_P[\rho(X,\theta + \frac{h}{\sqrt n})*q^{k_n}(Z)]\|_{\PSigma,p} \\
 \geq \|\WP(\theta)+ \sqrt n E_P[\rho(X,\theta + \frac{h}{\sqrt n})*q^{k_n}(Z)]\|_{\hat \Sigma_n,p}
\end{multline}
for any $\theta\in \IDsetRsieve$ and $h \in V_n(\theta,R|\delta_n)$.
Moreover, $\|\PSigma\|_{o,p}\leq C_0$, $0\in V_n(\theta,R|\delta_n)$ for any $\theta \in \Theta_n \cap R$, and Assumption \ref{ass:locrates}(ii) imply uniformly in $P\in \mathbf P$ that
\begin{multline}\label{lm:locsigma2}
\inf_{\theta \in \IDsetRsieve} \inf_{h\in V_{n}(\theta,R|\delta_n)} \|\WP(\theta) + \sqrt n E_P[\rho(X,\theta + \frac{h}{\sqrt n})*q^{k_n}(Z)]\|_{\PSigma,p} \\  \lesssim  \sup_{\theta \in \Theta_{n}\cap R} \|\WP(\theta)\|_p  + o(a_n) = O_P(k_n^{1/p}\sqrt{\log(1+k_n)}B_nJ_n) + o(a_n)
\end{multline}
where the final equality holds uniformly in $P\in \mathbf P_0$ by Lemma \ref{aux:tailbound} and Markov's inequality.
Therefore, results \eqref{lm:locsigma1}, \eqref{lm:locsigma2}, and Assumption \ref{ass:weights}(i) imply
\begin{multline}\label{lm:locsigma4}
\inf_{\theta \in \IDsetRsieve} \inf_{h \in V_{n}(\theta,R|\delta_n)} \|\WP(\theta) + \sqrt n E_P[\rho(X,\theta + \frac{h}{\sqrt n})*q^{k_n}(Z)]\|_{\PSigma,p} +  o_P(a_n) \\
\geq \inf_{\theta \in \IDsetRsieve} \inf_{h\in V_{n}(\theta,R|\delta_n)} \|\WP(\theta) + \sqrt n E_P[\rho(X,\theta + \frac{h}{\sqrt n})*q^{k_n}(Z)]\|_{\hat \Sigma_n,p}
\end{multline}
uniformly in $P\in \mathbf P_0$.
Next, note that Assumption \ref{ass:weights} implies Assumption \ref{app:ass:weights} and therefore Lemma \ref{aux:weights} yields that $\|\hat \Sigma_n\|_{o,p} \vee \|\hat \Sigma_n^{-1}\|_{o,p} = O_P(1)$ uniformly in $P\in \mathbf P$.
The lemma then follows from \eqref{lm:locsigma4} and noting that the reverse inequality also holds by identical arguments but relying on $\|\hat \Sigma_n\|_{o,p} \vee \|\hat \Sigma_n^{-1}\|_{o,p} = O_P(1)$ uniformly in $P\in \mathbf P$ rather than on $\|\PSigma\|_{o,p} \vee \|\PSigma^{-1}\|_{o,p} \leq C_0$. \qed

\begin{lemma}\label{aux:tailbound}
If Assumptions \ref{ass:startreg}(i) and \ref{ass:startreg}(iii) hold, then for some $C < \infty$ we have:
$$\sup_{P\in \mathbf P} E_P[\sup_{\theta \in \Theta_n \cap R} \|\WP(\theta)\|_p ] \leq C k_n^{1/p}\sqrt{\log(1+k_n)}B_nJ_n .$$
\end{lemma}

\noindent {\sc Proof:} Let $\mathcal G_n \equiv \{f q_{k,\jmath} : f\in \mathcal F_n, ~ 1\leq \jmath \leq \mathcal J, \text{ and } 1\leq k \leq k_{n,\jmath}\}$ and $\Iso$ be a Gaussian process on $\mathcal G_n$ satisfying $E[\Iso(g_1)] = 0$ and $E[\Iso (g_1)\Iso(g_2)] = \text{Cov}_P\{g_1(V),g_2(V)\}$ for any $g_1,g_2\in \mathcal G_n$.  
Then note $\|a\|_p \leq d^{1/p} \|a\|_\infty$ for any $a\in \mathbf R^d$ implies that 
\begin{multline}\label{aux:tailbound1}
E_P[\sup_{\theta \in \Theta_n \cap R} \|\WP(\theta)\|_p] \leq k_n^{1/p} E_P[\sup_{g\in \mathcal G_n} |\Iso(g)|] \\
\leq k_n^{1/p}\{E_P[|\Iso(g_0)|] + C_1 \int_0^\infty \sqrt {\log N_{[\hspace{0.02 in }]}(\epsilon,\mathcal G_n,\|\cdot\|_{P,2})}d\epsilon\} ,
\end{multline}
where the final inequality holds for any $g_0 \in \mathcal G_n$ and some $C_1<\infty$ by Corollary 2.2.8 in \cite{vandervaart:wellner:1996}.
Next, let $G_n \equiv B_n F_n$ for $F_n$ as in Assumption \ref{ass:startreg}(iii) and note Assumption \ref{ass:startreg}(i) implies $G_n$ is an envelope for $\mathcal G_n$.
Thus, $[-G_n,G_n]$ is a bracket of size $2\|G_n\|_{P,2}$ covering $\mathcal G_n$, and as a result we obtain
\begin{multline}\label{aux:tailbound2}
\int_0^\infty \sqrt {\log N_{[\hspace{0.02 in }]}(\epsilon,\mathcal G_n,\|\cdot\|_{P,2})}d\epsilon \\ \leq \int_0^{2\|G_n\|_{P,2}} \sqrt {1+\log N_{[\hspace{0.02 in}]}(\epsilon,\mathcal G_n,\|\cdot\|_{P,2})}d\epsilon \leq C_2 \sqrt{\log(1+k_n)}B_nJ_n ,
\end{multline}
where the final inequality holds for some $C_2 < \infty$ by result \eqref{aux:suppQ2} and $N_{[\hspace{0.02 in}]}(u,\mathcal G_n,\|\cdot\|_{P,2})$ being decreasing in $u$.
Furthermore, since $E_P[|\Iso(g_0)|] \leq \|g_0\|_{P,2} \leq \|G_n\|_{P,2}$ we have
\begin{equation}\label{aux:tailbound3}
E_P[|\Iso(g_0)|] \leq \|G_n\|_{P,2} \\ \leq \int_0^{\|G_n\|_{P,2}} \sqrt {1+\log N_{[\hspace{0.02 in}]}(u,\mathcal G_n,\|\cdot\|_{P,2})}du .
\end{equation}
Thus, the claim of the lemma follows from \eqref{aux:tailbound1}, \eqref{aux:tailbound2}, and \eqref{aux:tailbound3}. \qed

\begin{lemma}\label{aux:bessel}
Let Assumption \ref{ass:startreg}(ii) hold. It then follows that there exists a constant $C < \infty$ such that for all $P\in \mathbf P$, $n\geq 1$, $1 \leq \jmath \leq \mathcal J$, and functions $f\in L^2_P$ we have
\begin{equation}\label{aux:bessel:disp}
\sum_{k=1}^{k_{n,\jmath}} \langle f, q_{k,\jmath}\rangle_{L^2_P}^2 \leq CE_P[(E_P[f(V)|Z_{\jmath}])^2] .
\end{equation}
\end{lemma}

\noindent {\sc Proof:} Let $L^2_P(Z_{\jmath})$ denote the subspace of $L^2_P$ consisting of functions depending on $Z_{\jmath}$ only, and set $\ell^2(\mathbb N) \equiv \{\{c_k\}_{k=1}^\infty : c_k \in \mathbf R \text{ and } \|\{c_k\} \|_{\ell^2(\mathbb N)} < \infty\}$, where $\|\{c_k\}\|_{\ell^2(\mathbb N)}^2 \equiv \sum_k c_k^2$.
For any sequence $\{c_k\} \in \ell^2(\mathbb N)$, then define the map $J_{\jmath,n}:\ell^2(\mathbb N)\rightarrow L^2_P(Z_{\jmath})$ by
\begin{equation*} 
J_{\jmath,n}(\{c_k\}) = \sum_{k=1}^{k_{n,\jmath}} c_kq_{k,\jmath}.
\end{equation*}
Clearly, the maps $J_{\jmath,n}: \ell^2(\mathbb N)\rightarrow L^2_P(Z_{\jmath})$ are linear and, moreover, by Assumption \ref{ass:startreg}(ii) there is a $C < \infty$ such that the largest eigenvalue of $E_P[q^{k_{n,\jmath}}_{\jmath}(Z_{\jmath})q^{k_{n,\jmath}}_{\jmath}(Z_{\jmath})^\prime]$ is bounded by $C$ for all $n\geq 1$ and $P\in \mathbf P$.
Therefore, we can conclude that
\begin{multline}\label{aux:bessel2}
\sup_{P\in \mathbf P} \sup_{n\geq 1}\|J_{\jmath,n}\|_o^2 = \sup_{P\in \mathbf P} \sup_{n\geq 1} \sup_{\{c_k\} : \sum_k c_k^2 = 1} \|J_{\jmath,n}(\{c_k\})\|_{P,2}^2 \\ = \sup_{P\in \mathbf P} \sup_{n\geq 1} \sup_{\{c_k\} : \sum_k c_k^2 = 1}  E_P[(\sum_{k=1}^{k_{n,\jmath}} c_kq_{k,\jmath}(Z_{\jmath}))^2] \leq  \sup_{\{c_k\} : \sum_k c_k^2 = 1} C\sum_{k=1}^\infty c_k^2 = C
\end{multline}
which implies $J_{\jmath,n}$ is continuous.
Next, define $J^*_{\jmath,n} : L^2_P(Z_{\jmath}) \rightarrow \ell^2(\mathbb N)$ to be given by
\begin{equation*} 
J^*_{\jmath,n} (g) = \{a_k(g)\}_{k=1}^\infty \hspace{0.5 in} a_k(g) \equiv \left\{ \begin{array}{cl} \langle g,q_{k,\jmath} \rangle_{L^2_P} & \text{ if } k \leq k_{n,\jmath} \\ 0 & \text{ if } k > k_{n,\jmath} \end{array}\right. ,
\end{equation*}
and note $J_{\jmath,n}^*$ is the adjoint of $J_{\jmath,n}$.
Therefore, since $\|J_{\jmath,n}\|_o = \|J_{\jmath,n}^*\|_o$ by Theorem 6.5.1 in \cite{luenberger:1969}, we obtain for any $P\in \mathbf P$, $n \geq 1$, and $g\in L^2_P(Z_{\jmath})$ that
\begin{equation}\label{aux:bessel4}
\sum_{k=1}^{k_{n,\jmath}} \langle g,q_{k,\jmath}\rangle_{L^2_P}^2 = \|J^*_{\jmath,n} (g)\|^2_{\ell^2(\mathbb N)} \leq \|J^*_{\jmath,n}\|_o^2 \|g\|_{P,2}^2 = \|J_{\jmath,n}\|_o^2 \|g\|_{P,2}^2 .
\end{equation}
Therefore, since $E_P[ f(V)q_{k,\jmath}(Z_{\jmath})] = E_P[E_P[ f(V)|Z_{\jmath}] q_{k,\jmath}(Z_{\jmath})]$ for any $ f\in L^2_P$, setting $g(Z_{\jmath}) = E_P[f(V)|Z_{\jmath}]$ in \eqref{aux:bessel4} and employing \eqref{aux:bessel2} yields the lemma. \qed

\begin{lemma}\label{lm:geninf}
If $\Lambda$ is a set, $A: \Lambda\rightarrow \mathbf R^k$, $B:\Lambda\rightarrow \mathbf R^k$, and $W$ is a $k\times k$ matrix, then
$$|\inf_{\lambda \in \Lambda} \|W A(\lambda)\|_p - \inf_{\lambda \in \Lambda} \|W B(\lambda)\|_p | \leq \|W\|_{o,p}\times \sup_{\lambda \in \Lambda}\|A(\lambda) - B(\lambda)\|_p .$$
\end{lemma}

\noindent {\sc Proof:} Fix $\eta > 0$, and let $\lambda_a \in \Lambda$ satisfy $\|WA(\lambda_a)\|_p \leq \inf_{\lambda \in \Lambda} \|WA(\lambda)\|_p + \eta$.
Then,
\begin{multline}\label{lm:geninf1}
\inf_{\lambda \in \Lambda} \|WB(\lambda)\|_p - \inf_{\lambda \in \Lambda} \|WA(\lambda)\|_p \leq \|WB(\lambda_a)\|_p - \|WA(\lambda_a)\|_p + \eta \\ \leq \|W(B(\lambda_a) - A(\lambda_a))\|_p + \eta \leq \|W\|_{o,p} \times \sup_{\lambda \in \Lambda} \|A(\lambda) - B(\lambda)\|_p + \eta,
\end{multline}
where the second result follows from the triangle inequality, and the final result from $\|Wv\|_p \leq \|W\|_{o,p}\|v\|_p$ for any $v\in \mathbf R^k$.
By identical manipulations we also have
\begin{equation}\label{lm:geninf2}
\inf_{\lambda \in \Lambda} \|WA(\lambda)\|_p - \inf_{\lambda \in \Lambda} \|WB(\lambda)\|_p \leq \|W\|_{o,p} \times \sup_{\lambda \in \Lambda} \|A(\lambda) - B(\lambda)\|_p + \eta .
\end{equation}
Thus, since $\eta$ was arbitrary, the lemma follows from results \eqref{lm:geninf1} and \eqref{lm:geninf2}. \qed

%% file: Appendix/AppBoot.tex

\section{Bootstrap Approximation} \label{sec:appboot}

This appendix contains the proof of all results concerning the bootstrap approximation.
We first introduce two assumptions that generalize Assumption \ref{ass:bootrates} (at the cost of introducing additional notation) and deliver a stronger version of Theorem \ref{th:coupsmooth}.

\begin{assumption}\label{ass:genneigh}
There is an $\epsilon > 0$ and scalars $\mathcal D_n(\mathbf L, \mathbf E)$ and $\mathcal D_n(\mathbf B,\mathbf E)$ such that for any $P\in \mathbf P$, $\theta \in \IDsetRsieve$, and $\theta_1 \in \Theta_n \cap R$ satisfying $\|\theta_1 - \theta\|_{\mathbf E} \leq \epsilon$, there exists $\tilde \theta \in \IDsetRsieve$ such that $\|\theta - \tilde \theta\|_{\mathbf E} = 0$, $\|\tilde \theta - \theta_1\|_{\mathbf L} \leq \mathcal D_n(\mathbf L,\mathbf E)\|\tilde \theta - \theta_1\|_{\mathbf E}$, and $\|\tilde \theta - \theta_1\|_{\mathbf B} \leq \mathcal D_n(\mathbf B,\mathbf E)\|\tilde \theta - \theta_1\|_{\mathbf E}$.
\end{assumption}

\begin{assumption}\label{ass:genbootrates}
(i) Either $\Upsilon_F$ and $\Upsilon_G$ are affine or $(\mathcal R_n + \nu_n\tau_n) \mathcal D_n(\mathbf B,\mathbf E) = o(1)$;
(ii) $k_n^{1/p}\sqrt{\log(1+k_n)}B_n \sup_{P\in \mathbf P}J_{[\hspace{0.03 in}]}(\ell_n^{\kappa_\rho}\vee(\nu_n\tau_n)^{\kappa_\rho},\mathcal F_n,\|\cdot\|_{P,2}) =  o(a_n)$,
$K_m \ell_n^2 \mathcal S_n(\mathbf L,\mathbf E) = o(a_n n^{-\frac{1}{2}})$, $K_m\ell_n(\mathcal R_n+ \nu_n\tau_n) \mathcal D_n(\mathbf L, \mathbf E) =  o(a_nn^{-\frac{1}{2}})$, $\ell_n(\ell_n +  \{\mathcal R_n+\nu_n\tau_n\}  \mathcal D_n(\mathbf B,\mathbf E))1\{K_f > 0\} =  o(a_nn^{-\frac{1}{2}})$;
(iii)  $\limsup 1\{K_g>0\}\ell_n/r_n < 1/2$ and $(\mathcal R_n +\nu_n\tau_n)  \mathcal D_n(\mathbf B, \mathbf E) = o(r_n)$.
\end{assumption}

In particular, note Assumption \ref{ass:genneigh} holds with $\mathcal D_n(\mathbf L,\mathbf E) = \mathcal S_n(\mathbf L,\mathbf E)$, $\mathcal D_n(\mathbf B,\mathbf E) = \mathcal S_n(\mathbf B,\mathbf E)$, and $\tilde \theta = \theta$.
Hence, Assumption \ref{ass:bootrates} implies Assumptions \ref{ass:genneigh} and \ref{ass:genbootrates}.
In general, however, $\mathcal D_n(\mathbf L,\mathbf E)$ and $\mathcal D_n(\mathbf B,\mathbf E)$ can be smaller than $\mathcal S_n(\mathbf L,\mathbf E)$ and $\mathcal S_n(\mathbf B,\mathbf E)$ while the introduction of a $\tilde \theta \neq \theta$ eases requirements in partially identified models.

Our next theorem consists of two parts.
The first part, which replaces Assumption \ref{ass:bootrates} with \ref{ass:genneigh} and \ref{ass:genbootrates}, can by the preceding discussion be seen as a generalization of Theorem \ref{th:coupsmooth}.
The second part shows that, under additional restrictions, it is possible to replace the norm $\|\cdot\|_{\mathbf B}$ in the definition of $\hat V_n(\theta,R|\ell)$ (as in \eqref{eq:Vndef}) with the norm $\|\cdot\|_{\mathbf E}$ -- an observation that is sometimes helpful in easing rate restrictions.

\begin{theorem}\label{th:gencoupsmooth}
Let Assumptions \ref{ass:param}, \ref{ass:startreg}, \ref{ass:coupling}, \ref{ass:keycons}(i),  \ref{ass:driftlin}, \ref{ass:locrates}, \ref{ass:weights}, \ref{ass:locineq}, \ref{ass:loceq}, \ref{ass:ineqlindep}, \ref{ass:bootcoupling}, \ref{ass:extra}(i)(iii), \ref{ass:genneigh}, and \ref{ass:genbootrates} hold. Then, the following statements hold:
\vspace{-0.15 in}
\begin{packed_enum}
    \item[(i)] If Assumption \ref{ass:extra}(ii) holds, then there is a $\tilde \ell_n \asymp \ell_n$ such that uniformly in $P\in \mathbf P_0$ 
    \begin{equation*}
    \hat U_n(R|\ell_n) \geq \UpS(R| \tilde \ell_n) +  o_P(a_n).
    \end{equation*}
    \item[(ii)] In addition, suppose for some $\epsilon > 0$ and $\|\cdot\|_{\mathbf I}$ satisfying $\|h\|_{\mathbf E} \leq \|h\|_{\mathbf I}$ for all $h\in \mathbf B_n$, we have that for all $P\in \mathbf P_0$, $\{\theta \in \mathbf B_n: \overrightarrow{d}_H(\theta,\IDsetRsieve,\|\cdot\|_{\mathbf I})\leq \epsilon\} \subseteq \Theta_n$ and $P(\{\theta \in \mathbf B_n: \overrightarrow{d}_H(\theta,\hat \Theta_n^{\text{\rm r}},\|\cdot\|_{\mathbf I})\leq \epsilon\} \subseteq \Theta_n)$ tends to one uniformly in $P \in \mathbf P_0$. If $\Upsilon_F$ and $\Upsilon_G$ are affine, then part (i) holds with $\hat U_n(R|\ell_n)$ as in \eqref{eq:hatudef} but with 
        \begin{equation}\label{th:gencoupsmoothdisp2}
        \hat V_n(\theta,R|\ell) \equiv \{h \in \mathbf B_n : h \in G_n(\theta), ~\Upsilon_F(\theta + \frac{h}{\sqrt n}) = 0, \text{ and } \|\frac{h}{\sqrt n}\|_{\mathbf I} \leq \ell\}.
        \end{equation}
\end{packed_enum}
\end{theorem}

\noindent {\sc Proof:}  First note Assumptions \ref{ass:locrates}(i) and \ref{ass:genbootrates}(ii) imply $\mathcal R_n \vee \nu_n\tau_n = o(1)$.
Hence, by Assumption \ref{ass:genbootrates}(ii) we may apply Lemma \ref{lm:simcoup} to obtain uniformly in $P\in \mathbf P_0$
\begin{equation}\label{th:gencoupsmooth1}
\hat U_n(R|\ell_n) = \inf_{\theta \in \hat \Theta_n^{\text{r}}} \inf_{h \in \hat V_n(\theta,R|\ell_n)} \|\WPT(\theta) + \DerP(\theta)[h]\|_{\PSigma,p} + o_P(a_n) .
\end{equation}
Thus, we may select $\hat \theta_n \in \hat \Theta_n^{\text{r}}$ and $\hat h_n \in \hat V_n(\hat \theta_n,R| \ell_n)$ so that uniformly in $P\in \mathbf P_0$
\begin{equation}\label{th:gencoupsmooth2}
\hat U_n(R|\ell_n) = \|\WPT(\hat \theta_n) + \DerP(\hat \theta_n)[\hat h_n]\|_{\PSigma,p} + o_P(a_n) .
\end{equation}
Next note that by Assumptions \ref{ass:locrates}(i), \ref{ass:genneigh}, and \ref{ass:genbootrates} there is a $\delta_n$ so that $\delta_n \mathcal D_n(\mathbf B,\mathbf E) = o(r_n)$, $\delta_n\mathcal D_n(\mathbf B,\mathbf E) = o(1)$ if either $\Upsilon_F$ or $\Upsilon_G$ are not affine, $\mathcal R_n +\nu_n\tau_n= o(\delta_n)$, and
\begin{align}
& \ell_n \delta_n  \mathcal D_n(\mathbf B,\mathbf E)1\{K_f > 0\} =  o(a_nn^{-\frac{1}{2}})   \label{th:gencoupsmooth3p1} \\
& K_m\delta_n \ell_n \mathcal D_n(\mathbf L, \mathbf E) = o(a_nn^{-\frac{1}{2}})  \label{th:gencoupsmooth3p2} \\
& k_n^{1/p}\sqrt{\log(1+k_n)}B_n \times \sup_{P\in \mathbf P} J_{[\hspace{0.03 in}]}(\delta_n^{\kappa_\rho},\mathcal F_n,\|\cdot\|_{P,2}) =  o(a_n) . \label{th:gencoupsmooth3p3}
\end{align}

Next, notice that Corollary \ref{cor:bootsetrates}(i) implies that there exists a $\theta_{0n}\in \IDsetRsieve$ such that
\begin{equation}\label{th:gencoupsmooth4}
\|\hat \theta_n - \theta_{0n}\|_{\mathbf E} = o_P(\delta_n)
\end{equation}
uniformly in $P \in \mathbf P_0$ due to $(\mathcal R_n + \nu_n\tau_n) = o(\delta_n)$.
Furthermore, by Assumption \ref{ass:genneigh} we can assume without loss of generality that $\theta_{0n}$ in addition satisfies
\begin{equation}\label{th:gencoupsmooth5}
\|\hat \theta_n - \theta_{0n}\|_{\mathbf L} = o_P(\mathcal D_n(\mathbf L,\mathbf E)\delta_n) \hspace{0.5 in}
\|\hat \theta_n - \theta_{0n}\|_{\mathbf B} = o_P(\mathcal D_n(\mathbf B,\mathbf E)\delta_n)
\end{equation}
uniformly in $P\in \mathbf P_0$.
In addition note that since $\|q_{k,\jmath}\|_{\infty} \leq B_n$ for all $1\leq \jmath \leq \mathcal J$ and $1\leq k \leq k_{n,\jmath}$ by Assumption \ref{ass:startreg}(i), we obtain from Assumption \ref{ass:coupling}(ii) together with result \eqref{th:gencoupsmooth4} that with probability tending to one uniformly in $P\in \mathbf P_0$ we have
\begin{equation}\label{th:coupsmooth5}
E_P[\|\rho(X,\hat \theta_n) - \rho(X,\theta_{0n})\|^2_2q_{k,\jmath}^2(Z_{\jmath})] \leq B_n^2K_\rho^2\delta_n^{2\kappa_\rho} .
\end{equation}
Set $\mathcal G_n \equiv \{fq_{k,\jmath} : f\in \mathcal F_n,  1 \leq \jmath \leq \mathcal J, ~ 1\leq k \leq k_{n,\jmath}\}$ and let $\Iso$ be a Gaussian process on $\mathcal G_n$ satisfying $E[\Iso (g_1)\Iso(g_2)] = \text{Cov}_P\{g_1(V),g_2(V)\}$ and $E[\Iso(g_1)] = 0$ for any $g_1,g_2\in \mathcal G_n$. 
Since \eqref{th:coupsmooth5} holds with probability tending to one uniformly in $P\in \mathbf P_0$, Assumption \ref{ass:weights}(ii), result \eqref{lm:locprocess4}, and $\delta_n$ satisfying \eqref{th:gencoupsmooth3p3} imply for any $\epsilon > 0$ that
\begin{multline}\label{th:gencoupsmooth6}
\limsup_{n\rightarrow \infty} \sup_{P\in \mathbf P_0} P(\|\WPT(\hat \theta_n) - \WPT(\theta_{0n})\|_{\PSigma,p} > a_n \epsilon) \\ \leq  \limsup_{n\rightarrow \infty} \sup_{P\in \mathbf P_0} \frac{1}{a_n \epsilon} k_n^{1/p} E_P[\sup_{g_1,g_2\in \mathcal G_n : \|g_1-g_2\|_{P,2} \leq B_nK_\rho\delta_n^{\kappa_\rho}} |\Iso(g_1) - \Iso(g_2)|] = 0.
\end{multline}
Similarly, result \eqref{th:gencoupsmooth4} implies $\overrightarrow{d}_H(\hat \theta_n,\IDsetRsieve,\|\cdot\|_{\mathbf E}) \leq \epsilon$ with probability tending to one uniformly in $P\in \mathbf P_0$ for any $\epsilon >0$. Hence, Lemma \ref{aux:dererror} yields uniformly in $P\in \mathbf P_0$
\begin{multline}\label{th:gencoupsmooth7}
\|\DerP(\theta_{0n})[\hat h_n] - \DerP(\hat \theta_n)[\hat h_n] \|_{\PSigma,p} \lesssim \|\PSigma\|_{o,p} \times K_m\|\hat \theta_n - \theta_{0n}\|_{\mathbf L} \|\hat h_n\|_{\mathbf E} +  o_P(a_n) \\ \lesssim \|\PSigma\|_{o,p} \times K_m \mathcal D_n(\mathbf L, \mathbf E) \delta_n  \ell_n\sqrt n + o_P(a_n) = o_P(a_n),
\end{multline}
where the second inequality follows from $\|\hat h_n/\sqrt n\|_{\mathbf B} \leq  \ell_n$ due to $\hat h_n/\sqrt n \in \hat V_n(\hat \theta_n, R|\ell_n)$, Assumption \ref{ass:extra}(i), and \eqref{th:gencoupsmooth5}.
In turn, the final result in \eqref{th:gencoupsmooth7} follows from \eqref{th:gencoupsmooth3p2} and Assumption \ref{ass:weights}(ii).
Next, note the conditions of Theorem \ref{th:localsmooth}(i) hold because: Either $\Upsilon_F$ and $\Upsilon_G$ are affine (implying $K_f = K_g = 0$) or $\delta_n \mathcal D_n(\mathbf B,\mathbf E) = o(1)$, and $\delta_n \mathcal D_n(\mathbf B, \mathbf E) = o(r_n)$ and $\limsup  \ell_n /r_n1\{K_g > 0\} < 1/2$ by Assumption \ref{ass:genbootrates}(iii) imply
\begin{equation*}
r_n \geq  2(\ell_n + \delta_n\mathcal D_n(\mathbf B,\mathbf E))1\{K_g > 0\}
\end{equation*}
for $n$ sufficiently large.
Hence, Theorem \ref{th:localsmooth}(i), Assumption \ref{ass:extra}(ii), and $\|h\|_{\mathbf E} \lesssim \|h\|_{\mathbf B}$ for all $h \in \mathbf B_n$ by Assumption \ref{ass:extra}(i), imply that there is a constant $M <\infty$ for which with probability tending to one  uniformly in $P\in \mathbf P_0$ we have that
\begin{equation*}
\inf_{h\in V_{n}(\theta_{0n}, R|M  \ell_n)} \|\frac{\hat h_n}{\sqrt n} - \frac{h}{\sqrt n}\|_{\mathbf B} \leq M  \ell_n( \ell_n + \delta_n \mathcal D_n(\mathbf B, \mathbf E))1\{K_f > 0\} .
\end{equation*}
It follows from Assumption \ref{ass:genbootrates}(ii) and \eqref{th:gencoupsmooth3p1} that there is a $h_{0n} \in V_{n}(\theta_{0n},R|M \ell_n)$ such that $\|h_{0n} - \hat h_n\|_{\mathbf B} = o_P(a_n)$ uniformly in $P\in \mathbf P_0$, and hence Assumption \ref{ass:weights}(ii), Lemma \ref{aux:dererror}, and $\|h\|_{\mathbf E}\lesssim \|h\|_{\mathbf B}$ by Assumption \ref{ass:extra}(i) yield
\begin{equation}\label{th:gencoupsmooth11}
\|\DerP(\theta_{0n})[\hat h_n] - \DerP(\theta_{0n})[h_{0n}]\|_{\PSigma,p} \lesssim \|\PSigma\|_{o,p}\times \|\hat h - h_{0n}\|_{\mathbf E} = o_P(a_n)
\end{equation}
uniformly in $P\in \mathbf P_0$. Therefore, combining results \eqref{th:gencoupsmooth2}, \eqref{th:gencoupsmooth6}, \eqref{th:gencoupsmooth7}, and \eqref{th:gencoupsmooth11} together with $\theta_{0n} \in \IDsetRsieve$ and $h_{0n}\in V_{n}(\theta_{0n},R|M \ell_n)$ yields
\begin{multline*}
\hat U_n(R|\ell_n) = \|\WPT(\theta_{0n}) + \DerP(\theta_{0n})[h_{0n}]\|_{\PSigma,p} + o_P(a_n) \\
\geq \inf_{\theta \in \IDsetRsieve} \inf_{h \in V_{n}(\theta,R|M \ell_n)} \|\WPT(\theta) + \DerP(\theta)[h]\|_{\PSigma,p} + o_P(a_n)
\end{multline*}
uniformly in $P\in \mathbf P_0$. 
The first part of theorem then follows by setting $\tilde \ell_n = M \ell_n$.

In order to establish the second part of the theorem, note that the only assumptions that potentially require the norm $\|\cdot\|_{\mathbf B}$ to be stronger than $\|\cdot\|_{\mathbf I}$ are Assumptions \ref{ass:locineq}, \ref{ass:loceq}, \ref{ass:ineqlindep} (pertaining to the differentiability of $\Upsilon_F$ and $\Upsilon_G$) and Assumption \ref{ass:extra}(ii) (since a stronger norm $\|\cdot\|_{\mathbf B}$ makes $(\hat \Theta_n^{\text{r}})^\epsilon$ smaller).
We therefore establish part (ii) of the theorem by repeating the arguments employed in showing part (i) while carefully re-examining the role played by the norm $\|\cdot\|_{\mathbf B}$.
To this end, note that since
\begin{equation}\label{th:gencoupsmooth12}
\liminf_{n\rightarrow \infty} \inf_{P\in \mathbf P_0} P(\{\theta \in \mathbf B_n: \overrightarrow{d}_H(\theta,\hat \Theta_n^{\rm r},\|\cdot\|_{\mathbf I})\leq \epsilon\} \subseteq \Theta_n) = 1,
\end{equation}
we may apply Lemma \ref{lm:simcoup} with $\|\cdot\|_{\mathbf B}$ set to equal $\|\cdot\|_{\mathbf I}$ to still obtain that
\begin{equation}\label{th:gencoupsmooth13}
\hat U_n(R|\ell_n) = \inf_{\theta \in \hat \Theta_n^{\text{r}}} \inf_{h \in \hat V_n(\theta,R|\ell_n)} \|\WPT(\theta) + \DerP(\theta)[h]\|_{\PSigma,p} + o_P(a_n) .
\end{equation}
Letting $\hat \theta_n$ and $\hat h_n$ be defined as in \eqref{th:gencoupsmooth2} (but with $\hat V_n(\theta,R|\ell)$ as defined in \eqref{th:gencoupsmoothdisp2}), then observe that since results \eqref{th:gencoupsmooth6} and \eqref{th:gencoupsmooth7} do not rely on Assumptions \ref{ass:locineq}, \ref{ass:loceq}, \ref{ass:ineqlindep} or \ref{ass:extra}(ii), we can conclude from result \eqref{th:gencoupsmooth13} that uniformly in $P\in \mathbf P_0$
\begin{equation}\label{th:gencoupsmooth14}
\hat U_n(R|\ell_n) =  \|\WPT(\theta_{0n}) + \DerP(\theta_{0n})[\hat h_n]\|_{\PSigma,p} + o_P(a_n)
\end{equation}
for some $\theta_{0n} \in \IDsetRsieve$.
Next, note $\delta_n \mathcal D_n(\mathbf B,\mathbf E) = o(r_n)$ and $K_f = K_g = 0$ due to $\Upsilon_F$ and $\Upsilon_G$ being affine, together with Theorem \ref{th:localsmooth}(ii) imply that 
\begin{align*}
\hat V_n(\hat \theta_n,R|\ell_n) & \equiv \{h \in \mathbf B_n : h \in G_n(\hat \theta_n), ~ \Upsilon_F(\hat \theta_n + \frac{h}{\sqrt n}) = 0, ~ \|\frac{h}{\sqrt n}|\|_{\mathbf I} \leq \ell_n\} \notag \\
& \subseteq \{h \in \mathbf B_n : \Upsilon_G(\theta_{0n} + \frac{h}{\sqrt n}) \leq 0, ~ \Upsilon_F(\theta_{0n} + \frac{h}{\sqrt n}) = 0, ~ \|\frac{h}{\sqrt n}\|_{\mathbf I} \leq \ell_n\} \notag \\ & \subseteq V_n(\theta_{0n},R|\ell_n),
\end{align*}
with probability tending to one uniformly in $P\in \mathbf P_0$, and where the final inequality follows from $\ell_n \downarrow 0$, $\{\theta \in \mathbf B_n: \overrightarrow{d}_H(\theta, \IDsetRsieve, \|\cdot\|_{\mathbf I})\leq \epsilon\} \subseteq \Theta_n$ and $\|\cdot\|_{\mathbf E}\leq \|\cdot\|_{\mathbf I}$ by hypothesis.
Therefore, we can conclude that $\hat h_n \in V_n(\theta_{0n},R|\ell_n)$ with probability tending to one uniformly in $P\in \mathbf P_0$, which by \eqref{th:gencoupsmooth14} yields
\begin{equation*}
\hat U_n(R|\ell_n) \geq \inf_{\theta \in \IDsetRsieve}\inf_{h \in V_n(\theta,R|\ell_n)}  \|\WPT(\theta_{0n}) + \DerP(\theta_{0n})[h]\|_{\PSigma,p} + o_P(a_n),
\end{equation*}
and hence establishes the second claim of the theorem. \qed

\noindent {\sc Proof of Theorem \ref{th:coupsmooth}:} Follows from immediately from Theorem \ref{th:gencoupsmooth}(i) and Assumption \ref{ass:bootrates} implying Assumptions \ref{ass:genneigh} and \ref{ass:genbootrates} are satisfied by setting $\mathcal D_n(\mathbf B,\mathbf E) = \mathcal S_n(\mathbf B,\mathbf E)$, $\mathcal D_n(\mathbf L,\mathbf E) = \mathcal S_n(\mathbf L,\mathbf E)$ and $\theta = \tilde \theta$. \qed

\noindent {\sc Proof of Corollary \ref{th:critval}:} We establish the corollary by appealing to Lemmas \ref{lm:auxbootcrit} and \ref{lm:auxcrit}.
To this end, we first note Theorem \ref{th:coupsmooth} allows us to conclude that
\begin{equation}\label{th:critval1}
\hat U_n(R|\ell_n) \geq \UpS (R|\tilde \ell_n) + o_P(a_n)
\end{equation}
uniformly in $P\in \mathbf P_0$ with $\ell_n \asymp \tilde \ell_n$, while Assumption \ref{ass:bootrates}(ii) implies $K_m\tilde \ell_n^2 \mathcal S_n(\mathbf L,\mathbf E) = o(a_nn^{-\frac{1}{2}})$ and $k_n^{1/p}\sqrt{\log(1+k_n)}B_n\times \sup_{P\in \mathbf P} J_{[\hspace{0.03 in}]}(\tilde \ell_n^{\kappa_\rho}, \mathcal F_n,\|\cdot\|_{P,2}) = o(a_n)$, and hence
\begin{equation}\label{th:critval2}
I_n(R) \leq  \Up (R|\tilde \ell_n) + o_P(a_n)
\end{equation}
uniformly in $P\in \mathbf P_0$ by Theorem \ref{th:localdrift}(i).
Moreover, applying Lemma \ref{lm:auxbootcrit} with $B_n = \hat U_n(R|\ell_n)$, $D_n = \{V_i\}_{i=1}^n$, and $C_{P,n}^\star = \UpS(R|\tilde \ell_n)$ yields, for some $\delta_n = o(1)$, that
\begin{equation}\label{th:critval3}
    \liminf_{n\rightarrow \infty} \inf_{P\in \mathbf P_0}P(\hat c_n + \frac{a_n}{2} > q_{1-\alpha - \delta_n,P}(\UpS(R|\tilde \ell_n))) = 1,
\end{equation}
where $q_{\tau,P}(\UpS(R|\tilde \ell_n))$ denotes the $\tau$ quantile of $\UpS(R|\tilde \ell_n)$.
Since $\UpS(R|\tilde \ell_n)\stackrel{d}{=} \Up(R|\tilde \ell_n)$, results \eqref{th:critval2}, \eqref{th:critval3}, and Assumption \ref{ass:equic} verify the conditions of Lemma \ref{lm:auxcrit} (applied with $T_n = I_n(R)$ and $C_{P,n} = \Up(R|\tilde \ell_n)$) and therefore the corollary follows. \qed

\noindent {\sc Proof of Corollary \ref{th:incJcritval}:} In what follows, we use a ``$\rm u$" superscript for parameters associated with setting $R = \Theta$ -- e.g., $\mathbf B_n^{\rm u}$ denotes the vector subspace generated by $\Theta_n$. 
First note Theorem \ref{th:localdrift}(i) (for $R$ as in \eqref{def:Rset}) and Theorem \ref{th:localdrift}(ii) (for $R=\Theta$) imply that for any $\ell_n,\ell_n^{\rm u}\downarrow 0$ satisfying $k_n^{1/p}\sqrt{\log(1+k_n)}B_n\times \{\sup_{P\in \mathbf P} J_{[\hspace{0.02 in}]}(\ell_n^{\kappa_\rho},\mathcal F_n,\|\cdot\|_{P,2})+ \sup_{P\in \mathbf P} J_{[\hspace{0.02 in}]}((\ell_n^{\text{\rm u}})^{\kappa_\rho},\mathcal F_n^{\text{\rm u}},\|\cdot\|_{P,2})\} = o(a_n)$, $K_m (\ell_n^{\text{\rm u}})^2 \times \mathcal S_n^{\text{\rm u}}(\mathbf L,\mathbf E) = o(a_nn^{-1/2})$, $K_m \ell_n^2 \times \mathcal S_n(\mathbf L,\mathbf E) = o(a_n n^{-1/2})$ and $\mathcal R_n^{\text{\rm u}} = o(\ell_n^{\text{\rm u}})$ it follows uniformly in $P\in \mathbf P_0$ that
\begin{equation}\label{th:incJcritval1}
I_n(R)-I_n(\Theta) \leq \Up(R|\ell_n)-\Up(\Theta|\ell_n^{\rm u}) + o_P(a_n).
\end{equation}
Next note that we may apply Theorem \ref{th:coupsmooth} to obtain that uniformly in $P\in \mathbf P_0$ we have
\begin{equation}\label{cor:incJcritval2}
\hat U_n(R|\ell_n) \geq \UpS(R|\tilde \ell_n) + o_P(a_n)
\end{equation}
with $\tilde \ell_n \asymp \ell_n$.
Similarly, also note that Lemma \ref{lm:forincJboot} implies uniformly in $P\in \mathbf P_0$ that
\begin{equation}\label{cor:incJcritval3}
\hat U_n(\Theta|+\infty) \leq \UpS (\Theta|\tilde \ell_n^{\text{u}}) + o_P(a_n),
\end{equation}
for $\tilde \ell_n^{\text{u}}\downarrow 0$ satisfying Assumption \ref{ass:bootrates}(ii) and $\mathcal R_n^{\rm u} = o(\tilde \ell_n^{\rm u})$. 
In particular, it follows from results \eqref{cor:incJcritval2} and \eqref{cor:incJcritval3} that uniformly in $P\in \mathbf P_0$ we have
\begin{equation}\label{cor:incJcritval4}
\hat U_n(R|\ell_n) -\hat U_n(\Theta|+\infty) \geq \UpS(R|\tilde \ell_n) - \UpS(\Theta|\tilde \ell_n^{\rm u})+ o_P(a_n)
\end{equation}
for sequences $\tilde \ell_n,\tilde \ell_n^{\rm u}\downarrow 0$ satisfying the rate requirements needed for \eqref{th:incJcritval1} to hold (i.e.\ with $\ell_n,\ell_n^{\rm u}$ replaced by $\tilde \ell_n,\tilde \ell_n^{\rm u}$).
The corollary then follows by the same arguments as in Corollary \ref{th:critval} but employing \eqref{th:incJcritval1} and \eqref{cor:incJcritval4} in place of \eqref{th:critval1} and \eqref{th:critval2}. \qed

\begin{lemma}\label{lm:lnotbind}
Suppose there is a $\mathcal A_n(P) \subseteq \Theta_n \cap R$ such that $\|h\|_{\mathbf E} \leq \nu_n\|\DerP(\theta)[h]\|_p$ for all $\theta \in \mathcal A_n(P)$ and $h\in \sqrt n\{\mathbf B_n \cap R - \theta\}$.
If the estimator $\hat {\mathbb D}_n(\theta)$ satisfies
\begin{equation}\label{lm:lnotbinddisp1}
\sup_{\theta \in \mathcal A_n(P)} \sup_{h  \in \sqrt n\{\mathbf B_n\cap R - \theta\} :  \|\frac{h}{\sqrt n}\|_{\mathbf B} \geq \ell_n} \frac{\|\hat {\mathbb D}_n(\theta)[h] - \DerP(\theta)[h]\|_p}{\|h\|_{\mathbf E}} = o_P(\nu_n^{-1})
\end{equation}
and $\hat \Theta_n^{\text{\rm r}} \subseteq \mathcal A_n(P)$ with probability tending to one uniformly in $P\in \mathbf P_0$, Assumptions \ref{ass:startreg}(i)(iii), \ref{ass:weights}, \ref{ass:bootcoupling}  hold, and $\mathcal S_n(\mathbf B, \mathbf E)\mathcal R_n = o(\ell_n)$, then uniformly in $P\in \mathbf P_0$
\begin{equation}\label{lm:lnotbinddisp2}
\hat U_n(R|\ell_n) = \inf_{\theta \in \hat \Theta_n^{\text{\rm r}}}\inf_{h \in \hat V_n(\theta,R|+\infty)}\|\Bemp(\theta) + \hat {\mathbb D}_n(\theta)[h]\|_{\hat \Sigma_n,p} + o_P(a_n).
\end{equation}
\end{lemma}

\noindent {\sc Proof:} In the following arguments, we note that the only requirement on $\hat {\mathbb D}_n(\theta)$ is that it satisfy condition \eqref{lm:lnotbinddisp1}.
As a result, the lemma applies to estimators $\hat{\mathbb D}_n(\theta)$ besides the numerical derivative examined in the main text.

In order to establish the result, we first let $\hat \theta_n \in \hat \Theta_n^{\text{r}}$ and $\hat h_n \in \hat V_n(\hat \theta_n,R|+\infty)$ satisfy
\begin{equation*}
\inf_{\theta \in \hat \Theta_n^{\text{r}}}\inf_{h \in \hat V_n(\theta,R|+\infty)}\|\Bemp(\theta) + \hat {\mathbb D}_n(\theta)[h]\|_{\hat \Sigma_n,p} =
\|\Bemp(\hat \theta_n) + \hat {\mathbb D}_n(\hat \theta_n)[\hat h_n]\|_{\hat \Sigma_n,p} + o(a_n) .
\end{equation*}
Then note that in order to establish the claim of the lemma it suffices to show that
\begin{equation}\label{lm:lnotbind2}
\limsup_{n \rightarrow \infty} \sup_{P\in \mathbf P_0} P(\|\frac{\hat h_n}{\sqrt n}\|_{\mathbf B} \geq \ell_n) = 0 .
\end{equation}
To this end, note $0 \in \hat V_n(\theta,R|+\infty)$ for all $\theta \in \Theta_n \cap R$, the triangle inequality, $\|\hat \Sigma_n\|_{o,p} = O_P(1)$ uniformly in $P\in \mathbf P$ by Corollary \ref{cor:weights}, and Assumption \ref{ass:bootcoupling} yield
\begin{align}\label{lm:lnotbind3}
\|\hat {\mathbb D}_n(\hat \theta_n)[\hat h_n] \|_{\hat \Sigma_n,p} & \leq \|\Bemp(\hat \theta_n) + \hat {\mathbb D}_n(\hat \theta_n)[\hat h_n] \|_{\hat \Sigma_n,p} + \|\Bemp(\hat \theta_n) \|_{\hat \Sigma_n, p} \notag \\  & \leq  2 \|\hat \Sigma_n\|_{o,p} \|\WPT(\hat \theta_n)  \|_p + o_P(a_n)
\end{align}
uniformly in $P\in \mathbf P$.
Hence, since $\hat \theta_n \in \hat \Theta_n^{\text{r}} \subseteq \Theta_n \cap R$ almost surely, we obtain from result \eqref{lm:lnotbind3}, $\|\hat \Sigma_n\|_{o,p} = O_P(1)$ uniformly in $P\in \mathbf P$, and Lemma \ref{aux:tailbound} that 
\begin{equation}\label{lm:lnotbind4}
\|\hat {\mathbb D}_n(\hat \theta_n)[\hat h_n] \|_{\hat \Sigma_n,p} \leq 2 \|\hat \Sigma_n\|_{o,p}  \sup_{\theta \in \Theta_n \cap R} \|\WPT(\theta)  \|_p + o_P(a_n)  = O_P(k_n^{1/p} \sqrt{\log(1+k_n)} B_nJ_n)
\end{equation}
uniformly in $P\in \mathbf P$.
Since $\hat h_n\in \hat V_n(\hat \theta_n,R|+\infty)$ implies $\hat h_n \in \sqrt n\{\mathbf B_n \cap R - \hat \theta_n\}$ and $\hat \theta_n \in \hat \Theta_n^{\text{r}} \subseteq \mathcal A_n(P)$ with probability tending to one uniformly in $P\in \mathbf P_0$, we obtain from the first hypothesis of the lemma that $\|\hat h_n\|_{\mathbf E} \leq \nu_n \|\DerP(\hat \theta_n)[\hat h_n]\|_p$ with probability tending to one uniformly in $P\in \mathbf P_0$. Therefore, it follows that
\begin{align}\label{lm:lnotbind5}
\limsup_{n \rightarrow \infty} & \sup_{P\in \mathbf P_0} P(\ell_n \leq \|\frac{\hat h_n}{\sqrt n}\|_{\mathbf B}) \nonumber \\
& = \limsup_{n \rightarrow \infty} \sup_{P\in \mathbf P_0} P( \ell_n \leq \|\frac{\hat h_n}{\sqrt n}\|_{\mathbf B} \text{ and } \|\hat h_n\|_{\mathbf E} \leq \nu_n \|\DerP(\hat \theta_n)[\hat h_n]\|_{p}) \nonumber \\
& \leq \limsup_{n \rightarrow \infty} \sup_{P\in \mathbf P_0} P(\ell_n \leq \|\frac{\hat h_n}{\sqrt n}\|_{\mathbf B} \text{ and } \|\hat h_n\|_{\mathbf E} \leq 2\nu_n \|\hat{\mathbb D}_{n}(\hat \theta_n)[\hat h_n]\|_{p}) ,
\end{align}
where the inequality follows from condition \eqref{lm:lnotbinddisp1}.
Hence, results \eqref{lm:lnotbind4} and \eqref{lm:lnotbind5}, the definitions of $\mathcal S_n(\mathbf B,\mathbf E)$ and $\mathcal R_n$, and $\mathcal S_n(\mathbf B,\mathbf E) \mathcal R_n = o(\ell_n)$ by hypothesis yield
\begin{multline}\label{lm:notbind6}
\limsup_{n\rightarrow \infty}\sup_{P\in \mathbf P_0} P( \ell_n \leq \|\frac{\hat h_n}{\sqrt n}\|_{\mathbf B}) \\ \leq \limsup_{n \rightarrow \infty} \sup_{P\in \mathbf P_0} P(\ell_n \leq 2 \frac{\nu_n}{\sqrt n} \mathcal S_n(\mathbf B,\mathbf E) \|\hat{\mathbb D}_{n}(\hat \theta_n)[\hat h_n]\|_{p}) = 0 ,
\end{multline}
which establishes \eqref{lm:lnotbind2} and hence the claim of the lemma. \qed

\begin{lemma}\label{lm:simcoup}
Let Assumptions \ref{ass:param}(i), \ref{ass:startreg}, \ref{ass:coupling}, \ref{ass:keycons}(i),  \ref{ass:driftlin}(i), \ref{ass:locrates}(ii), \ref{ass:weights}, \ref{ass:bootcoupling}, \ref{ass:extra} hold and $\mathcal R_n \vee \nu_n\tau_n = o(1)$. 
If $\ell_n \downarrow 0$ satisfies $k_n^{1/p}\sqrt{\log(1+k_n)}B_n\sup_{P\in \mathbf P}J_{[\hspace{0.02 in}]}(\ell_n^{\kappa_\rho},\mathcal F_n,\|\cdot\|_{P,2}) = o(a_n)$ and $K_m \ell^2_n\mathcal S_n(\mathbf L,\mathbf E) = o(a_nn^{-\frac{1}{2}})$, then uniformly in $P\in \mathbf P_0$ we have
\begin{equation*}
\hat U_n(R|\ell_n) = \inf_{\theta \in \hat \Theta_n^{\text{\rm r}}} \inf_{h \in \hat V_n(\theta,R|\ell_n)} \|\WPT(\theta) + \DerP(\theta)[h]\|_{\PSigma,p} + o_P(a_n) .
\end{equation*}
\end{lemma}

\noindent {\sc Proof:} First note that Corollary \ref{cor:bootsetrates}(i) implies $\overrightarrow{d}_H(\hat \Theta_n^{\text{r}},\IDsetRsieve,\|\cdot\|_{\mathbf E}) = O_P(\mathcal R_n + \nu_n\tau_n)$ uniformly in $P\in \mathbf P_0$.
Hence, since $\mathcal R_n \vee \nu_n \tau_n = o(1)$, for any $\epsilon > 0$ it follows that
\begin{equation}\label{lm:simcoup1}
\liminf_{n\rightarrow \infty}\inf_{P\in \mathbf P_0} P(\hat \Theta_n^{\text{r}} \subseteq \{\theta \in \Theta_n \cap R : \overrightarrow{d}_H(\theta,\IDsetRsieve,\|\cdot\|_{\mathbf E}) \leq \epsilon\}) = 1.
\end{equation}
Furthermore, for any $\theta \in \hat\Theta_n^{\text{r}}$ and $h\in \hat V_n(\theta, R| \ell_n)$ note that $\Upsilon_G(\theta + h/\sqrt n) \leq 0$ and $\Upsilon_F(\theta + h/\sqrt n) = 0$ by definition of $\hat V_n(\theta,R|\ell_n)$.
Thus, $\theta + h/\sqrt n \in R$ for any $\theta \in \hat \Theta_n^{\text{r}}$ and $h \in \hat V_n(\theta,R|\ell_n)$, and hence Assumption \ref{ass:extra}(ii) allows us to conclude
\begin{multline}\label{lm:simcoup2}
\liminf_{n\rightarrow \infty}\inf_{P\in \mathbf P_0} P(\theta + \frac{h}{\sqrt n} \in \Theta_n \cap R\text{ for all } \theta \in \hat \Theta_n^{\text{r}} \text{ and } h\in \hat V_n(\theta, R|\ell_n)) \\
= \liminf_{n\rightarrow \infty}\inf_{P\in \mathbf P_0} P(\theta + \frac{h}{\sqrt n} \in \Theta_n \text{ for all } \theta \in \hat \Theta_n^{\text{r}} \text{ and } h \in \hat V_n(\theta,R| \ell_n))  = 1
\end{multline}
due to $\|h/\sqrt n\|_{\mathbf B} \leq \ell_n\downarrow 0 $ for any $h\in \hat V_n(\theta, R|\ell_n)$.
In particular, note that result \eqref{lm:simcoup2} and Assumption \ref{ass:extra}(i) imply that for some $ M <\infty$ we have $\hat V_n(\theta,R|\ell_n)\subseteq V_n(\theta,R|\ell_n/M)$ for all $\theta \in \hat \Theta_n^{\rm r}$ with probability tending to one uniformly in $P\in \mathbf P_0$.
Thus, \eqref{lm:simcoup1} and Lemma \ref{aux:derest} allow us to conclude that uniformly in $P\in \mathbf P_0$ we have
\begin{equation}\label{lm:simcoup3}
\sup_{\theta \in \hat \Theta_n^{\text{r}}} \sup_{h \in \hat V_n(\theta,R|\ell_n)}  \|\hat {\mathbb D}_n(\theta)[h] - \DerP(\theta)[h]\|_p = o_P(a_n) .
\end{equation}
Moreover, since $\hat \Theta_n^{\text{r}} \subseteq \Theta_n \cap R$ almost surely, we also have from Assumption \ref{ass:bootcoupling} that
\begin{equation}\label{lm:simcoup4}
\sup_{\theta \in \hat \Theta_n^{\text{r}}} \|\Bemp(\theta) - \WPT(\theta)\|_p   = o_P(a_n)
\end{equation}
uniformly in $P\in \mathbf P$.
Therefore, since $\|\hat \Sigma_n\|_{o,p} = O_P(1)$ uniformly in $P\in \mathbf P$ by Corollary \ref{cor:weights}, we obtain from results \eqref{lm:simcoup3} and \eqref{lm:simcoup4} and Lemma \ref{lm:geninf} that
\begin{equation}\label{lm:simcoup5}
\hat U_n(R|\ell_n) = \inf_{\theta \in \hat \Theta_n^{\text{r}}} \inf_{h \in \hat V_n(\theta,R|\ell_n)} \|\WPT(\theta) + \DerP(\theta)[h]\|_{\hat \Sigma_n,p} + o_P(a_n)
\end{equation}
uniformly in $P\in \mathbf P_0$.
Next, note that by Assumption \ref{ass:weights}(ii) there exists a constant $C_0 < \infty$ such that $\|\PSigma^{-1}\|_{o,p} \leq C_0$ for all $P\in \mathbf P$.
Thus, using that $\|\hat \Sigma_n a\|_p \leq \|\hat \Sigma_n\PSigma^{-1}\|_{o,p} \|\PSigma a\|_p$ for any $a\in \mathbf R^{k_n}$ and the triangle inequality we obtain
\begin{equation}\label{lm:simcoup6}
\|\WPT(\theta) + \DerP(\theta)[h]\|_{\hat \Sigma_n,p} \leq \{C_0\|\hat \Sigma_n - \PSigma\|_{o,p} + 1\}\|\WPT(\theta) + \DerP(\theta)[h]\|_{\PSigma,p}
\end{equation}
for any $\theta \in \Theta_n \cap R$, $h\in \mathbf B_n$, and $P\in \mathbf P$.
In particular, since $0 \in \hat V_n(\theta, R|\ell_n)$ for any $\theta \in \Theta_n \cap R$, Assumption \ref{ass:weights}, Markov's inequality, and Lemma \ref{aux:tailbound} yield
\begin{multline}\label{lm:simcoup7}
\|\hat \Sigma_n - \PSigma\|_{o,p} \times \inf_{\theta \in \hat \Theta_n^{\text{r}}} \inf_{h\in \hat V_n(\theta,R| \ell_n)} \|\WPT(\theta)  + \DerP(\theta)[h]\|_{\PSigma,p} \\ \leq \|\hat \Sigma_n - \PSigma\|_{o,p}\times \sup_{\theta \in \Theta_n \cap R} \|\WPT(\theta)\|_{\PSigma,p} = o_P(a_n)
\end{multline}
uniformly in $P\in \mathbf P$. It then follows from \eqref{lm:simcoup6} and \eqref{lm:simcoup7} that uniformly in $P\in \mathbf P$
\begin{multline}\label{lm:simcoup8}
\inf_{\theta \in \hat \Theta_n^{\text{r}}} \inf_{h \in \hat V_n(\theta, R|\ell_n)} \|\WPT(\theta) + \DerP(\theta)[h]\|_{\hat \Sigma_n,p} \\ \leq \inf_{\theta \in \hat \Theta_n^{\text{r}}} \inf_{h \in \hat V_n(\theta,R|\ell_n)} \|\WPT(\theta) + \DerP(\theta)[h]\|_{\PSigma,p} + o_P(a_n) .
\end{multline}
The reverse inequality to \eqref{lm:simcoup8} can be obtained by identical arguments but employing $\max\{\|\hat \Sigma_n\|_{o,p},\|\hat \Sigma_n^{-1}\|_{o,p}\} = O_P(1)$ uniformly in $P\in \mathbf P$ by Corollary \ref{cor:weights} instead of $\|\PSigma\|_{o,p}\vee\|\PSigma^{-1}\|_{o,p}$ being bouded uniformly in $P\in \mathbf P$. The claim of the Lemma then follows from \eqref{lm:simcoup5} and \eqref{lm:simcoup8} (and its reverse inequality). \qed

\begin{lemma}\label{aux:derest}
Let Assumptions \ref{ass:startreg}(i)(ii), \ref{ass:coupling}, and \ref{ass:driftlin}(i) hold, and define the sets
\begin{equation}\label{aux:derestdisp1}
V_{n}(\theta,R|\ell_n) \equiv\{h \in \mathbf B_n : \theta + \frac{h}{\sqrt n} \in \Theta_n \cap R \text{ and } \|\frac{h}{\sqrt n}\|_{\mathbf E} \leq \ell_n \} .
\end{equation}
If $\ell_n\downarrow 0$ satisfies $k_n^{1/p}\sqrt{\log(1+k_n)}B_n \times \sup_{P\in \mathbf P} J_{[\hspace{0.03 in}]}(\ell_n^{\kappa_\rho},\mathcal F_n,\|\cdot\|_{P,2}) = o(a_n)$ and $K_m \ell_n^2\times  \mathcal S_n(\mathbf L,\mathbf E) = o(a_nn^{-\frac{1}{2}})$, then there is an $\epsilon > 0$ such that uniformly in $P\in \mathbf P_0$:
\begin{equation}\label{aux:derestdisp}
\sup_{\theta \in \Theta_n \cap R : \overrightarrow{d}_H(\theta,\IDsetRsieve,\|\cdot\|_{\mathbf E}) \leq \epsilon} \sup_{h\in V_n(\theta,R|\ell_n)} \|\hat {\mathbb D}_{n}(\theta)[h] - \DerP(\theta)[h]\|_{p} = o_P(a_n) .
\end{equation}
\end{lemma}

\noindent {\sc Proof:} By definition of $V_{n}(\theta,R|\ell_n)$, we have $\theta + h/\sqrt n \in \Theta_n \cap R$ for any $\theta \in \Theta_n \cap R$, $h \in V_n(\theta,R|\ell_n)$.
Therefore, since $\|h/\sqrt n\|_{\mathbf E} \leq \ell_n$ for all $h \in V_n(\theta,R|\ell_n)$ we obtain
\begin{align}\label{aux:derest1}
\sup_{\theta \in \Theta_n \cap R}\sup_{h \in V_{n}(\theta,R|\ell_n)} & \|\hat {\mathbb D}_n(\theta)[h] - \sqrt nE_P[(\rho(X,\theta + \frac{h}{\sqrt n}) - \rho(X,\theta))*q^{k_n}(Z)]\|_{p} \notag \\
& \leq \sup_{\theta_1,\theta_2\in\Theta_n\cap R : \|\theta_1-\theta_2\|_{\mathbf E} \leq \ell_n}  \|\Gemp(\theta_1) - \Gemp(\theta_2)\|_{p} \notag \\
& \leq \sup_{\theta_1,\theta_2\in\Theta_n\cap R : \|\theta_1-\theta_2\|_{\mathbf E} \leq \ell_n}  \|\WP(\theta_1) - \WP(\theta_2)\|_p + o_P(a_n)
\end{align}
uniformly in $P\in \mathbf P$ by Assumption \ref{ass:coupling}(i).
Next note Assumptions \ref{ass:startreg}(i) and \ref{ass:coupling}(ii) imply that for any $1\leq \jmath \leq \mathcal J$ and $1\leq  k \leq k_{n,\jmath}$ we must have
\begin{equation}\label{aux:derest2}
\sup_{P\in \mathbf P} \sup_{\theta_1,\theta_2 \in \Theta_n \cap R : \|\theta_1-\theta_2\|_{\mathbf E} \leq \ell_n} E_P[\|\rho(X,\theta_1) - \rho(X,\theta_2)\|^2_2q_{k,\jmath}^2(Z_{\jmath})] \leq B_n^2K_\rho^2\ell_n^{2\kappa_\rho} .
\end{equation}
Define $\mathcal G_n \equiv \{fq_{k,\jmath} : f\in\mathcal F_n, 1\leq\jmath\leq \mathcal J \text{ and } 1\leq k \leq k_{n,\jmath}\}$ and let $\Iso$ be a Gaussian process on $\mathcal G_n$ satisfying $E[\Iso (g_1)\Iso(g_2)] = \text{Cov}_P\{g_1(V),g_2(V)\}$ and $E[\Iso(g_1)] = 0$ for any $g_1,g_2\in \mathcal G_n$. 
By result \eqref{aux:derest2} and $\|a\|_p \leq k_n^{1/p}\|a\|_\infty$ for any $a\in \mathbf R^{k_n}$ we obtain
\begin{multline}\label{aux:derest3}
E[\sup_{\theta_1,\theta_2\in\Theta_n\cap R : \|\theta_1-\theta_2\|_{\mathbf E} \leq \ell_n} \|\WP(\theta_1) - \WP(\theta_2) \|_{p}] \\
\leq k_n^{1/p} \times E[\sup_{g_1,g_2 \in \mathcal G_n : \|g_1 - g_2\|_{P,2} \leq B_nK_\rho \ell_n^{\kappa_\rho}} |\Iso(g_1) - \Iso(g_2)| ].
\end{multline}
Therefore, the calculations in \eqref{lm:locprocess4}, Markov's inequality, and $k_n^{1/p}\sqrt{\log(1+k_n)}B_n\times \sup_{P\in \mathbf P} J_{[\hspace{0.03 in}]}(\ell_n^{\kappa_\rho},\mathcal F_n,\|\cdot\|_{P,2}) = o(a_n)$ by hypothesis, yield that
\begin{equation}\label{aux:derest4}
\sup_{\theta \in \Theta_n \cap R}\sup_{h \in V_n(\theta,R|\ell_n)} \|\hat {\mathbb D}_n(\theta)[h] - \sqrt nE_P[(\rho(X,\theta + \frac{h}{\sqrt n}) - \rho(X,\theta))*q^{k_n}(Z)]\|_{p} = o_P(a_n)
\end{equation}
uniformly in $P\in \mathbf P$.
Next, let $\epsilon > 0$ be sufficiently small for Assumption \ref{ass:driftlin}(i) to hold and define the neighborhood $\mathcal N_n \equiv \{\theta \in \Theta_n \cap R : \overrightarrow d_H(\{\theta\},\IDsetRsieve,\|\cdot\|_{\mathbf E}) \leq \epsilon\}$.
We can then conclude from Lemmas \ref{lm:obviousineq} and \ref{aux:bessel}, and Assumption \ref{ass:driftlin}(i) that
\begin{multline}\label{aux:derest5}
\sup_{\theta \in \mathcal N_n}  \sup_{h \in V_n(\theta,R|\ell_n)}\|\sqrt nE_P[(\rho(X,\theta + \frac{h}{\sqrt n})-\rho(X,\theta))*q^{k_n}(Z)] - \DerP(\theta)[h]\|_{p} \\
\lesssim  \sup_{\theta \in \mathcal N_n} \sup_{h \in V_n(\theta,R|\ell_n)} \{K_m  \times \sqrt n \|\frac{h}{\sqrt n}\|_{\mathbf E} \|\frac{h}{\sqrt n}\|_{\mathbf L}\} = o(a_n),
\end{multline}
where the final equality follows from $K_m \ell_n^2\times \mathcal S_n(\mathbf L,\mathbf E) = o(a_nn^{-1/2})$ by hypothesis.
Hence, the Lemma follows from results \eqref{aux:derest4} and \eqref{aux:derest5}. \qed

\begin{lemma}\label{aux:dererror}
Let Assumptions \ref{ass:startreg}(ii) and \ref{ass:driftlin}(ii)(iii) hold. Then there are constants $\epsilon > 0$ and $C<\infty$ such that for all $n$, $P\in \mathbf P$, $\theta_0 \in \IDsetRsieve$, $\theta_1 \in \Theta_n \cap R$ satisfying $\overrightarrow{d}_H(\theta_1,\IDsetRsieve,\|\cdot\|_{\mathbf E}) \leq \epsilon$, and $h_0,h_1\in \mathbf B_n$ it follows that
$$\|\DerP(\theta_0)[h_0] - \DerP(\theta_1)[h_1]\|_p \leq C\{\|h_0 - h_1\|_{\mathbf E} + K_m\|\theta_0 - \theta_1\|_{\mathbf L}\|h_1\|_{\mathbf E}\} .$$
\end{lemma}

\noindent {\sc Proof:} We first fix $\epsilon > 0$ such that Assumptions \ref{ass:driftlin}(ii)(iii) are satisfied.
Then note that by Lemmas \ref{lm:obviousineq} and \ref{aux:bessel} it follows that there is a constant $C_0 < \infty$ with
\begin{equation*}\label{aux:derror1}
\|\DerP(\theta_0)[h_0] - \DerP(\theta_1)[h_1]\|_{p} \leq \{\sum_{\jmath = 1}^{\mathcal J} C_0 \|\nabla m_{P,\jmath}(\theta_0)[h_0] - \nabla m_{P,\jmath}(\theta_1)[h_1]\|_{P,2}^2 \}^{1/2} .
\end{equation*}
Moreover, since $(h_0 - h_1) \in \mathbf B_n$, we can also conclude from Assumptions \ref{ass:driftlin}(ii)(iii) that
\begin{align*}
\|\nabla m_{P,\jmath}& (\theta_0)[h_0]  - \nabla m_{P,\jmath}(\theta_1)[h_1]\|_{P,2} \notag \\
& \leq \|\nabla m_{P,\jmath}(\theta_0)[h_0-h_1]\|_{P,2} + \|\nabla m_{P,\jmath}(\theta_0)[h_1] - \nabla m_{P,\jmath}(\theta_1)[h_1]\|_{P,2} \notag \\
& \leq M \|h_0 - h_1\|_{\mathbf E} + K_m\|\theta_1 - \theta_0\|_{\mathbf L} \|h_1\|_{\mathbf E}
\end{align*}
for some $M < \infty$, and therefore the claim of the lemma follows. \qed

\begin{lemma}\label{lm:auxbootcrit}
Let $B_n$ and $D_n$ be observable random variables, $C^\star_{P,n}$ be a potentially unobservable random variable depending on $P\in \mathbf P$, and for any $\alpha \in (0,1)$ define
$$\hat q_\alpha \equiv \inf\{u : P(B_n \leq u |D_n) \geq \alpha \} \hspace{0.5 in} q_{\alpha,P} \equiv \inf\{u : P(C^\star_{P,n} \leq u) \geq \alpha\}.$$
If $B_n \geq C^\star_{P,n} + o_P(a_n)$ (with $a_n > 0$) uniformly in $P\in \mathbf P$ and $C^\star_{P,n}$ is independent of $D_n$, then there exists a $\delta_n \downarrow 0$ such that $\liminf_{n\rightarrow \infty}\inf_{P\in \mathbf P} P(\hat q_\alpha + a_n \geq q_{\alpha - \delta_n,P}) = 1.$
\end{lemma}

\noindent{Proof:} In the statement of the lemma, $\mathbf P$ and $a_n$ represent a generic set of distributions and positive sequence -- i.e.\ they need not be the same as in the main text. 
To establish the result, note Markov's inequality and the law of iterated expectations yield
\begin{equation*}
\limsup_{n\rightarrow \infty} \sup_{P\in \mathbf P} P( P(C^\star_{P,n} > B_n + a_n|D_n) > \epsilon) \\ \leq \limsup_{n\rightarrow \infty}\sup_{P\in \mathbf P} \frac{1}{\epsilon} P(C_{P,n}^\star > B_n + a_n) = 0 ,
\end{equation*}
where the final equality follows from $B_n \geq C^\star_{P,n} +o_P(a_n)$ uniformly in $P\in \mathbf P$ by hypothesis.
Thus, we conclude there exists some sequence $\delta_n \downarrow 0$ such that the event
\begin{equation*}
\Omega_n(P) \equiv \{D_n | P(C^\star_{P,n} > B_n + a_n|D_n) \leq \delta_n\}
\end{equation*}
satisfies $P(\Omega_n(P)^c) = o(1)$ uniformly in $P\in \mathbf P$.
Hence, for any $t\in \mathbf R$ we obtain that
\begin{align}\label{lm:auxbootcrit3}
P(B_n \leq t | D_n) 1\{D_n\in \Omega_n(P)\} & \leq P(B_n \leq t \text{ and } C^\star_{P,n} \leq B_n + a_n|D_n) + \delta_n \nonumber \\
& \leq P(C^\star_{P,n} \leq t + a_n) + \delta_n ,
\end{align}
where in the final inequality we employed that $C^\star_{P,n}$ is independent of $D_n$.
Therefore, setting $t = \hat q_{\alpha}$ in \eqref{lm:auxbootcrit3} implies that, under $\Omega_n(P)$, we have $\hat q_{\alpha} + a_n \geq q_{\alpha-\delta_n,P}$.
Since $\sup_{P\in \mathbf P} P(\Omega_n(P)^c) = o(1)$, the claim of the lemma follows. \qed

\begin{lemma}\label{lm:auxcrit}
Let $T_n \leq C_{P,n} + o_P(a_n)$ uniformly in $P\in \mathbf P_0$ with $0 < a_n = o(1)$, define $q_{\alpha,P} \equiv \inf\{u:P(C_{P,n} \leq u)\geq \alpha\}$, and suppose that, for some $\delta_n \downarrow 0$, $\hat c_n + a_n/2 \geq q_{1-\alpha - \delta_n,P}$ with probability tending to one uniformly in $P\in \mathbf P_0$. If for some $\eta_n \geq 0$
\begin{equation}\label{lm:auxcritdisp}
\limsup_{n\rightarrow \infty} \sup_{P\in \mathbf P_0} P(T_n > \hat c_n) = \limsup_{n\rightarrow \infty} \sup_{P\in \mathbf P_0} P(T_n > \hat c_n \vee \eta_n)
\end{equation}
and for some sequence $\varrho_n$ satisfying $\varrho_n a_n = o(1)$ we have $\sup_{P\in \mathbf P_0} P(|C_{P,n}- t| \leq \epsilon) \leq \varrho_n(\epsilon\wedge 1) + o(1)$  for all $t\in (\eta_n - a_n,+\infty)$, then it follows that
$$\limsup_{n\rightarrow \infty} \sup_{P\in \mathbf P_0} P(T_n > \hat c_n) \leq \alpha.$$
\end{lemma}

\noindent {\sc Proof:} First note that by condition \eqref{lm:auxcritdisp}, $T_n \leq C_{P,n}+o_P(a_n)$ uniformly in $P\in \mathbf P_0$ and the maintained hypothesis on $\hat c_n$ we can conclude that
\begin{align}\label{lm:auxcrit1}
\limsup_{n\rightarrow \infty} \sup_{P\in \mathbf P_0} P(T_n > \hat c_n) & = \limsup_{n\rightarrow \infty} \sup_{P\in \mathbf P_0} P(T_n > \hat c_n \vee \eta_n)  \notag
\\ & \leq \limsup_{n\rightarrow \infty} \sup_{P\in \mathbf P_0}  P(C_{P,n} + \frac{a_n}{2} > (q_{1-\alpha-\delta_n,P} -\frac{a_n}{2})\vee  \eta_n) \notag \\
& \leq \limsup_{n\rightarrow \infty} \sup_{P\in \mathbf P_0} P(C_{P,n} + a_n > q_{1-\alpha-\delta_n,P} \vee \eta_n).
\end{align}
Next observe that by direct calculation we also have the following inequalities
\begin{multline}\label{lm:auxcrit2}
P(C_{P,n} + a_n > q_{1-\alpha-\delta_n,P}\vee \eta_n) -P(C_{P,n} > q_{1-\alpha-\delta_n,P}) \\ \leq \left\{\begin{array}{cl} 0 & \text{ if } \eta_n-a_n \geq q_{1-\alpha-\delta_n,P} \\ P(|C_{P,n}-q_{1-\alpha-\delta_n,P}| \leq a_n) & \text{ if } \eta_n-a_n < q_{1-\alpha-\delta_n,P}\end{array}\right.  .
\end{multline}
Therefore, combining results \eqref{lm:auxcrit1} and \eqref{lm:auxcrit2} together with  $\sup_{P\in \mathbf P_0} P(|C_{P,n}- t| \leq \epsilon) \leq \varrho_n(\epsilon\wedge 1) + o(1)$ for all $t\in (\eta_n - a_n,+\infty)$ implies that
\begin{align}\label{lm:auxcrit3}
\limsup_{n\rightarrow \infty} & \sup_{P\in \mathbf P_0} P(T_n > \hat c_n) \notag \\
& \leq \limsup_{n\rightarrow \infty} \sup_{P\in \mathbf P_0} P(C_{P,n} > q_{1-\alpha-\delta_n,P})  + \limsup_{n\rightarrow \infty} \sup_{P\in \mathbf P_0} \sup_{t > \eta_n - a_n} P(|C_{P,n} -t| \leq a_n) \notag \\
& \leq \alpha + \delta_n + \varrho_n(a_n \wedge 1).
\end{align}
The claim of the lemma therefore follows from $\delta_n =o(1)$ and $\varrho_n a_n = o(1)$. \qed

\begin{lemma}\label{lm:forincJboot}
Let the conditions of Theorems \ref{th:localdrift}(ii) and \ref{th:coupsmooth} hold with $R = \Theta$ and suppose that $\ell_n^{\rm u}$ satisfies $k_n^{1/p}\sqrt{\log(1+k_n)}B_n\times \sup_{P\in \mathbf P}J_{[\hspace{0.03 in}]}((\ell_n^{\rm u})^{\kappa_\rho}\vee(\nu^{\rm u}_n\tau^{\rm u}_n)^{\kappa_\rho},\mathcal F_n^{\rm u},\|\cdot\|_{P,2}) =  o(a_n)$, $K_m \ell_n^{\rm u}( \ell_n^{\rm u} + \mathcal R_n^{\rm u}+ \nu_n^{\rm u}\tau_n^{\rm u})\times \mathcal S_n^{\rm u}(\mathbf L, \mathbf E) =  o(a_nn^{-1/2})$, and $\mathcal R_n^{\rm u} = o(\ell_n^{\rm u})$.
(i) If $\tau_n^{\rm u} \downarrow 0$ satisfies $J_n^{\text{\rm u}}B_n k_n^{1/p}\sqrt{\log(1+k_n)/n} = o(\tau_n^{\text{\rm u}})$ and $\nu_n^{\rm u}\tau_n^{\rm u}\times \mathcal S_n^{\rm u}(\mathbf B,\mathbf E) = o(1)$, then
$$\hat U_n(\Theta|+\infty) \leq \UpS (\Theta| \ell_n^{\text{\rm u}}) + o_P(a_n)$$
uniformly in $P\in \mathbf P_0$.
(ii) If $\mathcal S_n^{\rm u}(\mathbf B,\mathbf E)\times  \mathcal R_n^{\rm u} = o(1)$ and $\IDsetUsieve$ is a singleton for all $P\in \mathbf P_0$ and $n$ sufficiently large, then part (i) of the lemma continues to hold if $\tau_n^{\rm u} = 0$.
\end{lemma}

\noindent {\sc Proof:} First note that since we required $J_n^{\text{u}} B_n k_n^{1/p}\sqrt{\log(1+k_n)/n} = o(\tau_n^{\text{u}})$ and we assumed all other conditions of Corollary \ref{cor:bootsetrates}(ii) are satisfied when $\Theta = R$, it follows
\begin{equation}\label{lm:forincJbootfix1}
d_H(\hat \Theta_n^{\rm u},\IDsetUsieve,\|\cdot\|_{\mathbf E}) = O_P(\nu_n^{\rm u}\tau_n^{\rm u})    
\end{equation}
uniformly in $P\in \mathbf P_0$.
Therefore, Lemma \ref{aux:derest} yields, uniformly in $P\in \mathbf P_0$, that 
\begin{equation}\label{lm:forincJbootfix2}
\sup_{\theta \in \hat \Theta_n^{\rm u}} \sup_{h\in V_n(\theta,\Theta|\ell_n^{\rm u})} \|\hat {\mathbb D}_n(\theta)[h] - \mathbb D_P(\theta)[h]\|_p = o_P(a_n).    
\end{equation}
We further note that since $\hat \Theta_n^{\rm u} \subseteq \Theta_n$, Assumption \ref{ass:bootcoupling} holding with $R = \Theta$ implies
\begin{equation}\label{lm:forincJbootfix3}
\sup_{\theta \in \hat \Theta_n^{\rm u}} \|\Bemp(\theta) - \WPT(\theta)\|_p = o_P(a_n)
\end{equation}
uniformly in $P\in \mathbf P_0$.
Hence, by results \eqref{lm:forincJbootfix2} and \eqref{lm:forincJbootfix3}, $\|\hat \Sigma_n\|_{o,p} = O_P(1)$ uniformly in $P\in \mathbf P_0$ by Corollary \ref{cor:weights}, and $V_n(\theta,\Theta|\ell_n^{\rm u}) \subseteq \hat V_n(\theta,\Theta|+\infty)$ imply that
\begin{align}\label{lm:forincJbootfix4}
\hat U_n(\Theta|+\infty) & \leq \inf_{\theta \in \hat \Theta_n^{\rm u}} \inf_{h \in V_n(\theta,\Theta|\ell_n^{\rm u})} \|\WPT(\theta) + {\mathbb D}_P(\theta)[h]\|_{\hat \Sigma_n,p} + o_P(a_n) \notag \\
& = \inf_{\theta \in \hat \Theta_n^{\rm u}} \inf_{h \in V_n(\theta,\Theta|\ell_n^{\rm u})} \|\WPT(\theta) + {\mathbb D}_P(\theta)[h]\|_{\Sigma_P,p} + o_P(a_n) 
\end{align}
uniformly in $P\in \mathbf P_0$, and where the equality can be established by employing identical arguments to those used in Lemma \ref{lm:simcoup} (see, in particular, \eqref{lm:simcoup6}-\eqref{lm:simcoup8}).
Also note that, by hypothesis, there is an $\eta_n \downarrow 0$ satisfying $\nu_n^{\rm u}\tau_n^{\rm u}\times \mathcal S_n^{\rm u}(\mathbf B,\mathbf E) = o(\eta_n)$ and define
\begin{equation*}
\mathcal E_n(\theta) \equiv V_n(\theta,\Theta|\ell_n^{\rm u}) \cap \{h \in \mathbf B_n^{\rm u} : \|\frac{h}{\sqrt n}\|_{\mathbf B} \leq \eta_n\}.
\end{equation*}
Next, select $\theta_{0n} \in \IDsetUsieve$ and $h_{0n}\in \mathcal E_n(\theta_{0n})$ so that the following equality is satisfied
\begin{equation}\label{lm:forincJboot2}
\inf_{\theta \in \IDsetUsieve} \inf_{h\in \mathcal E_n(\theta)} \|\WPT(\theta) + \DerP(\theta)[h]\|_{\PSigma,p}
= \|\WPT(\theta_{0n}) + \DerP(\theta_{0n})[h_{0n}]\|_{\PSigma,p} + o(a_n).
\end{equation}
Assumption \ref{ass:bootrates} holding with $R = \Theta$ implies $K_m \ell_n^{\text{u}}(\nu_n^{\text{u}}\tau_n^{\text{u}})  \mathcal S_n^{\text{u}}(\mathbf L,\mathbf E) = o(a_nn^{-1/2})$ and $k_n^{1/p}\sqrt{\log(1+k_n)}B_n \sup_{P\in \mathbf P} J_{[\hspace{0.02 in}]}((\nu_n^{\rm u}\tau_n^{\rm u})^{\kappa_\rho},\mathcal F_n^{\rm u},\|\cdot\|_{P,2}) = o(a_n)$.
Hence, there is $\delta_n$ with
\begin{align}
& K_m \delta_n \ell_n^{\text{u}} \mathcal S_n^{\text{u}}(\mathbf L,\mathbf E)  = o(a_n n^{-1/2}) \label{lm:forincJboot3p1} \\
& k_n^{1/p} \sqrt{\log(1+k_n)}B_n \times \sup_{P\in \mathbf P} J_{[\hspace{0.02 in}]}(\delta_n^{\kappa_\rho},\mathcal F_n^{\text{u}},\|\cdot\|_{P,2})  = o(a_n), \label{lm:forincJboot3p2}
\end{align}
and $\nu_n^{\text{u}}\tau_n^{\text{u}} = o(\delta_n)$.
Moreover, note result \eqref{lm:forincJbootfix1} implies there is a $\hat \theta_{0n}$ in $\hat \Theta_n^{\text{u}}$ such that 
\begin{equation*}
\|\theta_{0n} - \hat \theta_{0n}\|_{\mathbf E} = O_P(\nu_n^{\text{u}}\tau_n^{\text{u}})
\end{equation*}
uniformly in $P\in \mathbf P_0$.
Thus, $\nu_n^{\text{u}}\tau_n^{\text{u}} = o(\delta_n)$ and $\hat \theta_{0n} \in \hat \Theta_n^{\text{u}} \subseteq \Theta_n$ implies that $\sqrt n(\hat \theta_{0n} - \theta_{0n}) \in V_n(\theta_{0n},\Theta|\delta_n)$ with probability tending to one uniformly in $P \in \mathbf P_0$.
Hence, applying Lemma \ref{lm:locprocess} with $\IDsetUsieve$ and $V_n(\theta,\Theta|\delta_n)$ in place of $\IDsetRsieve$ and $V_n(\theta,R|\delta_n)$, yields
\begin{equation}\label{lm:forincJboot5}
\|\WPT(\hat \theta_{0n})  - \WPT(\theta_{0n}) \|_p = o_P(a_n)
\end{equation}
 uniformly in $P\in \mathbf P_0$.
Furthermore, Lemma \ref{aux:dererror}, $h_{0n}\in \mathcal E_n(\theta_{0n})$ and result \eqref{lm:forincJboot3p1} imply that with probability tending to one uniformly in $P\in \mathbf P_0$ we must have
\begin{equation}\label{lm:forincJboot6}
\|\DerP(\hat \theta_{0n})[h_{0n}] - \DerP( \theta_{0n})[h_{0n}]\|_p \leq K_m\mathcal S_n^{\text{u}}(\mathbf L,\mathbf E) \delta_n \ell_n^{\text{u}} \sqrt n = o(a_n).
\end{equation}
Therefore, Assumption \ref{ass:weights}(ii), $\hat \theta_{0n} \in \hat \Theta_n^{\text{u}}$, $h_{0n} \in \mathcal E_n(\theta_{0n})$, $\mathcal E_n(\theta_{0n})\subseteq V_n(\hat \theta_{0n},\Theta|\ell_n^{\rm u})$ by Assumption \ref{ass:extra}(ii), and results \eqref{lm:forincJboot2}, \eqref{lm:forincJboot5}, and \eqref{lm:forincJboot6} yield that
\begin{multline}\label{lm:forincJboot7}
\inf_{\theta \in \IDsetUsieve} \inf_{h \in \mathcal E_n(\theta)} \|\WPT(\theta) + \DerP(\theta)[h]\|_{\PSigma,p} \\
\geq \inf_{\theta \in \hat \Theta_n^{\text{u}}} \inf_{h \in V_n(\theta,\Theta|\ell_n^{\rm u})} \|\WPT(\theta) + \DerP(\theta)[h]\|_{\PSigma,p} +o_P(a_n)
\end{multline}
uniformly in $P\in \mathbf P_0$.
To conclude, note that Assumption \ref{ass:keycons} holding with $R = \Theta$, Corollary \ref{cor:setrates}, and $\mathcal R^{\rm u}_n \times \mathcal S_n^{\rm u}(\mathbf B,\mathbf E) = o(\eta_n)$ due to $\mathcal R_n^{\rm u} = o(\tau_n^{\rm u} \nu_n^{\rm u})$ and $\nu_n^{\rm u}\tau_n^{\rm u}\times \mathcal S_n^{\rm u}(\mathbf B,\mathbf E) = o(\eta_n)$ allow us to conclude that uniformly in $P\in \mathbf P_0$ we have
\begin{align}\label{lm:forincJbootfix6}
I_n(\Theta) & = \inf_{\theta \in \IDsetUsieve} \inf_{ h\in \mathcal E_n(\theta)} \sqrt n Q_n(\theta + \frac{h}{\sqrt n}) + o_P(a_n) \notag \\
& =  \inf_{\theta \in \IDsetUsieve} \inf_{ h\in \mathcal E_n(\theta)} \|\WP(\theta) + \mathbb D_P(\theta)[h]\|_{\Sigma_P,p} + o_P(a_n),
\end{align}
where the second equality follows by identical arguments to those employed in Theorem \ref{th:localdrift}(ii).
Combining result \eqref{lm:forincJbootfix6} with Theorem \ref{th:localdrift}(ii) and employing the fact that $\WPT$ and $\WP$ share the same distribution we thus obtain, uniformly in $P\in \mathbf P_0$, that
\begin{multline}\label{lm:forincJbootfix7}
\inf_{\theta \in \IDsetUsieve}\inf_{h\in V_n(\theta,\Theta|\ell_n^{\rm u})} \|\WPT(\theta) + \DerP(\theta)[h]\|_{\PSigma,p} \\
= \inf_{\theta \in \IDsetUsieve}\inf_{h\in \mathcal E_n(\theta)} \|\WPT(\theta) + \DerP(\theta)[h]\|_{\PSigma,p} + o_P(a_n).
\end{multline}
The claim of part (i) of the lemma therefore follows \eqref{lm:forincJbootfix4}, \eqref{lm:forincJboot7}, and \eqref{lm:forincJbootfix7}.
To establish part (ii) note that if $\IDsetUsieve$ is a singleton, then $\overrightarrow d_H(\hat \Theta_n^{\rm u}, \IDsetUsieve,\|\cdot\|_{\mathbf E}) = d_H(\hat \Theta_n^{\rm u}, \IDsetUsieve,\|\cdot\|_{\mathbf E})$ and therefore Corollary \ref{cor:bootsetrates}(i) implies $d_H(\hat \Theta_n^{\rm u}, \IDsetUsieve,\|\cdot\|_{\mathbf E}) = O_P(\mathcal R_n^{\rm u})$ uniformly in $P\in\mathbf P_0$. 
Part (ii) of the lemma can then be established by replacing $\nu_n^{\rm u}\tau_n^{\rm u}$ with $\mathcal R_n^{\rm u}$ in the arguments employed in establishing part (i). \qed 

\begin{corollary}\label{cor:infactor}
Suppose that $I_n(R) \leq \Up (R|\tilde \ell_n) + o_P(a_n)$ and $\hat U_n(R|\ell_n) \geq \UpS (R|\tilde \ell_n) + o_P(a_n)$ uniformly in $P\in \mathbf P_0$ with $0<a_n =o(1)$, $\Up(R|\tilde \ell_n) \stackrel{d}{=} \UpS (R|\tilde \ell_n)$, and $\UpS (R|\tilde \ell_n)$ independent of $\{V_i\}_{i=1}^n$. Then for any constant $\eta \in (0,\alpha)$ it follows that
\begin{equation*}
\limsup_{n\rightarrow \infty} \sup_{P\in \mathbf P_0} P(I_n(R) > \hat q_{1-\alpha+\eta}(\hat  U_n(R|\ell_n)) + \eta) \leq \alpha .
\end{equation*}
\end{corollary}

\noindent {\sc Proof:} Since $I_n(R) \leq \Up (R|\tilde \ell_n) + o_P(a_n)$ uniformly in $P\in \mathbf P_0$ by hypothesis, we obtain
\begin{align}\label{cor:infactor1}
    \limsup_{n\rightarrow \infty} \sup_{P\in \mathbf P_0} & P(I_n(R) > \hat q_{1-\alpha+\eta}(\hat U_n(R|\ell_n))+\eta) \notag\\
    &\leq  \limsup_{n\rightarrow \infty} \sup_{P\in \mathbf P_0} P(\Up (R|\tilde \ell_n) + a_n > \hat q_{1-\alpha+\eta}(\hat U_n(R|\ell_n)) + \eta) \notag\\
    &\leq  \limsup_{n\rightarrow \infty} \sup_{P\in \mathbf P_0} P(\Up (R|\tilde \ell_n)  > q_{1-\alpha+\eta-\delta_n,P}(\UpS(R|\tilde \ell_n)) +\eta-2a_n) \notag \\
    & \leq \alpha,
\end{align}
where the second inequality holds for $q_{1-\alpha+\eta-\delta_n,P}(\UpS(R|\tilde \ell_n))$ the $1-\alpha+\eta-\delta_n$ quantile of $\UpS(R|\tilde \ell_n)$ and  some $\delta_n = o(1)$ by Lemma \ref{lm:auxbootcrit} applied with $B_n = \hat U_n(R|\ell_n)$, $C^\star_{P,n} = \UpS(R|\tilde \ell_n) $, and $D_n = \{V_i\}_{i=1}^n$. In turn, the final inequality in \eqref{cor:infactor1} follows from $\eta > 0$, $a_n = o(1)$, $\delta_n = o(1)$, and $\Up(R|\tilde \ell_n) \stackrel{d}{=} \UpS (R|\tilde \ell_n)$. \qed

%% file: Appendix/ExamplesProofs/AppExamplesProofs.tex
\section{Illustrative Examples} \label{sec:exproofs}

In this Section, we include the proofs for all the examples discussed in the main text and Supplemental Appendix I -- i.e., the results stated in Section \ref{sec:hetero} of the main text and in Section \ref{sec:examples} of Supplemental Appendix I.

\input{Appendix/ExamplesProofs/ExDemandProofs}  

\input{Appendix/ExamplesProofs/ExGMMProofs}     

\input{Appendix/ExamplesProofs/ExNPIVProofs}    

\input{Appendix/ExamplesProofs/ExQuantProofs}    

%% file: Appendix/ExamplesProofs/ExDemandProofs.tex

\subsection{Proofs for Section \ref{sec:hetero}}

\noindent {\sc Proof of Theorem \ref{th:heteroapprox}:} We establish the claim of the theorem by verifying the conditions of Theorem \ref{th:localdrift}(ii) for both $R$ as in \eqref{ex:hetero7paux} (to couple $I_n(R)$) and $R = \Theta$ (to couple $I_n(\Theta)$).
To this end, note that Assumption \ref{ass:param}(i) is imposed in Assumption \ref{ass:heterosieve}(i), Assumption \ref{ass:startreg}(i) holds with $B_n \asymp \sqrt{k_n}$ by Assumption \ref{ass:heteromoments}(i), Assumption \ref{ass:startreg}(ii) is directly imposed in Assumption \ref{ass:heteromoments}(ii), and Assumption \ref{ass:startreg}(iii) is satisfied with $J_n \asymp \sqrt{j_n\log(1+j_n)}$ by Lemma \ref{lm:heterocov} and $\|f\|_\infty \leq 3$ for any $f\in \mathcal F_n$.
The coupling requirement of Assumption \ref{ass:coupling}(i) is satisfied for $R = \Theta$, and hence also for $R$ as in \eqref{ex:hetero7paux}, with $a_n = (\log(n))^{-1/2}$ by Lemma \ref{lm:heterocoup} and Assumption \ref{ass:heteromoments}(iv).
Moreover, Assumptions \ref{ass:coupling}(ii), \ref{ass:keycons}, and \ref{ass:driftlin} also hold by Lemmas  \ref{lm:hetero3ver} and \ref{lm:heteroobviousver}.
To verify Assumption \ref{ass:locrates}, we first note that Assumption \ref{ass:locrates}(ii) is implied by Assumptions \ref{ass:heterosieve}(iv) and \ref{ass:heterosigma}(ii).
Furthermore, as argued, $B_n \asymp \sqrt{k_n}$, $J_n \asymp \sqrt{j_n\log(1+j_n)}$, and $\nu_n \asymp 1$ by Lemma \ref{lm:hetero3ver}, which yields that $\mathcal R_n \lesssim k_n \sqrt{j_n} \log(1+k_n)/\sqrt{n}$ since $k_n \geq j_n$ by Assumption \ref{ass:heteromoments}(iii).
Thus, $\kappa_\rho = 1$ by Lemma \ref{lm:heteroobviousver} and Lemma \ref{lm:heterocov} imply that Assumption \ref{ass:locrates}(i) holds by Assumption \ref{ass:heteromoments}(iv).
By similar arguments, it also follows that Assumption \ref{ass:weights} is implied by Assumption \ref{ass:heterosigma}, and that the requirements $k_n^{1/p}\sqrt{\log(1+k_n)}B_n \sup_{P\in \mathbf P} J_{[\hspace{0.02 in}]}(\ell_n^{\kappa_\rho},\mathcal F_n,\|\cdot\|_{P,2}) = o(a_n)$ and $\mathcal R_n = o(\ell_n)$ are implied by $k_n\sqrt{j_n}\log^2(n) \ell_n = o(1)$ and $k_n \sqrt{j_n}\log(n)/\sqrt n = o(\ell_n)$.
Since $K_m = 0$ in this application, it follows all the conditions of Theorem \ref{th:localdrift}(ii) hold for both $R = \Theta$ and $R$ as in \eqref{ex:hetero7paux}, and hence the theorem follows. \qed

\noindent {\sc Proof of Lemma \ref{lm:heterolipschitz}:} The result essentially follows from Theorem 1 in \cite{walkup1969lipschitzian}.
To map our problem into their setting, note that since $\{\delta_{s}\}_{s=1}^{s_n}$ are orthogonal, every $\mu \in \mathcal M_n$ can be identified with a unique $(\alpha_1,\ldots, \alpha_{s_n})\equiv \alpha \in \mathbf R^{s_n}$ through the relation $\mu = \sum_{s=1}^{s_n} \alpha_s \delta_{s}$ -- e.g., by $\alpha_{s} = \mu(S_{s})$ for $S_{s}$ the support of $\delta_{s}$.
With some abuse of notation, for the remaining of the proof we therefore employ $\alpha$ and $\mu$ interchangeably.
Further note that, for any $\theta = (\{F(c_\jmath|\cdot)\}_{\jmath =1}^{\mathcal J},\mu)$, the restrictions $\Upsilon_G(\theta) \leq 0$, $\Upsilon_F^{(\mu)}(\theta) = 0$, and $\Upsilon_F^{(\text{s})}(\theta) = 0$ depend only on $\mu$ and define a closed convex polyhedron on $\mathbf R^{s_n}$, which we denote by $K_n$.
Next, define the map $\Lambda_n : \mathbf R^{s_n} \to \mathbf R^{\mathcal J \mathcal L}$ to be given by
\begin{equation}\label{lm:heterolipschitz1}
\Lambda_n(\alpha) = \{\sum_{s=1}^{s_n} \alpha_s ( \int 1\{g(w_l,\eta)\leq c_\jmath \} \delta_{s}(d\eta))\}_{1\leq \jmath \leq \mathcal J, 1\leq l \leq \mathcal L}
\end{equation}
and note that for any $\theta = (\{F(c_\jmath|\cdot)\}_{\jmath =1 }^{\mathcal J},\mu) = \theta \in \Theta_n \cap R$, it follows by \eqref{eq:hetero23} that
\begin{equation}\label{lm:heterolipschitz2}
\Gamma_n(\theta) = K_n \cap \Lambda_n^{-1}(\{F(c_\jmath|w_l)\}_{1\leq \jmath \leq \mathcal J,1\leq l\leq \mathcal L}).
\end{equation}
Let $d_n$ denote the dimension of the null space of $\Lambda_n$, and note that if $d_n = s_n$, then $\Gamma_n(\theta_1) = \Gamma_n(\theta_2)$ for any $\theta_1,\theta_2\in \Theta_n \cap R$ by result \eqref{lm:heterolipschitz2}, and hence the conclusion of the lemma is immediate.
On the other hand, if $1\leq d_n \leq s_n - 1$, then Theorem 1 in \cite{walkup1969lipschitzian} implies there is a $C_n$ such that for any $\theta_1,\theta_2\in \Theta_n \cap R$ we have
\begin{align}\label{lm:heterolipschitz3}
d_H(\Gamma_n(\theta_1),\Gamma_n(\theta_2),\|\cdot\|_2) & \leq C_n\{\sum_{\jmath = 1}^{\mathcal J}\sum_{l= 1}^{\mathcal L} (F_1(c_\jmath|w_l) - F_2(c_\jmath|w_l))^2\}^{1/2} \notag \\
& \lesssim C_n \sum_{\jmath = 1}^{\mathcal J} \|F_1(c_\jmath|\cdot) - F_2(c_\jmath|\cdot)\|_\infty,
\end{align}
and where the norm $\|\cdot\|_2$ on $\Gamma_n(\theta)$ is understood as the usual Euclidean norm on the corresponding $\alpha \in \mathbf R^{s_n}$.
Similarly, we note that if $d_n = 0$, then $\Lambda_n$ is invertible and \eqref{lm:heterolipschitz3} holds with $C_n = \|\Lambda_{n}^{-1}\|_o$.
Also note that for any $\mu = \sum_{s=1}^{s_n} \alpha_{s}\delta_{s}$ and $\tilde \mu = \sum_{s=1}^{s_n} \tilde \alpha_{s} \delta_{s}$ we have $\|\mu - \tilde \mu\|_{TV} = \|\alpha - \tilde \alpha\|_1$ due to the measures $\{\delta_{s}\}_{s=1}^{s_n}$ being orthogonal. 
Hence, since $\|a\|_1 \leq \sqrt{s_n}\|a\|_2$ for any $a\in \mathbf R^{s_n}$, result \eqref{lm:heterolipschitz3} yields
\begin{equation*}\label{lm:heterolipschitz4}
d_H(\Gamma_n(\theta_1),\Gamma_n(\theta_2),\|\cdot\|_{TV}) \lesssim \sqrt{s_n}C_n \sum_{\jmath = 1}^{\mathcal J} \|F_1(c_\jmath|\cdot) - F_2(c_\jmath|\cdot)\|_\infty,
\end{equation*}
which establishes the claim of the lemma by setting $\zeta_n \asymp C_n \sqrt{s_n}$. \qed

\noindent {\sc Proof of Theorem \ref{th:heteroboot}:} Let $\hat V_n(\theta,R|\ell) \equiv \hat V_n(\theta,R|+\infty) \cap \{h \in \mathbf B_n : \|h/\sqrt n\|_{\mathbf E} \leq \ell\}$, recall $\|\theta\|_{\mathbf E} = \sum_{\jmath = 1}^{\mathcal J} \sup_{P\in \mathbf P} \|F(c_\jmath|\cdot)\|_{P,2}$ for any $\theta = (\{F(c_\jmath|\cdot)\}_{\jmath = 1}^{\mathcal J},\mu)$, define
\begin{equation*}
\hat E_n(R|\ell_n)  = \inf_{\theta \in \hat \Theta_n^{\text{\rm r}}}\inf_{h \in \hat V_n(\theta,R|\ell_n)} \{\sum_{\jmath = 1}^{\mathcal J} \|\hat {\mathbb W}_{\jmath, n}(\theta) + \hat {\mathbb D}_{\jmath,n}[h]\|_{\hat \Sigma_{\jmath,n},2}\}^{1/2}  ,
\end{equation*}
and note that for any $\ell_n$ satisfying the conditions of the theorem, Assumption \ref{ass:heteroboot}(iii) and Lemma \ref{lm:heteronol} imply $\hat U_n(R|+\infty) = \hat E_n(R|\ell_n) + o_P(a_n)$ uniformly in $P\in \mathbf P_0$.
Hence, to establish the theorem it suffices to show there are $\ell_n \asymp \tilde \ell_n$ and $\ell_n^{\text{\rm u}} \asymp \tilde \ell_n^{\text{\rm u}}$ such that
\begin{align}
\hat E_n(R|\ell_n) & \geq \UpS (R|\tilde \ell_n) + o_P(a_n) \notag \\
\hat E_n(R|\ell_n) - \hat U_n(\Theta|+\infty) & \geq \UpS (R|\tilde \ell_n) - \UpS (\Theta|\tilde \ell_n^{\text{u}}) + o_P(a_n) \label{th:heteroboot0}
\end{align}
uniformly in $P\in \mathbf P_0$.
To this end, we rely on Theorem \ref{th:gencoupsmooth}(ii) (for $\hat E_n(R|\ell_n)$) and Lemma \ref{lm:forincJboot}.
Also note that in the proof of Theorem \ref{th:heteroapprox} we showed Assumptions \ref{ass:heterosieve}, \ref{ass:heteromoments}, and \ref{ass:heterosigma} imply Assumptions \ref{ass:param}-\ref{ass:weights} hold with $B_n \asymp \sqrt {k_n}$, $J_n \asymp \sqrt{j_n\log(1+j_n)}$, $\nu_n \asymp 1$, $\mathcal R_n \asymp k_n \sqrt{j_n\log(1+k_n)\log(1+j_n)/n}$, $a_n = (\log(n))^{-1/2}$, $\kappa_\rho = 1$, $\|\theta\|_{\mathbf L} = \|\theta\|_{\mathbf E} = \sum_{\jmath =1 }^{\mathcal J}\sup_{P\in \mathbf P} \|F(c_\jmath|\cdot)\|_{P,2}$, and $\|\theta\|_{\mathbf B} = \sum_{\jmath=1}^{\mathcal J} \|F(c_\jmath|\cdot)\|_\infty + \|\mu\|_{TV}$ for $R = \Theta$ and $R$ as in \eqref{ex:hetero7paux}.

In order to apply Theorem \ref{th:gencoupsmooth}(ii), we set $\|\theta\|_{\mathbf I} =  \max_{1\leq \jmath \leq \mathcal J} \mathcal J\|F(c_\jmath|\cdot)\|_\infty$ for any $\theta = (\{F(c_\jmath|\cdot)\}_{\jmath = 1}^{\mathcal J},\mu)\in \mathbf B_n$ and note Assumption \ref{ass:heteroboot}(i) and Lemma \ref{lm:hetero4ver} verify Assumptions \ref{ass:locineq}, \ref{ass:loceq}, and \ref{ass:ineqlindep} are satisfied with $K_g=0$ and $K_f = 0$. 
Also note Assumption \ref{ass:heteroboot}(iii) and Lemma \ref{lm:heterobootcoup} verify Assumption \ref{ass:bootcoupling} and Assumptions \ref{ass:extra}(i)(iii) are immediate given the definitions of $\|\cdot\|_{\mathbf E}$ and $\|\cdot\|_{\mathbf B}$ and ${\mathcal V}_n(P) = \Theta_n \cap R$ by Lemma \ref{lm:hetero3ver}.
Also note $\{\theta \in \mathbf B_n :\overrightarrow{d}_H(\theta, \IDsetRsieve,\|\cdot\|_{\mathbf I}) \leq 1/2\} \subseteq \Theta_n$ for $n$ sufficiently large by Assumption \ref{ass:heteroboot}(iv) and the definitions of $\Theta_n$ and $\|\cdot\|_{\mathbf I}$.
Moreover, Assumptions \ref{ass:heterosieve}(ii)(iii) imply
\begin{equation}\label{heteroboot:0p5}
\sup_{h\in \mathbf B_n} \frac{\|h\|_{\mathbf I}}{\|h\|_{\mathbf E}} = \sup_{h\in \mathbf B_n}  \frac{\max_{1\leq \jmath \leq \mathcal J} \mathcal J \|F(c_\jmath|\cdot)\|_\infty}{\sum_{\jmath = 1}^{\mathcal J}\sup_{P\in \mathbf P} \|F(c_\jmath|\cdot)\|_{P,2}}\lesssim \sqrt{j_n},
\end{equation}
and hence Corollary \ref{cor:bootsetrates}(i), $\nu_n \asymp 1$, and $\mathcal R_n \asymp k_n \sqrt{j_n\log(1+k_n)\log(1+j_n)/n}$ yield 
\begin{equation*}
\overrightarrow{d}_H(\hat \Theta_n^{\rm r}, \IDsetRsieve,\|\cdot\|_{\mathbf I}) \lesssim \sqrt {j_n}\overrightarrow{d}_H(\hat \Theta_n^{\rm r}, \IDsetRsieve,\|\cdot\|_{\mathbf E}) = O_P(\frac{k_n j_n \log(n)}{\sqrt n} + \sqrt{j_n}\tau_n)
\end{equation*}
uniformly in $P\in \mathbf P_0$.
In particular, Assumptions \ref{ass:heteroboot}(iii)(v) imply $\overrightarrow{d}_H(\hat \Theta_n^{\rm r}, \IDsetRsieve,\|\cdot\|_{\mathbf I}) = o_P(1)$ uniformly in $P\in \mathbf P_0$, and therefore since, as argued, we have $\{\theta \in \mathbf B_n :\overrightarrow{d}_H(\theta, \IDsetRsieve,\|\cdot\|_{\mathbf I}) \leq 1/2\} \subseteq \Theta_n$ for $n$ sufficiently large, we obtain
\begin{equation}\label{th:heteroboot1}
\liminf_{n\rightarrow \infty} \inf_{P\in \mathbf P} P(\{\theta \in \mathbf B_n : \overrightarrow{d}_H(\{\theta\},\hat \Theta_n^{\text{\rm r}},\|\cdot\|_{\mathbf I}) \leq 1/4\} \subseteq \Theta_n) = 1.
\end{equation}
Next, observe Lemma \ref{lm:heterolipschitz}, Assumption \ref{ass:heterosieve}(ii) and the definitions of $\|\cdot\|_{\mathbf E}$, $\|\cdot\|_{\mathbf L}$, and $\|\cdot\|_{\mathbf B}$ imply Assumption \ref{ass:genneigh} holds with $\mathcal D_n(\mathbf B,\mathbf E) \asymp \zeta_n\sqrt{j_n}$ and $\mathcal D_n(\mathbf L,\mathbf E) = 1$.
Since $K_m = K_g = K_f = 0$ and $\Upsilon_F$ and $\Upsilon_G$ are affine, the only requirements imposed by Assumption \ref{ass:genbootrates} are that $k_n^{1/p}\sqrt{\log(1+k_n)}B_n \sup_{P\in \mathbf P}J_{[\hspace{0.03 in}]}(\ell_n^{\kappa_\rho}\vee(\nu_n\tau_n)^{\kappa_\rho},\mathcal F_n,\|\cdot\|_{P,2}) =  o(a_n)$ and $(\mathcal R_n +\nu_n\tau_n)  \mathcal D_n(\mathbf B, \mathbf E) = o(r_n)$, which are implied by Assumption \ref{ass:heteroboot}(v), Lemma \ref{lm:heterocov}, and $k_n\sqrt{j_n}\log^2(n)\ell_n = o(1)$ by hypothesis.
Hence, all the conditions of Theorem \ref{th:gencoupsmooth}(ii) hold, which implies there is a $\tilde \ell_n \asymp \ell_n$ such that uniformly in $P\in \mathbf P_0$
\begin{equation}\label{th:heteroboot2}
\hat E_n(R|\ell_n) \geq \UpS (R|\tilde \ell_n) + o_P(a_n).
\end{equation}

Finally, to apply Lemma \ref{lm:forincJboot} to $\hat U_n(\Theta|+\infty)$, note that we can set the norm $\|\cdot\|_{\mathbf B}$ to equal $\|\theta\|_{\mathbf B} = \max_{1\leq \jmath \leq \mathcal J} \|F(c_\jmath|\cdot)\|_\infty$ and interpret $\Upsilon_G$ and $\Upsilon_F$ as satisfying $\Upsilon_G(\theta) = \Upsilon_F(\theta) = 0$ for all $\theta \in \mathbf B$ (since $R=\Theta$).
Hence, Assumptions \ref{ass:locineq}, \ref{ass:loceq}, and \ref{ass:ineqlindep}, \ref{ass:extra}(i) are immediate, while Assumption \ref{ass:bootcoupling} is satisfied by Assumption \ref{ass:heteroboot}(iii) and Lemma \ref{lm:heterobootcoup}.
Further note since $\IDset$ is an equivalence class under $\|\cdot\|_{\mathbf E}$ and $\|\cdot\|_{\mathbf B}$, when studying the unconstrained statistic we can treat the model as identified.
As a result, we may set $\tau_n^{\text{\rm u}} = 0$ and Assumption \ref{ass:extra}(ii) holds by the same arguments employed in \eqref{th:heteroboot1}, while Assumption \ref{ass:extra}(iii) is immediate since $\mathcal V_n(P) = \Theta_n \cap R$.
In order to apply Lemma \ref{lm:forincJboot}(ii), it therefore only remains to verify that $k_n^{1/p}\sqrt{\log(1+k_n)}B_n \sup_{P\in \mathbf P}J_{[\hspace{0.03 in}]}(\ell_n^{\rm u},\mathcal F_n^{\rm u},\|\cdot\|_{P,2}) =  o(a_n),$ $\mathcal R_n^{\rm u} = o(\ell_n^{\rm u})$, and $\mathcal S_n^{\rm u}(\mathbf B,\mathbf E) \times \mathcal R_n^{\rm u} = o(1)$, which are implied by $k_n\sqrt{j_n}\log^2(n) \ell_n^{\rm u} = o(1)$, $k_n\sqrt{j_n}\log(n)/\sqrt n = o(\ell_n^{\rm u})$, and Assumption \ref{ass:heteroboot}(iii) respectively.
Thus, \eqref{th:heteroboot2} and Lemma \ref{lm:forincJboot}(ii) verify \eqref{th:heteroboot0} with $\tilde \ell_n^{\rm u} = \ell_n^{\rm u}$ and $\tilde \ell_n \asymp \ell_n$, which in turn establishes  the theorem. \qed

\begin{lemma}\label{lm:hetero3ver}
If Assumptions \ref{ass:heterosieve}(iii), \ref{ass:heteromoments}(iii), and \ref{ass:heterosigma}(ii) hold, then Assumption \ref{ass:keycons} holds with $R = \Theta$ and $R$ as in \eqref{ex:hetero7paux}, $\mathcal V_n(P) = \Theta_n \cap R$, $\|\theta\|_{\mathbf E} = \sum_{\jmath = 1}^{\mathcal J} \sup_{P\in \mathbf P}  \|F(c_\jmath|\cdot)\|_{P,2}$ for any $\theta = (\{F(c_\jmath|\cdot)\}_{\jmath = 1}^{\mathcal J},\mu) \in \mathbf B$, and $\nu_n^{-1} \asymp 1$.
\end{lemma}

\noindent {\sc Proof:} First note that since we are setting $\mathcal V_n(P) = \Theta_n \cap R$, Assumption \ref{ass:keycons}(ii) is immediate.
To verify Assumption \ref{ass:keycons}(i), let $\|\theta\|_{\mathbf E} = \sum_{\jmath = 1}^{\mathcal J}\sup_{P\in \mathbf P} \|F(c_\jmath|\cdot)\|_{P,2}$ for any $\theta = (\{F(c_\jmath|\cdot)\}_{\jmath = 1}^{\mathcal J},\mu) \in \mathbf B$.
Then note that any $(\{F(c_\jmath|\cdot)\}_{\jmath = 1}^{\mathcal J},\mu) = \theta \in \Theta_n$ must be such that $F(c_\jmath|\cdot) = p^{j_n\prime} \beta_{\jmath,\theta}$ for some $\beta_{\jmath,\theta} \in \mathbf R^{j_n}$ and, similarly, $\Pi_n\IDpoint = (\{F_n(c_\jmath|\cdot)\}_{\jmath =1}^{\mathcal J},\mu_n)$ must satisfy $F_n(c_\jmath|\cdot) = p^{j_n\prime} \beta_{\jmath,n}$.
The Cauchy Schwarz inequality, and Assumptions \ref{ass:heterosieve}(iii) and \ref{ass:heteromoments}(iii) then yield that uniformly in $P\in \mathbf P_0$ we must have
\begin{multline}\label{lm:hetero3ver2}
\|\theta - \Pi_n \IDpoint\|_{\mathbf E}  \lesssim  \sum_{\jmath = 1}^{\mathcal J} \|\beta_{\jmath,\theta} - \beta_{\jmath,n}\|_2 \lesssim \sum_{\jmath =1}^{\mathcal J} \|E_P[q^{k_n}(W) p^{j_n}(W)^\prime(\beta_{\jmath,\theta} - \beta_{\jmath,n})]\|_{2} \\
\lesssim \{\sum_{\jmath =1}^{\mathcal J} \|E_P[(F(c_\jmath|W) - F_n(c_\jmath|W))q^{k_n}(W)]\|_{\Sigma_{\jmath,P},2}^2\}^{1/2},
\end{multline}
where the final inequality holds due to $\|\Sigma_{\jmath,P}^{-1}\|_{o,2}$ being uniformly bounded by Assumption \ref{ass:heterosigma}(ii) and $\sum_{\jmath =1}^{\mathcal J} |a^{(\jmath)}| \leq \sqrt{\mathcal J}\|a\|_2$ for any $(a^{(1)},\ldots, a^{(\mathcal J)}) =a\in \mathbf R^{\mathcal J}$.
Result \eqref{lm:hetero3ver2} and the definition of $\rho_\jmath(X,\theta)$ in \eqref{ex:hetero7} verify Assumption \ref{ass:keycons}(i) holds with $\nu_n^{-1} \asymp 1$. \qed

\begin{lemma}\label{lm:heterocov}
Define the class $\mathcal F_n \equiv \{f: f(v) = (1\{y \leq c_\jmath\} - p^{j_n}(w)^\prime \beta) \text{ for some } 1\leq \jmath \leq \mathcal J \text{ and } \|p^{j_n\prime}\beta\|_\infty \leq 2\}$ and suppose that
Assumptions \ref{ass:heterosieve}(ii)(iii) hold.
Then, it follows that $\sup_{P\in \mathbf P} N_{[\hspace{0.02 in}]}(\epsilon,\mathcal F_n,\|\cdot\|_{P,2})\lesssim (1\vee (\sqrt{j_n}K/\epsilon)^{j_n})$  for some $K<\infty$, and in addition $\sup_{P\in \mathbf P}J_{[\hspace{0.02 in}]}(\epsilon,\mathcal F_n,\|\cdot\|_{P,2}) \lesssim \epsilon \sqrt{j_n} (1+\sqrt {\log(1 \vee (\sqrt{j_n}/\epsilon))})$.
\end{lemma}

\noindent {\sc Proof:} First note that for any $p^{j_n\prime}\beta_1$ and $p^{j_n\prime}\beta_2$, the Cauchy-Schwarz inequality yields
\begin{equation*}\label{lm:heterocov1}
|p^{j_n}(w)^\prime \beta_1 - p^{j_n}(w)^\prime \beta_2| \leq \sup_{w}\|p^{j_n}(w)\|_2 \|\beta_1-\beta_2\|_2 \lesssim \sqrt{j_n}  \|\beta_1-\beta_2\|_2,
\end{equation*}
where in the final inequality we employed Assumption \ref{ass:heterosieve}(ii).
Hence, Theorem 2.7.11 in \cite{vandervaart:wellner:1996}, $\|\beta\|_2 \asymp \sup_{P\in \mathbf P}\|p^{j_n\prime}\beta\|_{P,2}$ by Assumption \ref{ass:heterosieve}(iii), and $\sup_{P\in \mathbf P}\|p^{j_n\prime}\beta\|_{P,2} \leq \|p^{j_n\prime}\beta\|_{\infty}\leq 2$ for any $p^{j_n\prime}\beta \in \Theta_n$ imply
\begin{equation}\label{lm:heterocov2}
\sup_{P\in \mathbf P} N_{[\hspace{0.02 in}]}(\epsilon,\mathcal F_{n},\|\cdot\|_{P,2}) \lesssim  1\vee (\frac{K\sqrt{j_n}}{\epsilon })^{j_n},
\end{equation}
for some $K < \infty$, which establishes the first claim of the lemma.
For the second claim of the lemma, we employ \eqref{lm:heterocov2} and the change of variables $v = u/\epsilon$ to obtain
\begin{multline*}
\sup_{P\in \mathbf P} J_{[\hspace{0.02 in}]}(\epsilon,\mathcal F_n,\|\cdot\|_{P,2}) \lesssim \epsilon + \int_0^\epsilon (\log(1\vee (\frac{K\sqrt{j_n}}{u})^{j_n}))^{1/2}du  \\
= \epsilon(1 + \sqrt{j_n} \int_0^1(\log(1 \vee (\frac{K\sqrt{j_n}}{v\epsilon})))^{1/2}dv) \lesssim \sqrt{j_n}\epsilon(1+ \sqrt{\log(1\vee (\sqrt{j_n}/\epsilon))}),
\end{multline*}
where the final inequality follows from $(1\vee ab) \leq (1\vee a)(1\vee b)$ for any $a,b\in \mathbf R_+$. \qed

\begin{lemma}\label{lm:heteroobviousver}
Let $\rho_\jmath : \mathbf R\times \mathbf W\times \Theta$ be as defined in \eqref{ex:hetero7}. 
It then follows Assumptions \ref{ass:coupling}(ii) and \ref{ass:driftlin} hold with $\kappa_\rho = 1$, $K_\rho = 1$, $K_m = 0$, $M <\infty$, $\|\cdot\|_{\mathbf L} = \|\cdot\|_{\mathbf E}$, and $\|\theta\|_{\mathbf E} = \sum_{\jmath=1}^{\mathcal J} \sup_{P\in \mathbf P} \|F(c_\jmath|\cdot)\|_{P,2}$ for any $\theta = (\{F(c_\jmath|\cdot)\}_{\jmath = 1}^{\mathcal J},\mu)\in \mathbf B$.
\end{lemma}

\noindent {\sc Proof:} First note that for any $(\{F_1(c_\jmath|\cdot)\}_{\jmath = 1}^{\mathcal J},\mu_1) = \theta_1 \in \mathbf B$ and $(\{F_2(c_\jmath|\cdot)\}_{\jmath = 1}^{\mathcal J},\mu_2) = \theta_2 \in \mathbf B$, we obtain from \eqref{ex:hetero7} and the definition of $\|\cdot\|_{\mathbf E}$ that for all $P\in \mathbf P$
\begin{equation*}\label{eq:heteroobviousver1}
E_P[\|\rho(X,\theta_1) - \rho(X,\theta_2)\|_2^2 ] = \sum_{\jmath = 1}^{\mathcal J} E_P[(F_1(c_\jmath|W) - F_2(c_\jmath|W))^2] \leq \|\theta_1 - \theta_2\|_{\mathbf E}^2 ,
\end{equation*}
which verifies Assumption \ref{ass:coupling}(ii) holds with $\kappa_\rho = 1$ and $K_\rho =1$.
Next, for any $P\in \mathbf P$ define $\nabla m_{P,\jmath}(\theta)[h] = - F_h(c_\jmath|W)$ for all $\theta \in \mathbf B$ and $(\{F_h(c_\jmath|\cdot)\}_{\jmath = 1}^{\mathcal J},\mu_h) = h \in \mathbf B$.
Since $m_{P,\jmath}(\theta) = P(Y \leq c_\jmath|W) - F(c_\jmath|W)$ for any $\theta = (\{F(c_\jmath|\cdot)\}_{\jmath = 1}^{\mathcal J},\mu)\in \mathbf B$, direct calculation verifies Assumption \ref{ass:driftlin} holds with $K_m = 0$, $M = 1$, and $\|\cdot\|_{\mathbf L} = \|\cdot\|_{\mathbf E}$. \qed

\begin{lemma}\label{lm:heterocoup}
If $k_n^3 j_n^2 \log^2(n) = o(n)$, Assumptions \ref{ass:heterosieve}(i)-(iii) and \ref{ass:heteromoments}(i) hold, then Assumption \ref{ass:coupling}(i) holds with $R = \Theta$ for any $a_n$ with $k_n^3 j_n^2 \log^2(n)/n = o(a_n^2)$.
\end{lemma}

\noindent {\sc Proof:} We establish the result by applying Lemma \ref{lm:finitecoup}.
To this end, we let $\tilde j_n = \mathcal J + j_n$ set $\{r_{j}\}_{j=1}^{\tilde j_n} = \{1\{y \leq c_\jmath\}\}_{\jmath = 1}^{\mathcal J} \cup \{p_{j}\}_{j=1}^{j_n}$ and let $r^{\tilde j_n}(x) \equiv (r_{1}(x),\ldots, r_{\tilde j_n}(x))^\prime$.
Next note that any $f\in \mathcal F_n$ may be written as $r^{\tilde j_n\prime}\beta$ for some $\beta \in \mathbf R^{\tilde j_n}$.
Moreover, since $\sup_{P\in \mathbf P} \max_{1\leq \jmath \leq \mathcal J} \|F(c_\jmath|\cdot)\|_{P,2} \leq \max_{1\leq \jmath \leq \mathcal J} \|F(c_\jmath|\cdot)\|_\infty \leq 2$ for any $(\{F(c_\jmath|\cdot)\}_{\jmath = 1}^{\mathcal J},\mu) = \theta \in \Theta_n$, Assumption \ref{ass:heterosieve}(iii) implies that there exists a $C_0 < \infty$ (independent of $\tilde j_n$) such that $\|\beta\|_2 \leq C_0$ whenever $r^{\tilde j_n\prime} \beta \in \mathcal F_n$.
Hence, by Assumptions \ref{ass:heterosieve}(ii) and \ref{ass:heteromoments}(i), we may apply Lemma \ref{lm:finitecoup} with $b_{1n} \asymp \sqrt{j_n}$, $b_{2n} \asymp k_n$, and $C_n = O(1)$, from which the claim of the present lemma immediately follows. \qed

\begin{lemma}\label{lm:hetero4ver}
Let $\mathbf B = (\bigotimes_{\jmath =1}^{\mathcal J} C_B(\mathbf W))\times \mathcal M$ and $\Theta$, $\Upsilon_G$, and $\Upsilon_F$ be as defined in \eqref{ex:hetero6}, \eqref{ex:hetero9p0}, and \eqref{ex:hetero9p3}. If $\Psi(g,\cdot)$ is bounded on $\Omega$, then Assumptions \ref{ass:locineq}, \ref{ass:loceq}, and \ref{ass:ineqlindep} are satisfied with $K_g = 0$, $\nabla \Upsilon_G(\theta)[h] = \Upsilon_G(h)$, $K_f = 0$, and $\nabla \Upsilon_F(\theta)[h]$ equal to
\begin{equation}\label{lm:hetero4verdisp}
\nabla \Upsilon_F(\theta)[h] = (\Upsilon_F^{(\text{\rm e})}(h), \Upsilon_F^{(\mu)}(h) + 1, \Upsilon_F^{(\text{\rm s})}(h) + \lambda) .
\end{equation}
\end{lemma}

\noindent {\sc Proof:} For any measure $\mu \in \mathcal M$ let $\mu = \mu^+ - \mu^{-}$ denote its Jordan decomposition, $|\mu| = \mu^+ + \mu^{-}$, and recall the total variation of $\mu$ equals $\|\mu\|_{TV} = |\mu|(\Omega)$.
Since $\Upsilon_G : \mathbf B \to \ell^\infty(\mathcal B)$ is linear, in order to verify Assumption \ref{ass:locineq} we need only show that $\Upsilon_G$ is continuous.
To this end, recall that for any $(\{F(c_\jmath|\cdot)\}_{\jmath = 1}^{\mathcal J}, \mu) = \theta \in \mathbf B$ we had defined $\|\theta\|_{\mathbf B} = \sum_{\jmath = 1}^{\mathcal J}\|F(c_\jmath|\cdot)\|_\infty + \|\mu\|_{TV}$.
Hence, employing the definition of $\Upsilon_G$ we obtain
\begin{equation*}\label{lm:hetero4ver1}
\|\Upsilon_G\|_o = \sup_{\|\theta\|_{\mathbf B} = 1} \|\Upsilon_G(\theta)\|_\infty = \sup_{\mu : \|\mu\|_{TV} = 1} \sup_{B\in \mathcal B} |\mu(B)| \leq \sup_{\mu : \|\mu\|_{TV} = 1} |\mu|(\Omega) = 1,
\end{equation*}
which, by linearity of $\Upsilon_G$, implies Assumption \ref{ass:locineq} holds with $\nabla \Upsilon_G = \Upsilon_G$ and $K_g = 0$.
By similar arguments, note that $\Upsilon_F^{(\text{e})}:\mathbf B \to \mathbf R^{\mathcal J \mathcal L}$, as defined in \eqref{ex:hetero9p3}, is linear and 
\begin{multline}\label{lm:hetero4ver2}
\|\Upsilon_F^{\text{(\rm e)}}\|_o^2 = \sup_{\|\theta\|_{\mathbf B} = 1} \sum_{\jmath = 1}^{\mathcal J} \sum_{l = 1}^{\mathcal L} (F(c_\jmath|w_l) - \int 1\{g(w_l,\eta) \leq c_\jmath\}\mu(d\eta))^2 \\
\leq \sum_{\jmath = 1}^{\mathcal J} \sum_{l = 1}^{\mathcal L} \{2\sup_{\|F(c_\jmath|\cdot)\|_\infty = 1} (F(c_\jmath|w_l))^2 + 2 \sup_{\|\mu\|_{TV} = 1} (|\mu|(\Omega))^2\} = 4\mathcal J \mathcal L.
\end{multline}
Moreover, note that for any bounded $f : \Omega \to \mathbf R$ and $\mu_1,\mu_2\in \mathcal M$ it follows that
\begin{equation*}\label{lm:hetero4ver3}
\int_\Omega f(\eta)(\mu_1(d\eta) - \mu_2(d\eta)) \leq \|f\|_\infty |\mu_1 - \mu_2|(\Omega) = \|f\|_\infty \|\mu_1 - \mu_2\|_{TV},
\end{equation*}
which implies $\Upsilon_F^{(\mu)}$ and $\Upsilon_F^{\text{(\rm s)}}$ are Fr\'echet differentiable with $\nabla \Upsilon_F^{(\mu)} = \Upsilon_F^{(\mu)}+1$, $\nabla \Upsilon_F^{(\text{s})} = \Upsilon_F^{\text{(s)}} + \lambda$, $\|\nabla \Upsilon_F^{(\mu)}\|_o \leq 1$, and $\|\nabla \Upsilon_F^{\text{(s)}}\|_o \leq \|\Psi(g,\cdot)\|_\infty$.
By \eqref{lm:hetero4ver2} we may therefore conclude Assumptions \ref{ass:loceq}(i)(ii)(iii) are satisfied with $\nabla \Upsilon_F$ as in \eqref{lm:hetero4verdisp} and $K_f = 0$.
Furthermore, note that (provided $\Theta_n \cap R\neq \emptyset$) there is a $\theta^\star\in \mathbf B_n$ such that $\Upsilon_F(\theta^\star) = 0$, which together with \eqref{lm:hetero4verdisp} implies the range of $\nabla \Upsilon_F$ equals $\mathbf F_n$ and hence Assumption \ref{ass:loceq}(iv) holds.
Finally, we note Assumption \ref{ass:ineqlindep} is immediate due to $\Upsilon_F$ being affine. \qed

\begin{lemma}\label{lm:finitecoup}
Let $\{r_{j}\}_{j=1}^{j_n}$ be functions of $X$, $r^{j_n}(x) = (r_{1}(x),\ldots, r_{j_n}(x))^\prime$, define the class $\mathcal G_n = \{ r^{j_n\prime} \beta \text{ for some } \beta \text{ with } \|\beta\|_2 \leq  C_n\}$, and suppose $b_{1n} \equiv \sup_x \|r^{j_n}(x)\|_2$ and $b_{2n} \equiv \sup_z \|q^{k_n}(z)\|_2$ are finite. If $\{X_i,Z_i\}_{i=1}^n$ is i.i.d.\ with $(X,Z)\sim P\in \mathbf P$, then there is an isonormal Gaussian process $\Iso$ such that uniformly in $P\in \mathbf P$
\begin{multline}
\sup_{g \in \mathcal G_n} \|\frac{1}{\sqrt n}\sum_{i=1}^n(g(X_i)q^{k_n}(Z_i) - E_P[g(X)q^{k_n}(Z)]) - \Iso(gq^{k_n})\|_2 \\ = O_P(\frac{C_n\sqrt{k_nj_n}b_{1n} b_{2n}\log(n)}{\sqrt n}).
\end{multline}
\end{lemma}

\noindent {\sc Proof:} For notational simplicity, we first define a $k_n \times j_n$ matrix $\mathbb E^{(1)}_{n}$ to be given by
\begin{equation*}
\mathbb E^{(1)}_{n} \equiv \frac{1}{\sqrt n}\sum_{i=1}^n \{q^{k_n}(Z_i)r^{j_n}(X_i)^\prime - E_P[q^{k_n}(Z) r^{j_n}(X)^\prime]\}.
\end{equation*}
For any matrix $A$ let $\text{vec}\{A\}$ denote a column vector consisting of the unique elements of $A$ and set $\mathbb E_{n} \equiv \text{vec}\{\mathbb E^{(1)}_{n}\}$, noting that $\mathbb E_{n}$ has dimension (at most) $j_nk_n$.
Our first step is to couple $\mathbb E_{n}$ to a normal vector $\mathbb N_{P}$.
To this end, we note that
\begin{multline*}\label{lm:finitecoup2}
\sup_{z,x}\|\text{vec}\{q^{k_n}(z)r^{j_n}(x)^\prime - E_P[q^{k_n}(Z)r^{j_n}(X)^\prime]\}\|_2^2 \\ \leq \sup_{z,x} 4\text{trace}\{q^{k_n}(z)r^{j_n}(x)^\prime r^{j_n}(x) q^{k_n}(z)^\prime\} \leq 4 b_{1n}^2 b_{2n}^2
\end{multline*}
by definition of $b_{1n}$ and $b_{2n}$.
Since the dimension of $\mathbb E_{n}$ is at most $j_nk_n$, Theorem 1.1 in \cite{zhai2018high} and Markov's inequality imply, provided the underlying probability space is suitably rich, that there is a Gaussian vector $\mathbb N_{P}$ such that
\begin{equation}\label{lm:finitecoup3}
\|\mathbb E_{n} - \mathbb N_{P}\|_2 = O_P(\frac{\sqrt{k_nj_n} b_{1n} b_{2n}\log(n)}{\sqrt n})
\end{equation}
uniformly in $P\in \mathbf P$.
Next observe that for any $g \in \mathcal G_n$ there exists a $\beta \in \mathbf R^{j_n}$ such that
\begin{equation*}
\frac{1}{\sqrt n}\sum_{i=1}^n (g(X_i)q^{k_n}(Z_i) - E_P[g(X)q^{k_n}(Z)]) = \mathbb E_{n}^{(1)} \beta.
\end{equation*}
Hence, letting $\mathbb N_{P}^{(1)}$ denote the $k_n\times j_n$ matrix built from the corresponding entries of the normal vector $\mathbb N_{P}$, we define the Gaussian process $\Iso$ by setting
\begin{equation*}
\Iso(g q^{k_n}) = \mathbb N_{P}^{(1)} \beta
\end{equation*}
for any $r^{j_n\prime}\beta= g \in \mathcal G_n$.
Therefore, since $\|\beta\|_2 \leq C_n$ by definition of $\mathcal G_n$, and the operator norm is bounded by the Frobenius norm, we obtain from result \eqref{lm:finitecoup3} that
\begin{multline*}
\sup_{g\in \mathcal G_n} \|\frac{1}{\sqrt n}\sum_{i=1}^n(g(X_i)q^{k_n}(Z_i) - E_P[g(X)q^{k_n}(Z)]) - \Iso(gq^{k_n})\|_2\\
\leq \|\mathbb E_{n}^{(1)} - \mathbb N_{P}^{(1)}\|_{o,2} C_n = O_P(\frac{C_n\sqrt{k_nj_n}b_{1n} b_{2n}\log(n)}{\sqrt n})
\end{multline*}
uniformly in $P\in \mathbf P$, and hence the claim of the lemma follows. \qed

\begin{lemma}\label{lm:finiteder}
Let $\{r_{j}\}_{j=1}^{j_n}$ be a set of functions of $X$, $r^{j_n}(x) \equiv (r_{1}(x),\ldots, r_{j_n}(x))^\prime$, and suppose $\sup_x \|r^{j_n}(x)\|_2 \lesssim b_{1n}$, $\sup_{z}\|q^{k_n}(z)\|_2 \lesssim
b_{2n}$, and $E_P[q^{k_n}(Z)q^{k_n}(Z)^\prime]$ and $E_P[r^{j_n}(X) r^{j_n}(X)^\prime]$ have eigenvalues bounded uniformly in $P\in \mathbf P$, $j_n$, $k_n$.
If $\{X_i,Z_i\}_{i=1}^n$ is i.i.d. with $(X,Z) \sim P\in \mathbf P$, then there is a $K < \infty$ such that for all $\delta \geq 0$
\begin{multline*}
\sup_{P\in \mathbf P} P(\|\frac{1}{n}\sum_{i=1}^n q^{k_n}(Z_i)r^{j_n}(X_i)^\prime - E_P[q^{k_n}(Z)r^{j_n}(X)^\prime]\|_{o,2} > \delta) \\
\leq (j_n+k_n)\exp\{- \frac{n\delta^2 K}{b_{1n}^2\vee b^2_{2n}+ \delta b_{1n}b_{2n}}\}.
\end{multline*}
\end{lemma}

\noindent {\sc Proof:} We first define a $k_n\times j_n$ random matrix $\mathbb M_{i,n}$ satisfying $E_P[\mathbb M_{i,n}] = 0$ by
\begin{equation*}\label{lm:finiteder1}
\mathbb M_{i,n} \equiv \frac{1}{n}\{q^{k_n}(Z_i)r^{j_n}(X_i)^\prime - E_P[q^{k_n}(Z)r^{j_n}(X)^\prime]\}.
\end{equation*}
Since for any random matrix $A$ we have $\|E[A]\|_o \leq E[\|A\|_o]$ by Jensen's inequality, $\|A\|_o^2 \leq \text{trace}\{A^\prime A\}$, $\sup_x \|r^{j_n}(x)\|_2 \lesssim b_{1n}$, and $\sup_{z}\|q^{k_n}(z)\|_2 \lesssim b_{2n}$ imply
\begin{multline}\label{lm:finiteder2}
\|\mathbb M_{i,n}\|_o^2 \lesssim \|\frac{1}{n} q^{k_n}(Z_i)r^{j_n}(X_i)^\prime\|_o^2 + E_P[\|\frac{1}{n} q^{k_n}(Z)r^{j_n}(X)^\prime\|_o^2]\\
\lesssim \frac{\sup_{z}\|q^{k_n}(z)\|_2^2  \times \sup_{x} \|r^{j_n}(x)\|_2^2}{n^2} \lesssim \frac{b_{1n}^2 b_{2n}^2}{n^2}.
\end{multline}
Moreover, since the eigenvalues of $E_P[q^{k_n}(Z)q^{k_n}(Z)^\prime]$ are bounded uniformly in $P\in \mathbf P$ by assumption and $\sup_x \|r^{j_n}(x)\|_2 \lesssim b_{1n}$ it additionally follows that
\begin{equation}\label{lm:finiteder3}
\sup_{P\in \mathbf P} \|\sum_{i=1}^n E_P[\mathbb M_{i,n}\mathbb M_{i,n}^\prime]\|_o  \leq \sup_{P\in \mathbf P} \frac{2}{n} \|E_P[q^{k_n}(Z) q^{k_n}(Z)^\prime \|r^{j_n}(X)\|_2^2]\|_o \lesssim \frac{b_{1n}^2}{n}.
\end{equation}
Identical arguments but relying on the eigenvalues of $E_P[r^{j_n}(X)r^{j_n}(X)^\prime]$ being bounded uniformly in $P\in \mathbf P$ and $\sup_x \|q^{k_n}(x)\|_2 \lesssim b_{2n}$ by hypothesis further yield that
\begin{equation}\label{lm:finiteder4}
\sup_{P\in \mathbf P} \|\sum_{i=1}^n E_P[\mathbb M_{i,n}^\prime \mathbb M_{i,n}]\|_o  \lesssim  \frac{b_{2n}^2}{n}.
\end{equation}
The claim of the lemma then follows from results \eqref{lm:finiteder2}, \eqref{lm:finiteder3}, and \eqref{lm:finiteder4} allowing us to apply Theorem 1.6 in \cite{tropp:2012} with $\sigma^2 \asymp (b^2_{1n}\vee b_{2n}^2)/n$ and $R \asymp b_{1n}b_{2n}/n$.  \qed

\begin{lemma}\label{lm:heterobootcoup}
If Assumptions \ref{ass:heterosieve}(i)-(iii), \ref{ass:heteromoments}(i)(ii) hold, and $j_n^3 k_n^2 \log(1+j_nk_n) = o(n)$, then it follows that Assumption \ref{ass:bootcoupling} holds with $R = \Theta$ for any sequence $a_n$ satisfying $k_n^{1/p}(k_n^2j_n^5\log^3(1+k_nj_n)/n)^{1/4} = o(a_n)$.
\end{lemma}

\noindent {\sc Proof:} Let $\mathcal G_n \equiv \{g : g(x) = 1\{y \leq c_\jmath\} - p^{j_n}(w)^\prime \beta \text{ for some } 1\leq \jmath \leq \mathcal J \text{ and } \|p^{j_n\prime} \beta\|_\infty \leq 2\}$ and $\tilde {\mathcal F}_n \equiv \{gq_{k} :  g\in \mathcal G_n \text{ and } 1\leq k \leq k_n\}$.
Further let $\IsoW$ be a Gaussian process on $\tilde {\mathcal F}_n$ independent of $\{V_i\}_{i=1}^n$, satisfying $E[\IsoW(f_1)] =0$ and $E[\IsoW(f_1)\IsoW(f_2)] = \text{Cov}_P\{f_1,f_2\}$ for any $f_1,f_2\in \tilde{\mathcal F}_n$, and for any $f\in \tilde {\mathcal F}_n$ define $\Wboot(f)$ to be given by
$$\Wboot(f) \equiv \frac{1}{\sqrt n}\sum_{i=1}^n \omega_i\{f(V_i) - \frac{1}{n}\sum_{j=1}^n f(V_j)\}$$
where $\{\omega_i\}_{i=1}^n$ are the same weights used in building $\Bemp$.
Then note that when $R = \Theta$ and for $\WPT(\theta)\equiv (\IsoW(\rho_1(\cdot,\theta)q^{k_n})^\prime,\ldots \IsoW(\rho_{\mathcal J}(\cdot,\theta)q^{k_n})^\prime)^\prime$, we obtain
\begin{equation}\label{lm:heterobootcoup1}
\sup_{\theta \in \Theta_n} \|\Bemp(\theta) - \WPT(\theta)\|_p \lesssim k_n^{1/p} \sup_{f\in \tilde {\mathcal F}_n} |\Wboot(f)- \IsoW(f) |.
\end{equation}
We will therefore establish the lemma by employing \eqref{lm:heterobootcoup1} and applying Theorem \ref{th:mainbootcoup}(i) to the class $\tilde {\mathcal F}_n$.
To this end, define $f^{d_n}(V)$ to be given by
\begin{equation}\label{lm:heterobootcoup1p5}
f^{d_n}(V) \equiv g^{d_n}(V) - E_P[g^{d_n}(V)] \hspace{0.4 in} g^{d_n}(V) \equiv q^{k_n}(Z)\otimes\left(\begin{array}{c} p^{j_n}(W) \\ 1\{Y \leq c_1\} \\ \vdots \\ 1\{Y \leq c_{\mathcal J}\}\end{array}\right)
\end{equation}
and note $d_n = k_n(j_n + \mathcal J)$.
Next observe that applying Lemma \ref{lm:kronereig} with $D_1 \equiv (p^{j_n}(W)^\prime,1\{Y\leq c_1\},\ldots, 1\{Y \leq c_{\mathcal J}\})^\prime$ and $D_2 = q^{k_n}(Z)$ allows us to conclude
\begin{equation}\label{lm:heterobootcoup2}
\sup_{P\in \mathbf P} \overline{\text{eig}}\{E_P[g^{d_n}(V) g^{d_n}(V)^\prime] \}\leq \sup_{P\in \mathbf P} (\|\overline{\text{eig}}\{D_1D_1^\prime\}\|_{P,\infty} \times \overline{\text{eig}}\{E_P[D_2D_2^\prime]\} )\lesssim j_n,
\end{equation}
where the final inequality holds by Assumptions \ref{ass:heterosieve}(ii) and \ref{ass:heteromoments}(ii).
Hence, since in addition $\overline{\text{eig}}\{E_P[g^{d_n}(V)]E[g^{d_n}(V)^\prime]\} \leq \overline{\text{eig}}\{E_P[g^{d_n}(V)g^{d_n}(V)^\prime]\}$, results \eqref{lm:heterobootcoup1p5} and \eqref{lm:heterobootcoup2} imply Assumption \ref{ass:4seriescoup}(i) holds with $C_n \asymp j_n$.
Next note Assumption \ref{ass:4seriescoup}(ii) is satisfied with $K_n \asymp \sqrt{k_nj_n}$ by Assumptions \ref{ass:heterosieve}(ii) and \ref{ass:heteromoments}(i).
By Assumptions \ref{ass:heterosieve}(iii) it also follows that $\|\beta\|_2 \asymp \sup_{P\in \mathbf P} \|p^{j_n\prime}\beta\|_{P,2} \leq \|p^{j_n\prime}\beta\|_\infty$.
Hence, by definition of $\tilde {\mathcal F}_n$, there is a $C_0 < \infty$ such that any $f\in \tilde{\mathcal F}$ satisfies $f(V)-E_P[f(V)] = f^{d_n}(V)^\prime\beta$ for some $\beta$ in
\begin{equation*}
\mathcal B_n \equiv \{\beta \in \mathbf R^{d_n} : \beta = e_k \otimes \gamma \text{ for some } \gamma \in \mathbf R^{j_n+\mathcal J} \text{ with } \|\gamma\|_2 \leq C_0\},
\end{equation*}
where $e_k \in \mathbf R^{k_n}$ has its $k^{th}$ coordinate equal to one and all other coordinates equal to zero.
In particular, it follows that Assumption \ref{ass:4seriesreg}(i) is immediate with $G_{n,P}$ equal to the zero function and $J_{1n} = 0$.
Moreover, setting $\mathcal C_n \equiv \{\gamma \in \mathbf R^{j_n+\mathcal J} : \|\gamma\|_2 \leq C_0\}$, we can then conclude from the definition of $\mathcal B_n$ and $N(\epsilon,\mathcal C_n,\|\cdot\|_2) \lesssim 1 \vee (C_0/\epsilon)^{j_n}$ that
\begin{multline*}
\int_0^\infty \sqrt{\log(N(\epsilon,\mathcal B_n,\|\cdot\|_2))}d\epsilon \\
\lesssim \int_0^{C_0} \sqrt{\log(k_n) + \log(N(\epsilon,\mathcal C_n,\|\cdot\|_2))}d\epsilon \lesssim \sqrt{\log(k_n)} + \sqrt{j_n},
\end{multline*}
which verifies Assumption \ref{ass:4seriesreg}(ii) is satisfied with $J_{2n} \asymp \sqrt{\log(k_n)} + \sqrt{j_n}$.
Thus, applying Theorem \ref{th:mainbootcoup}(i) with $K_n \asymp \sqrt{k_n,j_n}$, $C_n \asymp j_n$, $d_n \lesssim k_nj_n$, $J_{1n} = 0$, and $J_{2n} \asymp \sqrt{\log(k_n)} + \sqrt{j_n}$ implies that uniformly in $P\in \mathbf P$ we have
\begin{equation}\label{lm:heterobootcoup5}
\sup_{f\in \tilde {\mathcal F}_n} |\Wboot(f)- \IsoW(f) | = O_P(\{\frac{k_n^2j_n^5\log^3(1+k_nj_n)}{n}\}^{1/4})
\end{equation}
provided that $j_n^3 k_n^2 \log(1+j_nk_n) = o(n)$.
Since the latter condition is satisfied by hypothesis, the claim of the lemma then follows from \eqref{lm:heterobootcoup1} and \eqref{lm:heterobootcoup5}. \qed

\begin{lemma}\label{lm:heteronol}
Define $\|\theta\|_{\mathbf E} = \sum_{\jmath =1}^{\mathcal J}\sup_{P\in \mathbf P} \|F(c_\jmath|\cdot)\|_{P,2}$ and for $\hat V_n(\theta,R|+\infty)$ as in \eqref{ex:hetero21} let $\hat V_n(\theta,R|\ell_n) = \hat V_n(\theta,R|+\infty) \cap \{h  : \|h/\sqrt n\|_{\mathbf E} \leq \ell_n\}$.
If Assumptions \ref{ass:heterosieve}, \ref{ass:heteromoments}, and \ref{ass:heterosigma} hold, then for any $a_n=o(1)$ and $\ell_n=o(1)$ satisfying $k_n^4j_n^5\log^3(1+k_nj_n)/n = o(a_n^4)$ and $k_n\sqrt{j_n}\log(n)/\sqrt n = o(\ell_n)$ it follows uniformly in $P\in \mathbf P_0$ that
\begin{equation*}
\hat U_n(R|+\infty) = \inf_{\theta \in \hat \Theta_n^{\text{\rm r}}}\inf_{h \in \hat V_n(\theta,R|\ell_n)} \{\sum_{\jmath = 1}^{\mathcal J} \|\hat {\mathbb W}_{\jmath,n}(\theta) + \hat {\mathbb D}_{\jmath,n}[h]\|_{\hat \Sigma_{\jmath,n},2}\}^{1/2} + o_P(a_n) .
\end{equation*}
\end{lemma}

\noindent {\sc Proof:} We establish the claim of the lemma by verifying the conditions of Lemma \ref{lm:lnotbind}.
To this end, recall that in the proof of Theorem \ref{th:heteroapprox} we argued that Assumptions \ref{ass:startreg}(i)(iii) and \ref{ass:weights} hold with $B_n \asymp \sqrt{k_n}$ and $J_n \asymp \sqrt{j_n\log(1+j_n)}$.
Moreover, Assumption \ref{ass:heterosieve}(iii) implies that for any $(\{p^{j_n\prime}\beta_{\jmath,h}\}_{\jmath = 1}^{\mathcal J},\mu) = h \in \mathbf B_n$ we have
\begin{equation}\label{lm:heteronol1}
\|h\|_{\mathbf E} \lesssim \sum_{\jmath =1}^{\mathcal J} \|\beta_{\jmath,h}\|_2 \lesssim  \{\sum_{\jmath = 1}^{\mathcal J} \|\mathbb D_{\jmath,P}[h]\|_2^2\}^{1/2} = \|\DerP[h] \|_2,
\end{equation}
where the second inequality follows from $\mathbb D_{\jmath,P}[h] = -E_P[q^{k_n}(Z)p^{j_n}(W)^\prime \beta_{\jmath,h}]$ and the smallest singular values of $E_P[q^{k_n}(Z)p^{j_n}(W)^\prime]$ being bounded away from zero uniformly in $P\in \mathbf P$ by Assumption \ref{ass:heteromoments}(iii).
Since $\nu_n \asymp 1$ by Lemma \ref{lm:hetero3ver} and the derivative $\DerP(\theta)$ does not depend on $\theta$, we conclude $\|h\|_{\mathbf E} \leq \nu_n \|\DerP[h]\|_2$ for all $h \in \mathbf B_n$ -- i.e., in verifying the conditions of Lemma \ref{lm:lnotbind} we may set $\mathcal A_n(P) = \Theta_n \cap R$.
In order to verify condition \eqref{lm:lnotbinddisp1} of Lemma \ref{lm:lnotbind} we note that since $\|h\|_{\mathbf E} \asymp \sum_{\jmath = 1}^{\mathcal J} \|\beta_{\jmath,h}\|_2$ by Assumption \ref{ass:heterosieve}(iii), the definitions of the operator norm $\|\cdot\|_{o,2}$, $\hat{\mathbb D}_{\jmath,n}$, and $\mathbb D_{\jmath,P}$ imply that
\begin{equation*}
\sup_{h \in \mathbf B_n} \frac{\|\hat {\mathbb D}_n[h] - \DerP[h]\|_2}{\|h\|_{\mathbf E}} \\ \lesssim \|\frac{1}{n}\sum_{i=1}^n q^{k_n}(Z_i)p^{j_n}(W_i)^\prime - E_P[q^{k_n}(Z)p^{j_n}(W)^\prime]\|_{o,2} = o_P(1),
\end{equation*}
where the final equality holds uniformly in $P\in \mathbf P$ by applying Lemma \ref{lm:finiteder} with $b_{1n} = \sqrt{j_n}$, $b_{2n} = k_n$ (by Assumptions \ref{ass:heterosieve}(ii) and \ref{ass:heteromoments}(i)) and employing that $k_n \geq j_n$ and $k_n^2\log(k_n)/n = o(1)$ by Assumptions \ref{ass:heteromoments}(iii)(iv).
Finally, we note that $j_n^5 k_n^4 \log^3(1+j_nk_n)/n = o(a_n^4)$ by hypothesis, and employing Lemma \ref{lm:heterobootcoup} with $p=2$ yields that Assumption \ref{ass:bootcoupling} holds for $R = \Theta$, and hence also for $R$ as in \eqref{ex:hetero7paux}.
The only condition of Lemma \ref{lm:lnotbind} that remains to be verified is that $\mathcal S_n(\mathbf B,\mathbf E) \mathcal R_n = o(\ell_n)$.
To this end, we observe that since $\hat V_n(\theta,R|\ell_n)$ is defined through the constraint $\|h\|_{\mathbf E} \leq \ell_n$ (instead of $\|\cdot\|_{\mathbf B} \leq \ell_n$), it suffices to verify $\mathcal R_n = o(\ell_n)$ -- i.e. for the purposes of this lemma we may set $\|\cdot\|_{\mathbf B} = \|\cdot\|_{\mathbf E}$. However, since as argued $J_n \asymp \sqrt{j_n\log(1+j_n)}$, $B_n = \sqrt{k_n}$, and $\nu_n \asymp 1$, we have $\mathcal R_n \asymp k_n \sqrt{j_n\log(1+k_n)\log(1+j_n)}/\sqrt n$, and the requirement $\mathcal R_n = o(\ell_n)$ is implied by $k_n\sqrt{j_n}\log(n)/\sqrt n = o(\ell_n)$.
Thus, the claim of the lemma follow from Lemma \ref{lm:lnotbind}. \qed  

%% file: Appendix/ExamplesProofs/ExGMMProofs.tex

\subsection{Proofs for Section \ref{sec:gmmex}}\label{sec:gmmexproofs}

\noindent {\sc Proof of Theorem \ref{th:gmmapprox}:} We establish the theorem by simply applying Theorem \ref{th:localdrift}(ii) to both $R$ as in \eqref{eq:appgmm0p5} (to couple $I_n(R)$) and to $R = \Theta$ (to couple $I_n(\Theta)$).
To this end, note that as discussed Assumption \ref{ass:param}(ii)(iii) holds, while Assumption \ref{ass:param}(i) is directly imposed in \ref{ass:gmmbasic}(i).
Since $q^{k_n}(Z)$ equals the vector $(1,\ldots, 1)^\prime \in \mathbf R^{\mathcal J}$, it further follows Assumption \ref{ass:startreg}(i) holds with $B_n = 1$, while Assumption \ref{ass:startreg}(ii) is automatically satisfied.
We further note that Assumption \ref{ass:startreg}(iii) holds for $R = \Theta$ (and hence also for $R$ as in \eqref{eq:appgmm0p5}) with $J_n = C_0$ for some $C_0 < \infty$ by Assumption \ref{ass:gmmrho}(ii) and Lemma \ref{lm:gmmcov}.
Also note Assumption \ref{ass:coupling}(i) is satisfied for $R = \Theta$, and hence also for $R$ as in \eqref{eq:appgmm0p5}, by Lemma \ref{lm:gmmcoupling}.
Additionally, since $\Theta$ is convex by Assumption \ref{ass:gmmbasic}(iii), the mean value theorem and Assumption \ref{ass:gmmrho}(ii) imply that
$$E_P[\|\rho(X,\theta_1) - \rho(X,\theta_2)\|_2^2] \leq E_P[\sup_{\theta\in \Theta} \|\nabla_\theta \rho(X,\theta)\|_{o,2}^2] \|\theta_1 - \theta_2\|_2^2$$
for all $\theta_1,\theta_2\in \Theta$, which verifies Assumption \ref{ass:coupling}(ii) holds with $\kappa_\rho = 1$ and $\|\cdot\|_{\mathbf E} = \|\cdot\|_2$.
Lemma \ref{lm:gmm3ver} additionally verifies that Assumption \ref{ass:keycons} holds with $\|\cdot\|_{\mathbf E} = \|\cdot\|_2$ and $\nu_n^{-1} = \eta$ for some $\eta > 0$ when $R = \Theta$ and hence also when $R$ is as in \eqref{eq:appgmm0p5}.
Furthermore, we note that in this problem $\mathcal R_n \asymp n^{-1/2}$ because $\nu_n \asymp 1$, $J_n = O(1)$, $k_n = \mathcal J$, and $B_n = 1$.
To verify Assumption \ref{ass:driftlin}, note that in this application $\nabla m_{P,\jmath}(\theta) = E_P[\nabla_\theta \rho_\jmath(X,\theta)]$.
Hence, Assumptions \ref{ass:driftlin}(i)(ii) hold with $\|\cdot\|_{\mathbf L} = \|\cdot\|_2$ due to $E_P[\sup_{\theta \in \Theta} \|\nabla^2_\theta \rho_\jmath(X,\theta)\|_{o,2}]$ being bounded in $P\in \mathbf P$ by Assumption \ref{ass:gmmrho}(ii). Similarly, Assumption \ref{ass:driftlin}(iii) is satisfied due to $E_P[\sup_{\theta \in \Theta} \|\nabla_\theta \rho(X,\theta)\|_{o,2}]$ being bounded by Assumption \ref{ass:gmmrho}(ii).
Finally, we note that since $\mathcal R_n \asymp n^{-1/2}$ and $\kappa_\rho = 1$, Lemma \ref{lm:gmmcov} verifies Assumption \ref{ass:locrates}(i).
Assumption \ref{ass:locrates}(ii) is immediate since $E_P[\rho(X,\IDpoint)] = 0$, while Assumption \ref{ass:weights} holds by Assumption \ref{ass:gmmsigma}.
To conclude, simply note that the condition $k_n^{1/p} \sqrt{\log(1+k_n)}B_n\times \sup_{P\in \mathbf P} J_{[\hspace{0.02 in}]}(\ell_n^{\kappa_\rho},\mathcal F_n,\|\cdot\|_{P,2}) = o(a_n)$ is implied by $\ell_n\sqrt{\log(1/\ell_n)} = o(a_n)$ by Lemma \ref{lm:gmmcov}, and $K_m\mathcal R_n^2 = o(a_n/\sqrt n)$ is implied by $n^{-1/2} = o(a_n)$. \qed

\noindent {\sc Proof of Theorem \ref{th:gmmboot}:} We first define a variable $\hat E_n(R|\ell_n)$ to be given by
\begin{equation*}
\hat E_n(R|\ell_n) \equiv \inf_{h \in \hat V_n(\hat \theta_n,R|\ell_n)} \|\hat {\mathbb W}_n(\hat \theta_n) + \hat{\mathbb D}_n(\hat \theta_n) [h]\|_{\hat \Sigma_n,2}
\end{equation*}
and note Lemma \ref{lm:gmmnoell} implies $\hat U_n(R|+\infty) = \hat E_n(R|\ell_n) + o_P(a_n)$ uniformly in $P\in \mathbf P_0$ for any $\ell_n\downarrow 0$ satisfying the conditions of the theorem.
Therefore, to establish the theorem it suffices to show that uniformly in $P\in \mathbf P_0$ we have
\begin{align*}
\hat E_n(R|\ell_n) & \geq \UpS (R|\tilde \ell_n) + o_P(a_n) \\
\hat E_n(R|\ell_n) - \hat U_n(\Theta|+\infty) & \geq \UpS (R|\tilde \ell_n) - \UpS (\Theta|\tilde \ell_n^{\text{u}}) + o_P(a_n).
\end{align*}
with $\ell_n \asymp \tilde \ell_n$ and $\tilde \ell_n^{\rm u}$ satisfying the conditions of the theorem.
To this end we rely on Theorem \ref{th:coupsmooth} (for $\hat E_n(R|\ell_n)$) and Lemma \ref{lm:forincJboot}.
Next note that in the proof of Theorem \ref{th:gmmapprox} we established that Assumptions \ref{ass:gmmbasic}, \ref{ass:gmmrho}, \ref{ass:gmmid}, and \ref{ass:gmmsigma}, imply Assumptions \ref{ass:param}, \ref{ass:startreg}, \ref{ass:coupling}, \ref{ass:keycons},  \ref{ass:driftlin}, \ref{ass:locrates}, and  \ref{ass:weights} hold with $\mathcal R_n \asymp n^{-1/2}$, $\nu_n \asymp 1$,  $\|\cdot\|_{\mathbf B} = \|\cdot\|_{\mathbf E} = \|\cdot\|_{\mathbf L} = \|\cdot\|_2$, $\kappa_\rho = 1$, and $a_n = \sqrt{\log(n)}/n^{\frac{1}{10 + 5d_\theta}}$ for $R = \Theta$ and $R$ as in \eqref{eq:appgmm0p5}.
We thus avoid repeating the arguments, and verify only that Assumptions \ref{ass:locineq}, \ref{ass:loceq}, \ref{ass:ineqlindep}, \ref{ass:bootcoupling}, \ref{ass:extra}, and \ref{ass:bootrates} hold for $R = \Theta$ and $R$ as in \eqref{eq:appgmm0p5}.

Next note Lemma \ref{lm:gmm4ver} implies Assumptions \ref{ass:locineq}, \ref{ass:loceq}, and \ref{ass:ineqlindep} are satisfied, while Lemma \ref{lm:bootgmmcoupling} verifies Assumption \ref{ass:bootcoupling} with $a_n = \sqrt{\log(n)}/n^{\frac{1}{10 + 5d_\theta}}$ for $R = \Theta$, and hence also for $R$ as in \eqref{eq:appgmm0p5}.
Assumption \ref{ass:extra}(i) is immediate since $\|\cdot\|_{\mathbf E} = \|\cdot\|_{\mathbf B} = \|\cdot\|_2$, while Assumptions \ref{ass:extra}(ii)(iii) are implied by Assumption \ref{ass:gmmbootp1}(i), $\|\hat \theta_n - \IDpoint\|_2 = o_P(1)$ uniformly in $P\in \mathbf P_0$ (which we showed in establishing Theorem \ref{th:gmmapprox}), and $\mathcal V_n(P) \equiv \{\theta \in \Theta :\|\theta - \IDpoint\|_2 \leq \epsilon\}$ for some $\epsilon > 0$ by Lemma \ref{lm:gmm3ver}.
Assumption \ref{ass:bootrates}(i) is immediate since $\mathcal S_n(\mathbf B,\mathbf E) = 1$ and the choices of $\hat \theta_n$ and $\hat \theta_n^{\text{u}}$ correspond to setting $\tau_n = o(n^{-1/2})$.
Similarly,  Lemma \ref{lm:gmmcov}, $\mathcal S_n(\mathbf L,\mathbf E) = 1$, and $n^{-1/2} = o(\ell_n)$ imply that the condition $\ell_n^2 \sqrt{\log(1/\ell_n)} = o(a_nn^{-\frac{1}{2}})$ verifies Assumption \ref{ass:bootrates}(ii).
Moreover, since $\ell_n = o(r_n)$ and $n^{-1/2} = o(r_n)$ Assumption \ref{ass:bootrates}(iii) holds.
Hence, Theorem \ref{th:coupsmooth} implies
\begin{equation}\label{eq:gmmboot1}
\hat E_n(R|\ell_n)  \geq  \UpS (R|\tilde \ell_n) + o_P(a_n)
\end{equation}
uniformly in $P\in \mathbf P_0$ for some $\ell_n \asymp \tilde \ell_n$.
Similarly, since $\mathcal R_n^{\text{u}} \asymp n^{-1/2}$, the conditions of Lemma \ref{lm:forincJboot}(ii) are immediate and hence by \eqref{eq:gmmboot1} there are $\ell_n\asymp \tilde \ell_n$ and $\ell_n^{\text{u}}\asymp \tilde \ell_n^{\text{u}}$ with
\begin{equation}\label{eq:gmmboot2}
\hat E_n(R|\ell_n) - \hat U_n(\Theta|+\infty) \geq \UpS (R|\tilde \ell_n) - \UpS (\Theta|\tilde \ell_n^{\text{u}}) + o_P(a_n).
\end{equation}
The theorem therefore follows from \eqref{eq:gmmboot1}, \eqref{eq:gmmboot2} and Lemma \ref{lm:gmmnoell}. \qed

\begin{lemma}\label{lm:gmm3ver}
If Assumptions \ref{ass:gmmbasic}, \ref{ass:gmmrho}, \ref{ass:gmmid}, and \ref{ass:gmmsigma}(ii) hold, then Assumption \ref{ass:keycons} is satisfied with $R = \Theta$ and $R$ as in \eqref{eq:appgmm0p5}, $\|\cdot\|_{\mathbf E} =\|\cdot\|_2$, $\nu_n^{-1} = \eta$ for some $\eta > 0$, and $\mathcal V_n(P) \equiv \{\theta \in \Theta : \|\theta - \IDpoint\|_2 \leq \epsilon\}$ for some $\epsilon > 0$.
\end{lemma}

\noindent {\sc Proof:} To verify Assumption \ref{ass:keycons}(ii), note Assumptions \ref{ass:gmmbasic}(i), \ref{ass:gmmrho}(ii), \ref{ass:gmmsigma}, and \ref{ass:gmmid}(i) and Lemma \ref{lm:gmmcov} allow us to apply Lemma \ref{lm:setcons}(i) with $\|\cdot\|_{\mathbf A} = \|\cdot\|_2$, $J_n = O(1)$ and $S_n(\epsilon) > 0$ to conclude $\hat \theta_n \in \mathcal V_n(P) \equiv \{\theta \in \Theta_n : \|\theta - \IDpoint\|_2 \leq \epsilon\}$ with probability tending to one uniformly in $P\in \mathbf P_0$ for any $\epsilon > 0$ and for both $R = \Theta$ and $R$ as in \eqref{eq:appgmm0p5}.
In order to verify Assumption \ref{ass:keycons}(i), next note that $\Theta$ being convex and Assumption \ref{ass:gmmrho}(ii) imply that for some $C_0 < \infty$ we have
\begin{equation*}\label{lm:gmm3ver1}
\|E_P[\rho(X,\theta)] - E_P[\rho(X,\IDpoint)] - E_P[\nabla_\theta \rho(X,\IDpoint)](\theta - \IDpoint)\|_2 \leq C_0 \|\theta - \IDpoint\|_2^2
\end{equation*}
for all $\theta \in \Theta$.
Hence, since the smallest singular value of $E_P[\nabla_\theta \rho(X,\IDpoint)]$ is bounded away from zero uniformly in $P\in \mathbf P_0$ by Assumption \ref{ass:gmmid}(ii), we obtain for some $C_1 < \infty$
\begin{align}\label{lm:gmm3ver2}
\|\theta - \IDpoint\|_2 & \leq C_1 \|E_P[\nabla_\theta \rho(X,\IDpoint)](\theta - \IDpoint)\|_2 \notag \\
& \leq C_1\{\|E_P[\rho(X,\theta)] - E_P[\rho(X,\IDpoint)]\|_2 + C_0\|\theta - \IDpoint\|_2^2\}
\end{align}
for all $\theta \in \Theta$ and $P\in \mathbf P_0$.
Therefore, provided $\epsilon > 0$ is set sufficiently small in defining $\mathcal V_n(P) \equiv \{\theta \in \Theta_n : \|\theta - \IDpoint\|_2 \leq \epsilon\}$, it follows that Assumption \ref{ass:keycons}(i) holds with $\|\cdot\|_{\mathbf E} = \|\cdot\|_2$ and $\nu_n^{-1} = \eta$ for some $\eta > 0$ due to \eqref{lm:gmm3ver2} and Assumption \ref{ass:gmmsigma}(ii). \qed

\begin{lemma}\label{lm:gmmcov}
Let $\mathcal F \equiv \{\rho_{\jmath}(\cdot,\theta) : \text{ for some } \theta\in \Theta \text{ and } 1\leq \jmath \leq \mathcal J\}$.
If Assumptions \ref{ass:gmmbasic}(iii) and \ref{ass:gmmrho} hold, then it follows that $\sup_{P\in \mathbf P} N_{[\hspace{0.02 in}]}(\epsilon,\mathcal F,\|\cdot\|_{P,2})\lesssim 1\vee \epsilon^{-d_\theta}$ and $\sup_{P\in \mathbf P}J_{[\hspace{0.02 in}]}(\epsilon,\mathcal F,\|\cdot\|_{P,2}) \lesssim \epsilon (1+\sqrt {\log(1 \vee \epsilon^{-1})})$.
\end{lemma}

\noindent {\sc Proof:} Since $\Theta$ is convex by Assumption \ref{ass:gmmbasic}(iii), the mean value theorem and Assumption \ref{ass:gmmrho}(i) imply for any $\theta_1,\theta_2\in \Theta$ and $1\leq \jmath \leq \mathcal J$ that
\begin{equation}\label{lm:gmmcov1}
|\rho_\jmath(x,\theta_1) - \rho_\jmath(x,\theta_2)| \leq \sup_{\theta \in \Theta} \|\nabla_\theta \rho(x,\theta)\|_{o,2} \|\theta_1 - \theta_2\|_2.
\end{equation}
Setting $D(x)\equiv \sup_{\theta \in \Theta} \|\nabla_\theta \rho(x,\theta)\|_{o,2}$, then note that Theorem 2.7.11 in \cite{vandervaart:wellner:1996} and the right hand side of \eqref{lm:gmmcov1} not depending on $\jmath$ imply
\begin{equation}\label{lm:gmmcov2}
N_{[\hspace{0.02 in}]}(\epsilon,\mathcal F,\|\cdot\|_{P,2}) \leq \mathcal J \times N(\frac{\epsilon}{2\|D\|_{P,2}}, \Theta, \|\cdot\|_2) \lesssim 1\vee \epsilon^{-d_\theta},
\end{equation}
where we employed that $N(\epsilon,\Theta,\|\cdot\|_2) \lesssim 1\vee \epsilon^{-d_\theta}$ due to $\Theta$ being bounded by Assumption \ref{ass:gmmbasic}(iii) and $\sup_{P\in \mathbf P} \|D\|_{P,2} < \infty$ by Assumption \ref{ass:gmmrho}(ii).

For the second claim of the Lemma we employ the bound in \eqref{lm:gmmcov2} to obtain
\begin{multline*}\label{lm:gmmcov3}
\sup_{P\in \mathbf P} J_{[\hspace{0.02 in}]}(\epsilon,\mathcal F,\|\cdot\|_{P,2}) \lesssim \int_0^\epsilon (1+ \log (1\vee u^{-d_\theta}))^{1/2}du\\
= \epsilon \int_0^1 (1+ \log(1\vee (\epsilon v)^{-d_\theta}))^{1/2} dv \lesssim \epsilon (1+ \sqrt{\log(1\vee \epsilon^{-1})}),
\end{multline*}
where the first equality follows from the change of variables $v = u/\epsilon$ and the final inequality is implied by the inequality $1\vee (ab) \leq (1\vee a)(1\vee b)$. \qed

\begin{lemma}\label{lm:gmmcoupling}
If Assumptions \ref{ass:gmmbasic}(i)(iii) and \ref{ass:gmmrho} hold, then it follows that Assumption \ref{ass:coupling}(i) is satisfied with $R = \Theta$ and $a_n = \sqrt{\log(n)}/n^{\frac{1}{6+5d_\theta}}$.
\end{lemma}

\noindent {\sc Proof:} Let $\epsilon_n = \sqrt{\log(n)}/n^{\frac{1}{6+5d_\theta}}$ and set $\delta_{n} \equiv 1\wedge (\epsilon_n^2 \sqrt n)^{-\frac{2}{2+5d_\theta}}$, which note satisfies $1\geq \delta_n = o(1)$.
Further define $N_n \equiv N(\delta_n,\Theta,\|\cdot\|_2)$ and set $\{\theta_k\}_{k=1}^{N_n}$ to be the center of the $N_n$ balls covering $\Theta$.
For notational simplicity, we also let
\begin{equation*}
r_{n,P}(x) \equiv ((\rho(x,\theta_1) - E_P[\rho(X,\theta_1)])^\prime,\ldots , (\rho(x,\theta_{N_n}) - E_P[\rho(X,\theta_{N_n})])^\prime)^\prime
\end{equation*}
and note $r_{n,P}(x) \in \mathbf R^{\mathcal J N_n}$.
For any $P\in \mathbf P$ and $\eta > 0$ further define $C_{n,P}(\eta)$ to equal
\begin{equation}\label{lm:gmmcoupling2}
C_{n,P}(\eta) \equiv \frac{(\mathcal J N_n) E_P[\|r_{n,P}(X)\|_2^3]}{\eta^3\epsilon_n^3 \sqrt n}.
\end{equation}
It then follows by Yurinskii's coupling (see, e.g., Theorem 10.10 in \cite{pollard2002user}) that there exists a Gaussian vector $\mathbb N_{n,P}\in \mathbf R^{\mathcal J N_n}$ and universal constant $K_0$ such that
\begin{equation}\label{lm:gmmcoupling3}
 P(\|\frac{1}{\sqrt n}\sum_{i=1}^n r_{n,P}(X_i)  - \mathbb N_{n,P}\|_2 > 3\eta\epsilon_n) \leq K_0 C_{n,P}(\eta)(1 + \frac{|\log(1/C_{n,P}(\eta))|}{\mathcal J N_n}).
\end{equation}
Next note Assumption \ref{ass:gmmrho}(ii), Jensen's inequality, and the convexity of $u\mapsto |u|^{\frac{3}{2}}$ yield
\begin{equation}\label{lm:gmmcoupling4}
\sup_{P\in \mathbf P} E_P[\|r_{n,P}(X)\|_2^3] \lesssim (\mathcal J N_n)^{\frac{3}{2}} \times \sup_{P\in \mathbf P} \frac{1}{\mathcal J N_n} \sum_{k=1}^{N_n} \sum_{\jmath = 1}^{\mathcal J} E_P[ |\rho_{\jmath}(X,\theta_k)|^3] \lesssim N_n^{\frac{3}{2}}.
\end{equation}
In particular, since $N(\epsilon, \Theta,\|\cdot\|_2)\lesssim 1\vee \epsilon^{-d_\theta}$, it follows from $\delta_{n}\leq 1$ that $N_n \lesssim \delta_{n}^{-d_\theta}$, and hence by \eqref{lm:gmmcoupling4} and the definition of $C_{n,P}(\eta)$ in \eqref{lm:gmmcoupling2} we obtain
\begin{equation}\label{lm:gmmcoupling5}
\sup_{P\in \mathbf P} C_{n,P}(\eta) \lesssim \frac{N_n^{\frac{5}{2}}}{\eta^3\epsilon_n^3 \sqrt n} \lesssim \frac{1}{\eta^3\epsilon_n^3(n\delta_{n}^{5d_\theta})^{\frac{1}{2}}}.
\end{equation}
Moreover, since the function $u\mapsto u(1 + |\log(1/u)|/A)$ with $A\geq 1$ is increasing in $u$ on the interval $(0,1]$ and $\epsilon_n^3(n\delta_n^{5d_\theta})^{\frac{1}{2}} \to \infty$, we obtain from results \eqref{lm:gmmcoupling3} and \eqref{lm:gmmcoupling5} that
\begin{multline}\label{lm:gmmcoupling6}
\limsup_{n\rightarrow \infty} \sup_{P\in \mathbf P} P(\|\frac{1}{\sqrt n}\sum_{i=1}^n (r_{n,P}(X_i)  - \mathbb N_{n,P}\|_2 > 3\eta \epsilon_n) \\
\lesssim  \limsup_{n\rightarrow \infty} \frac{1}{\eta^3\epsilon_n^3(n\delta_{n}^{5d_\theta})^{\frac{1}{2}}}(1 + \frac{|\log(\eta^3\epsilon_n^3(n\delta_{n}^{5d_\theta})^{\frac{1}{2}})|}{\mathcal J N_n}) = 0,
\end{multline}
where the final result follows by direct calculation.
Letting $\mathbb S_{n,P}$ denote the linear span of $r_{n,P}$ in $L^2_P$ we then employ $\mathbb N_{n,P}$ to define a Gaussian process $\Iso^{(1)}$ on $\mathbb S_{n,P}$ by setting
\begin{equation}\label{lm:gmmcoupling7}
\Iso^{(1)}(\sum_{k=1}^{N_n} \lambda_{k}^\prime \rho (\cdot,\theta_k)) \equiv (\lambda_{1}^\prime,\ldots, \lambda_{N_n}^\prime) \mathbb N_{n,P}
\end{equation}
for any $\{\lambda_k\}_{k=1}^{N_n}$ with $\lambda_k \in \mathbf R^{\mathcal J}$.
Letting $\text{Proj}\{f|\mathbb S_{n,P}\}$ denote the projection of $f $ onto $\mathbb S_{n,P}$ under $\|\cdot\|_{P,2}$, and assuming the probability space is suitably large to carry an isonormal process $\Iso^{(2)}$ on $\{(f -\int fdP) - \text{Proj}\{f-\int f dP|\mathbb S_{n,P}\} : f\in \mathcal F\}$ that is independent of $\Iso^{(1)}$, we then define the isonormal process $\Iso$ to be given by
\begin{equation}\label{lm:gmmcoupling8}
\Iso(f) \equiv \Iso^{(1)}(\text{Proj}\{f|\mathbb S_{n,P}\}) + \Iso^{(2)}(f - \text{Proj}\{f|\mathbb S_{n,P}\}).
\end{equation}

Next let $\Pi_n \theta$ denote the projection of any $\theta \in \Theta$ onto $\{\theta_k\}_{k=1}^{N_n}$ under $\|\cdot\|_2$ and define
\begin{equation}\label{lm:gmmcoupling9}
\mathcal G_{n,P} \equiv \{(\rho_\jmath(\cdot,\theta) - \rho_\jmath(\cdot,\Pi_n\theta))- E_P[(\rho_\jmath(X,\theta)-\rho_\jmath(X,\Pi_n\theta))] : \theta \in \Theta, 1\leq \jmath \leq \mathcal J\}.
\end{equation}
By the mean value theorem, $\Theta$ being convex by Assumption \ref{ass:gmmbasic}(iii), and $\|\theta - \Pi_n\theta\|_2 \leq \delta_n$ for every $\theta \in \Theta$ due to $\delta_n$-balls around $\{\theta_k\}_{k=1}^{N_n}$ covering $\Theta$, it follows that
\begin{multline*}
\sup_{\theta \in \Theta} |(\rho_\jmath(x,\theta) - \rho_\jmath(x,\Pi_n\theta)) - E_P[(\rho_\jmath(X,\theta)-\rho_\jmath(X,\Pi_n\theta))]| \\
\leq \{\sup_{\theta \in \Theta} \|\nabla_\theta \rho(x,\theta)\|_{o,2} + \sup_{P\in \mathbf P} E_P[\sup_{\theta \in \Theta} \|\nabla_\theta \rho(X,\theta)\|_{o,2}]\}\times \delta_n.
\end{multline*}
Hence, setting $G(x) \equiv 1\vee \{\sup_{\theta \in \Theta} \|\nabla_\theta \rho(x,\theta)\|_{o,2} + \sup_{P\in \mathbf P} E_P[\sup_{\theta \in \Theta} \|\nabla_\theta \rho(X,\theta)\|_{o,2}]\}$ it follows that $G\delta_n$ is an envelope for $\mathcal G_{n,P}$, which by Assumption \ref{ass:gmmrho}(ii) satisfies $\sup_{P\in \mathbf P}\|G\delta_n\|_{P,2} \lesssim \delta_n$.
Further note that if $[f_l,f_u]$ is a bracket containing a function $f$, then $[f_l - E_P[f_u(X)], f_u - E_P[f_l(X)]]$ contains $f - E_P[f(X)]$ and satisfies
\begin{equation*}
\|f_u - f_l - E_P[f_l(X) - f_u(X)]\|_{P,2} \leq 2\|f_u - f_l\|_{P,2}
\end{equation*}
by Jensen's inequality and the triangle inequality.
Therefore, Lemma \ref{lm:gmmcov} implies
\begin{equation*}\label{lm:gmmcoupling12}
\sup_{P\in \mathbf P} N_{[\hspace{0.02 in}]}(\epsilon,\mathcal G_{n,P},\|\cdot\|_{P,2}) \lesssim N_n\times (1\vee \epsilon^{-d_\theta}),
\end{equation*}
and hence Theorem 2.14.2 in \cite{vandervaart:wellner:1996} together with $\mathcal G_{n,P}$ having envelope $\delta_n G$ with $G \geq 1$, $\sup_{P\in \mathbf P} \|G\|_{P,2} < \infty$, and $N_n \lesssim \delta_n^{-d_\theta}$ yield
\begin{align}\label{lm:gmmcoupling13}
\sup_{P\in \mathbf P} E_P[\sup_{g\in \mathcal G_{n,P}}& |\frac{1}{\sqrt n}\sum_{i=1}^n(g(X_i)-E_P[g(X)])|] \notag \\
& \lesssim  \sup_{P\in \mathbf P}\{\delta_n \|G\|_{P,2} \int_0^{1} (1 + \log N_{[\hspace{0.02 in}]}(\epsilon\delta_n \|G\|_{P,2},\mathcal G_{n,P},\|\cdot\|_{P,2}))^{\frac{1}{2}}d\epsilon\} \notag \\
& \lesssim \delta_n  \int_0^{1}( 1 + \log(N_n) + \log(1\vee (\epsilon\delta_n)^{-d_\theta}))^{1/2}d\epsilon\notag \\  & \lesssim \delta_n (1+ \log(\delta_n^{-d_\theta}))^{1/2}.
\end{align}
Therefore, the definitions of $\delta_n$ and $\epsilon_n$, result \eqref{lm:gmmcoupling13} and Markov's inequality imply
\begin{multline}\label{lm:gmmcoupling14}
\limsup_{n\rightarrow \infty}\sup_{P\in \mathbf P} P(\sup_{g\in \mathcal G_{n,P}} |\frac{1}{\sqrt n}\sum_{i=1}^n (g(X_i) - E_P[g(X)])| > \eta \epsilon_n) \\
\lesssim \limsup_{n\rightarrow \infty} \frac{\delta_n(1+\log(\delta_n^{-d_\theta}))^{1/2}}{\eta \epsilon_n} = 0.
\end{multline}
Similarly, since $\Iso$ is Gaussian and $0 \in \mathcal G_{n,P}$, Corollary 2.2.8 in \cite{vandervaart:wellner:1996} and packing numbers being bounded by bracketing numbers imply
\begin{multline}\label{lm:gmmcoupling15}
\sup_{P\in \mathbf P} E_P[\sup_{g\in \mathcal G_{n,P}}|\Iso(g)|] \lesssim \sup_{P\in \mathbf P} \int_0^{\infty} (\log N_{[\hspace{0.02 in}]}(\epsilon,\mathcal G_{n,P},\|\cdot\|_{P,2}))^{\frac{1}{2}} d\epsilon \\
\lesssim \sup_{P\in \mathbf P} \int_0^{2\delta_n\|G\|_{P,2}}(\log N_{[\hspace{0.02 in}]}(\epsilon,\mathcal G_{n,P},\|\cdot\|_{P,2}))^{\frac{1}{2}} d\epsilon\lesssim \delta_n(1+ \log(\delta_n^{-d_\theta}))^{1/2},
\end{multline}
where in the second inequality we employed that the bracket $[-\delta_nG, \delta_n G]$ covers $\mathcal G_{n,P}$ due to $\delta_n G$ being an envelope for $\mathcal G_{n,P}$, and the final inequality follows from the change of variables $u = \epsilon/(2\delta_n\|G\|_{P,2})$ and the same manipulations as in \eqref{lm:gmmcoupling13}.
Hence,
\begin{equation}\label{lm:gmmcoupling16}
\limsup_{n\rightarrow \infty}\sup_{P\in \mathbf P} P(\sup_{g\in \mathcal G_{n,P}} |\Iso(g)| > \eta \epsilon_n) \lesssim \limsup_{n\rightarrow \infty} \frac{\delta_n(1+\log(\delta_n^{-d_\theta}))^{1/2}}{\eta \epsilon_n} = 0,
\end{equation}
by result \eqref{lm:gmmcoupling15} and Markov's inequality.
To conclude, for any $\theta \in \Theta$ set $\mathbb W_P(\theta)$ to be
$$\mathbb W_P(\theta) \equiv (\mathbb G_P(\rho_1(\cdot,\theta)),\ldots, \mathbb G_P(\rho_{\mathcal J}(\cdot,\theta))^\prime$$ 
and note that the definitions of $\mathbb G_{P}$ in \eqref{lm:gmmcoupling7} and \eqref{lm:gmmcoupling8}, and of $\mathcal G_{n,P}$ in \eqref{lm:gmmcoupling9}, yield
\begin{multline*}
\sup_{\theta \in \Theta} \|\mathbb G_{n}(\theta) - \mathbb W_{P}(\theta)\|_2 \leq \|\frac{1}{\sqrt n}\sum_{i=1}^n r_{n,P}(X_i) -\mathbb N_{n,P}\|_2 \\+ \sup_{g\in \mathcal G_{n,P}} \sqrt{\mathcal J} |\frac{1}{\sqrt n}\sum_{i=1}^n (g(X_i) - E_P[g(X)]) | + \sup_{g\in \mathcal G_{n,P}} \sqrt{\mathcal J}|\mathbb G_{P}(g)|.
\end{multline*}
Thus the lemma follows from \eqref{lm:gmmcoupling6}, \eqref{lm:gmmcoupling14}, and \eqref{lm:gmmcoupling16}. \qed

\begin{lemma}\label{lm:bootgmmcoupling}
If Assumptions \ref{ass:gmmbasic}(i)(iii) and \ref{ass:gmmrho} hold, then it follows that Assumption \ref{ass:bootcoupling} is satisfied with $R = \Theta$ and $a_n = \log^{3/4}(n)/n^{\frac{1}{12+2d_\theta}}$.
\end{lemma}

\noindent {\sc Proof:}
We establish the lemma by relying on Theorem \ref{th:mainbootcoup}(i) in Section \ref{sec:bootcoup}.
To this end set $\zeta_n = n^{-\frac{1}{2(6+d_\theta)}}$, $M_n = n^{\frac{1}{6+ d_\theta}}$, and $N_n \equiv N(\zeta_n,\Theta,\|\cdot\|_2)$.
By Assumption \ref{ass:gmmrho}(ii) the function $F(x) \equiv (1+\sup_{\theta \in \Theta} \|\rho(x,\theta)\|_2)$ is integrable, and for any $\theta \in \Theta$ we let
\begin{equation*}\label{lm:bootgmmcoupling1}
\tilde \rho(x,\theta) \equiv (\rho_1(x,\theta)1\{F(x) \leq M_n\},\ldots, \rho_{\mathcal J}(x,\theta)1\{F(x) \leq M_n\})^\prime.
\end{equation*}
Defining $d_n = \mathcal J N_n$ and $\{\theta_k\}_{k=1}^{N_n}$ to be the centers of the $\zeta_n$-balls covering $\Theta$ we then let
\begin{equation*}\label{lm:bootgmmcoupling2}
f_{n}^{d_n}(X) \equiv (\tilde \rho(X,\theta_1)^\prime - E_P[\tilde \rho(X,\theta_1)^\prime] ,\ldots, \tilde \rho(X,\theta_{N_n})^\prime - E_P[\tilde \rho(X,\theta_{N_n})^\prime])^\prime.
\end{equation*}
Next note that since each entry of the matrix $f_n^{d_n}(X)f_{n}^{d_n}(X)^\prime$ is almost surely bounded by $2M_n^2$ it follows that $\|E_P[f_n^{d_n}(X)f_n^{d_n}(X)^\prime]\|_{o,2} \leq 2d_n M_n^2$, and hence Assumption \ref{ass:4seriescoup} in Section \ref{sec:bootcoup} holds with $C_n \asymp d_n M_n^2$ and $K_n \asymp M_n$.
For every $\theta \in \Theta$ let $\Pi_n \theta$ denote its projection (under $\|\cdot\|_2$) onto $\{\theta_k\}_{k=1}^{N_n}$ and define the class $\mathcal G_{n,P} \equiv \{(\rho_j(\cdot,\theta) - \tilde \rho_j(\cdot,\Pi_n \theta)) - E_P[\rho_j(X,\theta) - \tilde \rho_j(X,\Pi_n \theta)] : \theta \in \Theta \text{ and } 1\leq \jmath \leq \mathcal J\}$.
Further observe that
\begin{align}\label{lm:bootgmmcoupling2p5}
 &\sup_{g \in \mathcal G_{n,P}} |g(x)| \notag \\ & \leq \max_{1\leq \jmath \leq \mathcal J} \sup_{\theta \in \Theta} 2|\rho_\jmath(x,\theta) - \rho_\jmath(x,\Pi_n\theta)| + F(x)1\{F(x) > M_n\} + E_P[F(X)1\{F(X)>M_n\}]\notag \\
& \leq \sup_{\theta \in \Theta} 2\|\nabla_\theta \rho(x,\theta)\|_{o,2} \|\theta - \Pi_n\theta\|_2 + F(x)1\{F(x) > M_n\} + E_P[F(X)1\{F(X)>M_n\}]
\end{align}
where in the second inequality we employed the mean value theorem and $\Theta$ being convex by Assumption \ref{ass:gmmbasic}(iii).
In particular, since the $\zeta_n$-balls centered around $\{\theta_k\}_{k=1}^{N_n}$ cover $\Theta$ and $\zeta_n \leq 1$, result \eqref{lm:bootgmmcoupling2p5} implies that the function
\begin{equation*}
G(x) \equiv 2\sup_{\theta \in \Theta} \|\nabla_\theta \rho(x,\theta)\|_{o,2}  + F(x) + \sup_{P\in \mathbf P}E_P[F(X)]
\end{equation*}
is an envelope for $\mathcal G_{n,P}$, while Assumption \ref{ass:gmmrho}(ii) implies $\sup_{P\in \mathbf P} E_P[G^2(X)] < \infty$.
Moreover, result \eqref{lm:bootgmmcoupling2p5} and Markov's, Jensen's, and Holder's inequalities  yield that
\begin{align}\label{lm:bootgmmcoupling3}
\sup_{g\in \mathcal G_{n,P}} \|g\|_{P,2}  & \leq \zeta_n \|G\|_{P,2} +  2\{E_P[F^3(X)]\}^{\frac{2}{3}}\{P(F(X) > M_n)\}^{\frac{1}{3}} \notag \\ & \leq (\zeta_n + M_n^{-1/2}\times 2\sup_{P\in \mathbf P} (E_P[F^3(X)])^{1/2})\times \|G\|_{P,2},
\end{align}
where in the final equality we employed that $\|G\|_{P,2} \geq 1$ because $F(X) \geq 1$.
Thus, by result \eqref{lm:bootgmmcoupling3} and Assumption \ref{ass:gmmrho}(ii), we may set $\delta_n \equiv C(\zeta_n +M_n^{-1/2})$ and obtain $\|g\|_{P,2} \leq \delta_n \|G\|_{P,2}$ for all $g\in \mathcal G_{n,P}$ and $P\in \mathbf P$ provided $C$ is chosen large enough.
Next note that since $\Theta$ being bounded by Assumption \ref{ass:gmmbasic}(iii) implies $N_n \lesssim \zeta_n^{-d_\theta}$, we obtain
\begin{equation}\label{lm:bootgmmcoupling5}
\sup_{P\in \mathbf P}N_{[\hspace{0.02 in}]}(\epsilon,\mathcal G_{n,P},\|\cdot\|_{P,2}) \lesssim N_n\times (1\vee \epsilon)^{-d_\theta} \lesssim \zeta_n^{-d_\theta}\times (1\vee \epsilon)^{-d_\theta}
\end{equation}
due to Lemma \ref{lm:gmmcov}.
Hence, the change of variables $u = \epsilon/(\delta_n \|G\|_{P,2})$ implies that
\begin{align}\label{lm:bootgmmcoupling6}
\sup_{P\in \mathbf P} \int_0^{\delta_n\|G\|_{P,2}} (1+ & \log N_{[\hspace{0.02 in}]}(\epsilon,\mathcal G_{n,P},\|\cdot\|_{P,2}))^{1/2}d\epsilon \notag\\
& \lesssim \sup_{P\in \mathbf P} \delta_n\|G\|_{P,2} \int_0^1 (1+ \log(\zeta_n^{-d_\theta}) + \log(1\vee(u \delta_n \|G\|_{P,2})^{-d_\theta}))^{1/2}du \notag \\ &\lesssim \delta_n (1+\log(\zeta_n^{-1}))^{1/2}
\end{align}
where in the inequalities we employed result \eqref{lm:bootgmmcoupling5}, $\zeta_n \lesssim \delta_n$, and $\sup_{P\in \mathbf P} \|G\|_{P,2} < \infty$.
In particular, results \eqref{lm:bootgmmcoupling5} and \eqref{lm:bootgmmcoupling6} together with Lemma \ref{lm:auxbootcalc} imply that Assumption \ref{ass:4seriesreg}(i) in Section \ref{sec:bootcoup} is satisfied with $J_{1n} \lesssim \delta_n(1+\log(\zeta_n^{-1}))^{1/2}$.
Similarly, note that in this application, the set $\mathcal B_n$ in Assumption \ref{ass:4seriesreg}(ii) consists of $0 \in \mathbf R^{d_n}$ and the set of vectors in $\mathbf R^{d_n}$ with one coordinate equal to one and all other coordinates equal to zero.
Thus, Assumption \ref{ass:4seriesreg}(ii) holds with $J_{2n} = (\log(1+ d_n))^{1/2} \lesssim (1+\log(\zeta_n^{-1}))^{1/2}$.

We have so far verified Assumptions \ref{ass:4seriescoup} and \ref{ass:4seriesreg} in Section \ref{sec:bootcoup} hold with $d_n \lesssim \zeta_n^{-d_\theta}$, $K_n = M_n$, $C_n \lesssim M_n^2\zeta_n^{-d_\theta}$, $J_{1n} \lesssim \delta_n(1+\log(\zeta_n^{-1}) )^{1/2}$, and $J_{2n} \lesssim (1+\log(\zeta_n^{-1}))^{1/2}$.
Since we had set $\zeta_n = n^{-1/(2(6+d_\theta))}$, $M_n = n^{1/(6+d_\theta)}$, and $\delta_n = C(\zeta_n + M_n^{-1/2})$ the requirement that $d_n\log(1+d_n)K_n^2 = o(n)$ imposed by Theorem \ref{th:mainbootcoup}(i) holds as well.
Therefore, Assumption \ref{ass:gmmbasic}(i) and Theorem \ref{th:mainbootcoup}(i) finally enable us to conclude that there exists a process $\WPT$ that is independent of the data
$\{X_i\}_{i=1}^n$ and such that
\begin{equation*}
\sup_{\theta \in \Theta} \|\hat{\mathbb W}_{n} (\theta) - \WPT(\theta)\|_2  = O_P(\log^{3/4}(n)n^{-1/(12+2d_\theta)})
\end{equation*}
uniformly in $P\in \mathbf P$, which conclude the proof of the lemma. \qed

\begin{lemma}\label{lm:gmm4ver}
If Assumption \ref{ass:gmmbasic}(ii), \ref{ass:gmmbootp1}(ii)-(vi), and \ref{ass:gmmbootp2} hold, then it follows that Assumptions \ref{ass:locineq}, \ref{ass:loceq}, and \ref{ass:ineqlindep} are satisfied.
\end{lemma}

\noindent {\sc Proof:} Recall that in this setting $\mathbf G = \mathbf R^{d_G}$ and $\|\cdot\|_{\mathbf G} = \|\cdot\|_\infty$.
For $\epsilon$ and $B^\epsilon$ as in Assumption \ref{ass:gmmbootp1}, let $N_\epsilon(\IDpoint) \equiv \{\theta \in \Theta : \|\theta - \IDpoint\|_2\leq \epsilon\}$ noting that $N_\epsilon(\IDpoint) \subseteq B^\epsilon$ and that $N_\epsilon(\IDpoint)$ implicitly depends on $P$ through $\IDpoint$ (which depends on $P$ through \eqref{eq:appgmm0}).
For any $\theta_1,\theta_2\in N_\epsilon(\IDpoint)$, $N_\epsilon(\theta_0)\subseteq B^\epsilon$ and Proposition 7.3.3 in \cite{luenberger:1969} imply
\begin{equation}\label{eq:gmm4ver1}
\|\Upsilon_G(\theta_1) - \Upsilon_G(\theta_2) - \nabla \Upsilon_G(\theta_1)[\theta_1-\theta_2]\|_{\mathbf G}
\leq \{\sup_{\theta \in B^\epsilon} \max_{1\leq j \leq d_G} \|\nabla^2 \Upsilon_{G,j}(\theta)\|_{o,2}\} \frac{\|\theta_1 - \theta_2\|_2^2}{2}.
\end{equation}
Similarly, for any $\theta_1,\theta_2\in N_\epsilon(\theta_0)$, Proposition 7.3.2 in \cite{luenberger:1969} yields
\begin{align}\label{eq:gmm4ver2}
\|\nabla \Upsilon_G(\theta_1) - \nabla \Upsilon_G(\theta_2)\|_o & = \sup_{\|h\|_2 = 1} \max_{1\leq j \leq d_{G}} |(\nabla \Upsilon_{G,j}(\theta_1) - \nabla \Upsilon_{G,j}(\theta_2))[h]| \notag
\\ & \leq  \{\sup_{\theta \in B^\epsilon} \max_{1\leq j \leq d_G} \|\nabla^2 \Upsilon_{G,j}(\theta)\|_{o,2}\} \|\theta_1 -\theta_2\|_2.
\end{align}
Since $\|\nabla^2 \Upsilon_{G,j}(\theta)\|_{o,2}$ is uniformly bounded on $B^\epsilon$ by Assumption \ref{ass:gmmbootp1}(v), it follows from results \eqref{eq:gmm4ver1} and \eqref{eq:gmm4ver2} that Assumptions \ref{ass:locineq}(i)(ii) are satisfied with
$$K_g \equiv \sup_{\theta \in B^\epsilon} \max_{1\leq j \leq d_G} \|\nabla^2 \Upsilon_{G,j}(\theta)\|_{o,2}.$$
Assumption \ref{ass:gmmbootp1}(iii) additionally implies $\sup_{\theta \in B^\epsilon} \|\nabla \Upsilon_G(\theta)\|_{o,2} < \infty$, and hence verifies Assumption \ref{ass:locineq}(iii).
By identical arguments, but recalling $\mathbf F = \mathbf R^{d_F}$ and $\|\cdot\|_{\mathbf F} = \|\cdot\|_2$, it follows Assumptions \ref{ass:gmmbootp1}(iii)-(iv) imply Assumptions \ref{ass:loceq}(i)-(iii) hold with
\begin{equation}\label{eq:gmm4ver3}
K_f \equiv \sqrt {d_F} \sup_{\theta \in B^\epsilon} \max_{1\leq j \leq d_F} \|\nabla^2 \Upsilon_{F,j}(\theta)\|_{o,2} .
\end{equation}

To conclude, note that since Assumption \ref{ass:gmmbootp1}(vi) implies the range of $\nabla \Upsilon_F(\theta)$ equals $\mathbf R^{d_F}$ for all $\theta \in B^\epsilon$, it follows that $\nabla \Upsilon_F(\theta)$ admits a right inverse.
Moreover, if $\Upsilon_F$ is affine, then $K_f = 0$ and hence Assumption \ref{ass:loceq}(iv) holds.
On the other hand, if $\Upsilon_F$ is nonlinear, then note $\nabla \Upsilon_F(\theta)^{-} = \nabla \Upsilon_F(\theta)^\prime(\nabla \Upsilon_F(\theta) \nabla \Upsilon_F(\theta)^\prime)^{-1}$ and therefore $\|\nabla \Upsilon_F(\theta)^{-}\|_{o,2}$ is bounded for all $\theta \in B^\epsilon$ due to $\|\nabla \Upsilon_F(\theta)\|_{o,2}$ being bounded on $B^\epsilon$ by Assumption \ref{ass:gmmbootp1}(ii), and the smallest singular value of $\nabla \Upsilon_F(\theta)^\prime$ being bounded away from zero on $B^\epsilon$ by Assumption \ref{ass:gmmbootp2}(ii). It follows Assumption \ref{ass:loceq}(iv) holds as well.
Since Assumption \ref{ass:gmmbootp2} directly implies Assumption \ref{ass:ineqlindep}, the lemma follows. \qed

\begin{lemma}\label{lm:gmmnoell}
Let Assumptions \ref{ass:gmmbasic}, \ref{ass:gmmrho}, \ref{ass:gmmid}, and \ref{ass:gmmsigma} hold, and set $a_n = \sqrt{\log(n)}/n^{\frac{1}{10 + 5d_\theta}}$.
For any $\ell_n$ with $n^{-1/2} = o(\ell_n)$, it  follows uniformly in $P\in \mathbf P_0$ that
$$\hat U_n(R|+\infty)  = \inf_{h \in \hat V_n(\hat \theta_n,R|\ell_n)} \|\hat {\mathbb W}_n(\hat \theta_n) + \hat{\mathbb D}_n(\hat \theta_n)[h]\|_{\hat \Sigma_n,2} + o_P(a_n).$$
\end{lemma}

\noindent {\sc Proof:}
We establish the lemma by relying on Lemma \ref{lm:lnotbind}.
To this end note that in the proof of Theorem \ref{th:gmmapprox}, Assumptions \ref{ass:param}(i), \ref{ass:startreg}(i)(iii), and \ref{ass:weights} were verified and $\hat \theta_n$ was shown to be consistent for $\IDpoint$ uniformly in $P\in \mathbf P_0$ with $\|\cdot\|_{\mathbf B} = \|\cdot\|_{\mathbf E} = \|\cdot\|_2$, $\mathcal R_n = n^{-1/2}$, and $\nu_n \asymp 1$ for both $R$ as in \eqref{eq:appgmm0p5}.
Next, note Lemma \ref{lm:bootgmmcoupling} verifies Assumption \ref{ass:bootcoupling} holds with $a_n = \sqrt{\log(n)}/n^{\frac{1}{10 + 5d_\theta}}$ for $R = \Theta$ and hence also for $R$ as in \eqref{eq:appgmm0p5}.
Moreover, the mean value theorem and $\Theta$ being convex imply that
\begin{equation}\label{eq:gmmnoell1}
|\frac{\partial}{\partial \theta_k} \rho_\jmath(x,\theta_1) - \frac{\partial}{\partial \theta_k} \rho_\jmath(x,\theta_2)| \leq \max_{1\leq \jmath \leq \mathcal J} \sup_{\theta \in \Theta} \|\nabla_\theta^2 \rho_\jmath (x,\theta)\|_{o,2} \|\theta_1 - \theta_2\|_2
\end{equation}
for any $\theta_1,\theta_2 \in \Theta$, $1\leq \jmath \leq \mathcal J$, and $1\leq k \leq d_\theta$.
Thus, Assumption \ref{ass:gmmrho}(ii) implies there exists a $C_0 < \infty$ such that for all $P\in \mathbf P$ and $\theta_1,\theta_2\in \Theta$ it follows that
\begin{equation*}\label{eq:gmmnoell2}
\|E_P[\nabla_\theta \rho(X,\theta_1) - \nabla_\theta \rho(X,\theta_2)]\|_{o,2}  \leq C_0 \|\theta_1-\theta_2\|_2.
\end{equation*}
In particular, the function $\theta \mapsto E_P[\nabla_\theta \rho(X,\theta)]$ is uniformly continuous in $\theta$ and $P\in \mathbf P$, which implies by Assumption \ref{ass:gmmid}(ii) that there is an $\epsilon_0 > 0$ such that the smallest singular value of $E_P[\nabla_\theta \rho(X,\theta)]$ is bounded away from zero on $\{\theta \in \Theta :\|\theta - \IDpoint\|_2 \leq \epsilon_0 \text{ for some } P\in \mathbf P_0\}$ (where recall $\IDpoint$ implicitly depends on $P$ through \eqref{eq:appgmm0}). Since $\|\DerP(\theta)[h]\|_2 \equiv \|E_P[\nabla_\theta \rho(X,\theta)h]\|_2$, $\nu_n \asymp 1$, $p = 2$, and $\|\cdot\|_{\mathbf E} = \|\cdot\|_2$, the Lemma \ref{lm:lnotbind} requirement that $\|h\|_{\mathbf E} \leq \nu_n\|\DerP(\theta)[h]\|_2$ for all $\theta \in \mathcal A_n(P)$, $P\in \mathbf P_0$, and $h\in \sqrt n \{\mathbf B_n \cap R - \theta\}$ holds with $\mathcal A_n(P) = (\IDpoint)^{\epsilon_0}$ and $R = \Theta$ (and hence also for $R$ as in \eqref{eq:appgmm0p5}).
Moreover, by uniform consistency (in $P\in \mathbf P_0$) of $\hat \theta_n$ it follows that $\hat \theta_n \in \mathcal A_n(P)$ with probability tending to one uniformly in $P\in \mathbf P_0$.

To conclude, define $\mathcal F \equiv \{ \frac{\partial}{\partial \theta_k} \rho_\jmath(\cdot,\theta) : \text{ for some } \theta \in \Theta,~ 1 \leq \jmath \leq \mathcal J, ~ 1\leq k \leq d_\theta\}$ and let $F(x) \equiv  \max_{1\leq \jmath \leq \mathcal J} \sup_{\theta \in \Theta} \|\nabla_\theta^2 \rho_\jmath (x,\theta)\|_{o,2}$.
Then note that if $\epsilon$-balls around $\{\theta_i\}_{i=1}^{N_\epsilon}$ cover $\Theta$, then result \eqref{eq:gmmnoell1} implies that the brackets $[\frac{\partial}{\partial \theta_k} \rho_\jmath(\cdot,\theta_i) - \epsilon F, \frac{\partial}{\partial \theta_k} \rho_\jmath(\cdot,\theta_i) + \epsilon F]$ cover $\mathcal F$. 
Writing these brackets as $\{[f_{l,k},f_{u,k}]\}_{k=1}^{K_\epsilon}$ for conciseness, further note that $K_\epsilon = \mathcal J d_\theta N_\epsilon < \infty$ since $N_\epsilon < \infty$ due to $\Theta$ being compact, and $C_1 \equiv \sup_{P\in \mathbf P} \|F\|_{P,1} < \infty$ by Assumption \ref{ass:gmmrho}(ii).
Moreover, by definition of $[f_{l,k},f_{u,k}]$ it further follows that
\begin{equation}\label{eq:gmmnoell3}
E_P[f_{u,k}(X) - f_{l,k}(X)] = \|f_{u,k} - f_{l,k}\|_{P,1} \leq 2\epsilon C_1
\end{equation}
for all $P\in \mathbf P$.
Hence, employing the bound $f(x) - E_P[f(X)] \leq f_{u,k}(x)  - E_P[f_{l,k}(X)]$ for $[f_{l,k}, f_{u,k}]$ the bracket containing $f$, we obtain from result \eqref{eq:gmmnoell3} that
\begin{multline}\label{eq:gmmnoell4}
\sup_{f \in \mathcal F} \{\frac{1}{n}\sum_{i=1}^n f(X_i) - E_P[f(X)]\} \\ \leq \max_{1\leq k \leq K_\epsilon}  |\frac{1}{n} \sum_{i=1}^n f_{u,k}(X_i) - E_P[f_{u,k}(X)]| +  2\epsilon C_1 = 2\epsilon C_1 + o_P(1),
\end{multline}
where the equality holds uniformly in $P\in \mathbf P$ by Assumption \ref{ass:gmmrho}(ii), $K_\epsilon < \infty$, and Theorem 2.8.1 in \cite{vandervaart:wellner:1996}.
By identical arguments, we have
\begin{multline}\label{eq:gmmnoell5}
\inf_{f \in \mathcal F} \{\frac{1}{n}\sum_{i=1}^n f(X_i) - E_P[f(X)]\} \\ \geq - \max_{1\leq k \leq K_\epsilon}  |\frac{1}{n} \sum_{i=1}^n f_{l,k}(X_i) - E_P[f_{l,k}(X)]| - 2\epsilon C_1 = - 2\epsilon C_1 + o_P(1),
\end{multline}
uniformly in $P\in \mathbf P$.
We thus conclude from results \eqref{eq:gmmnoell4} and \eqref{eq:gmmnoell5} that $\mathcal F$ is Glivenko-Cantelli uniformly in $P\in \mathbf P$.
Since by Assumption \ref{ass:gmmbootp1}(i) there exists an $\epsilon > 0$ such that $\{\theta : \|\theta - \IDpoint\|_2 \leq \epsilon\} \subseteq \Theta$ for all $P\in \mathbf P_0$, we can conclude
\begin{multline}\label{eq:gmmnoell6}
\sup_{\theta : \|\theta - \IDpoint\|_2\leq \epsilon} \sup_{h \in \mathbf R^{d_\theta} : \|\frac{h}{\sqrt n}\|_{2} \geq \ell_n} \frac{\|\hat {\mathbb D}_n(\theta)[h] - \DerP(\theta)[h]\|_2}{\|h\|_2} \\
\leq \sup_{\theta \in \Theta}  \|\frac{1}{n}\sum_{i=1}^n \nabla_\theta \rho(X_i,\theta) - E_P[\nabla_\theta \rho(X,\theta)]\|_{o,2} = o_P(1)
\end{multline}
uniformly in $P\in \mathbf P$, and where the equality follows from $\mathcal F$ being Glivenko-Cantelli uniformly in $P\in \mathbf P$.
Since $\nu_n \asymp 1$, result \eqref{eq:gmmnoell6} verifies condition \eqref{lm:lnotbinddisp1} in Lemma \ref{lm:lnotbind}.
This concludes verifying the requirements of Lemma \ref{lm:lnotbind} and hence the present Lemma follows for any $\ell_n$ satisfying $\mathcal S_n(\mathbf B,\mathbf E)\mathcal R_n = o(\ell_n)$, which in this application is equivalent to $n^{-1/2} = o(\ell_n)$ due to $\|\cdot\|_{\mathbf B} = \|\cdot\|_{\mathbf E} = \|\cdot\|_2$ and $\mathcal R_n \asymp n^{-1/2}$. \qed

%% file: Appendix/ExamplesProofs/ExNPIVProofs.tex

\subsection{Proofs for Section \ref{sec:exconsumer}}\label{sec:exconsumerproofs}

\noindent {\sc Proof of Theorem \ref{th:slutskyapprox}:}
We establish the theorem by simply verifying the conditions of Theorem \ref{th:localdrift}(ii) for both $R$ as corresponding to \eqref{eq:appslutsky3} and \eqref{eq:appslutsky4} (to couple $I_n(R)$) and to $R = \Theta$ (to couple $I_n(\Theta)$).
To this end, note Assumption \ref{ass:param}(i) is imposed in Assumption \ref{ass:slutskyproc}(i), Assumption \ref{ass:startreg}(i) holds with $B_n \asymp \sqrt{k_n}$ by Assumption \ref{ass:slutskymoments}(i), and Assumption \ref{ass:startreg}(ii) is satisfied by Assumption \ref{ass:slutskymoments}(ii).
Further note that for $R = \Theta$, the class $\mathcal F_n$ has bounded envelope $F_n$ by Assumption \ref{ass:slutskyproc}(iii) and $\|g\|_\infty \leq C_0$ for any $(g,\gamma)\in \Theta$.
Hence, by Lemma \ref{lm:slutskycov} it follows that Assumption \ref{ass:startreg}(iii) holds with $J_n \asymp \sqrt{j_n\log(1+j_n)}$ when $R = \Theta$, and therefore also for $R$ as corresponding to \eqref{eq:appslutsky3} and \eqref{eq:appslutsky4}.
Next, we note Lemma \ref{lm:slutskycoup} and $j_n^2 k_n^3 \log^3(n)/n = o(1)$ by Assumption \ref{ass:slutskymoments}(iv) imply Assumption \ref{ass:coupling}(i) for both specifications of $R$ under consideration and any $a_n$ satisfying $a_n = O((\log(n))^{-1})$.
To verify Assumption \ref{ass:coupling}(ii), we observe that for any $(g_1,\gamma_1) \in \Theta_n$ and $(g_2,\gamma_2)\in \Theta_n$, Assumption \ref{ass:slutskyproc}(iii) implies
\begin{multline}\label{lm:slutskycoup1}
E_P[((Q - g_1(S,Y) - W^\prime \gamma_1) - (Q - g_2(S,Y) - W^\prime \gamma_2))^2] \\ \lesssim \sup_{P\in \mathbf P} \|g_1 -g_2\|_{P,2}^2 + \|\gamma_1 - \gamma_2\|_2^2.
\end{multline}
Hence, since $\|\cdot\|_{\mathbf E} \equiv \sup_{P\in \mathbf P} \|\cdot\|_{P,2} + \|\cdot\|_2$ it follows Assumption \ref{ass:coupling}(ii) holds with $\kappa_\rho = 1$ and some $K_\rho < \infty$.
Lemma \ref{lm:slutsky3ver} additionally verifies that Assumption \ref{ass:keycons} holds with $\|\cdot\|_{\mathbf E} \equiv \sup_{P\in \mathbf P} \|\cdot\|_{P,2} + \|\cdot\|_2$, $\mathcal V_n(P) = \Theta_n \cap R$, and $\nu_n^{-1} \asymp \text{s}_n$ for both $R = \Theta$ and $R$ as corresponding to \eqref{eq:appslutsky3} and \eqref{eq:appslutsky4}.
Further note that in this application
$$\nabla m_{P}(\theta)[h] = -E_P[g_h(S,Y) + W^\prime \gamma_h|Z]$$
for any $\theta \in \mathbf B$ and $(g_h,\gamma_h) = h \in \mathbf B$.
By direct calculation it then follows Assumptions \ref{ass:driftlin}(i)(ii) hold with $K_m = 0$ for $\|\cdot\|_{\mathbf L} = \|\cdot\|_{\mathbf E}$, and Assumption \ref{ass:driftlin}(iii) is satisfied for some $M < \infty$ by result \eqref{lm:slutskycoup1} and Jensen's inequality.
To verify Assumption \ref{ass:locrates}(i), note that since as argued $\nu_n \asymp 1/\text{s}_n$, $p = 2$, $J_n \asymp \sqrt{j_n\log(1+j_n)}$, and $B_n \asymp \sqrt{k_n}$, it follows $\mathcal R_n \asymp k_n\sqrt{j_n}\log(1+k_n)/\text{s}_n\sqrt n$ due to $j_n \leq k_n$ by Assumption \ref{ass:slutskymoments}(iii).
Therefore, since $\kappa_\rho = 1$, Lemma \ref{lm:slutskycov} implies Assumption \ref{ass:locrates}(i) demands $k_n\sqrt{j_n \log(1+k_n)}\mathcal R_n(1+\sqrt{\log(1\vee (\sqrt{j_n}/\mathcal R_n))}) = o(a_n)$, which is satisfied with $a_n = 1/\sqrt{\log(n)}$ by Assumption \ref{ass:slutskymoments}(iv).
In turn, Assumption \ref{ass:locrates}(ii) holds with $a_n = (\log(n))^{-1/2}$ by Assumptions \ref{ass:slutskysieve}(iv) and \ref{ass:slutskysigma}(ii).
Finally, we note that Assumption \ref{ass:slutskysigma} implies Assumption \ref{ass:weights} due $B_n \lesssim \sqrt{k_n}$, $p = 2$, $J_n \asymp \sqrt{j_n\log(1+j_n)}$, and $a_n = (\log(n))^{-1/2}$.
To conclude, note that since $K_m = 0$, the condition $K_m \mathcal R_n^2 \mathcal S_n(\mathbf L,\mathbf E) = o(a_nn^{-1/2})$ is automatically satisfied, the requirement $\mathcal R_n = o(\ell_n)$ is equivalent to $k_n\sqrt{j_n}\log(1+k_n)/\text{s}_n\sqrt{n} = o(\ell_n)$, and the condition $k_n^{1/p}\sqrt{\log(1+k_n)}B_n\times \sup_{P\in \mathbf P} J_{[\hspace{0.02 in}]}(\ell_n^{\kappa_\rho},\mathcal F_n,\|\cdot\|_{P,2}) =o(a_n)$ is implied by $k_n\sqrt{j_n\log(1 + k_n)}\ell_n\sqrt{\log(\sqrt{j_n}/\ell_n)} = o((\log(n))^{-1/2})$.
Thus, all the conditions of Theorem \ref{th:localdrift}(ii) hold for both $R = \Theta$ and $R$ corresponding to \eqref{eq:appslutsky3} and \eqref{eq:appslutsky4}, and thus the claim of the theorem follows. \qed

\noindent {\sc Proof of Theorem \ref{th:slutskyboot}:} We first define the variable $\hat E_n(R|\ell_n)$ to be given by
\begin{equation*}
\hat E_n(R|\ell_n) \equiv \inf_{h \in \hat V_n(\hat \theta_n,R|\ell_n)} \|\Bemp(\hat \theta_n) + \hat{\mathbb D}_n[h]\|_{\hat \Sigma_n,2}
\end{equation*}
and note that for any sequence $\ell_n$ satisfying the conditions of the theorem, Lemma \ref{lm:slutskynoell} implies $\hat U_n(R|+\infty) = \hat E_n(R|\ell_n) + o_P(a_n)$ uniformly in $P\in \mathbf P_0$.
To establish the theorem, it therefore suffices to show that uniformly in $P\in \mathbf P_0$
\begin{align*}
\hat E_n(R|\ell_n) & \geq \UpS (R|\tilde \ell_n) + o_P(a_n) \\
\hat E_n(R|\ell_n) - \hat U_n(\Theta|+\infty) & \geq \UpS (R|\tilde \ell_n) - \UpS (\Theta|\tilde \ell_n^{\text{u}}) + o_P(a_n)
\end{align*}
with $\ell_n \asymp \tilde \ell_n$ and $ \tilde \ell_n^{\text{\rm u}}$ satisfying the requirements of the theorem.
To this end we rely on Theorem \ref{th:coupsmooth} (for $\hat E_n(R|\ell_n)$) and Lemma \ref{lm:forincJboot}(ii).
We note that in the proof of Theorem \ref{th:slutskyapprox} we established that Assumptions \ref{ass:slutskyproc}, \ref{ass:slutskysieve}, \ref{ass:slutskymoments}, and \ref{ass:slutskysigma} imply Assumptions \ref{ass:param}, \ref{ass:startreg}, \ref{ass:coupling}, \ref{ass:keycons},  \ref{ass:driftlin},  \ref{ass:locrates} and \ref{ass:weights} hold with $\mathcal R_n \asymp k_n\sqrt{j_n}\log(1+k_n)/\text{s}_n\sqrt n$, $B_n \asymp \sqrt{k_n}$, $\nu_n \asymp 1/\text{s}_n$,  $\|\theta\|_{\mathbf B} = \|g\|_{1,\infty} \vee \|\gamma\|_2$ and $\|\theta\|_{\mathbf L} = \|\theta\|_{\mathbf E} = \sup_{P\in \mathbf P} \|g\|_{P,2} + \|\gamma\|_2$ for $\theta = (g,\gamma)$, $\kappa_\rho = 1$, and $a_n = (\log(n))^{-1/2}$ for $R = \Theta$ and $R$ as corresponding to \eqref{eq:appslutsky3} and \eqref{eq:appslutsky4}.
We thus only verify Assumptions \ref{ass:locineq}, \ref{ass:loceq}, \ref{ass:ineqlindep}, \ref{ass:bootcoupling}, \ref{ass:extra}, and \ref{ass:bootrates} for $R = \Theta$ and $R$ as corresponding to \eqref{eq:appslutsky3} and \eqref{eq:appslutsky4}.

Next note that Lemma \ref{lm:slutsky4ver} implies Assumptions \ref{ass:locineq}, \ref{ass:loceq}, and \ref{ass:ineqlindep} are satisfied with $K_g =2$ and $K_f =0$, while Lemma \ref{lm:slutskzybootcoup} and Assumption \ref{ass:slutskyboot}(iv) imply Assumption \ref{ass:bootcoupling} holds with $R = \Theta$ (and hence for $R$ corresponding to \eqref{eq:appslutsky3} and \eqref{eq:appslutsky4}) with $a_n = (\log(n))^{-1/2}$.
Further note that since $\sup_{P\in \mathbf P} \|g\|_{P,2} + \|\gamma\|_2 \leq 2 (\|g\|_{1,\infty} \vee \|\gamma\|_2)$ for any $(g,\gamma) \in C^1_B(\Omega)\times \mathbf R^{d_w}$, it follows that Assumption \ref{ass:extra}(i) holds with $M = 2$.
To verify Assumption \ref{ass:extra}(ii) note that by the definitions of $\|\cdot\|_{\mathbf B}$ and $\|\cdot\|_{\mathbf E}$ in this application and the eigenvalues of $E_P[p^{j_n}(S,Y)p^{j_n}(S,Y)^\prime]$ being bounded away from zero uniformly in $P\in \mathbf P$ by Assumption \ref{ass:slutskysieve}(iii) we obtain
\begin{multline}\label{eq:slutskyboot0}
\mathcal S_n(\mathbf B,\mathbf E) = \sup_{(\beta,\gamma)} \frac{\|p^{j_n\prime} \beta\|_{1,\infty} \vee \|\gamma\|_2}{\sup_{P\in \mathbf P} \|p^{j_n\prime}\beta\|_{P,2} + \|\gamma\|_2} \\ \leq 1 \vee \sup_{\beta} \frac{\|p^{j_n\prime}\beta\|_{1,\infty}}{\sup_{P\in \mathbf P} \|p^{j_n\prime}\beta\|_{P,2}}  \lesssim 1 \vee \sup_{\beta} \frac{\|p^{j_n\prime}\beta\|_{1,\infty}}{\|\beta\|_2} \lesssim j_n^{3/2}
\end{multline}
where the final equality follows from Assumptions \ref{ass:slutskysieve}(i)(ii).
In particular, note that result \eqref{eq:slutskyboot0}, $\mathcal R_n \asymp k_n\sqrt{j_n}\log(1+k_n)/\text{s}_n\sqrt n$, and Assumption \ref{ass:slutskymoments}(iv) imply that $\mathcal R_n \mathcal S_n(\mathbf B,\mathbf E) = o(1)$.
Thus, since setting $\hat \theta_n$ and $\hat \theta_n^{\text{u}}$ to be the minimizers of $Q_n$ (respectively over $\Theta_n \cap R$ and $\Theta_n$) corresponds to setting $\tau_n = 0$, Assumption \ref{ass:slutskyboot}(ii) and $\mathcal R_n \mathcal S_n(\mathbf B,\mathbf E) = o(1)$ implies Assumption \ref{ass:extra}(ii) holds.
We also note Assumption \ref{ass:extra}(iii) is immediate since $\mathcal V_n(P) = \Theta_n \cap R$ by Lemma \ref{lm:slutsky3ver}.
To conclude, note Assumption \ref{ass:bootrates}(i) holds since we showed $\mathcal R_n \mathcal S_n(\mathbf B,\mathbf E) = o(1)$.
Moreover, since $B_n \asymp \sqrt{k_n}$ and $K_f = K_m = 0$, Lemma \ref{lm:slutskycov} implies Assumption \ref{ass:bootrates}(ii) holds for any $\ell_n, \ell_n^{\rm u}$ satisfying $k_n\sqrt{j_n\log(1+k_n)}(\ell_n\vee \ell_n^{\rm u})(1+\sqrt{\log(\sqrt{j_n}/(\ell_n\vee \ell_n^{\rm u})}) = o(a_n)$.
Similarly, we obtain that Assumption \ref{ass:bootrates}(iii) is satisfied provided $\ell_n = o(r_n)$ (imposed in the theorem) and $k_nj_n^2\log(1+k_n)/\text{s}_n\sqrt n = o(r_n)$ (implied by Assumption \ref{ass:slutskyboot}(iii)), while the requirement $\mathcal R_n^{\rm u} = o(\ell_n^{\rm u})$ is implied by $k_nj_n^2\log(1+k_n)/{\rm s}_n\sqrt{n} = o(\ell_n^{\rm u})$.
Hence, the conditions of Theorem \ref{th:coupsmooth} and Lemma \ref{lm:forincJboot}(ii) hold, and the theorem follows. \qed

\begin{lemma}\label{lm:slutsky3ver}
If Assumptions \ref{ass:slutskyproc}(ii), \ref{ass:slutskysieve}(iii), \ref{ass:slutskymoments}(iii), and \ref{ass:slutskysigma}(ii) hold, then Assumption  \ref{ass:keycons} is satisfied with both $R = \Theta$ and $R$ corresponding to \eqref{eq:appslutsky3} and \eqref{eq:appslutsky4}, $\|\cdot\|_{\mathbf E} = \sup_{P\in \mathbf P} \|\cdot\|_{P,2} + \|\cdot\|_2$, $\mathcal V_n(P) = \Theta_n \cap R$, and $\nu_n^{-1} \asymp \text{\rm s}_n$.
\end{lemma}

\noindent {\sc Proof:} By Assumption \ref{ass:slutskyproc}(ii) there is a unique $\IDpoint \equiv (g_0,\gamma_0)\in \Theta \cap R$ for which  \eqref{eq:appslutsky0} holds, and we let $\Pi_n \IDpoint = (g_n,\gamma_0) \in \Theta_n \cap R$ for $g_n = p^{j_n\prime}\beta_n$.
To verify Assumption \ref{ass:keycons}(i) is satisfied we set $\|\theta\|_{\mathbf E} \equiv \sup_{P\in \mathbf P} \|g\|_{P,2} + \|\gamma\|_2$ for any $(g,\gamma) = \theta \in \mathbf B$.
Since the eigenvalues of $E_P[p^{j_n}(S,Y)p^{j_n}(S,Y)^\prime]$ are bounded uniformly in $j_n$ and $P\in \mathbf P$ by Assumption \ref{ass:slutskysieve}(iii) we can conclude for any $\theta = (p^{j_n\prime}\beta,\gamma)$ that
\begin{align}\label{lm:slutsky3ver1}
\text{s}_n \|\theta - \Pi_n\IDpoint\|_{\mathbf E} & \lesssim \text{s}_n \{\|\beta - \beta_n\|_2 + \|\gamma - \gamma_0\|_2\} \notag \\
& \lesssim \|E_P[(p^{j_n}(S,Y)^\prime(\beta - \beta_n) + W^\prime(\gamma - \gamma_0))q^{k_n}(Z)]\|_{\PSigma,2} \notag \\ & = \|E_P[(\rho(X,\theta) - \rho(X,\Pi_n\theta_0))q^{k_n}(Z)]\|_{\PSigma,2} 
\end{align}
where the second inequality holds by Assumptions \ref{ass:slutskymoments}(iii) and \ref{ass:slutskysigma}(ii), while the final equality holds by definition of $\rho(X,\theta)$ (see \eqref{eq:appslutsky1p5}). 
Thus, we conclude from \eqref{lm:slutsky3ver1} that Assumption \ref{ass:keycons}(i) holds with $\nu_n^{-1} \asymp \text{s}_n$ and $\mathcal V_n(P) = \Theta_n \cap R$.
Finally, note Assumption \ref{ass:keycons}(ii) is immediate since $\mathcal V_n(P) = \Theta_n \cap R$. \qed

\begin{lemma}\label{lm:slutskycov}
Define the class $\mathcal F_n \equiv \{f: f(v) = (q - g(s,y) - w^\prime \gamma) \text{ for some } (g,\gamma) \in \Theta_n\}$ and suppose that Assumptions \ref{ass:slutskyproc}(iii) and \ref{ass:slutskysieve}(i)(iii) hold.
Then, it follows that $\sup_{P\in \mathbf P} N_{[\hspace{0.02 in}]}(\epsilon,\mathcal F_n,\|\cdot\|_{P,2})\lesssim 1\vee (\sqrt{j_n}K/\epsilon)^{j_n+d_w}$  for some $K<\infty$, and in addition $\sup_{P\in \mathbf P}J_{[\hspace{0.02 in}]}(\epsilon,\mathcal F_n,\|\cdot\|_{P,2}) \lesssim \epsilon \sqrt{j_n} (1+\sqrt {\log(1 \vee (\sqrt{j_n}/\epsilon))})$.
\end{lemma}

\noindent {\sc Proof:} Define the classes $\mathcal F_{1n} \equiv \{f : f(v) = q - w^\prime \gamma \text{ with } \|\gamma\|_2 \leq C_0\}$ and $\mathcal F_{2n} \equiv \{ p^{j_n \prime} \beta : \|p^{j_n \prime} \beta\|_{1,\infty}\leq C_0\}$, and then note that by definition of $\mathcal F_n$ we have
\begin{equation}\label{lm:slutskycov1}
\sup_{P\in \mathbf P} N_{[\hspace{0.02 in}]}(\epsilon,\mathcal F_n,\|\cdot\|_{P,2}) \leq \sup_{P\in \mathbf P} N_{[\hspace{0.02 in}]}(\frac{\epsilon}{2},\mathcal F_{1n},\|\cdot\|_{P,2}) \times \sup_{P\in \mathbf P} N_{[\hspace{0.02 in}]}(\frac{\epsilon}{2},\mathcal F_{2n},\|\cdot\|_{P,2}).
\end{equation}
Next observe that since the support of $W$ is bounded uniformly in $P\in \mathbf P$ by Assumption \ref{ass:slutskyproc}(iii), the Cauchy-Schwarz inequality, the covering numbers of $\{\gamma \in \mathbf R^{d_w} : \|\gamma\|_2 \leq C_0\}$ being bounded (up to a multiplicative constant) by $1\vee \epsilon^{-d_w}$, and Theorem 2.7.11 in \cite{vandervaart:wellner:1996} allow us to conclude that
\begin{equation}\label{lm:slutskycov2}
\sup_{P\in \mathbf P} N_{[\hspace{0.02 in}]}(\frac{\epsilon}{2},\mathcal F_{1n},\|\cdot\|_{P,2}) \lesssim 1\vee \epsilon^{-d_w}.
\end{equation}
Similarly, for any $p^{j_n\prime}\beta_1, p^{j_n\prime}\beta_2 \in \mathcal F_{2n}$, the Cauchy-Schwarz inequality implies that
\begin{equation*}
|p^{j_n}(s,y)^\prime \beta_1 - p^{j_n}(s,y)^\prime \beta_2| \leq \sup_{(s,y)}\|p^{j_n}(s,y)\|_2 \|\beta_1-\beta_2\|_2 \lesssim \sqrt{j_n}  \|\beta_1-\beta_2\|_2,
\end{equation*}
where in the final inequality we employed Assumption \ref{ass:slutskysieve}(i).
Hence, Theorem 2.7.11 in \cite{vandervaart:wellner:1996}, $\|\beta\|_2 \asymp \|p^{j_n\prime}\beta\|_{P,2}$ uniformly in $P\in \mathbf P$ by Assumption \ref{ass:slutskysieve}(iii), and $\|p^{j_n\prime}\beta\|_{P,2} \leq \|p^{j_n\prime}\beta\|_{\infty}\leq C_0$ for any $p^{j_n\prime}\beta \in \Theta_n$ imply that
\begin{equation}\label{lm:slutskycov4}
\sup_{P\in \mathbf P} N_{[\hspace{0.02 in}]}(\frac{\epsilon}{2},\mathcal F_{2n},\|\cdot\|_{P,2}) \lesssim  1\vee (\frac{K\sqrt{j_n}}{\epsilon })^{j_n}
\end{equation}
for some $K <\infty$. Thus, the first claim of the lemma follows from results \eqref{lm:slutskycov1}, \eqref{lm:slutskycov2}, and \eqref{lm:slutskycov4}.
For the second claim of the lemma, we employ the first claim of the lemma and the change of variables $v = u/\epsilon$ to obtain the bound
\begin{align*}
\sup_{P\in \mathbf P} &J_{[\hspace{0.02 in}]}(\epsilon,\mathcal F_n,\|\cdot\|_{P,2}) \lesssim \epsilon + \int_0^\epsilon (\log(1\vee (\frac{K\sqrt{j_n}}{u})^{j_n+d_w}))^{1/2}du  \\
& = \epsilon(1 + \sqrt{j_n + d_w} \int_0^1(\log(1 \vee (\frac{K\sqrt{j_n}}{v\epsilon})))^{1/2}dv) \lesssim \sqrt{j_n}\epsilon(1+ \sqrt{\log(1\vee (\sqrt{j_n}/\epsilon))}),
\end{align*}
where we used that $(1\vee ab) \leq (1\vee a)(1\vee b)$ whenever $a$ and $b$ are positive. \qed

\begin{lemma}\label{lm:slutskycoup}
Let Assumptions \ref{ass:slutskyproc}(i)(iii), \ref{ass:slutskysieve}(i)(iii), and \ref{ass:slutskymoments}(i) hold. If $a_n \downarrow 0$ and $k_n^3 j_n^2 \log^2(n)/n = o(a_n)$, then Assumption \ref{ass:coupling}(i) holds with $R = \Theta$.
\end{lemma}

\noindent {\sc Proof:} We establish the claim of the lemma by relying on Lemma \ref{lm:finitecoup}.
To this end, let $\tilde j_n = (1+d_w)+j_n$, set $r^{j_n}(x) \equiv (q,w^\prime, p^{j_n}(x)^\prime)^\prime$, and observe any $f\in \mathcal F_n$ can be written as $f = r^{j_n\prime}\delta$ for some $\delta \in \mathbf R^{\tilde j_n}$. Moreover, by Assumption \ref{ass:slutskysieve}(iii) and definition of $\Theta_n$, it follows that there exists an $M<\infty$ such that $\mathcal F_n \subseteq \{r^{\tilde j_n\prime} \delta : \|\delta\|_2 \leq M\}$, while Assumptions \ref{ass:slutskyproc}(iii), \ref{ass:slutskysieve}(i), and \ref{ass:slutskymoments}(i) imply $\sup_{x} \|r^{j_n}(x)\|_2 \lesssim \sqrt{j_n}$ and $\sup_{z} \|q^{k_n}(z)\|_2 \leq \sqrt{k_n} \max_{1\leq k \leq k_n} \|q_k\|_\infty \lesssim k_n$. The claim of the lemma therefore follows from applying Lemma \ref{lm:finitecoup} with $b_{1n} \asymp \sqrt{j_n}$, $b_{2n}\asymp k_n$, and $C_n = M$. \qed

\begin{lemma}\label{lm:slutskzybootcoup}
Suppose Assumptions \ref{ass:slutskyproc}(i)(iii), \ref{ass:slutskysieve}(i)(iii), \ref{ass:slutskymoments}(i)(ii) hold and let $\mathcal C_n \equiv \{\beta \in \mathbf R^{j_n} : \|p^{j_n\prime}\beta\|_{1,\infty} \leq C_0\}$ and $\mathcal E_n \equiv \int_0^\infty \sqrt{\log(N(\epsilon,\mathcal C_n,\|\cdot\|_2))}d\epsilon$.
If $j_n^2 k_n^2 \log(1+k_nj_n) = o(n)$, then it follows that Assumption \ref{ass:bootcoupling} holds with $R = \Theta$ for any sequence $a_n$ satisfying $k_n^{1/p}(\sqrt{\log(k_n)} + \mathcal E_n)j_n^{3/4}k_n^{1/2}\log^{1/4}(1+j_nk_n)/n^{1/4} = o(a_n)$.
\end{lemma}

\noindent {\sc Proof:}
Recall that in this application $X \equiv (Q,S,Y,W^\prime)^\prime$ and, when $R = \Theta$, we have $\mathcal F_n \equiv \{f : f(x) = (q- g(s,y) - w^\prime \gamma) \text{ for some } (g,\gamma) \in \Theta_n\}$.
Also define $\tilde {\mathcal F}_n \equiv \{f q_{k} : f \in \mathcal F_n \text{ and } 1 \leq k \leq k_n\}$,  for $\{\omega_i\}_{i=1}^n$ the weights used in building $\Bemp$ set
$$\hat {\mathbb G}_n(f) \equiv \frac{1}{\sqrt n} \sum_{i=1}^n \omega_i \{f(V_i) - \frac{1}{n}\sum_{j=1}^n f(V_j)\}$$
for any $f\in \tilde{\mathcal F_n}$, and let $\IsoW$ denote an isonormal process on $\mathcal F_n$ independent of $\{V_i\}_{i=1}^n$. Setting $\WPT(\theta) \equiv (\IsoW(\rho(\cdot,\theta)q_1),\ldots, \IsoW(\rho(\cdot,\theta)q_{k_n}))^\prime$, then note that
\begin{equation}\label{lm:slutskybootcoup1}
\sup_{\theta \in \Theta} \|\Bemp(\theta) - \WPT(\theta)\|_p \leq k_n^{1/p} \sup_{f\in \tilde {\mathcal F}_n} |\hat {\mathbb G}_n(f)- \IsoW(f) |.
\end{equation}
We will establish the lemma by relying on \eqref{lm:slutskybootcoup1} and applying Theorem \ref{th:mainbootcoup} to couple $\hat {\mathbb W}_n$ and $\IsoW$ on $\tilde {\mathcal F}_n$.
To this end, define $d_n = k_n(j_n + d_w+1)$ and let 
\begin{equation}\label{lm:slutskybootcoup1p5}
f_n^{d_n}(V) \equiv g^{d_n}(V) - E_P[g^{d_n}(V)] \hspace{0.3 in} g^{d_n}(V) \equiv q^{k_n}(Z)\otimes (p^{j_n}(S,Y)^\prime, Q, W^\prime)^\prime .
\end{equation}
Next, we set $D_1 \equiv (Q, W^\prime, p^{j_n}(S,Y)^\prime)^\prime$ and $D_2 = q^{k_n}(Z)$, and for $\overline{\text{eig}}\{D_1D_1^\prime\}$ the largest eigenvalue of the matrix $D_1D_1^\prime$, then note that we must have
\begin{equation}\label{lm:slutskybootcoup2}
\sup_{P\in \mathbf P}\|\overline{\text{eig}}\{D_1D_1^\prime\}\|_{P,\infty} \leq \sup_{P\in \mathbf P} \|\text{trace}\{D_1D_1^\prime\}\|_{P,\infty} \lesssim j_n,
\end{equation}
where the final inequality follows from Assumptions \ref{ass:slutskyproc}(iii) and \ref{ass:slutskysieve}(i).
Hence, since $\overline{\text{eig}}\{E_P[q^{k_n}(Z)q^{k_n}(Z)^\prime]\}$ is bounded uniformly in $P\in \mathbf P$ by Assumption \ref{ass:slutskymoments}(ii), result \eqref{lm:slutskybootcoup2} and Lemma \ref{lm:kronereig} imply $\overline{\text{eig}}\{E_P[g^{d_n}(V)g^{d_n}(V)^\prime] \lesssim j_n$.
It thus follows from $\overline{\text{eig}}\{E_P[g^{d_n}(V)]E_P[g^{d_n}(V)^\prime]\} \leq \overline{\text{eig}}\{E_P[g^{d_n}(V)g^{d_n}(V)^\prime]\}$ and definition \eqref{lm:slutskybootcoup1p5} that Assumption \ref{ass:4seriescoup}(i) is satisfied with $C_n \asymp j_n$.
Similarly, note that Assumptions \ref{ass:slutskyproc}(iii), \ref{ass:slutskysieve}(i), and  \ref{ass:slutskymoments}(i) imply Assumption \ref{ass:4seriescoup}(ii) holds with $K_n \asymp \sqrt{j_nk_n}$.
Moreover, Assumption \ref{ass:4seriesreg}(i) is immediate with $G_{n,P}$ equal to the zero function and $J_{1n} = 0$.
Finally, note that any function $ f \in \tilde{\mathcal F_n}$ has the structure
\begin{equation}\label{lm:slutskybootcoup3}
 f(v) = q_{k}(z)(q - p^{j_n}(s,y)^\prime \beta - w^\prime \gamma) \text{ for some } (p^{j_n\prime}\beta,\gamma) \in \Theta_n.
\end{equation}
Therefore, for $\mathcal B_n$ as defined in Assumption \ref{ass:4seriesreg}(ii), $\mathcal C_n \equiv \{\beta \in \mathbf R^{j_n} : \|p^{j_n\prime}\beta\|_{1,\infty} \leq C_0\}$, and $\mathcal G_n \equiv \{\gamma \in \mathbf R^{d_w} : \|\gamma\|_2 \leq C_0\}$, we can conclude that
\begin{multline}\label{lm:slutskybootcoup4}
N(\epsilon,\mathcal B_n,\|\cdot\|_2) \leq k_n \times N(\frac{\epsilon}{2},\mathcal G_n,\|\cdot\|_2)  \times N(\frac{\epsilon}{2},\mathcal C_n,\|\cdot\|_2)  \\
\lesssim k_n \times ((\frac{1}{\epsilon})^{d_w} \vee 1) \times N(\frac{\epsilon}{2},\mathcal C_n,\|\cdot\|_2) ,
\end{multline}
where in the second inequality we employed that $N(\epsilon,\mathcal G_n,\|\cdot\|_2) \lesssim (1/\epsilon)^{d_w} \vee 1$.
Furthermore, note that Assumption \ref{ass:slutskysieve}(iii) implies that $\|\beta\|_2 \asymp \|p^{j_n\prime}\beta\|_{P,2}$ uniformly in $j_n$ and $P\in \mathbf P$, and hence since $\|p^{j_n\prime}\beta \|_{P,2} \leq \|p^{j_n\prime}\beta \|_{1,\infty}$, the definition of $\Theta_n$ and \eqref{lm:slutskybootcoup3} implies that the radius of $\mathcal B_n$ under $\|\cdot\|_2$ is uniformly bounded in $n$.
Thus, the bound in \eqref{lm:slutskybootcoup4} yields that for some $M < \infty$ we must have
\begin{align*}
\int_0^\infty & \sqrt{\log(N(\epsilon,\mathcal B_n,\|\cdot\|_2))}d\epsilon \notag\\
& \lesssim \int_0^M\sqrt{\log(k_n)} d\epsilon +  \int_0^{1} \sqrt{\log(1/\epsilon)}d\epsilon + \int_0^M \sqrt{\log(N(\epsilon/2,\mathcal C_n,\|\cdot\|_2))}d\epsilon \notag\\
& \lesssim \sqrt{\log(k_n)} + \int_0^\infty \sqrt{\log(N(u,\mathcal C_n,\|\cdot\|_2))}du,
\end{align*}
where the final inequality follows from $N(\epsilon,\mathcal C_n,\|\cdot\|_2)$ being (weakly) larger than one for all $\epsilon$ and the change of variables $u = \epsilon/2$.
Hence, Assumption \ref{ass:4seriesreg}(ii) holds with $J_{2n} = \sqrt{\log(k_n)} + \mathcal E_n$, and as a result Theorem \ref{th:mainbootcoup} implies uniformly in $P\in \mathbf P$
\begin{equation}\label{lm:slutskybootcoup6}
\sup_{f\in \tilde{\mathcal F}_n} |\hat {\mathbb G}_n(f) - \IsoW(f)| = O_P((\sqrt{\log(k_n)} + \mathcal E_n)\{\frac{j_n^3 k_n^2 \log(1+j_nk_n)}{n}\}^{1/4}).
\end{equation}
The claim of the lemma therefore follows from \eqref{lm:slutskybootcoup1} and \eqref{lm:slutskybootcoup6}. \qed

\begin{lemma}\label{lm:slutsky4ver}
If $\mathbf B = C^1_B(\Omega)\times \mathbf R^{d_w}$ and $\Upsilon_G$, $\Upsilon_F$, and $\Theta$ are as defined in \eqref{eq:appslutsky3}, \eqref{eq:appslutsky4}, and \eqref{eq:slutskyTheta}, then it follows that Assumptions \ref{ass:locineq}, \ref{ass:loceq}, and \ref{ass:ineqlindep} are satisfied with $K_g = 2$, $K_f = 0$, and for any $\theta = (g,\gamma)$ and $h = (g_h,\gamma_h)$, $\nabla \Upsilon_G(\theta)[h]$ equals
\begin{equation*}\label{lm:slutsky4verdisp}
\nabla \Upsilon_G(\theta)[h](s,y) = \frac{\partial}{\partial s} g_h(s,y) + g(s,y) \frac{\partial}{\partial y} g_h(s,y) + g_h(s,y)\frac{\partial}{\partial y} g(s,y).
\end{equation*}
\end{lemma}

\noindent {\sc Proof:} Recall that in this application $\mathbf G = C_B^0(\Omega)$ and $\|\theta\|_{\mathbf B} = \max\{\|g\|_{1,\infty},\|\gamma\|_2\}$.
Hence, for any $\theta_1 = (g_1,\gamma_1) \in \mathbf B$ and $\theta_2 = (g_2,\gamma_2) \in \mathbf B$ we obtain that
\begin{multline*}
\|\Upsilon_G(\theta_1) - \Upsilon_G(\theta_2) - \nabla \Upsilon_G(\theta_1)[\theta_1-\theta_2]\|_{\mathbf G} \\
\leq \sup_{(s,y)\in \Omega} |g_1(s,y) - g_2(s,y)| \times \sup_{(s,y)\in \Omega} |\frac{\partial}{\partial y} (g_1(s,y) - g_2(s,y))| \leq \|g_1 -g_2\|_{1,\infty}^2,
\end{multline*}
which verifies Assumption \ref{ass:locineq}(i) holds with $K_g = 2$.
Similarly, we additionally conclude
\begin{align}\label{lm:slutsky4ver2}
\|\nabla & \Upsilon_G(\theta_1) - \nabla \Upsilon_G(\theta_2)\|_o \notag \\
 & = \sup_{g_h:\|g_h\|_{1,\infty} \leq 1} \sup_{(s,y)\in \Omega} |(g_1(s,y)-g_2(s,y))\frac{\partial}{\partial y} g_h(s,y) + g_h(s,y)\frac{\partial}{\partial y}(g_1(s,y) - g_2(s,y))| \notag \\
 & \leq 2 \|g_1 - g_2\|_{1,\infty},
\end{align}
which verifies Assumption \ref{ass:locineq}(ii) holds with $K_g = 2$ as well.
Moreover, note that since any $\theta = (g,\gamma) \in \Theta$ satisfies $\|g\|_{1,\infty} \leq C_0$, it follows that $\|\tilde g\|_{1,\infty} \leq C_0 + \epsilon$ for any $\tilde g\in \Theta^\epsilon$.
Thus, by identical arguments to those in \eqref{lm:slutsky4ver2} we obtain
\begin{equation*}
\|\nabla \Upsilon_G(\theta)\|_o \leq 2\|g\|_{1,\infty} \leq 2(C_0 + \epsilon),
\end{equation*}
which thus verifies Assumption \ref{ass:locineq}(iii) holds with $M = 2(C_0+\epsilon)$.

Next note $\Upsilon_F: \mathbf B \to \mathbf F$ is affine and continuous, and hence $\nabla \Upsilon_F(\theta)[h] = \Upsilon_F(h) - c_0$ for all $\theta,h \in \mathbf B$.
Therefore, Assumptions \ref{ass:loceq}(i)(ii) hold with $K_f = 0$, while
\begin{equation*}
\sup_{g_h : \|g_h\|_{1,\infty}\leq 1} |g_h(s_0,y_0)| \leq 1
\end{equation*}
implies Assumption \ref{ass:loceq}(iii) is satisfied with $M = 1$. Since $\Upsilon_F$ being affine and $K_f = 0$ further imply that Assumptions \ref{ass:loceq}(iv) and \ref{ass:ineqlindep} hold, the lemma follows. \qed

\begin{lemma}\label{lm:slutskynoell}
Let $a_n = (\log(n))^{-1/2}$ and Assumptions \ref{ass:slutskyproc}, \ref{ass:slutskysieve}, \ref{ass:slutskymoments}, \ref{ass:slutskysigma} hold.
If $\ell_n$ satisfies $k_nj_n^2\log(1+k_n)/\text{\rm s}_n\sqrt n = o(\ell_n)$, then uniformly in $P\in \mathbf P_0$:
$$\hat U_n(R|+\infty)  = \inf_{h \in \hat V_n(\hat \theta_n,R|\ell_n)} \|\Bemp(\hat \theta_n) + \hat{\mathbb D}_n[h]\|_{\hat \Sigma_n,2} + o_P(a_n).$$
\end{lemma}

\noindent {\sc Proof:} We establish the lemma by applying Lemma \ref{lm:lnotbind}.
To this end, recall that in the proof of Theorem \ref{th:slutskyapprox}, Assumptions \ref{ass:startreg}(i)(iii) and \ref{ass:weights} were verified to hold with $B_n \asymp \sqrt{k_n}$ and $J_n \asymp \sqrt{j_n\log(1+j_n)}$.
Since the eigenvalues of $E_P[p^{j_n}(S,Y)p^{j_n}(S,Y)^\prime]$ are bounded uniformly in $P\in \mathbf P$ by Assumption \ref{ass:slutskysieve}(iii) and $\|\theta\|_{\mathbf E} = \sup_{P\in \mathbf P} \|g\|_{P,2} + \|\gamma\|_2$ for any $\theta = (g,\gamma)$, it also follows that for any $h = (p^{j_n\prime}\beta_h,\gamma_h)$ we have
\begin{multline*}
\|h\|_{\mathbf E} = \sup_{P\in \mathbf P} \|p^{j_n\prime}\beta_h\|_{P,2} + \|\gamma\|_2 \lesssim \|\beta_h\|_2 + \|\gamma_h\|_2 \\
\lesssim \frac{1}{\text{s}_n} \|E_P[q^{k_n}(Z)(p^{j_n}(S,Y)^\prime \beta_h + W^\prime \gamma_h)]\|_2 = \frac{1}{\text{s}_n} \|\DerP[h]\|_2,
\end{multline*}
where the second inequality holds by Assumption \ref{ass:slutskymoments}(iii) and the final equality follows from the definition of $\DerP[h]$.
Hence, since $\nu_n \asymp 1/\text{s}_n$ by Lemma \ref{lm:slutsky3ver} and $p=2$, we conclude the Lemma \ref{lm:lnotbind} requirement that $\|h\|_{\mathbf E} \leq \nu_n \|\DerP(\theta)[h]\|_p$ for all $\theta \in \mathcal A_n(P)$ holds with $\mathcal A_n(P) = \Theta_n \cap R$.
Next, define the $k_n\times (j_n+d_w)$ matrix 
\begin{equation}\label{lm:slutskynoell3}
\mathbb M_{i,n} \equiv \frac{1}{n}\{q^{k_n}(Z_i)(p^{j_n}(S_i,Y_i)^\prime ~ W_i^\prime) - E_P[q^{k_n}(Z)(p^{j_n}(S,Y)^\prime ~ W^\prime)]\},
\end{equation}
which satisfies $E_P[\mathbb M_{i,n}] = 0$.
Since $\|(p^{j_n\prime}\beta,\gamma)\|_{\mathbf E} \asymp \|\beta\|_2 + \|\gamma\|_2$ by Assumption \ref{ass:slutskysieve}(iii), we then conclude from \eqref{lm:slutskynoell3} that for some $C_1 < \infty$ we must have
\begin{multline}\label{lm:slutskynoell7}
\sup_{P\in \mathbf P} P(\sup_{h \in \mathbf B_n} \frac{\|\hat {\mathbb D}_n[h] - \DerP[h]\|_2}{\|h\|_{\mathbf E}} > \text{s}_n)
\leq \sup_{P\in \mathbf P} P( \|\frac{1}{n}\sum_{i=1}^n \mathbb M_{i,n}\|_{o,2} >  C_1\text{s}_n)\\
\leq (j_n+d_w+k_n) \exp\{-\frac{n\text{s}_n^2C_2}{(k_n^2\vee j_n) + \text{s}_n k_n\sqrt{j_n}}\} = o(1),
\end{multline}
where the final inequality follows by applying Lemma \ref{lm:finiteder} with $b_{1n} = \sqrt{j_n}$ (by Assumptions \ref{ass:slutskyproc}(iii) and \ref{ass:slutskysieve}(i)) and $b_{2n} = k_n$ (by Assumption \ref{ass:slutskymoments}(i)), while the final equality results from $\log(k_n)k_n^2/\text{s}_n^2 n = o(1)$ by Assumption \ref{ass:slutskymoments}(iv) and $k_n \geq j_n$ by Assumption \ref{ass:slutskysieve}(iii).
Hence, $\nu_n \asymp 1/\text{s}_n$ and \eqref{lm:slutskynoell7} imply condition \eqref{lm:lnotbinddisp1} in Lemma \ref{lm:lnotbind} holds.
Finally, we note that by Assumption \ref{ass:slutskyboot}(iv), we may apply Lemma \ref{lm:slutskzybootcoup} with $p = 2$ to conclude that Assumption \ref{ass:bootcoupling} holds with $R = \Theta$ (and hence for $R$ as corresponding to \eqref{eq:appslutsky3} and \eqref{eq:appslutsky4})  with $a_n = (\log(n))^{-1/2}$.
This concludes verifying the requirements of Lemma \ref{lm:lnotbind} and therefore the present Lemma follows for any $\ell_n$ satisfying $\mathcal S_n(\mathbf B,\mathbf E)\mathcal R_n = o(\ell_n)$, which in this application is equivalent to $k_nj_n^2\log(1+k_n)/\text{s}_n\sqrt n = o(\ell_n)$ due to $\mathcal S_n(\mathbf B,\mathbf E) \lesssim j_n^{3/2}$ and $\mathcal R_n \asymp k_n\sqrt{j_n}\log(1+k_n)/\text{s}_n\sqrt n$. \qed

\begin{lemma}\label{lm:kronereig}
Let $D_1\in \mathbf R^{d_1}$, $D_2\in \mathbf R^{d_2}$ be distributed according to $Q$, and for any matrix $A$ let $\overline{\text{\rm eig}}\{A\}$ denote its largest eigenvalue.
Then it follows that
$$\overline{\text{\rm eig}}\{E_Q[(D_1\otimes D_2)(D_1 \otimes D_2)^\prime]\} \leq \|\overline{\text{\rm eig}}\{D_1D_1^\prime\}\|_{Q,\infty} \times \overline{\text{\rm eig}}\{E_Q[D_2D_2^\prime]\}. $$
\end{lemma}

\noindent {\sc Proof:} Let $\mathcal A \equiv \{\{a_i\}_{i=1}^{d_1} : a_i \in \mathbf R^{d_2} \text{ and } \sum_{i=1}^{d_1} \|a_i\|^2_2 \leq 1\}$, set $(D_1^{(1)},\ldots, D_1^{(d_1)}) = D_1 \in \mathbf R^{d_1}$, and then note that by direct calculation we obtain that
\begin{align*}
\overline{\text{eig}}\{E_Q & [(D_1\otimes D_2)(D_1 \otimes D_2)^\prime]\} \\
& = \sup_{\{a_i\}_{i=1}^{d_1}\in \mathcal A} (a_1^\prime,\ldots, a_{d_1}^\prime)E_Q[(D_1\otimes D_2)(D_1 \otimes D_2)^\prime](a_1^\prime,\ldots, a_{d_1}^\prime)^\prime \\
& = \sup_{\{a_i\}_{i=1}^{d_1}\in \mathcal A}  E_Q[(\sum_{i=1}^{d_1} (a_i^\prime D_2) D_1^{(i)})^2] \\ & \leq \|\overline{\text{\rm eig}}\{D_1D_1^\prime\}\|_{Q,\infty} \sup_{\{a_i\}_{i=1}^{d_1}\in \mathcal A} \sum_{i=1}^{d_1} E_Q[(a_i^\prime D_2)^2].
\end{align*}
However, since $\sum_{i=1}^{d_1} \|a_i\|_2^2 \leq 1$ for all $\{a_i\}_{i=1}^{d_1} \in \mathcal A$, we additionally have the inequality
$$\sup_{\{a_i\}_{i=1}^{d_1}\in \mathcal A} \sum_{i=1}^{d_1} E_Q[(a_i^\prime D_2)^2] \leq \sup_{\{a_i\}_{i=1}^{d_1}\in \mathcal A}  \sum_{i=1}^{d_1} \overline{\text{eig}}\{E_Q[D_2D_2^\prime]\} \|a_i\|^2_2 =  \overline{\text{eig}}\{E_Q[D_2D_2^\prime]\},$$
and therefore the claim of the lemma follows. \qed

\begin{lemma}\label{lm:auxsplines}
Let $\lambda$ be the Lebesgue measure, $\{B^{(1)}_{b}\}_{b=1}^{j_{1n}}$ and $\{B^{(2)}_{b}\}_{b=1}^{j_{2n}}$ be B-splines on $[0,1]$ of order $r \geq 3$ with no interior knot multiplicity, mesh ratio bounded in $n$, and $\|\cdot\|_{\lambda,2}$ normalized to have norm one. If $\{p_{j}\}_{j=1}^{j_n}$ is the tensor product of $\{B^{(1)}_{b}\}_{b=1}^{j_{1n}}$ and $\{B^{(2)}_{b}\}_{b=1}^{j_{2n}}$ and $\mathcal C_n \equiv \{ \beta \in \mathbf R^{j_n} : \|p^{j_n\prime}\beta\|_{1,\infty} \leq C_0\}$, then it follows that
$$\int_0^\infty \sqrt{\log(N(\epsilon,\mathcal C_n,\|\cdot\|_2))}d\epsilon \lesssim \sqrt{j_{1n}\wedge j_{2n}}\log(j_n+1).$$
\end{lemma}

\noindent {\sc Proof:} We rely heavily on Chapter 5 in \cite{devore1993constructive}, and note that $B_{j}$ corresponds to $N_j/\|N_j\|_{\lambda,2}$ in their notation. 
Throughout, for two sequences $a_n$ and $b_n$ we employ $a_n \asymp b_n$ to mean that there exist constants $\underline c$ and $\bar c$ such that $\underline c a_n \leq b_n \leq \bar c a_n$ for all $n$.
In what follows it will also prove convenient to  index the elements of $\beta \in \mathbf R^{j_n}$ by $\beta_{b_1,b_2}$ with $1\leq b_1 \leq j_{1n}$ and $1\leq b_2 \leq j_{2n}$.
Then note that the mesh ratios corresponding to $\{B^{(1)}_{b}\}_{b=1}^{j_{1n}}$ and $\{B_{b}^{(2)}\}_{b=1}^{j_{2n}}$ being bounded uniformly in $j_{1n}$, $j_{2n}$ and two applications of Theorem 5.4.2 in \cite{devore1993constructive} imply that
\begin{multline}\label{lm:auxsplines1}
\|p^{j_n \prime}\beta\|_\infty =  \sup_{u_1 \in [0,1]} \sup_{u_2\in [0,1]} |\sum_{b_2 = 1}^{j_{2n}} B_{b_2}^{(2)}(u_2)  \sum_{b_1 = 1}^{j_{1n}} \beta_{b_1, b_2} B_{b_1}^{(1)}(u_1) | \\
\asymp \sup_{u_1\in [0,1]} \max_{1\leq b_2 \leq j_{2n}} \sqrt{j_{2n}} | \sum_{b_1 = 1}^{j_{1n}} \beta_{b_1,b_2} B_{b_1}^{(1)}(u_1)| \asymp \sqrt{j_{1n} j_{2n}} \|\beta\|_\infty
\end{multline}
uniformly in $\beta \in \mathbf R^{j_n}$.
By similar arguments we also obtain uniformly in $\beta \in \mathbf R^{j_n}$ that
\begin{align}\label{lm:auxsplines2}
\sup_{u_1 \in (0,1)} \sup_{u_2\in [0,1]} |\sum_{b_2 = 1}^{j_{2n}} & B_{b_2}^{(2)}(u_2) \sum_{b_1 = 1}^{j_{1n}} \frac{\partial}{\partial u_1} \{ \beta_{b_1, b_2} B_{b_1}^{(1)}(u_1) \}| \notag \\
& \asymp \sup_{u_1\in (0,1)} \max_{1\leq b_2 \leq j_{2n}} \sqrt{j_{2n}} | \sum_{b_1 = 1}^{j_{1n}} \frac{\partial}{\partial u_1} \{\beta_{b_1,b_2}B_{b_1}^{(1)}(u_1)\}| \notag \\
& \asymp \max_{1\leq b_2 \leq j_{2n}} \max_{2\leq b_1 \leq j_{1n}} \sqrt{j_{2n}} j_{1n}^{3/2} |\beta_{b_1,b_2} - \beta_{b_1-1,b_2}|,
\end{align}
where the second result follows by employing equation (3.11) and Theorem 5.4.2 in Chapter 5 of \cite{devore1993constructive} and the mesh ratio of $\{B_{b}^{(1)}\}_{b=1}^{j_{1n}}$ being bounded.
Since by identical arguments we can also derive the symmetric (to \eqref{lm:auxsplines2}) relationship
\begin{multline}\label{lm:auxsplines3}
\sup_{u_1 \in [0,1]} \sup_{u_2\in (0,1)} |\sum_{b_1 = 1}^{j_{1n}} B_{b_1}^{(1)}(u_1) \sum_{b_2 = 1}^{j_{2n}} \frac{\partial}{\partial u_2} \{ \beta_{b_1, b_2} B_{b_2}^{(2)}(u_2) \}| \\
 \asymp \max_{1\leq b_1 \leq j_{1n}} \max_{2\leq b_2 \leq j_{2n}} \sqrt{j_{1n}} j_{2n}^{3/2} |\beta_{b_1,b_2} - \beta_{b_1,b_2-1}|,
 \end{multline}
it follows from results \eqref{lm:auxsplines1}, \eqref{lm:auxsplines2}, and \eqref{lm:auxsplines3} that there is an $M_0 < \infty$ such that
\begin{align}\label{lm:auxsplines4}
\max_{1\leq b_1 \leq j_{1n}} \max_{1\leq b_2 \leq j_{2n}} |\beta_{b_1,b_2}| & \leq M_0/\sqrt{j_{n}} \notag\\
\max_{1\leq b_2 \leq j_{2n}} \max_{2\leq b_1 \leq j_{1n}} |\beta_{b_1,b_2} - \beta_{b_1-1,b_2}| & \leq M_0/(j_{1n}\sqrt{j_{n}}) \notag \\
\max_{1\leq b_1 \leq j_{1n}} \max_{2\leq b_2 \leq j_{2n}} |\beta_{b_1,b_2} - \beta_{b_1,b_2-1}| & \leq M_0/(j_{2n}\sqrt{j_{n}})
\end{align}
for all $\beta \in \mathcal C_n$.
Hence, in order to establish the claim of the lemma, it suffices to bound the covering numbers for the set defined by \eqref{lm:auxsplines4}.

We proceed by combining two bounds, one for ``small" $\epsilon$ and one for ``large" $\epsilon$.
First, assume without loss of generality  $j_{1n} \geq j_{2n}$, let $c_n \equiv \lceil \log(j_{1n}+1) \rceil$, and define the sets
\begin{align}\label{lm:auxsplines5}
\frac{\epsilon}{3\sqrt {j_n}} k_1  & \leq \beta_{b_1,b_2} \leq \frac{\epsilon }{3\sqrt {j_n}} (k_{1} + 1) \text { for all } b_1 = mc_n + 1 \text{ with } 0 \leq m \leq \lceil j_{1n}/c_n \rceil -1 \notag \\
\frac{\epsilon}{3c_n\sqrt{j_n}} k_{2} & \leq \beta_{b_1,b_2} - \beta_{b_1-1,b_2} \leq \frac{\epsilon}{3c_n\sqrt {j_n}} (k_{2}+1) \text{ otherwise }
\end{align}
where $k_1,k_2$ are non-zero integers -- i.e.\ the sets (in $\mathbf R^{j_n}$) defined in \eqref{lm:auxsplines5} consist of ``chains" along the $b_1$ dimension that reset every $c_n$ integers.
To compute the diameter of the sets in \eqref{lm:auxsplines5}, then note that since all ``chains" have the same structure
\begin{align}\label{lm:auxsplines6}
\sup \|\beta - \tilde \beta\|_2^2 & \text{ s.t. } \beta,\tilde \beta \text{ satisfying \eqref{lm:auxsplines5}} \notag \\
& \leq  \sup j_{2n}\lceil \frac{j_{1n}}{c_n}\rceil\sum_{b_1 = 1}^{c_n} (\beta_{b_1,j_{2n}} - \tilde \beta_{b_1,j_{2n}})^2  \text{ s.t. } \beta,\tilde \beta \text{ satisfying \eqref{lm:auxsplines5}} \notag \\
& \leq j_{2n}\lceil \frac{j_{1n}}{c_n}\rceil  \sum_{b_1 = 1}^{c_n} \frac{\epsilon^2}{9j_n}\{1 + \frac{2(b_1-1)}{c_n}\}^2,
\end{align}
where the final inequality follows from \eqref{lm:auxsplines5}.
Since $\lceil j_{1n}/c_n\rceil c_n \leq (j_{1n} + c_n) \leq 2j_{1n}$ due to $\lceil j_{1n}/c_n\rceil \leq 1+j_{1n}/c_n$ and $c_n \leq j_{1n}$, it follows from $j_n = j_{1n}j_{2n}$ that every set in \eqref{lm:auxsplines5} is contained in a ball of radius $\epsilon$.
Moreover, by \eqref{lm:auxsplines4} the total number of sets with the structure in \eqref{lm:auxsplines5} needed to cover the set $\mathcal C_n$ is bounded by
\begin{equation}\label{lm:auxsplines7}
(\lceil \frac{6M_0}{\epsilon} \rceil )^{j_{2n}\lceil \frac{j_{1n}}{c_n}\rceil}(\lceil \frac{6M_0c_n}{\epsilon j_{1n}} \rceil )^{j_{2n}\lceil \frac{j_{1n}}{c_n}\rceil c_n} .
\end{equation}
Next, we employ again the bound $\lceil j_{1n}/c_n\rceil c_n \leq 2j_{1n}$ and $\lceil a \rceil \leq 2a$ whenever $a \geq 1$, to obtain from \eqref{lm:auxsplines6} and \eqref{lm:auxsplines7} that whenever $\epsilon \leq 6M_0c_n/j_{1n}$ we have
\begin{multline}\label{lm:auxsplines8}
N(\epsilon,\mathcal C_n,\|\cdot\|_2) \leq (\frac{12M_0}{\epsilon})^{\frac{2j_{n}}{c_n}}(\frac{12M_0c_n}{\epsilon j_{1n}})^{2j_{n}} \\
= (\frac{12M_0c_n}{\epsilon j_{1n}}(\frac{j_{1n}}{c_n})^{\frac{1}{c_n+1}})^{\frac{2j_n(c_n+1)}{c_n}} \leq  (\frac{M_1 \log(1+j_{1n})}{\epsilon j_{1n}})^{4j_n},
\end{multline}
where the final equality holds for some $M_1 < \infty$ due to $(c_n+1)/c_n \leq 2$ and $(j_{1n}/c_n)^{\frac{1}{c_n+1}} \leq j_{1n}^{\frac{1}{\log(1+j_{1n})}} =O(1)$ because $c_n = \lceil \log(1+j_{1n})\rceil$.

The bound in \eqref{lm:auxsplines8} is valid only for $\epsilon \leq 6M_0c_n/j_{1n}$.
To obtain a bound for $\epsilon \geq 6M_0 c_n/j_{1n}$, let $\{\mathbb Z_{b_1,b_2}\}_{b_1,b_2}$ be independent standard normal random variables.
By Sudakov's inquality (see, e.g., Proposition A.2.5 in \cite{vandervaart:wellner:1996}), it then follows that for some $M_2 < \infty$ independent of $n$ we have that
\begin{equation}\label{lm:auxsplines9}
\sqrt{\log(N(\epsilon,\mathcal C_n,\|\cdot\|_2))} \leq \frac{M_2}{\epsilon} E[\sup_{\beta \in \mathcal C_n} \sum_{b_1 = 1}^{j_{1n}} \sum_{b_2 = 1}^{j_{2n}} \beta_{b_1,b_2} \mathbb Z_{b_1,b_2}].
\end{equation}
Next, for notational convenience define $\Delta_{b_1} \beta_{b_1,b_2} = (\beta_{b_1,b_2} - \beta_{b_1-1,b_2})$ and $\Delta_{b_2} \beta_{b_1,b_2} = (\beta_{b_1,b_2} - \beta_{b_1,b_2-1})$, and then note that by \eqref{lm:auxsplines4} it follows that
\begin{align*}
&\sup_{\beta \in \mathcal C_n} \sum_{b_1 =1}^{j_{1n}} \sum_{b_2 = 1}^{j_{2n}} \beta_{b_1,b_2} \mathbb Z_{b_1,b_2} \notag \\
& = \sup_{\beta \in \mathcal C_n} \sum_{b_2 =1}^{j_{2n}}\sum_{b_1=2}^{j_{1n}} \Delta_{b_1} \beta_{b_1,b_2} \sum_{\tilde b_1 = b_1}^{j_{1n}}  \mathbb Z_{\tilde b_1,b_2} + \sum_{b_2 = 2}^{j_{2n}} \Delta_{b_2}\beta_{1,b_2}\sum_{b_1 = 1}^{j_{1n}} \sum_{\tilde b_2 = b_2}^{j_{2n}} \mathbb Z_{b_1,\tilde b_2} + \beta_{1,1}\sum_{b_1 = 1}^{j_{1n}}\sum_{b_2=1}^{j_{2n}} \mathbb Z_{b_1,b_2} \notag \\
& \leq \sum_{b_2 =1}^{j_{2n}}\sum_{b_1=2}^{j_{1n}} \frac{M_0}{j_{1n}\sqrt{j_n}}|\sum_{\tilde b_1 = b_1}^{j_{1n}}  \mathbb Z_{\tilde b_1,b_2}| + \sum_{b_2 = 2}^{j_{2n}} \frac{M_0}{j_{2n}\sqrt {j_n}}|\sum_{b_1 = 1}^{j_{1n}} \sum_{\tilde b_2 = b_2}^{j_{2n}} \mathbb Z_{b_1,\tilde b_2}| + \frac{M_0}{\sqrt{j_n}}|\sum_{b_1 = 1}^{j_{1n}}\sum_{b_2=1}^{j_{2n}} \mathbb Z_{b_1,b_2}|.
\end{align*}
Hence, employing that if $\mathbb W \sim N(0,\sigma^2)$ then $E[|\mathbb W|] \lesssim \sigma$, we can conclude that
\begin{multline*}
E[\sup_{\beta \in \mathcal C_n} \sum_{b_1 = 1}^{j_{1n}} \sum_{b_2 = 1}^{j_{2n}} \beta_{b_1,b_2} \mathbb Z_{b_1,b_2}]
\lesssim \sum_{b_2 =1}^{j_{2n}}\sum_{b_1=2}^{j_{1n}} \frac{\sqrt{j_{1n}-b_1}}{j_{1n}\sqrt{j_n}} + \sum_{b_2 = 2}^{j_{2n}} \frac{\sqrt{j_{1n}(j_{2n} - b_2)}}{j_{2n}\sqrt {j_n}} + 1 \\ \leq \frac{j_n\sqrt{j_{1n}}}{j_{1n}\sqrt{j_n}} + \frac{j_{2n}\sqrt{j_{1n}j_{2n}}}{j_{2n}\sqrt{j_n}} + 1 \leq 3\sqrt{j_{2n}}
\end{multline*}
where in the final inequality we employed that $j_n = j_{1n}j_{2n}$. Hence, by \eqref{lm:auxsplines9} we have
\begin{equation}\label{lm:auxsplines11}
\sqrt{\log(N(\epsilon,\mathcal C_n,\|\cdot\|_2))} \lesssim \frac{\sqrt{j_{2n}}}{\epsilon}.
\end{equation}

To conclude the proof, we combine the bounds in \eqref{lm:auxsplines8} and \eqref{lm:auxsplines11}.
In particular, setting $\delta_n \equiv 6M_0\lceil \log(j_{1n} + 1)\rceil/j_{1n}$ and observing that $\|\beta\|_2 \asymp \|p^{j_n\prime}\beta\|_{\lambda,2}\leq C_0$ for all $\beta \in \mathcal C_n$ allows us to conclude that for some $M_2 <\infty$ we must have
\begin{multline*}
\int_0^\infty \sqrt{\log(N(\epsilon,\mathcal C_n,\|\cdot\|_2))} d\epsilon \lesssim \int_{\delta_n}^{M_2} \frac{\sqrt{j_{2n}}}{\epsilon} + \sqrt{j_n} \int_0^{\delta_n}(\log(\frac{M_1\log(j_{1n})}{\epsilon j_{1n}}))^{1/2}d\epsilon \\
\lesssim \sqrt{j_{2n}}\log(1+j_{1n}) + \frac{\sqrt{j_n} \log(1+j_{1n})}{j_{1n}} \int_0^1 (\log(\frac{1}{u}))^{1/2}du \lesssim \sqrt{j_{2n}}\log(1+j_{1n})
\end{multline*}
where the second inequality follows from the change of variables $u = \epsilon/\delta_n$ and the final inequality employed that $j_n = j_{1n} j_{2n}$ and $j_{2n} \leq j_{1n}$.
Substituting $j_{2n} = j_{1n} \wedge j_{2n}$ and employing $j_{1n} \leq j_n$ establishes the lemma. \qed 

%% file: Appendix/ExamplesProofs/ExQuantProofs.tex

\subsection{Proofs for Section \ref{sec:exquant}} \label{sec:exquantproofs}

\noindent {\sc Proof of Theorem \ref{th:quantapprox}:} We establish the result by applying Theorem \ref{th:localdrift} to both $R = \Theta$ and $R$ corresponding to \eqref{ex:quant4} and \eqref{ex:quant5}.
To this end, note that Assumption \ref{ass:param}(i) is directly imposed in Assumption \ref{ass:quantproc}(i), Assumption \ref{ass:startreg}(i) holds with $B_n = O(1)$ by Assumption \ref{ass:quantmoments}(i), Assumption \ref{ass:startreg}(ii) is directly imposed by Assumption \ref{ass:quantmoments}(iii), and Assumption \ref{ass:startreg}(iii) holds with $F_n = 1$ and $J_n = O(1)$ by Lemma \ref{lm:quantcov}.
Next, we apply Lemma \ref{lm:quantcoup} with $\pi_{0n} = O(1)$ and $\pi_{1n} = O(k_n)$ (which is possible by Assumptions \ref{ass:quantmoments}(i)(ii)) to obtain that Assumption \ref{ass:coupling}(i) holds for both $R = \Theta$ and $R$ corresponding to \eqref{ex:quant4} and \eqref{ex:quant5} for any $a_n$ satisfying $k_n^{1/p}\sqrt{j_n}\log(n)(n^{1/6}\vee k_n)/n^{1/3} = o(a_n)$ and in particular it holds for $a_n \asymp (\log(n))^{-1/2}$ by Assumption \ref{ass:quantmoments}(v).
We also note Assumptions \ref{ass:coupling}(ii), \ref{ass:keycons}, and \ref{ass:driftlin} hold with $\|\cdot\|_{\mathbf L} = \|\cdot\|_\infty$, $\|\cdot\|_{\mathbf E} = \sup_{P\in \mathbf P} \|\cdot\|_{P,2}$, and $\kappa_\rho =1/2$ by Lemmas \ref{lm:quant3ver} and \ref{lm:quantrhover}.
To verify Assumption \ref{ass:locrates}, note $J_n = O(1)$, $B_n = O(1)$, and $\nu_n \asymp \sqrt{k_n}/\text{\rm s}_n k_n^{1/p}$ by Lemma \ref{lm:quant3ver} imply in this application we have $\mathcal R_n \asymp \sqrt{k_n\log(1+k_n)}/\text{s}_n\sqrt n$.
Thus, Lemma \ref{lm:quantcov} and Assumption \ref{ass:quantmoments}(v) verify Assumption \ref{ass:locrates}(i), while Assumption \ref{ass:quantsieve}(iii) implies Assumption \ref{ass:locrates}(ii), and Assumption \ref{ass:quantsigma} implies Assumption \ref{ass:weights}.
Finally, since $\|\cdot\|_{\mathbf L} = \|\cdot\|_\infty$ by Lemma \ref{lm:quantrhover}, Assumptions \ref{ass:quantsieve}(i)(ii) yield
\begin{equation}\label{th:quantapprox1}
\sup_{\beta} \frac{\|p^{j_n\prime} \beta \|_\infty}{\sup_{P\in \mathbf P} \|p^{j_n\prime}\beta\|_{P,2}} \lesssim \frac{\sqrt{j_n}\|\beta\|_2}{\|\beta\|_2} = \sqrt{j_n}.
\end{equation}
Therefore, the condition $K_m\ell_n^2 \times \mathcal S_n(\mathbf L,\mathbf E) = o(a_n n^{-1/2})$ is equivalent to $\ell_n^2\sqrt{n j_n\log(n)} = o(1)$.
Moreover, by Lemma \ref{lm:quantcov}, Assumption \ref{ass:quantsieve}(ii), $\kappa_\rho = 1/2$, and $B_n = O(1)$ the condition $k_n^{1/p}\sqrt{\log(1+k_n)}B_n\times \sup_{P\in \mathbf P} J_{[\hspace{0.02 in}]}(\ell_n^{\kappa_\rho},\mathcal F_n,\|\cdot\|_{P,2}) = o(a_n)$ is implied by the restriction $k_n^{1/p} \sqrt{j_n\ell_n\log(1+k_n)\log(1/\ell_n)} = o((\log(n))^{-1/2})$.
Thus, the first claim of the theorem follows from Theorem \ref{th:localdrift}(i) applied to $I_n(R)$.

The second claim of the Theorem follows from applying Theorem \ref{th:localdrift}(ii) to $I_n(\Theta)$.
To this end, note that the only remaining condition to verify is that $\mathcal R_n^2 \times \mathcal S_n(\mathbf L,\mathbf E) = o(a_n n^{-1/2})$.
Using that, as already argued, $\mathcal R_n \asymp \sqrt{k_n\log(1+k_n)}/\text{s}_n\sqrt n$ and result \eqref{th:quantapprox1} we note that a sufficient condition for this final requirement is that $k_n\log(1+k_n)\sqrt{j_n\log(n)}/\text{s}_n^2 \sqrt n = o(1)$ and therefore the theorem follows. \qed

\noindent {\sc Proof of Theorem \ref{th:quantboot}:} We proceed by relying on Theorem \ref{th:coupsmooth} and Lemma \ref{lm:forincJboot}(ii).
To this end, we note that, for both $R = \Theta$ and $R$ corresponding to \eqref{ex:quant4} and \eqref{ex:quant5}, the proof of Theorem \ref{th:quantapprox} established that Assumptions \ref{ass:param}(i), \ref{ass:startreg},  \ref{ass:keycons}, \ref{ass:coupling}, \ref{ass:driftlin}, \ref{ass:locrates}, and \ref{ass:weights} are satisfied with $B_n = O(1)$, $J_n = O(1)$, $\|\cdot\|_{\mathbf E} = \sup_{P\in \mathbf P} \|\cdot\|_{P,2}$, $\|\cdot\|_{\mathbf L} = \|\cdot\|_\infty$, $\nu_n \asymp \sqrt{k_n}/\text{s}_n k_n^{1/p}$, $\mathcal R_n \asymp \sqrt{k_n\log(1+k_n)}/\text{s}_n\sqrt n$, $\kappa_\rho = 1/2$, and $\mathcal S_n(\mathbf L,\mathbf E) \lesssim \sqrt{j_n}$.

Next, note that Assumptions \ref{ass:locineq}, \ref{ass:loceq}, and \ref{ass:ineqlindep} are satisfied by Lemma \ref{lm:quantile4ver} with $K_g = 0$ and $K_f > 0$.
To verify Assumption \ref{ass:bootcoupling} we apply Lemma \ref{lm:quantbootcoup} with $\pi_{0n} = O(1)$ and $\pi_{1n} \lesssim k_n$, which is possible by Assumptions \ref{ass:quantmoments}(i)(ii).
In particular, Lemma \ref{lm:quantbootcoup} evaluated at $d_n \asymp (nk_n)^{3/13}$ and Assumption \ref{ass:4bootquant}(iii) yield
\begin{equation}\label{th:quantboot1}
\sup_{\theta \in \Theta_n} \|\Bemp(\theta) - \WPT(\theta)\|_p = o_P((\log(n))^{-1/2})
\end{equation}
uniformly in $P\in \mathbf P$, which implies Assumption \ref{ass:bootcoupling} is satisfied for both $R = \Theta$ and $R$ corresponding to \eqref{ex:quant4} and \eqref{ex:quant5}.
Next note that Assumption \ref{ass:extra}(i) is immediate since $\|\cdot\|_{\mathbf E} = \sup_{P\in \mathbf P_0} \|\cdot\|_{P,2}$ and $\|\cdot\|_{\mathbf B} = \|\cdot\|_{2,\infty}$, while Assumption \ref{ass:extra}(iii) follows from Lemma \ref{lm:quantcons} and $\mathcal V_n(P) \equiv \{\theta \in \Theta_n : \|\theta - \Pi_n \theta_0\|_{1,\infty} \leq \epsilon\}$ for some $\epsilon > 0$ by Lemma \ref{lm:quant3ver}.
In order to verify Assumption \ref{ass:extra}(ii) and the rate conditions of Assumption \ref{ass:bootrates}, note that the eigenvalues of $E_P[p^{j_n}(D)p^{j_n}(D)^\prime]$ being bounded away from zero uniformly in $P\in \mathbf P$ by Assumption \ref{ass:quantsieve}(ii), Assumptions \ref{ass:quantsieve}(i) and \ref{ass:4bootquant}(iv) together with $\|\cdot\|_{\mathbf B} = \|\cdot\|_{2,\infty}$ and $\|\cdot\|_{\mathbf E} = \sup_{P\in \mathbf P} \|\cdot\|_{P,2}$ deliver the bound
\begin{equation}\label{th:quantboot2}
\mathcal S_n(\mathbf B,\mathbf E) = \sup_{\beta} \frac{\|p^{j_n\prime} \beta \|_{2,\infty}}{\sup_{P\in \mathbf P} \|p^{j_n\prime}\beta\|_{P,2}} \lesssim \frac{j_n^{5/2}\|\beta\|_2}{\|\beta\|_2} = j_n^{5/2}.
\end{equation}
Thus, we note Assumption \ref{ass:bootrates}(i) follows from Assumption \ref{ass:4bootquant}(v), result \eqref{th:quantboot2}, and $\mathcal R_n \asymp \sqrt{k_n\log(1+k_n)}/\text{s}_n\sqrt n$ implying $\mathcal R_n \mathcal S_n(\mathbf B,\mathbf E) = o(1)$.
Furthermore, we note $\mathcal R_n \mathcal S_n(\mathbf B,\mathbf E) = o(1)$, Assumption \ref{ass:4bootquant}(ii), and the definition of $\Theta$ in \eqref{eq:quantTheta} imply assumption \ref{ass:extra}(ii) is satisfied as well.
Assumption \ref{ass:bootrates}(ii) similarly follows from Assumption \ref{ass:4bootquant}(v), result \eqref{th:quantboot2}, $\mathcal R_n \asymp \sqrt{k_n\log(1+k_n)}/\text{s}_n\sqrt n$, $\kappa_\rho =1/2$, and Lemma \ref{lm:quantcov}. Finally, Assumption \ref{ass:bootrates}(iii) also follows by Assumption \ref{ass:4bootquant}(v), result \eqref{th:quantboot2}, $\mathcal R_n \asymp \sqrt{k_n(\log(1+k_n))}/\text{s}_n\sqrt n$, and $K_g = 0$.
We have thus verified the conditions of Theorem \ref{th:coupsmooth} for $R$ corresponding to \eqref{ex:quant4} and \eqref{ex:quant5}, and hence 
$$\hat U_n(R|\ell_n) \geq \UpS(R|\tilde \ell_n) +o_P(a_n)$$  
uniformly in $P\in \mathbf P_0$ for some $\tilde \ell_n \asymp \ell_n$.
Similarly, under the additional conditions imposed on the second part of this theorem, Lemma \ref{lm:forincJboot}(ii) implies that for any $\tilde \ell_n^{\rm u}$ satisfying the conditions of Theorem \ref{th:quantapprox}(ii) it follows that uniformly in $P\in \mathbf P_0$
$$\hat U_n(\Theta|+\infty) \leq \UpS(\Theta|\tilde \ell_n^{\rm u}) + o_P(a_n),$$
which implies the second claim of the theorem also holds. \qed

\begin{lemma}\label{lm:quantcons}
Let Assumptions \ref{ass:quantproc}(i)(iii), \ref{ass:quantsieve}(ii)(iii),  \ref{ass:quantid}(i), \ref{ass:quantsigma}, and \ref{ass:quantmoments}(i)(iii)(iv) hold, $k_n\log(1+k_n)/ n = o(1)$, and suppose that $$Q_n(\hat \theta_n) \leq \inf_{\theta \in \Theta_n \cap R} Q_n(\theta) + o(n^{-1/2}) \hspace{0.5 in} Q_n(\hat \theta_n^{\rm u}) \leq \inf_{\theta \in \Theta_n } Q_n(\theta) + o(n^{-1/2}).$$
Then: $\|\hat \theta_n - \Pi_n \theta_0\|_{1,\infty} \vee \|\hat \theta_n^{\rm u} - \Pi_n \theta_0\|_{1,\infty} = o_P(1)$ uniformly in $P\in \mathbf P_0$.
\end{lemma}

\noindent {\sc Proof:} We establish the result by verifying the conditions of Lemma \ref{lm:setcons} with $\tau_n = o(n^{-1/2})$.
First note that, for both $R = \Theta$ and $R$ corresponding to \eqref{ex:quant4}, Assumption \ref{ass:param}(i) is directly imposed in Assumption \ref{ass:quantproc}(i), Assumption \ref{ass:startreg}(i) holds with $B_n = O(1)$ by Assumption \ref{ass:quantmoments}(i), Assumption \ref{ass:startreg}(iii) holds with $F_n = 1$ and $J_n = O(1)$ by Lemma \ref{lm:quantcov}, Assumption \ref{app:ass:weights} is implied by Assumption \ref{ass:quantsigma}, and Assumption \ref{app:ass:rates}(i) follows from Assumption \ref{ass:quantsieve}(iii).
Next, define the following neighborhood
\begin{equation}\label{lm:quantcons1}
{\mathcal V}_n(P) \equiv \{\theta \in \Theta_n : \|\theta - \Pi_n\IDpoint\|_{1,\infty} \leq \epsilon\}
\end{equation}
for any $\epsilon > 0$ and $P\in \mathbf P_0$ (where recall $\theta_0$ is implicitly a function of $P$ through \eqref{ex:quant1}).
Further set $Q_P(\theta) \equiv \|E_P[\rho(X,\theta)q^{k_n}(Z)]\|_{\PSigma ,p}$ and note that since for any $a\in \mathbf R^{k_n}$ we have $\|a\|_p \leq \|\PSigma^{-1}\|_{o,p} \|\PSigma a \|_{p}$ and $\|\PSigma^{-1}\|_{o,p}$ is bounded uniformly in $k_n$ and $P\in \mathbf P$ by Assumption \ref{ass:quantsigma}(ii), we obtain from Lemma \ref{lm:obviousineq} and Assumption \ref{ass:quantsieve}(iii) that
\begin{align}\label{lm:quantcons2}
S_n(\epsilon) & \equiv \inf_{P\in \mathbf P_0} \{\inf_{\theta \in (\Theta_n \cap R) \setminus {\mathcal V}_n(P)} Q_{P}(\theta) - \inf_{\theta \in \Theta_n \cap R} Q_{P}(\theta)\} \notag \\
& \gtrsim \inf_{P\in \mathbf P_0} \inf_{\theta \in (\Theta_n \cap R) \setminus {\mathcal V}_n(P)} \frac{k_n^{1/p}}{\sqrt{k_n}} \|E_P[q^{k_n}(Z)\rho(X,\theta)]\|_2 + O((n\log(n))^{-1/2}).
\end{align}
We further note that the eigenvalues of $E_P[q^{k_n}(Z)q^{k_n}(Z)^\prime]$ being bounded uniformly in $k_n$ and $P\in \mathbf P$ together with Lemma \ref{aux:bessel} and Assumption \ref{ass:quantmoments}(iv) yield
\begin{multline}\label{lm:quantcons3}
\sup_{P\in \mathbf P}\sup_{\theta \in \Theta}\|E_P[q^{k_n}(Z)(E_P[\rho(X,\theta)|Z] -q^{k_n}(Z)^\prime \pi_n(\theta))]\|_2^2 \\
\lesssim \sup_{P\in \mathbf P}\sup_{\theta \in \Theta} E_P[(E_P[\rho(X,\theta)|Z] - q^{k_n}(Z)^\prime \pi_n(\theta))^2] = o(1).
\end{multline}
Therefore, since the eigenvalues of $E_P[q^{k_n}(Z)q^{k_n}(Z)^\prime]$ are bounded away from zero by Assumption \ref{ass:quantmoments}(iii), we obtain from result \eqref{lm:quantcons3} that
\begin{align}\label{lm:quantcons4}
\inf_{P\in \mathbf P_0} \inf_{\theta \in (\Theta_n \cap R) \setminus {\mathcal V}_n(P)}  &\|E_P[q^{k_n}(Z)\rho(X,\theta)]\|_2 \notag \\
& \geq \inf_{P\in \mathbf P_0} \inf_{\theta \in (\Theta_n \cap R) \setminus {\mathcal V}_n(P)} \|E_P[q^{k_n}(Z)q^{k_n}(Z)^\prime \pi_n(\theta)]\|_2 + o(1) \notag \\
& \gtrsim \inf_{P\in \mathbf P_0} \inf_{\theta \in (\Theta_n \cap R) \setminus {\mathcal V}_n(P)} (E_P[(q^{k_n}(Z)^\prime \pi_n(\theta))^2])^{1/2} + o(1).
\end{align}
Also note that Assumption \ref{ass:quantsieve}(iii) implies that for $n$ sufficiently large, we have the set inclusion $\Theta_n \cap R \setminus {\mathcal V}_n(P) \subseteq \{\theta \in \Theta : \|\theta - \Pi_n\IDpoint\|_{1,\infty} \geq \epsilon\} \subseteq \{\theta \in \Theta : \|\theta - \IDpoint\|_{1,\infty} \geq \epsilon/2\}$ holding for all $P\in \mathbf P_0$.
Hence, \eqref{lm:quantcons2}, \eqref{lm:quantcons4}, and Assumption \ref{ass:quantmoments}(iv) yield
\begin{equation}\label{lm:quantcons5}
S_n(\epsilon) \gtrsim \inf_{P\in \mathbf P_0} \inf_{\theta \in \Theta :\|\theta - \IDpoint\|_{1,\infty} \geq \epsilon} \frac{k_n^{1/p}}{\sqrt{k_n}}(E_P[(P(Y\leq \theta(D)|Z) - \tau)^2])^{1/2} + o(\frac{k_n^{1/p}}{\sqrt{k_n}}).
\end{equation}
Since $J_n \asymp 1$ and $B_n \asymp 1$, Assumption \ref{ass:quantid}(i), result \eqref{lm:quantcons5} and $k_n\log(1+k_n)/ n = o(1)$ by hypothesis imply that $k_n^{1/p}\sqrt{\log(1+k_n)}J_nB_n/\sqrt n = o(S_n(\epsilon))$ as required by lemma \ref{lm:setcons}.
The preceding arguments apply for both $R=\Theta$ and $R$ corresponding to \eqref{ex:quant4} and \eqref{ex:quant5}, and therefore the claim of the lemma follows. \qed

\begin{lemma}\label{lm:quant3ver}
Let $k_n\log(1+k_n)/n = o(1)$, and Assumptions \ref{ass:quantproc}(i)(iii), \ref{ass:quantid}, \ref{ass:quantsieve}(ii)(iii),  \ref{ass:quantsigma},  and \ref{ass:quantmoments}(i)(iii)(iv) hold.
For both $R = \Theta$ and $R$ corresponding to \eqref{ex:quant4} and \eqref{ex:quant5}, it  follows that Assumption \ref{ass:keycons} holds with $\|\cdot\|_{\mathbf E} \equiv \sup_{P\in \mathbf P} \|\cdot\|_{P,2}$, $\nu_n \asymp \sqrt{k_n}/\text{\rm s}_n k_n^{1/p}$, and $\mathcal V_n(P) \equiv \{\theta \in \Theta_n : \|\theta - \Pi_n \theta_0\|_{1,\infty} \leq \epsilon\}$ for some $\epsilon > 0$.
\end{lemma}

\noindent {\sc Proof:} For either $R = \Theta$ or $R$ corresponding to \eqref{ex:quant4} and \eqref{ex:quant5} set $\mathcal V_n(P) \equiv \{\theta \in \Theta_n \cap R : \|\theta - \Pi_n \IDpoint\|_{1,\infty} \leq \epsilon\}$ and note that Assumption \ref{ass:keycons}(ii) is then satisfied by Lemma \ref{lm:quantcons}.
Further observe that since $\Pi_n \IDpoint \in \Theta_n$, there exists a $\beta_{0n}$ such that $\Pi_n \IDpoint = p^{j_n\prime}\beta_{0n}$.
For any $\theta = p^{j_n\prime}\beta \in \mathcal V_n(P)$, it then follows by Assumptions \ref{ass:quantsieve}(ii) and \ref{ass:quantid}(ii) that
\begin{multline*}
\sup_{P\in \mathbf P} \|p^{j_n\prime}(\beta - \beta_{0n})\|_{P,2} \lesssim \|\beta - \beta_{0n}\|_2 \\
\leq \frac{1}{\text{s}_n} \times \inf_{P\in \mathbf P_0} \inf_{\|\theta-\Pi_n \IDpoint\|_{1,\infty} \leq \epsilon} \|E_P[f_{Y|DZ,P}(\theta(D)|D,Z)q^{k_n}(Z)p^{j_n}(D)^\prime(\beta_n - \beta_{0n})]\|_2.
\end{multline*}
Hence, since $\|\cdot\|_{\mathbf E} = \sup_{P\in \mathbf P} \|\cdot\|_{P,2}$, the mean value theorem allows us to conclude that
\begin{align}\label{lm:quant3ver3}
&\frac{\text{s}_n k_n^{1/p}}{\sqrt{k_n}}  \|p^{j_n\prime}(\beta - \beta_{0n})\|_{\mathbf E} \notag \\
&  \lesssim \frac{k_n^{1/p}}{\sqrt{k_n}}\|E_P[(P(Y \leq p^{j_n}(D)^\prime \beta|D,Z) - P(Y\leq p^{j_n}(D)^\prime \beta_{0n}|D,Z))q^{k_n}(Z)]\|_2 \notag \\
&  \lesssim \|E_P[(P(Y \leq p^{j_n}(D)^\prime \beta|D,Z) - P(Y\leq p^{j_n}(D)^\prime \beta_{0n}|D,Z))q^{k_n}(Z)]\|_{\PSigma,p}
\end{align}
for any $\theta = p^{j_n\prime}\beta \in \mathcal V_n(P)$ and $P\in \mathbf P_0$, and where the final inequality follows from Lemma \ref{lm:obviousineq}, $\|a\|_p \leq \|\PSigma^{-1}\|_{o,p} \|\PSigma a\|_{o,p}$ for any $a\in \mathbf R^{k_n}$, and Assumption \ref{ass:quantsigma}(ii).
Since the preceding arguments apply to both $R = \Theta$ and $R$ corresponding to \eqref{ex:quant4} and \eqref{ex:quant5}, result \eqref{lm:quant3ver3} verifies Assumption \ref{ass:keycons}(i) is satisfied with $\nu_n \asymp \sqrt{k_n}/\text{s}_nk_n^{1/p}$ for both choices of $R$ and therefore the claim of the lemma follows. \qed

\begin{lemma}\label{lm:quantmodulus}
Let Assumption \ref{ass:quantproc}(iii) hold. If $f(y,d) = 1\{y \leq \theta(d)\} - \tau$ for some $\theta \in \Theta$ ($\Theta$ as in \eqref{eq:quantTheta}) and $z\mapsto q(z)$ is differentiable with bounded level and derivative, then there exists a $K < \infty$ independent of $f$, such that for all $P\in \mathbf P$ we have
$$\varpi(fq,h,P) \leq K\times \{\|q\|_\infty \sqrt{h} + \|q\|_{1,\infty} h\}.$$
\end{lemma}

\noindent {\sc Proof:} First note that since $\|f\|_\infty \leq 1$, we can obtain by direct calculation and the definition of the integral modulus of continuity in \eqref{appcoup:eq1} the upper bound
\begin{equation}\label{lm:quantmodulus1}
\varpi^2(fq,h,P) \leq 2\|q\|_\infty^2 \varpi^2(f,h,P) + 2\varpi^2(q,h,P).
\end{equation}
For $\Omega_Z(P)$ the support of $Z$ under $P$, the mean value theorem then implies that
\begin{equation}\label{lm:quantmodulus2}
\varpi^2(q,h,P) \equiv \sup_{\|s\|_2\leq h} E_P[(q(Z+s) - q(Z))^21\{Z+s \in \Omega_Z(P)\}] \leq \|q\|_{1,\infty}^2 h^2.
\end{equation}
Furthermore, for any $(s_y,s_d)\in \mathbf R^2$ and $d \in [0,1]$ such that $d+ s_d \in [0,1]$, we also note that the mean value theorem and Assumption \ref{ass:quantproc}(iii) imply that
\begin{align}\label{lm:quantmodulus3}
E_P[(1\{Y + s_y & \leq \theta(D + s_d)\}  - 1\{Y \leq \theta(D)\})^2|D = d] \notag\\
 & = | P(Y\leq \theta(D + s_d) - s_y|D =d ) - P(Y \leq \theta(D)|D = d)| \notag \\
 & \lesssim |\theta(d + s_d) - s_y - \theta(d)|.
\end{align}
Hence, by the law of iterated expectations, a second application of the mean value theorem, and employing that $\|\theta\|_{1,\infty} \leq C_0$ by definition of $\Theta$, we can conclude
\begin{multline}\label{lm:quantmodulus4}
\varpi^2(f,h,P) \leq \sup_{\|(s_y,s_d)\|_2 \leq h} E_P[(1\{Y + s_y \leq \theta(D + s_d)\} - 1\{Y \leq \theta(D)\})^21\{D + s_d \in [0,1]\}]\\
\lesssim \sup_{\|(s_y,s_d)\|_2 \leq h} E_P[|\theta(D + s_d) - s_y - \theta(D)|1\{D + s_d \in [0,1]\}] \lesssim h.
\end{multline}
The claim of the lemma then follows from \eqref{lm:quantmodulus1}, \eqref{lm:quantmodulus2}, and \eqref{lm:quantmodulus4}. \qed

\begin{lemma}\label{lm:quantcov}
Define the class $\mathcal F_n \equiv \{f: f(v) = 1\{y \leq \theta(d)\} - \tau \text{ for some } \theta \in \Theta_n\}$ for $\Theta_n$ as in \eqref{eq:quantThetaN}, and suppose that Assumption \ref{ass:quantproc}(iii) and \ref{ass:quantsieve}(ii) hold.
For $\zeta_{j_n} \geq (1\wedge \sup_{d\in[0,1]}\|p^{j_n}(d)\|_2)$, it then follows that for all $\epsilon \leq 1$ and some $1\leq K<\infty$
$$ \sup_{P\in \mathbf P} N_{[\hspace{0.02 in}]}(\epsilon,\mathcal F_n,\|\cdot\|_{P,2}) \leq \exp\{\frac{K}{\epsilon}\} \wedge (\frac{K\sqrt{\zeta_{j_n}}}{\epsilon})^{2j_n},$$
and $\sup_{P\in \mathbf P} J_{[\hspace{0.02 in}]}(\epsilon,\mathcal F_n,\|\cdot\|_{P,2}) \lesssim \sqrt{1\wedge \epsilon} \wedge \sqrt{j_n(\log(\zeta_{j_n})+\log(1\vee \epsilon^{-1}))}(1\wedge \epsilon)$.
\end{lemma}

\noindent {\sc Proof:} We first note that if $\theta_l(d) \leq \theta(d) \leq \theta_u(d)$, then it immediately follows that
\begin{equation}\label{lm:quantcov1}
1\{y \leq \theta_l(d)\} - \tau \leq 1\{y \leq \theta(d)\} - \tau \leq 1\{y \leq \theta_u(d)\},
\end{equation}
which implies brackets for $\Theta_n$ readily yield brackets in $\mathcal F_n$. Moreover, by the mean value theorem and Assumption \ref{ass:quantproc}(iii) we can in addition conclude that
\begin{multline}\label{lm:quantcov2}
E_P[(1\{Y \leq \theta_l(D)\} - 1\{Y \leq \theta_u(D)\})^2] \\
= E_P[P(Y \leq \theta_u(D)|D) - P(Y \leq \theta_l(D)|D)] \lesssim E_P[|\theta_u(D) - \theta_l(D)|] .
\end{multline}
Hence, combining results \eqref{lm:quantcov1} and \eqref{lm:quantcov2} it follows that for some $M_0 < \infty$ we have
\begin{equation}\label{lm:quantcov3}
N_{[\hspace{0.02 in}]}(\epsilon,\mathcal F_n,\|\cdot\|_{P,2}) \leq N_{[\hspace{0.02 in}]}(\frac{\epsilon^2}{M_0},\Theta_n,\|\cdot\|_{P,1}).
\end{equation}
On the other hand, since $\Theta_n \subseteq \Theta$, we also obtain by Corollary 2.7.2 in \cite{vandervaart:wellner:1996}, $\|\cdot\|_{P,2}\leq \|\cdot\|_\infty$ and inequality \eqref{lm:quantcov3} that
\begin{equation}\label{lm:quantcov4}
\sup_{P\in \mathbf P} N_{[\hspace{0.02 in}]}(\epsilon,\mathcal F_n,\|\cdot\|_{P,2}) \leq N_{[\hspace{0.02 in}]}(\frac{\epsilon^2}{M_0},\Theta_n,\|\cdot\|_{\infty}) \leq \exp\{\frac{M_1}{\epsilon}\}.
\end{equation}
In addition, the Cauchy-Schwarz inequality implies $\sup_{P\in \mathbf P} \|p^{j_n\prime}(\beta_1-\beta_2)\|_{P,1} \leq \zeta_{j_n} \|\beta_1 - \beta_2\|_2$
for any $\beta_1,\beta_2\in \mathbf R^{j_n}$.
Therefore, defining $\mathcal B_n \equiv \{\beta \in \mathbf R^{j_n}: p^{j_n\prime}\beta \in \Theta_n\}$ Theorem 2.7.11 in \cite{vandervaart:wellner:1996} allows us to conclude
\begin{equation}\label{lm:quantcov5}
\sup_{P\in \mathbf P} N_{[\hspace{0.02 in}]}(\epsilon,\mathcal F_n,\|\cdot\|_{P,1}) \leq N(\frac{\epsilon^2}{2M_0 \zeta_{j_n}}, \mathcal B_n, \|\cdot\|_2).
\end{equation}
Further note that by Assumption \ref{ass:quantsieve}(ii), we have $\|p^{j_n\prime}\beta\|_{P,2} \asymp \|\beta\|_2$ uniformly in $P\in \mathbf P$ and $n$, and hence since $\sup_{P\in \mathbf P} \|p^{j_n\prime}\beta\|_{P,2}\leq \|p^{j_n\prime}\beta\|_\infty \leq C_0$, it follows
\begin{equation}\label{lm:quantcov6}
\sup_{P\in \mathbf P} N_{[\hspace{0.02 in}]}(\epsilon,\mathcal F_n,\|\cdot\|_{P,1}) \leq  (\frac{M_2 \sqrt{\zeta_{j_n}}}{\epsilon})^{2j_n}
\end{equation}
for some $M_2 < \infty$ due to result \eqref{lm:quantcov5}.
The first claim of the lemma therefore follows from \eqref{lm:quantcov4} and \eqref{lm:quantcov6}.
Moreover, noting $N_{[\hspace{0.02 in}]}(1,\mathcal F_n,\|\cdot\|_\infty) = 1$ we also obtain
\begin{multline}\label{lm:quantcov7}
\sup_{P\in \mathbf P} J_{[\hspace{0.02 in}]}(\epsilon,\mathcal F_n,\|\cdot\|_{P,2}) \leq (\int_0^{1\wedge \epsilon}(1+  \frac{K}{u})^{1/2}du ) \wedge (\int_0^{1\wedge \epsilon} (1+2j_n\log(\frac{K\sqrt{\zeta_{j_n}}}{u}))^{1/2}du)  \\
 \lesssim \sqrt{1\wedge \epsilon} \wedge \{\sqrt{j_n\log(\zeta_{j_n})} + \int_0^1 (j_n\log(\frac{1}{v(1\wedge \epsilon)}))^{1/2}dv\}(1\wedge \epsilon),
\end{multline}
where the second inequality follows by the change of variables $v = u/(1\wedge \epsilon)$. The claim of the lemma thus follows from \eqref{lm:quantcov7} and direct calculation. \qed

\begin{lemma}\label{lm:quantcoup}
Let Assumptions \ref{ass:quantproc}(i)(iii)(iv) and \ref{ass:quantsieve}(ii) hold, $\Theta_n$ be as in \eqref{eq:quantThetaN}, set $\pi_{0n} \equiv \max_{1 \leq k \leq k_n} \|q_k\|_\infty$ and $\pi_{1n} \equiv \max_{1\leq k \leq k_n} \|q_k\|_{1,\infty}$, and suppose $\log(k_n \vee \pi_{0n} \vee\sup_{d}\|p^{j_n}(d)\|_2) = O(\log(n))$. If $j_n/n = o(1)$, then Assumption \ref{ass:coupling}(i) holds with $R = \Theta$ for any $a_n$ with $(\pi_{0n} k_n^{1/p} \log(n)\sqrt{j_n} / \sqrt {n})(\sqrt{j_n} + \pi_{0n} n^{1/3} + \pi_{1n} n^{1/6}) = o(a_n)$.
\end{lemma}

\noindent {\sc Proof:} We establish the lemma by applying Theorem \ref{th:coup}.
To this end, define the class $\tilde {\mathcal F}_n  \equiv \{ fq_k \text{ for some } f\in \mathcal F_n, 1\leq k \leq k_n\}$ and let $\Iso$ be a Gaussian process on $\tilde {\mathcal F}_n$ satisfying $E[\Iso(f_1)] = 0$ and $E[\Iso (f_1)\Iso(f_2)] = \text{Cov}_P\{f_1(V),f_2(V)\}$ for any $f_1,f_2\in \tilde{\mathcal F_n}$. 
For any $\theta \in \Theta_n$, set $\WP(\theta) \equiv (\Iso(\rho(\cdot,\theta)q_1),\ldots, \Iso(\rho(\cdot,\theta)q_{k_n})^\prime$ and note
\begin{equation}\label{lm:quantcoup1}
\sup_{\theta \in \Theta_n} \|\Gemp(\theta)  - \WP(\theta)\|_p \leq \sup_{f\in \tilde {\mathcal F}_n} k_n^{1/p} |\frac{1}{\sqrt n}\sum_{i=1}^n (f(V_i) - E_P[f(V)]) -\Iso(f)|.
\end{equation}
We proceed by applying Theorem \ref{th:coup} to the class $\tilde {\mathcal F}_n$ with $\delta_n \asymp \sqrt{j_n/n}$.
Note that Assumptions \ref{ass:coupreg} and \ref{ass:supportreg} are directly imposed in Assumption \ref{ass:quantproc}(iv), while Assumption \ref{ass:coupprocessreg}(i) is satisfied by Lemma \ref{lm:quantmodulus}, and Assumption \ref{ass:coupprocessreg}(ii) holds with $K_n = \pi_{0n}$ since $\mathcal F$ has envelope $1$.
Furthermore, for $S_n$ as in \eqref{th:coupdisp1}, we have
\begin{equation*}\label{lm:quantcoup2}
S_n^2 \lesssim  \sum_{i=0}^{\lceil \log_2 n\rceil} 2^i \{\frac{\pi_{0n}^2}{2^{\frac{i}{3}}} + \frac{\pi_{1n}^2}{2^{\frac{2i}{3}}}\} \lesssim \pi_{0n}^2 n^{2/3} + \pi_{1n}^2 n^{1/3}
\end{equation*}
due to Lemma \ref{lm:quantmodulus}.
Also note that Lemmas \ref{aux:simpbracket} and \ref{lm:quantcov} together imply
\begin{equation*}\label{lm:quantcoup3}
\sup_{P\in \mathbf P} \log(N_{[\hspace{0.02 in}]}(\delta_n,\tilde {\mathcal F}_n,\|\cdot\|_{P,2})) \leq  \sup_{P\in \mathbf P} \log(k_n N_{[\hspace{0.02 in}]}(\frac{\delta_n}{\pi_{0n}},\mathcal F_n,\|\cdot\|_{P,2})) \lesssim j_n \log(n),
\end{equation*}
where we employed that $\log(k_n \vee \pi_{0n} \vee \sup_d \|p^{j_n}(d)\|_2) = O(\log(n))$  and $\delta_n = \sqrt{j_n/n}$.
Similarly, Lemmas \ref{aux:simpbracket} and \ref{lm:quantcov} and the change of variables $v = u/\pi_{0n}$ yield
\begin{multline*}\label{lm:quantcoup4}
\sup_{P\in \mathbf P} J_{[\hspace{0.02 in}]}(\delta_n,\tilde {\mathcal F}_n,\|\cdot\|_{P,2}) \leq \sup_{P\in \mathbf P} \int_0^{\delta_n}(1+ \log(k_n) + \log (N_{[\hspace{0.02 in}]}(\frac{u}{\pi_{0n}},\mathcal F_n,\|\cdot\|_{P,2})))^{1/2}du\\
\leq \delta_n \sqrt{\log(k_n)} + \sup_{P\in \mathbf P} J_{[\hspace{0.02 in}]}(\frac{\delta_n}{\pi_{0n}},\mathcal F_n,\|\cdot\|_{P,2})\pi_{0n} \lesssim \frac{j_n \sqrt{\log(n)}}{\sqrt n},
\end{multline*}
where the final inequality follow from $\log(\pi_{0n}) = O(\log(n))$, $\log(k_n) = O(\log(n))$, and $\log(\sup_{d} \|p^{j_n}(d)\|_2) = O(\log(n))$.
Hence, Theorem \ref{th:coup} implies that
\begin{equation*}\label{lm:quantcoup4}
\sup_{f\in \tilde {\mathcal F}_n}  |\frac{1}{\sqrt n}\sum_{i=1}^n (f(V_i) - E_P[f(V)]) -\Iso(f)|
= O_P(\frac{\pi_{0n} \log(n)\sqrt{j_n}}{\sqrt {n}}\{\sqrt{j_n} + \pi_{0n} n^{1/3} + \pi_{1n} n^{1/6}\})
\end{equation*}
uniformly in $P\in \mathbf P$, which together with \eqref{lm:quantcoup1} establishes the Lemma. \qed

\begin{lemma}\label{lm:quantrhover}
If Assumption \ref{ass:quantproc}(iii) holds, then it follows that Assumptions \ref{ass:coupling}(ii) and \ref{ass:driftlin} are satisfied when $R = \Theta$ ($\Theta$ as in \eqref{eq:quantTheta}) with $\|\cdot\|_{\mathbf E} = \sup_{P\in \mathbf P} \|\cdot\|_{P,2}$, $\|\cdot\|_{\mathbf L} = \|\cdot\|_\infty$, $\kappa_\rho = 1/2$, $m_P(\theta)(Z) \equiv P(Y\leq \theta(D)|Z)$, and
\begin{equation}\label{lm:quantrhoverdisp}
\nabla m_P(\theta)[h](Z) \equiv E_P[f_{Y|DZ,P}(\theta(D)|D,Z)h(D)|Z].
\end{equation}
\end{lemma}

\noindent {\sc Proof:} For $\theta_1 \vee \theta_2$ and $\theta_1 \wedge \theta_2$ the pointwise minimum and maximum of $\theta_1$ and $\theta_2$, note that the conditional density $f_{Y|DZ,P}$ being bounded in $(D,Z)$ and $P\in \mathbf P$ by Assumption \ref{ass:quantproc}(iii) together with the mean value theorem imply that
\begin{multline*}
E_P[(\rho(X,\theta_1) - \rho(X,\theta_2))^2] = E_P[P(Y \leq \theta_1(D) \vee \theta_2(D)|D) -P(Y \leq \theta_1(D) \wedge \theta_2(D)|D)] \\
\lesssim E_P[\theta_1(D) \vee \theta_2(D) -  \theta_1(D) \wedge \theta_2(D)] \leq \sup_{P\in \mathbf P} \|\theta_1 - \theta_2\|_{P,2},
\end{multline*}
where in the final inequality we employed Jensen's inequality and that $\theta_1(d) \vee \theta_2(d) - \theta_1(d) \wedge \theta_2(d) = |\theta_1(d) - \theta_2(d)|$.
It thus follows Assumption \ref{ass:coupling}(ii) holds with $\|\cdot\|_{\mathbf E} = \sup_{P\in \mathbf P} \|\cdot\|_{P,2}$ and $\kappa_\rho = 1/2$.
Moreover, Jensen's inequality and the mean value theorem imply for some $\bar \theta$ such that $\bar \theta(d)$ is a convex combination of $\theta_1(d)$ and $\theta_2(d)$ that
\begin{align*}\label{lm:quantrhover2}
E_P[(P(Y&\leq\theta_1(D)|Z) - P(Y \leq \theta_2(D)|Z) - \nabla m_P(\theta_2)[\theta_1 - \theta_2](Z))^2] \notag \\
& \leq E_P[(\{f_{Y|DZ,P}(\bar \theta(D)|D,Z) - f_{Y|DZ,P}(\theta_2(D)|D,Z)\}\{\theta_1(D) - \theta_2(D)\})^2] \notag \\
& \lesssim \|\theta_1 - \theta_2\|_\infty^2 \times \sup_{P\in \mathbf P} E_P[(\theta_1(D) - \theta_2(D))^2],
\end{align*}
where the final inequality follows from $f_{Y|DZ,P}$ being Lipschitz uniformly in $(D,Z)$ and $P\in \mathbf P$.
Hence, we may conclude Assumption \ref{ass:driftlin}(i) is satisfied with $\|\cdot\|_{\mathbf L}= \|\cdot\|_\infty$ and $\|\cdot\|_{\mathbf E} = \sup_{P\in \mathbf P} \|\cdot\|_{P,2}$.
Furthermore, once again employing Jensen's inequality  and that $f_{Y|DZ,P}$ is Lipschitz uniformly in $(D,Z)$ and $P\in \mathbf P$ yields
\begin{multline}
E_P[(E_P[\{f_{Y|DZ,P}(\theta_1(D)|D,Z) - f_{Y|DZ,P}(\theta_2(D)|D,Z)\}h(D)|Z])^2]  \\ \lesssim \|\theta_1 - \theta_2\|_\infty^2 \times  \sup_{P\in \mathbf P} \|h\|_{P,2}^2
\end{multline}
which implies Assumption \ref{ass:driftlin}(ii) is also satisfied under the stated choices of $\|\cdot\|_{\mathbf L}$ and $\|\cdot\|_{\mathbf E}$.
Finally, we note Assumption \ref{ass:driftlin}(iii) is immediate due to Jensen's inequality and $f_{Y|DZ,P}$ being bounded uniformly in $(D,Z)$ and $P\in \mathbf P$. \qed

\begin{lemma}\label{lm:quantile4ver}
If Assumption \ref{ass:4bootquant}(i) holds, $\mathbf B = C^2_B([0,1])$ and $\Upsilon_G$, $\Upsilon_F$, and $\Theta$ are as defined in \eqref{ex:quant4}, \eqref{ex:quant5}  (with $\lambda\neq 0$), and \eqref{eq:quantTheta}, then it follows that Assumptions \ref{ass:locineq}, \ref{ass:loceq}, and \ref{ass:ineqlindep} are satisfied with $K_g = 0$, $\nabla \Upsilon_G(\theta)[h] = -\nabla^2 h$, and
\begin{equation}\label{lm:quantile4verdisp}
\nabla \Upsilon_F(\theta)[h] = 2\int_0^1 \theta(u)h(u)du - 2 (\int_0^1 \theta(u)du)(\int_0^1 h(u)du).
\end{equation}
\end{lemma}

\noindent {\sc Proof:} Note that since $\Upsilon_G$ is linear and continuous, it immediately follows that Assumptions \ref{ass:locineq}(i) and \ref{ass:locineq}(ii) hold with $\nabla \Upsilon_G = \Upsilon_G$ and $K_g = 0$.
It further follows from $\nabla \Upsilon_G = \Upsilon_G$ and the definitions of the operator norm $\|\cdot\|_o$ and $\|\cdot\|_{m,\infty}$ that
\begin{equation}\label{lm:quantile4ver1}
\|\nabla \Upsilon_G(\theta)\|_o = \sup_{\|h\|_{2,\infty} =1} \|-\nabla^2 h\|_\infty \leq 1,
\end{equation}
which implies Assumption \ref{ass:locineq}(iii) holds with $M = 1$. Moreover, by direct calculation
\begin{multline}\label{lm:quantile4ver2}
|\Upsilon_F(\theta_1) - \Upsilon_F(\theta_2) - \nabla \Upsilon_F(\theta_1)[\theta_1 - \theta_2]|  \\ = |\int_0^1(\theta_1(u) - \theta_2(u))^2du - (\int_0^1 (\theta_1(u) - \theta_2(u))du)^2| \leq \|\theta_1 - \theta_2\|_{2,\infty}^2,
\end{multline}
which implies $\Upsilon_F$ is indeed Fr\'echet differentiable and its derivative is equal to $\nabla \Upsilon_F$ as defined in \eqref{lm:quantile4verdisp}.
In addition, by \eqref{lm:quantile4verdisp} and Jensen's inequality we have
\begin{multline}\label{lm:quantile4ver3}
\|\nabla \Upsilon_F(\theta_1) - \nabla \Upsilon_F(\theta_2)\|_o \\
= \sup_{\|h\|_{2,\infty} = 1} 2|\int_0^1(\theta_1(u)-\theta_2(u))(h(u) - \int_0^1 h(\tilde u) d\tilde u) du| \leq 2\|\theta_1 - \theta_2\|_{2,\infty},
\end{multline}
which together with \eqref{lm:quantile4ver2} implies Assumptions \ref{ass:loceq}(i) and \ref{ass:loceq}(ii) hold with $K_f = 2$.
Next, note that since $\lambda \neq 0$ it follows that $\mathbf F_n = \mathbf R$.
For any $\theta \in \mathbf B_n$ such that $\Upsilon_F(\theta) \neq 0$, we then define $\nabla \Upsilon_F(\theta)^{-} : \mathbf F_n \to \mathbf B_n$ to be given (for any $c\in \mathbf R$) by
\begin{equation}\label{lm:quantile4ver4}
\nabla \Upsilon_F(\theta)^{-}[c](d) \equiv c\times \frac{ \theta(d) - \int_0^1  \theta(u)du}{2\Upsilon_F(\theta)},
\end{equation}
and note that since $\theta \in \mathbf B_n$ and the constant function is in $\mathbf B_n$ by Assumption \ref{ass:4bootquant}(i), it follows that $\nabla \Upsilon_F(\theta)^{-}[c] \in \mathbf B_n$.
Moreover, by direct calculation we obtain
\begin{equation}\label{lm:quantile4ver5}
\nabla \Upsilon_F(\theta) \nabla \Upsilon_F(\theta)^{-}[c] = 2\int_0^1 \theta(u)\{c\times \frac{ \theta(u) - \int_0^1  \theta(\tilde u)d\tilde u}{2\Upsilon_F(\theta)} \}du = c\times \frac{2\Upsilon_F(\theta)}{2\Upsilon_F(\theta)} = c,
\end{equation}
which verifies $\nabla \Upsilon_F(\theta)^{-}$ is indeed the right inverse of $\nabla \Upsilon_F(\theta)$. In addition note that
\begin{equation}\label{lm:quantile4ver6}
\|\nabla \Upsilon_F(\theta)^{-}\|_o = \sup_{|c| = 1} \|c\times \frac{\theta - \int_0^1  \theta(u)du}{2\Upsilon_F(\theta)}\|_{2,\infty} \leq \frac{\|\theta\|_{2,\infty}}{|\Upsilon_F(\theta)|},
\end{equation}
and hence, since $\|\theta\|_{2,\infty} \leq C_0$ and $\Upsilon_F(\theta) = \lambda$ for any $\theta \in \IDsetRsieve$, it follows that we may select an $\epsilon >0$ such that Assumption \ref{ass:loceq}(iv) holds with $M = 4C_0/\lambda$.

Next, let $\theta_2$ be the function given by $\theta_2(d) = d^2$ and note that by Assumption \ref{ass:4bootquant}(i) it follows that $\theta_2 \in \mathbf B_n$.
For any $\theta \in \IDsetRsieve$ we may then set $h$ to equal
\begin{equation}\label{lm:quantile4ver7}
h \equiv \frac{2\lambda}{C_0} \theta_2 - \frac{\nabla \Upsilon_F(\theta)[\theta_2]}{C_0} \theta,
\end{equation}
which belongs to $\mathbf B_n$ since $\theta_2,\theta \in \mathbf B_n$.
Further observe $\nabla \Upsilon_F(\theta)[\theta] = 2\Upsilon_F(\theta) = 2\lambda$ due to $\theta \in R$, and hence by linearity of $\nabla \Upsilon_F(\theta)$ and \eqref{lm:quantile4ver7} we can conclude that $h \in \mathbf B_n \cap \mathcal N(\nabla \Upsilon_F(\theta))$.
In addition, it also follows from $\Upsilon_G  = \nabla \Upsilon_G$ that
\begin{equation}\label{lm:quantile4ver8}
\Upsilon_G(\theta)(u) + \nabla \Upsilon_G(\theta)[h](u) = - \nabla^2 \theta(u)(1 - \frac{\nabla \Upsilon_F(\theta)[\theta_2]}{C_0}) - \frac{4\lambda}{C_0} \leq -\frac{4\lambda}{C_0},
\end{equation}
where the inequality results from $-\nabla^2\theta(u) \leq 0$ due to $\theta \in \IDsetRsieve$, $ \theta_2(d) =  d^2$, and $|\nabla \Upsilon_F(\theta)[\theta_2]| \leq  C_0$ because $\|\theta\|_{2,\infty} \leq C_0$ since $\theta \in \IDsetRsieve \subseteq \Theta_n$.
By similar arguments and the triangle inequality we also have $\|h\|_{2,\infty} \leq 4\lambda/C_0+ C_0$ and hence by \eqref{lm:quantile4ver8} we conclude Assumption \ref{ass:ineqlindep} is satisfied. \qed

\begin{lemma}\label{lm:quantbootcoup}
Suppose Assumptions \ref{ass:quantproc}(i)(iii)(iv) and \ref{ass:quantsieve}(ii) hold, $\Theta_n$ be as in \eqref{eq:quantThetaN}, and let $\pi_{0n} \equiv \max_{1\leq k \leq k_n} \|q_k\|_\infty$ and $\pi_{1n} \equiv \max_{1\leq k \leq k_n}\|q_k\|_{1,\infty}$.
For any sequence $d_n\uparrow \infty$ such that $d_n^4\log(1+d_n) = o(n)$ and $\delta_n \asymp d_n^{-1/6} + \pi_{1n}/(\pi_{0n}d_n^{1/3})$ satisfies $\delta_n \log(1+k_n) = o(1)$ it follows that uniformly in $P\in \mathbf P$ we have
\begin{multline*}
\sup_{\theta \in \Theta_n} \|\Bemp(\theta) - \WPT(\theta)\|_p \\
= O_P(\frac{k_n^{1/p} d_n^2 \pi_{0n} \sqrt{\log(1+d_n)}}{\sqrt n} + \pi_{0n}k_n^{1/p}(\sqrt{\delta_n} + \sqrt n\exp\{-n\delta_n^3\})).
\end{multline*}
\end{lemma}

\noindent {\sc Proof:} We first define the class $\tilde{\mathcal F}_n \equiv \{fq_k \text{ for some } f \in \mathcal F_n \text{ and } 1\leq k \leq k_n\}$, let $\IsoW$ be an isonormal Gaussian process on $\tilde {\mathcal F}_n$ independent of $\{V_i\}_{i=1}^n$, set $\WPT(\theta) \equiv (\IsoW(\rho(\cdot,\theta)q_1),\ldots, \IsoW(\rho(\cdot,\theta)q_{k_n}))^\prime$, and for any $f\in \mathcal F_n$ define $\Wboot(f)$ to equal
$$\Wboot(f) \equiv \frac{1}{\sqrt n}\sum_{i=1}^n \omega_i(f(V_i) - \frac{1}{n}\sum_{j=1}^n f(V_j))$$
where $\{\omega_i\}_{i=1}^n$ are the same weights employed in  $\Bemp$. 
These definitions then imply 
\begin{equation}\label{lm:quantbootcoup1}
\sup_{\theta \in \Theta_n} \|\Bemp(\theta) - \WPT(\theta)\|_p \leq \sup_{f\in \tilde {\mathcal F}_n} k_n^{1/p} |\Wboot(f) - \IsoW(f)|.
\end{equation}
In what follows, we aim to establish the lemma by applying Theorem \ref{th:mainbootcoup} to the class $\tilde {\mathcal F}_n$ by relying on a Haar basis expansion as in Lemmas \ref{aux:coup} and \ref{aux:haar}.
Specifically, note that by Assumption \ref{ass:quantproc}(iv) and Lemma \ref{aux:coup}, there exists a sequence of partitions $\Delta_{n}(P) = \{\Delta_{d,n}(P) : d = 1,\ldots, d_n\}$ of the support of $V \equiv (Y,D,Z)$ such that $P(\Delta_{d,n}(P)) = 1/d_n$. For any $1\leq d \leq d_n-1$ we then set $\{f_{d,n,P}\}_{d=1}^{d_n-1}$ to be given by
\begin{equation}\label{lm:quantbootcoup2}
f_{d,n,P}(V) \equiv \frac{(d_n1\{V \in \Delta_{d,n}(P)\}-1)}{\sqrt{d_n-1}}
\end{equation}
and let $f_{n,P}^{d_n}(v) \equiv (f_{1,n,P}(v),\ldots, f_{d_n-1,n,P}(v))^\prime$.
Then note that $E_P[f_{n,P}^{d_n}(V)] = 0$ and 
\begin{equation}\label{lm:quantbootcoup2p5}
E_P[f_{d,n,P}(V)f_{\tilde d,n,P}(V)] = \left\{\begin{array}{cl} 1 & \text{ if } d = \tilde d \\ -\frac{1}{d_n-1} & \text{ if } d \neq \tilde d\end{array}\right. ~.
\end{equation}
By result \eqref{lm:quantbootcoup2p5} and direct calculation it follows Assumption \ref{ass:4seriescoup}(i) holds with $C_n \asymp 1$, while \eqref{lm:quantbootcoup2} implies Assumption \ref{ass:4seriescoup}(ii) holds with $K_n \asymp \sqrt{d_n}$.
Also note that $\text{Var}_P\{f^{d_n}_{n,P}(V)\} = E_P[f^{d_n}_{n,P}(V)f^{d_n}_{n,P}(V)^\prime]$ and therefore using that by \eqref{lm:quantbootcoup2p5} the smallest eigenvalue of $E_P[f^{d_n}_{n,P}(V)f^{d_n}_{n,P}(V)^\prime]$ is of order $1/d_n$, we obtain uniformly in $P\in \mathbf P$ that
\begin{equation}\label{lm:quantbootcoup2p6}
\|\text{Var}_P^{-1}\{f^{d_n}_{n,P}(V)\}\|_{o,2} \lesssim d_n .
\end{equation}

We next aim to verify that Assumption \ref{ass:4seriesreg} is satisfied by setting $\beta_{n,P}(f)$ to be
\begin{equation}\label{lm:quantbootcoup3}
\beta_{n,P}(f) \equiv \left(\begin{array}{c} \frac{\sqrt{d_n-1}}{d_n}(E_P[f(V)|V \in \Delta_{1,n}(P)]-E_P[f(V)|V \in \Delta_{d_n,n}(P)])\\ \vdots \\ \frac{\sqrt{d_n-1}}{d_n}(E_P[f(V)|V \in \Delta_{d_n-1,n}(P)]-E_P[f(V)|V \in \Delta_{d_n,n}(P)])\end{array}\right)
\end{equation}
for any $f\in \tilde{\mathcal F}_n$.
Then observe that, by direct calculation, for any $f \in \tilde {\mathcal F}_n$ we have that
\begin{align}\label{lm:quantbootcoup4}
&f_{n,P}^{d_n}(V)^\prime \beta_{n,P}(f)\notag \\ & = \sum_{d=1}^{d_n-1} (E_P[f(V)|V \in \Delta_{d,n}(P)]-E_P[f(V)|V \in \Delta_{d_n,n}(P)])(1\{V \in \Delta_{d,n}(P)\}-1/d_n) \notag \\
& = \sum_{d=1}^{d_n}E_P[f(V)|V \in \Delta_{d,n}(P)] 1\{V \in \Delta_{d,n}(P)\} - E_P[f(V)],
\end{align}
where the final equality follows from $\{\Delta_{d,n}(P)\}_{d=1}^{d_n}$ being a partition of the support of $V$ that satisfies $P(V\in \Delta_{d,n}(P)) = 1/d_n$.
Defining  $\mathcal G_{n,P} \equiv \{(f -\int f dP) - f_{n,P}^{d_n\prime}\beta_{n,P}(f) : f \in \tilde{\mathcal F}_n\}$, then observe that since $\mathcal F_n$ has envelope 1, it follows the class $\tilde {\mathcal F}_n$ has envelope $\pi_{0n}$ and hence by \eqref{lm:quantbootcoup4} and Jensen's inequality, the class $\mathcal G_{n,P}$ has envelope $G_{n,P} \equiv 2\pi_{0n}$.
Moreover, by Lemmas \ref{aux:haar} and \ref{lm:quantmodulus} we can in addition conclude that
\begin{equation*}
\sup_{P\in \mathbf P}\|(f - \int f dP) -  f_{n,P}^{d_n\prime} \beta_{n,P}(f)\|_{P,2}^2  \lesssim \frac{\pi_{0n}^2}{d_n^{1/3}} + \frac{\pi_{1n}^2}{d_n^{2/3}},
\end{equation*}
and hence it follows that $\|g\|_{P,2} \leq \delta_n \|G_{n,P}\|_{P,2}$ for all $g \in \mathcal G_{n,P},$ $P\in \mathbf P$, and $\delta_n$ satisfying
\begin{equation*}
\delta_n \asymp \frac{1}{d_n^{1/6}} + \frac{\pi_{1n}}{\pi_{0n} d_n^{1/3}}.
\end{equation*}
Next, note that if $\underline f(V) \leq f(V) \leq \bar f(V)$ almost surely, then result \eqref{lm:quantbootcoup4} yields that
\begin{multline}\label{lm:quantbootcoup7}
\underline f(V) -  \sum_{d=1}^{d_n}E_P[\bar f(V)|V \in \Delta_{d,n}(P)] 1\{V \in \Delta_{d,n}(P)\}  \leq (f(V)-\int f dP) - f_n^{d_n}(V)^\prime \beta_{n,P}( f)  \\  \leq \bar f(V) - \sum_{d=1}^{d_n}E_P[\underline f(V)|V \in \Delta_{d,n}(P)] 1\{V \in \Delta_{d,n}(P)\}
\end{multline}
which implies brackets for $\tilde {\mathcal F}_n$ can be employed to obtain brackets for $\mathcal G_{n,P}$.
Moreover, by the triangle and Jensen's inequality, the width of the brackets built in \eqref{lm:quantbootcoup7} is bounded by $2\|\bar f - \underline f\|_{P,2}$.
Thus, Lemma \ref{aux:simpbracket}, and $\tilde {\mathcal F}_n$ having envelope $\pi_{0n}$ yields
\begin{multline}\label{lm:quantbootcoup9}
\sup_{P\in \mathbf P} \log(N_{[\hspace{0.02 in}]}(\epsilon, \mathcal G_{n,P},\|\cdot\|_{P,2})) \leq \sup_{P\in \mathbf P} \log(N_{[\hspace{0.02 in}]}(\frac{\epsilon}{2},\tilde {\mathcal F}_{n},\|\cdot\|_{P,2})) \\
\leq \log(k_n) + \sup_{P\in \mathbf P} \log(N_{[\hspace{0.02 in}]}(\frac{\epsilon}{2\pi_{0n}},\mathcal F_n,\|\cdot\|_{P,2}))  \lesssim \log(k_n) +  \frac{\pi_{0n}}{\epsilon}1\{\epsilon \leq 2\pi_{0n}\}
\end{multline}
where the final inequality follows for any $\epsilon \leq 2\pi_{0n}$ by Lemma \ref{lm:quantcov}, and for any $\epsilon > 2\pi_{0n}$ by observing that $\mathcal F_n$ is contained in the brackets $[-\tau,1-\tau]$ which has width 1, and hence $N_{[\hspace{0.02 in}]}(1,\mathcal F_n,\|\cdot\|_{P,2}) = 1$.
Recalling that $\mathcal G_{n,P}$ has envelope $G_{n,P} \equiv 2\pi_{0n}$, we can then use result \eqref{lm:quantbootcoup9} to obtain the following upper bound
\begin{multline}\label{lm:quantbootcoup10}
\sup_{P\in \mathbf P} J_{[\hspace{0.02 in}]}(\delta_n \|G_{n,P}\|_{P,2},\mathcal G_{n,P},\|\cdot\|_{P,2}) \lesssim  \int_0^{2\delta_n \pi_{0n}}\sqrt{1+ \log(k_n) + \frac{\pi_{0n}}{\epsilon}1\{\epsilon \leq 2\pi_{0n}\}}d\epsilon \\
\lesssim \delta_n\pi_{0n}\sqrt{\log(1+k_n)} + \int_0^{2\delta_n\pi_{0n}}\sqrt{\frac{\pi_{0n}}{\epsilon}} d\epsilon  \lesssim \pi_{0n} \sqrt{\delta_n},
\end{multline}
where in the final inequality we employed that $\delta_n\log(1+k_n) = o(1)$ by hypothesis.
Similarly, for $\eta_{n,P} \equiv 1 + \log N_{[\hspace{0.02 in}]}(\delta_n\|G_{n,P}\|_{P,2},\mathcal G_{n,P},\|\cdot\|_{P,2})$ we can conclude that
\begin{align}\label{lm:quantbootcoup11}
\sqrt{n} E_P[G_{n,P}(V) \exp\{-\frac{n\delta_n^2 \|G_{n,P}\|^2_{P,2}}{G_{n,P}^2(V) \eta_{n,P}} \}] & \lesssim \sqrt{n}\pi_{0n}\exp\{-\frac{n\delta_n^2}{1+\log(k_n)+ \frac{1}{2\delta_n}1\{\delta_n\leq 1\}}\} \notag \\
& \leq \sqrt n \pi_{0n} \exp\{ - n \delta_n^3\}
\end{align}
where the second inequality holds for $n$ sufficiently large due to $\delta_n \log(1+k_n) = o(1)$.
Together, results \eqref{lm:quantbootcoup10} and \eqref{lm:quantbootcoup11} verify Assumption \ref{ass:4seriesreg}(i) is satisfied with $J_{1n} \asymp \pi_{0n}(\sqrt{\delta_n} + \sqrt{n}\exp\{-n\delta_n^3\})$.
Finally, let $\mathcal B_{n} \equiv \{\beta_{n,P}(f) : f\in \tilde {\mathcal F}_n, P\in \mathbf P\}$ and note \eqref{lm:quantbootcoup3}, $P(\Delta_{i,n}(P)) = 1/d_n$, Jensen's inequality, and $\|f\|_\infty \leq \pi_{0n}$ for any $f\in \tilde {\mathcal F}_n$ imply $\|\beta_{n,P}(f)\|_2 \lesssim \pi_{0n}$ for all $f\in \tilde {\mathcal F}_n$ and $P\in \mathbf P$.
It thus follows that $\mathcal B_n$ is contained in a ball of radius $M \pi_{0n}$ for some $M <\infty$, which allows us to conclude
\begin{multline}\label{lm:quantbootcoup12}
\int_0^\infty \sqrt{N(\epsilon,\mathcal B_n,\|\cdot\|_2)}d\epsilon \lesssim \int_0^{M\pi_{0n}} \sqrt{d_n\log(\frac{M\pi_{0n}}{\epsilon})}d\epsilon \\
= \sqrt{d_n}M\pi_{0n}\int_0^1\sqrt{\log(\frac{1}{u})}du = O(\sqrt{d_n}\pi_{0n}),
\end{multline}
where the first equality follows from the change of variables $u = \epsilon/M\pi_{0n}$.
Result \eqref{lm:quantbootcoup12} verifies Assumption \ref{ass:4seriesreg}(ii) is satisfied with $J_{2n} \asymp \sqrt{d_n}\pi_{0n}$.
In summary, since $d_n^4\log(1+d_n) = o(n)$ by hypothesis, it follows that the conditions of Theorem \ref{th:mainbootcoup}(ii) hold with $C_n \asymp 1$, $K_n \asymp \sqrt{d_n}$, $\xi_n \asymp d_n$, $J_{1n}\asymp \pi_{0n}(\sqrt{\delta_n} + \sqrt n \exp\{-n\delta_n^3\})$, and $J_{2n} \asymp \sqrt{d_n}\pi_{0n}$.
Therefore, Theorem \ref{th:mainbootcoup}(ii) allows us to conclude, uniformly in $P\in \mathbf P$, that
\begin{equation*}
\sup_{f\in \tilde {\mathcal F}_n} |\Wboot(f) - \IsoW(f)| = O_P(\frac{\pi_{0n} d_n^2 \sqrt{\log(1+d_n)}}{\sqrt n} + \pi_{0n}(\sqrt{\delta_n} + \sqrt n\exp\{-n\delta_n^3\})),
\end{equation*}
which together with \eqref{lm:quantbootcoup1} establishes the claim of the Lemma. \qed 

%% file: Appendix/AppLocParam.tex

\section{Local Parameter Space}\label{sec:localparam}

This section contains analytical results concerning our approximation to the local parameter space.
The main result of this section is the following theorem.

\begin{theorem}\label{th:localsmooth}
Let Assumptions \ref{ass:param}(ii)(iii), \ref{ass:locineq}, \ref{ass:loceq}, and \ref{ass:ineqlindep} hold, $\{\ell_n,\delta_n,r_n\}_{n=1}^\infty$ satisfy $\ell_n \downarrow 0$, $\delta_n1\{K_f > 0\} \downarrow 0$, $r_n \geq  2(\ell_n + \delta_n)1\{K_g > 0\}$, $r_n/\delta_n \downarrow 0$, and define
\begin{align*}
G_n(\theta) & \equiv  \{h \in \mathbf B_n : \Upsilon_G(\theta + \frac{h}{\sqrt n}) \leq (\Upsilon_G(\theta) - K_gr_n\|\frac{h}{\sqrt n}\|_{\mathbf B}\mathbf {1_G}) \vee (-r_n\mathbf {1_G}) \} \\ 
A_n(\theta) & \equiv  \{h \in \mathbf B_n : h \in G_n(\theta), ~ \Upsilon_F(\theta + \frac{h}{\sqrt n}) = 0 \text{ and } \|\frac{h}{\sqrt n}\|_{\mathbf B} \leq \ell_n\} \\ 
T_n(\theta) & \equiv \{h \in \mathbf B_n : \Upsilon_F(\theta + \frac{h}{\sqrt n}) = 0, ~ \Upsilon_G(\theta + \frac{h}{\sqrt n}) \leq 0 \text{ and } \|\frac{h}{\sqrt n}\|_{\mathbf B} \leq 2\ell_n\} . 
\end{align*}
(i) Then, there exist $M < \infty$, $\epsilon > 0$, and $n_0 < \infty$ such that for all $n > n_0$, $P\in \mathbf P_0$, $\theta_0 \in \IDsetRsieve$, and $\theta_1\in (\IDsetRsieve)^\epsilon\cap R$ satisfying $\|\theta_0-\theta_1\|_{\mathbf B} \leq \delta_{n}$ we have
\begin{equation}\label{th:localsmoothdisp4}
\sup_{h_1 \in A_n(\theta_1)} \inf_{h_0 \in T_{n}(\theta_0)} \|\frac{h_1}{\sqrt n} - \frac{h_0}{\sqrt n}\|_{\mathbf B} \leq M\times \ell_n(\ell_n+\delta_n) 1\{K_f > 0\} .
\end{equation}
(ii) If in addition $\Upsilon_G$ and $\Upsilon_F$ are affine, then for any $\theta_0,\theta_1\in \mathbf B_n \cap R$ with $\|\theta_0 - \theta_1\|_{\mathbf B} \leq \delta_n$
\begin{multline*}
\{h \in \mathbf B_n : h \in G_n(\theta_1) \text{ and }  \Upsilon_F(\theta_1 + \frac{h}{\sqrt n}) = 0\} \\
\subseteq \{h \in \mathbf B_n : \Upsilon_G(\theta_0 + \frac{h}{\sqrt n}) \leq 0 \text{ and }  \Upsilon_F(\theta_0 + \frac{h}{\sqrt n}) =0\}.
\end{multline*}
\end{theorem}

\noindent {\sc Proof:} We begin by establishing part (ii).
First note that if $\Upsilon_G$ is affine, then $K_g = 0$ and since $r_n/\delta_n = o(1)$, Lemma \ref{lm:localineq}(ii) implies that for $n$ sufficiently large
\begin{equation}\label{th:localsmooth0}
G_n(\theta_1) \subseteq \{h \in \mathbf B_n : \Upsilon_G(\theta_0 + \frac{h}{\sqrt n}) \leq 0\}
\end{equation}
for any $\theta_0,\theta_1\in \mathbf B_n$ with $\|\theta_0-\theta_1\|_{\mathbf B} \leq \delta_n$.
Moreover, if $\Upsilon_F$ is affine and continuous, then $\Upsilon_F(\theta) = L(\theta) + c_0$ for some continuous linear map $L : \mathbf B \to \mathbf F$ and $c_0 \in \mathbf F$.
It follows that $\nabla \Upsilon_F(\theta)[h] = L(h)$, which does not depend on $\theta$, and since any $\theta \in R$ must satisfy $L(\theta) = -c_0$ (since $\Upsilon_F(\theta) = 0$), we can conclude that $\{h : \Upsilon_F(\theta + h) = 0\} = \{h : L(h) = 0\}$ whenever $\theta \in R$.
Therefore part (ii) follows from result \eqref{th:localsmooth0} and $\theta_1,\theta_2\in R$.

We next turn to the proof of part (i).
Throughout, let $\tilde \epsilon$ be such that Assumptions \ref{ass:locineq} and \ref{ass:loceq} hold and set $\epsilon = \tilde \epsilon/2$.
We break up the proof into four steps.

\noindent \underline{\sc Step 1:} (Decompose $h/\sqrt n$). For any $P\in \mathbf P_0$, $\theta_0 \in \IDsetRsieve$, and $h\in \mathbf B_n$ set
\begin{equation}\label{th:localsmooth1}
h^{\perp_{\theta_0}} \equiv \nabla \Upsilon_F(\theta_0)^{-} \nabla \Upsilon_F(\theta_0)[h] \hspace{0.8 in} h^{\mathcal N_{\theta_0}} \equiv h - h^{\perp_{\theta_0}} ,
\end{equation}
where recall $\nabla \Upsilon_F(\theta_0)^-:\mathbf F_n\rightarrow \mathbf B_n$ denotes the right inverse of $\nabla \Upsilon_F(\theta_0):\mathbf B_n \rightarrow \mathbf F_n$.
Further note that $h^{\mathcal N_{\theta_0}} \in \mathcal N(\nabla \Upsilon_F(\theta_0))$ since $\nabla \Upsilon_F(\theta_0)\nabla \Upsilon_F(\theta_0)^{-} =I$ implies that
\begin{equation}\label{th:localsmooth2}
\nabla \Upsilon_F(\theta_0)[h^{\mathcal N_{\theta_0}}] = \nabla \Upsilon_F(\theta_0)[h] - \nabla\Upsilon_F(\theta_0)\nabla \Upsilon_F(\theta_0)^{-} \nabla \Upsilon_F(\theta_0)[h] = 0 ,
\end{equation}
by definition of $h^{\perp_{\theta_0}}$ in \eqref{th:localsmooth1}.
Next, observe that if $\theta_1 \in (\IDsetRsieve)^\epsilon \cap R$ and $h \in \mathbf B_n$ satisfies $\|h/\sqrt n\|_{\mathbf B}\leq \ell_n$ and $\Upsilon_F(\theta_1 + h/\sqrt n) = 0$, then $\theta_1 + h/\sqrt n \in (\IDsetRsieve)^{\tilde \epsilon}$ for $n$ sufficiently large, and hence by Assumption \ref{ass:loceq}(i) and $\Upsilon_F(\theta_1) = 0$ due to $\theta_1 \in R$
\begin{equation}\label{th:localsmooth3}
\|\nabla \Upsilon_F(\theta_1)[\frac{h}{\sqrt n}]\|_{\mathbf F} = \|\Upsilon_F(\theta_1+ \frac{h}{\sqrt n}) - \Upsilon_F(\theta_1) - \nabla \Upsilon_F(\theta_1)[\frac{h}{\sqrt n}]\|_{\mathbf F} \leq K_f\|\frac{h}{\sqrt n}\|_{\mathbf B}^2 .
\end{equation}
Hence, Assumption \ref{ass:loceq}(ii), result \eqref{th:localsmooth3}, $\|\theta_0 - \theta_1\|_{\mathbf B} \leq \delta_n$, and $\|h/\sqrt n\|_{\mathbf B}\leq \ell_n$ imply
\begin{multline}\label{th:localsmooth4}
\|\nabla \Upsilon_F(\theta_0)[\frac{h}{\sqrt n}]\|_{\mathbf F} \\ \leq \|\nabla \Upsilon_F(\theta_0) - \nabla \Upsilon_F(\theta_1)\|_o \|\frac{h}{\sqrt n}\|_{\mathbf B} + K_f\|\frac{h}{\sqrt n}\|_{\mathbf B}^2 \leq K_f \ell_n(\delta_{n} + \ell_n) .
\end{multline}
Moreover, since $\nabla \Upsilon_F(\theta_0):\mathbf F_n \rightarrow \mathbf B_n$ satisfies Assumption \ref{ass:loceq}(iv), we also have that
\begin{multline}\label{th:localsmooth5}
K_f\|h^{\perp_{\theta_0}}\|_{\mathbf B} = K_f\|\nabla \Upsilon_F(\theta_0)^{-} \nabla\Upsilon(\theta_0)[h]\|_{\mathbf B} \\ \leq K_f\|\nabla \Upsilon_F(\theta_0)^{-}\|_o \|\nabla \Upsilon_F(\theta_0)[h]\|_{\mathbf F} \leq M_f\|\nabla \Upsilon_F(\theta_0)[h]\|_{\mathbf F} 
\end{multline}
for some $M_f < \infty$. 
Further note that if $K_f = 0$, then \eqref{th:localsmooth1} and \eqref{th:localsmooth4} imply $h^{\perp_{\theta_0}} = 0$.
Thus, combining results \eqref{th:localsmooth4} and \eqref{th:localsmooth5} to handle the case $K_f > 0$ we conclude that for any $P\in \mathbf P_0$, $\theta_0\in \IDsetRsieve$, $\theta_1\in (\IDsetRsieve)^\epsilon\cap R$ satisfying $\|\theta_0 - \theta_1\|_{\mathbf B}\leq \delta_n$ and any $h\in \mathbf B_n$ such that $\Upsilon_F(\theta_1 + h/\sqrt n) = 0$ and $\|h/\sqrt n\|_{\mathbf B} \leq \ell_n$ we must have 
\begin{equation}\label{th:localsmooth6}
\|\frac{h^{\perp_{\theta_0}}}{\sqrt n}\|_{\mathbf B} \leq  M_f\ell_n(\delta_{n} + \ell_n)1\{K_f > 0\} .
\end{equation}

\noindent \underline{\sc Step 2:} (Inequality Constraints).
In what follows, it is convenient to define the set
\begin{equation*}
S_n(\theta_0,\theta_1) \equiv \{ h \in \mathbf B_n : \Upsilon_G(\theta_0 + \frac{h}{\sqrt n}) \leq 0, ~ \Upsilon_F(\theta_1 + \frac{h}{\sqrt n}) = 0, ~ \|\frac{h}{\sqrt n}\|_{\mathbf B} \leq \ell_n\} .
\end{equation*}
Then note $r_n \geq  2(\ell_n + \delta_n)1\{K_g > 0\}$, $r_n/\delta_n = o(1)$, and Lemma \ref{lm:localineq}(i) imply that
\begin{equation}\label{th:localsmooth8}
A_n(\theta_1) \subseteq S_n(\theta_0,\theta_1)
\end{equation}
for $n$ sufficiently large, all $P\in \mathbf P_0$, $\theta_0 \in \IDsetRsieve$, and $\theta_1 \in (\IDsetRsieve)^\epsilon$ satisfying $\|\theta_0 - \theta_1\|_{\mathbf B} \leq \delta_{n}$.
The proof will proceed by verifying \eqref{th:localsmoothdisp4} holds with $S_n(\theta_0,\theta_1)$ in place of $A_n(\theta_1)$.
In particular, if $K_f = 0$, then $\Upsilon_F(\theta_0) = \Upsilon_F(\theta_1)$ due to $\theta_0,\theta_1 \in R$, and Assumptions \ref{ass:loceq}(i)(ii) together with \eqref{th:localsmooth8} imply $A_n(\theta_1)\subseteq S_n(\theta_0,\theta_1) \subseteq T_n(\theta_0)$.
Hence, result \eqref{th:localsmoothdisp4}  holds for the case $K_f = 0$.

For the rest of the proof we therefore assume $K_f > 0$. We further note that Lemma \ref{lm:intpert} implies that for any $\eta_n \downarrow 0$, there is an $n_0 <\infty$ and $1\leq C < \infty$ (independent of $\eta_n$) such that for all $P\in \mathbf P_0$, $n > n_0$, and $\theta_0 \in \IDsetRsieve$ there exists a $h_{\theta_0,n} \in \mathbf B_n \cap \mathcal N(\nabla \Upsilon_F(\theta_0))$ such that for any $\tilde h \in \mathbf B_n$ for which there exists a $h\in S_n(\theta_0,\theta_1)$ satisfying $\|(\tilde h - h)/\sqrt n\|_{\mathbf B} \leq \eta_n$ the following inequalities hold
\begin{equation}\label{th:localsmooth9}
\Upsilon_G(\theta_0 + \frac{h_{\theta_0,n}}{\sqrt n} + \frac{\tilde h}{\sqrt n}) \leq 0 \hspace{0.5 in} \|\frac{h_{\theta_0,n}}{\sqrt n}\|_{\mathbf B} \leq C\eta_n.
\end{equation}

\noindent \underline{\sc Step 3:} (Equality Constraints). The results in this step allow us to address the challenge that $h\in S_n(\theta_0,\theta_1)$ satisfies $\Upsilon_F(\theta_1 + h/\sqrt n) = 0$ but not necessarily $\Upsilon_F(\theta_0 + h/\sqrt n) = 0$.
To this end, let $\mathcal R(\nabla \Upsilon_F(\theta_0)^{-}\nabla \Upsilon_F(\theta_0))$ denote the range of the operator $\nabla \Upsilon_F(\theta_0)^{-}\nabla \Upsilon_F(\theta_0):\mathbf B_n \rightarrow \mathbf B_n$ and define the vector subspaces
\begin{equation}\label{th:localsmooth10}
\mathbf B_n^{\mathcal N_{\theta_0}} \equiv \mathbf B_n\cap \mathcal N(\nabla \Upsilon_F(\theta_0)) \hspace{0.5 in} \mathbf B_n^{\perp_{\theta_0}} \equiv \mathcal R(\nabla \Upsilon_F(\theta_0)^{-}\nabla \Upsilon_F(\theta_0)) .
\end{equation}
Since $h^{\mathcal N_{\theta_0}} \in \mathbf B_n^{\mathcal N_{\theta_0}}$ by \eqref{th:localsmooth2}, the definitions in \eqref{th:localsmooth1} and \eqref{th:localsmooth10} imply that $\mathbf B_n = \mathbf B_n^{\mathcal N_{\theta_0}} + \mathbf B_n^{\perp_{\theta_0}}$. Furthermore, since $\nabla \Upsilon_F(\theta_0)\nabla \Upsilon_F(\theta_0)^{-} = I$, we also have
\begin{equation}\label{th:localsmooth11}
\nabla \Upsilon_F(\theta_0)^{-} \nabla \Upsilon_F(\theta_0)[h] = h
\end{equation}
for any $h \in \mathbf B_n^{\perp_{\theta_0}}$, and thus that $\mathbf B_n^{\mathcal N_{\theta_0}} \cap \mathbf B_n^{\perp_{\theta_0}} = \{0\}$.
Since $\mathbf B_n = \mathbf B_n^{\mathcal N_{\theta_0}} + \mathbf B_n^{\perp_{\theta_0}}$, it then follows that $\mathbf B_n = \mathbf B_n^{\mathcal N_{\theta_0}} \oplus \mathbf B_n^{\perp_{\theta_0}}$ -- i.e. the decomposition in \eqref{th:localsmooth1} is unique.
Moreover, we observe that $\mathbf B_n^{\mathcal N_{\theta_0}} \cap \mathbf B_n^{\perp_{\theta_0}} = \{0\}$ further implies the restricted map $\nabla \Upsilon_F(\theta_0) : \mathbf B_n^{\perp_{\theta_0}} \rightarrow \mathbf F_n$ is in fact bijective, and by \eqref{th:localsmooth11} its inverse is $\nabla \Upsilon_F(\theta_0)^{-} : \mathbf F_n \rightarrow \mathbf B_n^{\perp_{\theta_0}}$.

Recall $\Upsilon_F$ is Fr\'echet differentiable on $ (\IDsetRsieve)^{\tilde \epsilon}$ by Assumption \ref{ass:loceq}(i).
The preceding discussion and Assumption \ref{ass:loceq} imply we may apply Lemma \ref{lm:auximplicit} with $\mathbf A_1 = \mathbf B_n^{\mathcal N_{\theta_0}}$, $\mathbf A_2 = \mathbf B_n^{\perp_{\theta_0}}$, and some $K_0 < \infty$ to obtain that for any $P \in \mathbf P_0$, $\theta_0 \in \IDsetRsieve$ and $h^{\mathcal N_{\theta_0}} \in \mathbf B^{\mathcal N_{\theta_0}}_n$ satisfying $\|h^{\mathcal N_{\theta_0}}\|_{\mathbf B} \leq \{\tilde \epsilon/2 \wedge (2K_0)^{-2}\wedge 1\}^2$, there is a $h^\star(h^{\mathcal N_{\theta_0}}) \in \mathbf B_{n}^{\perp_{\theta_0}}$ such that
\begin{equation}\label{th:localsmooth12}
\Upsilon_F(\theta_0 + h^{\mathcal N_{\theta_0}} + h^\star(h^{\mathcal N_{\theta_0}})) = 0 \hspace{0.5 in} \|h^\star(h^{\mathcal N_{\theta_0}})\|_{\mathbf B} \leq 2K_0^2\|h^{\mathcal N_{\theta_0}}\|_{\mathbf B}^2 .
\end{equation}
Moreover, for any $P \in \mathbf P_0$, $\theta_0 \in \IDsetRsieve$, $\theta_1 \in (\IDsetRsieve)^\epsilon \cap R$ with $\|\theta_0 - \theta_1\|_{\mathbf B}\leq \delta_n$, and any $h\in \mathbf B_n$ such that $\Upsilon_F(\theta_1 + h/\sqrt n) = 0$ and $\|h/\sqrt n\|_{\mathbf B}\leq \ell_n$, result \eqref{th:localsmooth6}, the decomposition in \eqref{th:localsmooth1}, $\delta_n \downarrow 0$ (since $K_f > 0$), and $\ell_n \downarrow 0$ imply that for $n$ large
\begin{equation}\label{th:localsmooth13}
\|\frac{h^{\mathcal N_{\theta_0}}}{\sqrt n}\|_{\mathbf B} \leq \|\frac{h}{\sqrt n}\|_{\mathbf B} + \|\frac{h^{\perp_{\theta_0}}}{\sqrt n}\|_{\mathbf B} \leq 2 \ell_n .
\end{equation}
Thus, for $h_{\theta_{0,n}} \in \mathbf B_n^{\mathcal N_{\theta_0}}$ as in \eqref{th:localsmooth9}, $C\geq 1$, and results \eqref{th:localsmooth12} and \eqref{th:localsmooth13} imply that for $n$ sufficiently large we must have for all $P\in \mathbf P_0$, $\theta_0 \in \IDsetRsieve$, $\theta_1\in \mathbf B_n \cap R$ with $\|\theta_0 - \theta_1\|_{\mathbf B}\leq \delta_n$ and $h \in \mathbf B_n$ satisfying $\Upsilon_F(\theta_1 + h/\sqrt n) = 0$ that
\begin{align}
\Upsilon_F(\theta_0 + \frac{h_{\theta_0,n}}{\sqrt n}+ \frac{h^{\mathcal N_{\theta_0}}}{\sqrt n}+ h^\star(\frac{h_{\theta_0,n}}{\sqrt n}+ \frac{h^{\mathcal N_{\theta_0}}}{\sqrt n})) & = 0 \notag \\ \|h^\star(\frac{h_{\theta_0,n}}{\sqrt n}+ \frac{h^{\mathcal N_{\theta_0}}}{\sqrt n})\|_{\mathbf B} - 16K_0^2C^2(\ell_n^2 + \eta_n^2) & \leq 0 \label{th:localsmooth15}.
\end{align}

\noindent \underline{Step 4:} (Build Approximation). In order to employ Steps 2 and 3, we now set $\eta_n$ to
\begin{equation}\label{th:localsmooth16}
\eta_n = 32(M_f + C^2K_0^2)\ell_n(\ell_n + \delta_n).
\end{equation}
In addition, for any $P\in \mathbf P_0$, $\theta_0 \in \IDsetRsieve$, $\theta_1 \in \mathbf B_n\cap R$ satisfying $\|\theta_0 - \theta_1\|_{\mathbf B} \leq \delta_n$, and any $h \in S_n(\theta_0,\theta_1)$, we let $h^{\mathcal N_{\theta_0}}$ be as in \eqref{th:localsmooth1} and define
\begin{equation}\label{th:localsmooth17}
\frac{\hat h }{\sqrt n} \equiv \frac{h_{\theta_0,n}}{\sqrt n} + \frac{h^{\mathcal N_{\theta_0}}}{\sqrt n} + h^\star(\frac{h_{\theta_0,n}}{\sqrt n} + \frac{h^{\mathcal N_{\theta_0}}}{\sqrt n}) .
\end{equation}
From Steps 2 and 3 it then follows that for $n$ sufficiently large (independent of $P\in \mathbf P_0$, $\theta_0 \in \IDsetRsieve$, $\theta_1 \in \mathbf B_n \cap R$ with $\|\theta_0-\theta_1\|_{\mathbf B}\leq \delta_n$ or $h \in S_n(\theta_0,\theta_1)$) we have
\begin{equation}\label{th:localsmooth18}
\Upsilon_F(\theta_0 + \frac{\hat h}{\sqrt n}) = 0 .
\end{equation}
Moreover, from results \eqref{th:localsmooth15} and \eqref{th:localsmooth16} we also obtain that for $n$ sufficiently large
\begin{equation}\label{th:localsmooth19}
\|h^\star(\frac{h_{\theta_0,n}}{\sqrt n} + \frac{h^{\mathcal N_{\theta_0}}}{\sqrt n})\|_{\mathbf B }\leq 16C^2K_0^2(\ell_n^2 + \eta_n^2) \leq \frac{\eta_n}{2} + 16C^2K^2_0 \eta_n^2 \leq \frac{3}{4} \eta_n .
\end{equation}
Thus, $h = h^{\mathcal N_{\theta_0}} + h^{\perp_{\theta_0}}$, \eqref{th:localsmooth6}, \eqref{th:localsmooth16}, \eqref{th:localsmooth17} and \eqref{th:localsmooth19} imply for large $n$ that $\|(\hat h - h - h_{\theta_0,n})/\sqrt n\|_{\mathbf B} \leq \eta_n$, and employing \eqref{th:localsmooth9} with $\tilde h = (\hat h - h_{\theta_0,n})/\sqrt n)$ yields
\begin{equation}\label{th:localsmooth20}
\Upsilon_G(\theta_0 + \frac{\hat h}{\sqrt n}) \leq 0  .
\end{equation}
Since $\|h_{\theta_0,n}/\sqrt n\|_{\mathbf B} \leq C\eta_n$ by \eqref{th:localsmooth9}, results \eqref{th:localsmooth6}, \eqref{th:localsmooth15}, and $\|h/\sqrt n\|_{\mathbf B} \leq \ell_n$ for any $h \in S_n(\theta_0,\theta_1)$ imply by \eqref{th:localsmooth16} and $\ell_n \downarrow 0$, $\delta_n\downarrow 0$ that
\begin{multline}\label{th:localsmooth21}
\|\frac{\hat h}{\sqrt n}\|_{\mathbf B} \leq \|\frac{h_{\theta_0,n}}{\sqrt n}\|_{\mathbf B} + \|h^\star(\frac{h_{\theta_0,n}}{\sqrt n} + \frac{h^{\mathcal N_{\theta_0}}}{\sqrt n})\|_{\mathbf B} + \|\frac{h^{\perp_{\theta_0}}}{\sqrt n}\|_{\mathbf B} + \|\frac{h}{\sqrt n}\|_{\mathbf B} \\ \leq C\eta_n + 16C^2K_0^2(\ell_n^2 + \eta_n^2) + M_f\ell_n(\delta_n + \ell_n) +\ell_n \leq 2\ell_n
\end{multline}
for $n$ sufficiently large.
Therefore, we conclude from \eqref{th:localsmooth18}, \eqref{th:localsmooth20}, and \eqref{th:localsmooth21} that $\hat h \in T_n(\theta_0)$.
Similarly, \eqref{th:localsmooth6}, \eqref{th:localsmooth9}, \eqref{th:localsmooth15}, and \eqref{th:localsmooth16} yield for some $M < \infty$
\begin{multline*} 
\|\frac{\hat h}{\sqrt n} - \frac{h}{\sqrt n}\|_{\mathbf B} \leq \|\frac{h_{\theta_0,n}}{\sqrt n}\|_{\mathbf B} + \|h^\star(\frac{h_{\theta_0,n}}{\sqrt n} + \frac{h^{\mathcal N_{\theta_0}}}{\sqrt n})\|_{\mathbf B} + \|\frac{h^{\perp_{\theta_0}}}{\sqrt n}\|_{\mathbf B} \\ \leq C\eta_n + 16C^2K_0^2(\ell_n^2 + \eta_n^2) + M_f\ell_n(\ell_n + \delta_n) \leq M\ell_n(\ell_n + \delta_n) ,
\end{multline*}
which establishes the \eqref{th:localsmoothdisp4} for the case $K_f > 0$. \qed

\begin{lemma}\label{lm:localineq}
Let Assumptions \ref{ass:param}(ii)(iii), \ref{ass:locineq} hold, and $\ell_n \downarrow 0$ be given.
(i) Then, there are $n_0, M_g < \infty$ and $\epsilon > 0$ such that for all $n > n_0$, $P\in \mathbf P_0$, $\theta_0 \in \IDsetRsieve$,  $\theta_1 \in (\IDsetRsieve)^\epsilon$: 
\begin{multline*}
\{h \in \mathbf B_n : \Upsilon_G(\theta_1 + \frac{h}{\sqrt n}) \leq (\Upsilon_G(\theta_1) - K_gr\|\frac{h}{\sqrt n}\|_{\mathbf B}\mathbf {1_G}) \vee (-r\mathbf {1_G}) \text{ and } \|\frac{h}{\sqrt n}\|_{\mathbf B} \leq \ell_n \} \\
\subseteq \{ h \in \mathbf B_n : \Upsilon_G(\theta_0 + \frac{h}{\sqrt n}) \leq 0 \text{ and } \|\frac{h}{\sqrt n}\|_{\mathbf B} \leq \ell_n \}
\end{multline*}
for any $r \geq \{M_g\|\theta_0 - \theta_1\|_{\mathbf B} + K_g \|\theta_0 - \theta_1\|_{\mathbf B}^2\}\vee 2\{\ell_n + \|\theta_0 - \theta_1\|_{\mathbf B}\}1\{K_g > 0\}$.
(ii) If in addition $\Upsilon_G$ is affine, then for any $n$, $\theta_0,\theta_1\in \mathbf B_n$, and $r\geq M_g\|\theta_0-\theta_1\|_{\mathbf B}$ we have
$$\{h \in \mathbf B_n : \Upsilon_G(\theta_1+\frac{h}{\sqrt n}) \leq \Upsilon_G(\theta_1)\vee (-r\mathbf{1_G}) \} \subseteq \{h \in \mathbf B_n : \Upsilon_G(\theta_0 + \frac{h}{\sqrt n}) \leq 0\}.$$
\end{lemma}

\noindent {\sc Proof:} Let $\tilde \epsilon > 0$ be such that Assumption \ref{ass:locineq} holds and set $M_g <\infty$ to satisfy
\begin{equation}\label{lm:localineq0}
\|\nabla \Upsilon_G(\theta)\|_o \leq M_g
\end{equation}
for any $\theta \in (\IDsetRsieve)^{\tilde \epsilon}$, which is possible by Assumption \ref{ass:locineq}(iii).
Next, set $\epsilon = \tilde \epsilon /2$ and define $N(\delta) \equiv \{\theta \in \mathbf B_n : \overrightarrow d_H(\{\theta\},\IDsetRsieve,\|\cdot\|_{\mathbf B}) < \delta\}$ for any $\delta > 0$.
Then note that for any $\theta_1 \in N(\epsilon)$ and $\|h/\sqrt n\|_{\mathbf B} \leq \ell_n$ we have $\theta_1 + h/\sqrt n \in N(\tilde \epsilon)$ for $n$ sufficiently large.
Therefore, Assumption \ref{ass:locineq}(i) allows us to conclude that
\begin{equation}\label{lm:localineq1}
\|\Upsilon_G(\theta_1 + \frac{h}{\sqrt n}) - \Upsilon_G(\theta_1) -\nabla \Upsilon_G(\theta_1)[\frac{h}{\sqrt n}]\|_{\mathbf G} \leq K_g\|\frac{h}{\sqrt n}\|_{\mathbf B}^2 .
\end{equation}
Similarly, Assumption \ref{ass:locineq}(ii) implies that if $\theta_0 \in \IDsetRsieve$ and $\theta_1 \in N(\epsilon)$, then we have
\begin{multline}\label{lm:localineq2}
\|\nabla \Upsilon_G(\theta_0)[\frac{h}{\sqrt n}] - \nabla \Upsilon_G(\theta_1)[\frac{h}{\sqrt n}]\|_{\mathbf G} \\ \leq \|\nabla \Upsilon_G(\theta_0) - \nabla \Upsilon_G(\theta_1)\|_o \|\frac{h}{\sqrt n}\|_{\mathbf B} \leq K_g\|\theta_0 - \theta_1\|_{\mathbf B} \|\frac{h}{\sqrt n}\|_{\mathbf B}
\end{multline}
for any $h\in \mathbf B_n$.
Hence, since $\Upsilon_G(\theta_0) \leq 0$ due to $\theta_0 \in \IDsetRsieve\subseteq \Theta_n \cap R$ we obtain that
\begin{align}\label{lm:localineq3}
\Upsilon_G(\theta_0 + \frac{h}{\sqrt n}) & + \{\Upsilon_G(\theta_1) - \Upsilon_G(\theta_1 + \frac{h}{\sqrt n})\} \nonumber \\
& \leq \{\Upsilon_G(\theta_0 + \frac{h}{\sqrt n}) - \Upsilon_G(\theta_0)\} + \{\Upsilon_G(\theta_1) - \Upsilon_G(\theta_1 + \frac{h}{\sqrt n})\} \nonumber \\ & \leq K_g\|\frac{h}{\sqrt n}\|_{\mathbf B}\{2\|\frac{h}{\sqrt n}\|_{\mathbf B} + \|\theta_0 -\theta_1\|_{\mathbf B} \}\mathbf {1_G} ,
\end{align}
by \eqref{lm:localineq1}, \eqref{lm:localineq2}, and Lemma \ref{lm:obviousam}.
Also note for any $\theta_0 \in \IDsetRsieve$, $\theta_1 \in N(\epsilon)$, and $h\in \mathbf B_n$ with $\|h/\sqrt n\|_{\mathbf B} \leq \ell_n$ we have $\theta_0 + h/\sqrt n \in N(\tilde \epsilon)$ and $\theta_1 + h/\sqrt n \in N(\tilde \epsilon)$ for $n$ sufficiently large.
Thus, by Assumptions \ref{ass:locineq}(i), result \eqref{lm:localineq0}, and Lemma \ref{lm:obviousam}
\begin{align}\label{lm:localineq4}
\Upsilon_G(\theta_0 + \frac{h}{\sqrt n}) - \Upsilon_G(\theta_1 + \frac{h}{\sqrt n})
& \leq \nabla \Upsilon_G(\theta_0 + \frac{h}{\sqrt n})[\theta_0 - \theta_1] + K_g\|\theta_0 - \theta_1\|^2_{\mathbf B} \mathbf {1_G} \nonumber \\ & \leq \{M_g\|\theta_0 - \theta_1\|_{\mathbf B} + K_g\|\theta_0 - \theta_1\|_{\mathbf B}^2\}\mathbf {1_G} .
\end{align}
Hence, \eqref{lm:localineq3} and \eqref{lm:localineq4} yield for $r \geq \{M_g\|\theta_0 - \theta_1\|_{\mathbf B} + K_g \|\theta_0 - \theta_1\|_{\mathbf B}^2\}\vee 2\{\ell_n + \|\theta_0 - \theta_1\|_{\mathbf B}\}1\{K_g >0\}$, $\theta_0\in \IDsetRsieve$, $\theta_1\in N(\epsilon)$, $\|h/\sqrt n\|_{\mathbf B} \leq \ell_n$, and $n$ large
\begin{align}\label{lm:localineq5}
\Upsilon_G (\theta_0 + \frac{h}{\sqrt n}) & \leq \Upsilon_G(\theta_1 + \frac{h}{\sqrt n}) + (K_gr\|\frac{h}{\sqrt n}\|_{\mathbf B} - \Upsilon_G(\theta_1))\mathbf {1_G} \wedge r \mathbf {1_G} \nonumber \\
& = \Upsilon_G(\theta_1 + \frac{h}{\sqrt n}) - (\Upsilon_G(\theta_1)- K_gr\|\frac{h}{\sqrt n}\|_{\mathbf B})\mathbf {1_G}\vee (-r \mathbf {1_G})
\end{align}
where the equality follows from $(-a)\vee(-b) = - (a\wedge b)$ by Theorem 8.6 in \cite{aliprantis:border:2006}. Since $a_1\leq a_2$ and $b_1 \leq b_2$ implies $a_1 \wedge b_1 \leq a_2 \wedge b_2$ in $\mathbf G$ by Corollary 8.7 in \cite{aliprantis:border:2006}, \eqref{lm:localineq5} implies that for $n$ sufficiently large and any $\theta_0\in \IDsetRsieve$, $\theta_1 \in N(\epsilon)$ and $h \in \mathbf B_n$ satisfying $\|h/\sqrt n\|_{\mathbf B} \leq \ell_n$ and
\begin{equation*}
\Upsilon_G(\theta_1 + \frac{h}{\sqrt n}) \leq (\Upsilon_G(\theta_1) - K_gr\|\frac{h}{\sqrt n}\|_{\mathbf B}\mathbf {1_G}) \vee (-r\mathbf {1_G})
\end{equation*}
we must have $\Upsilon_G(\theta_0 + h/\sqrt n) \leq 0$, which verifies the first claim of the lemma.
For the second claim, just note that if $\Upsilon_G$ is affine, then we may set $K_g = 0$ and $\epsilon = +\infty$ in Assumption \ref{ass:locineq}, which leads to the desired simplification. \qed

\begin{lemma}\label{lm:intpert}
If Assumptions \ref{ass:param}(ii)(iii), \ref{ass:locineq}, \ref{ass:ineqlindep}(ii) hold, and $\eta_n \downarrow 0$, $\ell_n \downarrow 0$, then there is a $n_0$ (depending on $\eta_n,\ell_n$) and a $C < \infty$ (independent of $\eta_n,\ell_n$) such that for all $n > n_0$, $P\in \mathbf P_0$, and $\theta \in \IDsetRsieve$ there is $h_{\theta,n} \in \mathbf B_n\cap \mathcal N(\nabla \Upsilon_F(\theta))$ with
\begin{equation}\label{lm:intpertdisp}
\Upsilon_G(\theta + \frac{h_{\theta,n}}{\sqrt n} + \frac{\tilde h}{\sqrt n}) \leq 0 \hspace{0.7 in} \|\frac{h_{\theta,n}}{\sqrt n}\|_{\mathbf B} \leq C\eta_n
\end{equation}
for all $\tilde h \in \mathbf B_n$ for which there is a $h \in \mathbf B_n$ satisfying $\|(\tilde h-h)/\sqrt n\|_{\mathbf B} \leq \eta_n$, $\|h/\sqrt n\|_{\mathbf B} \leq \ell_n$, and the inequality $\Upsilon_G(\theta + h/\sqrt n)\leq 0$.
\end{lemma}

\noindent {\sc Proof:} By Assumption \ref{ass:ineqlindep}(ii) there are $\epsilon > 0$ and $M_d < \infty$ such that for every $P\in \mathbf P_0$, $n$, and $\theta \in \IDsetRsieve$ there exists a $\bar h_{\theta,n} \in \mathbf B_n \cap \mathcal N(\nabla \Upsilon_F(\theta))$ satisfying
\begin{equation}\label{lm:intpert1}
\Upsilon_G(\theta) + \nabla \Upsilon_G(\theta)[\bar h_{\theta,n}] \leq -\epsilon \mathbf {1_G} \hspace{0.5 in} \|\bar h_{\theta,n}\|_{\mathbf B} \leq M_d.
\end{equation}
Also note Assumption \ref{ass:locineq}(iii) and $\ell_n = o(1)$ imply that there is an $M_g < \infty$ such that for $n$ sufficiently large and any $h\in \mathbf B_n$ satisfying $\|h/\sqrt n\|_{\mathbf B}\leq \ell_n$ we must have
\begin{equation}\label{lm:intpert1p5}
\|\nabla \Upsilon_G(\theta + \frac{h}{\sqrt n})\|_o \leq M_g.
\end{equation}
Moreover, result \eqref{lm:intpert1p5}, Assumption \ref{ass:locineq}(i), Lemma \ref{lm:obviousam}, and $\ell_n = o(1)$ imply that for $n$ sufficiently large and any $h \in \mathbf B_n$ with $\|h/\sqrt n\|_{\mathbf B} \leq \ell_n$ we must have
\begin{multline}\label{lm:intpert2}
\Upsilon_G(\theta + \frac{h}{\sqrt n}) \leq \Upsilon_G(\theta) + \nabla \Upsilon_G(\theta)[\frac{h}{\sqrt n}] + K_g\|\frac{h}{\sqrt n}\|_{\mathbf B}^2 \mathbf {1_G} \\
\leq \Upsilon_G(\theta) + \{\|\nabla \Upsilon_G(\theta)\|_o\ell_n + K_g\ell_n^2\} \mathbf {1_G}
\leq \Upsilon_G(\theta) + 2M_g\ell_n\mathbf {1_G}  .
\end{multline}
Hence, \eqref{lm:intpert1} and \eqref{lm:intpert2} imply for any $h \in \mathbf B_n$ with $\|h/\sqrt n\|_{\mathbf B}\leq \ell_n$ we must have
\begin{equation}\label{lm:intpert3}
\Upsilon_G(\theta + \frac{h}{\sqrt n}) + \nabla \Upsilon_G(\theta)[\bar h_{\theta,n}] \leq \{2M_g\ell_n -\epsilon\} \mathbf {1_G} .
\end{equation}

Next, we let $C_0 > 8M_g/\epsilon$ and aim to show \eqref{lm:intpertdisp} holds with $C = C_0M_d$ by setting
\begin{equation}\label{lm:intpert4}
\frac{h_{\theta,n}}{\sqrt n} \equiv C_0\eta_n\bar h_{\theta,n} .
\end{equation}
To this end, we first note that if $\theta \in \IDsetRsieve$, $h \in \mathbf B_n$ satisfies $\|h/\sqrt n\|_{\mathbf B} \leq \ell_n$ and $\Upsilon_G(\theta + h/\sqrt n) \leq 0$, and $\tilde h\in \mathbf B_n$ is such that $\|(h-\tilde h)/\sqrt n\|_{\mathbf B} \leq \eta_n$, then definition \eqref{lm:intpert4} implies that $\|(h_{\theta,n} + \tilde h)/\sqrt n\|_{\mathbf B} = o(1)$.
Therefore, Assumption \ref{ass:locineq}(i), Lemma \ref{lm:obviousam}, and $\|(\tilde h - h)/\sqrt n\|_{\mathbf B}\leq \eta_n$ together allow us to conclude that
\begin{align}\label{lm:intpert5}
\Upsilon_G& (\theta + \frac{h_{\theta,n}}{\sqrt n} + \frac{\tilde h}{\sqrt n})\nonumber  \\
&\leq \Upsilon_G(\theta + \frac{h}{\sqrt n}) + \nabla \Upsilon_G(\theta + \frac{h}{\sqrt n})[\frac{h_{\theta,n}}{\sqrt n} + \frac{(\tilde h - h)}{\sqrt n}] + 2K_g(\|\frac{h_{\theta,n}}{\sqrt n}\|_{\mathbf B}^2 + \eta_n^2)\mathbf {1_G} \nonumber \\ & \leq \Upsilon_G(\theta + \frac{h}{\sqrt n}) + \nabla \Upsilon_G(\theta + \frac{h}{\sqrt n})[\frac{h_{\theta,n}}{\sqrt n}] + \{2K_g\|\frac{h_{\theta,n}}{\sqrt n}\|_{\mathbf B}^2 + 2M_g\eta_n\}\mathbf {1_G} ,
\end{align}
where the final result follows from result \eqref{lm:intpert1p5} and $2K_g\eta_n^2 \leq M_g\eta_n$ for $n$ sufficiently large. Similarly, Assumption \ref{ass:locineq}(ii) and Lemma \ref{lm:obviousam} yield
\begin{align}\label{lm:intpert6}
\nabla \Upsilon_G(\theta + \frac{h}{\sqrt n})[\frac{h_{\theta,n}}{\sqrt n}] & \leq \nabla \Upsilon_G(\theta)[\frac{h_{\theta,n}}{\sqrt n}] + \|\nabla \Upsilon_G(\theta + \frac{h}{\sqrt n}) - \nabla \Upsilon_G(\theta)\|_o \|\frac{h_{\theta,n}}{\sqrt n}\|_{\mathbf B} \mathbf {1_G} \nonumber \\ & \leq \nabla \Upsilon_G(\theta)[\frac{h_{\theta,n}}{\sqrt n}] + K_g\ell_n\|\frac{h_{\theta,n}}{\sqrt n}\|_{\mathbf B} \mathbf {1_G} .
\end{align}
Hence,  results \eqref{lm:intpert5} and \eqref{lm:intpert6}, $\|h_{\theta,n}/\sqrt n\|_{\mathbf B} \leq C_0M_d\eta_n$ due to $\|\bar h_{\theta,n}\|_{\mathbf B} \leq M_d$ by \eqref{lm:intpert1}, and $\eta_n \downarrow 0$, $\ell_n \downarrow 0$, imply that for $n$ sufficiently large we have
\begin{equation}\label{lm:intpert7}
\Upsilon_G(\theta + \frac{h_{\theta,n}}{\sqrt n} + \frac{\tilde h}{\sqrt n}) \leq \Upsilon_G(\theta + \frac{h}{\sqrt n})  + \nabla \Upsilon_G(\theta)[\frac{h_{\theta,n}}{\sqrt n}] + 4M_g\eta_n\mathbf {1_G} .
\end{equation}
In addition, since $C_0\eta_n \downarrow 0$, we have $C_0\eta_n \leq 1$ eventually, and hence $\Upsilon_G(\theta + h/\sqrt n)\leq 0$, $2M_g\ell_n \leq \epsilon/2$ for $n$ sufficiently large due to $\ell_n \downarrow 0$, and result \eqref{lm:intpert3} imply that
\begin{align}\label{lm:intpert8}
\Upsilon_G(\theta + \frac{h}{\sqrt n})  + C_0\eta_n \nabla \Upsilon_G(\theta)[\bar h_{\theta,n}] & \leq C_0\eta_n\{\Upsilon_G(\theta + \frac{h}{\sqrt n})  + \nabla \Upsilon_G(\theta)[\bar h_{\theta,n}]\} \notag \\ & \leq C_0\eta_n\{2M_g\ell_n - \epsilon\} \mathbf {1_G} \leq -\frac{C_0\eta_n\epsilon}{2} \mathbf {1_G} .
\end{align}
Thus, we can conclude from results \eqref{lm:intpert4}, \eqref{lm:intpert7}, \eqref{lm:intpert8}, and $C_0 > 8M_g/\epsilon$ that
\begin{equation*}
\Upsilon_G(\theta + \frac{h_{\theta,n}}{\sqrt n} + \frac{\tilde h}{\sqrt n}) \leq \{4M_g - \frac{C_0\epsilon}{2}\}\eta_n \mathbf {1_G} \leq 0 ,
\end{equation*}
for $n$ sufficiently large, which establishes the claim of the Lemma. \qed

\begin{lemma}\label{lm:obviousam}
If $\mathbf A$ is an AM space with norm $\|\cdot\|_{\mathbf A}$ and unit $\mathbf {1_A}$, and $a_1,a_2\in \mathbf A$, then it follows that $a_1 \leq a_2 + C\mathbf1_{\mathbf A}$ for any $a_1,a_2 \in \mathbf A$ satisfying $\|a_1- a_2\|_{\mathbf A} \leq C$.
\end{lemma}

\noindent {\sc Proof:} Since $\mathbf A$ is an AM space with unit $\mathbf {1_A}$ we have that $\|a_1 - a_2\|_{\mathbf A} \leq C$ implies $|a_1 - a_2| \leq C\mathbf {1_A}$, and hence the claim follows trivially from $a_1 - a_2 \leq |a_1 - a_2|$. \qed

\begin{lemma}\label{lm:auximplicit}
Let $\mathbf A$ and $\mathbf C$ be Banach spaces with norms $\|\cdot\|_{\mathbf A}$ and $\|\cdot\|_{\mathbf C}$, $\mathbf A = \mathbf A_{1} \oplus \mathbf A_{2}$ and $F:\mathbf A \rightarrow \mathbf C$. Suppose $F(a_0) = 0$ and that there are $\epsilon_0 > 0$ and $K_0 < \infty$ such that:
\vspace{-0.15 in}
\begin{packed_enum}
    \item[(i)] $F:\mathbf A \rightarrow \mathbf C$ is Fr\'echet differentiable at all $a\in \mathcal B_{\epsilon_0}(a_0) \equiv \{a \in \mathbf A : \|a - a_0\|_{\mathbf A} \leq \epsilon_0\}$.
    \item[(ii)] $\|F(a + h) - F(a) - \nabla F(a)[h]\|_{\mathbf C} \leq K_0\|h\|^2_{\mathbf A}$ for all $a,a+h \in \mathcal B_{\epsilon_0}(a_0)$.
    \item[(iii)] $\|\nabla F(a_1) - \nabla F(a_2)\|_o \leq K_0\|a_1 - a_2\|_{\mathbf A}$ for all $a_1,a_2\in \mathcal B_{\epsilon_0}(a_0)$.
    \item[(iv)] $\nabla F(a_0) : \mathbf A \rightarrow \mathbf C$ has $\|\nabla F(a_0)\|_o \leq K_0$.
    \item[(v)] $\nabla F(a_0) : \mathbf A_2\rightarrow \mathbf C$ is bijective and $\|\nabla F(a_0)^{-1}\|_o \leq K_0$.
\end{packed_enum}
\vspace{-0.1 in} Then, for all $h_1 \in \mathbf A_1$ with $\|h_1\|_{\mathbf A} \leq (\epsilon_0/2 \wedge (4K_0^2)^{-1} \wedge 1)^2$ there is a unique $h_2^\star(h_1) \in \mathbf A_2$ with $F(a_0 + h_1 + h_2^\star(h_1)) = 0$. In addition, $h_2^\star(h_1)$ satisfies $\|h^\star_2(h_1)\|_{\mathbf A} \leq 4K_0^2\|h_1\|_{\mathbf A}$ for arbitrary $\mathbf A_1$, and $\|h^\star_2(h_1)\|_{\mathbf A} \leq 2K_0^2\|h_1\|_{\mathbf A}^2$ when $\mathbf A_1 = \mathcal N(\nabla F(a_0))$.
\end{lemma}

\noindent {\sc Proof:} We closely follow the arguments in the proof of Theorems 4.B in \cite{zeidler:1985}.
First, we define $g:\mathbf A_1 \times \mathbf A_2 \rightarrow \mathbf C$ pointwise for any $h_1 \in \mathbf A_1$ and $h_2\in \mathbf A_2$ by
\begin{equation}\label{lm:auximplicit1}
g(h_1,h_2) \equiv \nabla F(a_0)[h_2] - F(a_0 + h_1 + h_2) .
\end{equation}
Since $\nabla F(a_0) : \mathbf A_2 \rightarrow \mathbf C$ is bijective by hypothesis, $F(a_0 + h_1 + h_2) = 0$ if and only if
\begin{equation}\label{lm:auximplicit2}
h_2 = \nabla F(a_0)^{-1}[g(h_1,h_2)] .
\end{equation}
Letting $T_{h_1} : \mathbf A_2 \rightarrow \mathbf A_2$ be given by $T_{h_1}(h_2) \equiv \nabla F(a_0)^{-1}[g(h_1,h_2)]$, we see from \eqref{lm:auximplicit2} that the desired $h_2^\star(h_1)$ must be a fixed point of $T_{h_1}$. Next, define the set
\begin{equation*}
M_0 \equiv \{h_2 \in \mathbf A_2 : \|h_2\|_{\mathbf A} \leq \delta_0\}
\end{equation*}
for $\delta_0 \equiv (\epsilon_0/2) \wedge (4K_0^2)^{-1}\wedge 1$, and consider an arbitrary $h_1 \in \mathbf A_1$ with $\|h_1\|_{\mathbf A} \leq \delta_0^2$.
Notice that then $a_0 +h_1 +h_2 \in \mathcal B_{\epsilon_0}(a_0)$ for any $h_2 \in M_0$ and hence $g$ is differentiable with respect to $h_2$ with derivative $\nabla_2 g(h_1,h_2) \equiv \nabla F(a_0) - \nabla F(a_0 + h_1 + h_2)$.
Thus, if $h_2,\tilde h_2 \in M_0$, then Proposition 7.3.2 in \cite{luenberger:1969} implies that
\begin{align}\label{lm:auximplicit4}
\|g(h_1,h_2) - g(h_1,\tilde h_2)\|_{\mathbf C} & \leq \sup_{0 < \tau < 1} \|\nabla_2g(h_1, h_2 + \tau(\tilde h_2 - h_2))\|_o \|h_2 - \tilde h_2\|_{\mathbf A} \nonumber \\
& \leq \frac{1}{2K_0} \|h_2 - \tilde h_2\|_{\mathbf A} ,
\end{align}
where the final inequality follows by Condition (iii) and $\delta_0^2 \leq \delta_0 \leq (4K_0^2)^{-1}$.
Moreover,
\begin{multline}\label{lm:auximplicit5}
\|\nabla F(a_0)[h_2] - \nabla F(a_0 + h_1)[h_2]\|_{\mathbf C} \\ \leq \|\nabla F(a_0) - \nabla F(a_0+h_1)\|_o\|h_2\|_{\mathbf A} \leq K_0\|h_1\|_{\mathbf A} \|h_2\|_{\mathbf A} \leq \frac{\|h_2\|_{\mathbf A}}{4K_0}
\end{multline}
by Condition (iii) and $\|h_1\|_{\mathbf A} \leq \delta_0 \leq (4K_0^2)^{-1}$.
Similarly, for any $h_2 \in M_0$ we have
\begin{equation}\label{lm:auximplicit6}
\|F(a_0 + h_1+h_2) - F(a_0 + h_1) - \nabla F(a_0 + h_1)[h_2]\|_{\mathbf C} \leq K_0\|h_2\|_{\mathbf A}^2 \leq \frac{\|h_2\|_{\mathbf A}}{4K_0}
\end{equation}
due to $a_0 + h_1 \in \mathcal B_{\epsilon_0}(a_0)$ and Condition (ii).
Moreover, since $F(a_0) = 0$ by hypothesis, Conditions (ii) and (iv), $\|h_1\|_{\mathbf A} \leq \delta^2_0$, and $\delta_0\leq (4K_0^2)^{-1}$ yield that
\begin{equation}\label{lm:auximplicit7}
\|F(a_0 + h_1)\|_{\mathbf C} = \|F(a_0 + h_1) - F(a_0)\|_{\mathbf C} \leq K_0\|h_1\|^2_{\mathbf A} + \|\nabla F(a_0)\|_o \|h_1\|_{\mathbf A} \leq \frac{\delta_0}{2K_0} .
\end{equation}
Hence, by \eqref{lm:auximplicit1} and \eqref{lm:auximplicit5}-\eqref{lm:auximplicit7} we obtain for any $h_2 \in M_0$ and $h_1$ with $\|h_1\|_{\mathbf A} \leq \delta_0^2$
\begin{equation}\label{lm:auximplicit8}
\|g (h_1,h_2) \|_{\mathbf C} \leq \frac{\|h_2\|_{\mathbf A}}{2K_0} + \frac{\delta_0}{2K_0} \leq \frac{\delta_0}{K_0} .
\end{equation}
Thus, since $\|\nabla F(a_0)^{-1}\|_o \leq K_0$ by Condition (v), result \eqref{lm:auximplicit8} implies $T_{h_1} : M_0 \rightarrow M_0$, and \eqref{lm:auximplicit4} yields $\|T_{h_1}(h_2) - T_{h_1}(\tilde h_2)\|_{\mathbf A} \leq 2^{-1}\|h_2 - \tilde h_2\|_{\mathbf A}$ for any $h_2,\tilde h_2 \in M_0$.
By Theorem 1.1.1.A in \cite{zeidler:1985} we then conclude $T_{h_1}$ has a unique fixed point $h^\star_2(h_1) \in M_0$, and the first claim of the lemma follows from \eqref{lm:auximplicit1} and \eqref{lm:auximplicit2}.

Next, we note that since $h_2^\star(h_1)$ is a fixed point of $T_{h_1}$, we can conclude that
\begin{equation}\label{lm:auximplicit9}
\|h^\star_2(h_1)\|_{\mathbf A} = \|T_{h_1}(h^\star_2(h_1))\|_{\mathbf A} \leq \|T_{h_1}(h^\star_2(h_1)) - T_{h_1}(0)\|_{\mathbf A} + \|T_{h_1}(0)\|_{\mathbf A} .
\end{equation}
Thus, since $\|T_{h_1}(h^\star_2(h_1)) - T_{h_1}(0)\|_{\mathbf A} \leq 2^{-1}\|h^\star_2(h_1)\|_{\mathbf A}$ by \eqref{lm:auximplicit4} and $\|\nabla F(a_0)^{-1}\|_o \leq K_0$, it follows from result \eqref{lm:auximplicit9} and $T_{h_1}(0) \equiv -\nabla F(a_0)^{-1}F(a_0 + h_1)$ that
\begin{multline}\label{lm:auximplicit10}
\frac{1}{2} \|h^\star_2(h_1)\|_{\mathbf A} \leq \|T_{h_1}(0)\|_{\mathbf A} \leq K_0\|F(a_0 + h_1)\|_{\mathbf C} \\ \leq K_0\{K_0\|h_1\|^2_{\mathbf A} + \|\nabla F(a_0)\|_o \|h_1\|_{\mathbf A}\} \leq 2K_0^2\|h_1\|_{\mathbf A} ,
\end{multline}
where in the second inequality we employed $\|\nabla F(a_0)^{-1}\|_o \leq K_0$, in the third inequality we used \eqref{lm:auximplicit7}, and in the final inequality we exploited $\|h_1\|_{\mathbf A} \leq 1$.
While the estimate in \eqref{lm:auximplicit10} applies for generic $\mathbf A_1$, we note that if in addition $\mathbf A_1 = \mathcal N(\nabla F(a_0))$, then
\begin{equation*}
\frac{1}{2} \|h^\star_2(h_1)\|_{\mathbf A} \leq \|T_{h_1}(0)\|_{\mathbf A} \leq K_0\|F(a_0 + h_1)\|_{\mathbf C} \leq K_0^2\|h_1\|^2_{\mathbf A} ~,
\end{equation*}
due to $F(a_0) = 0$ and $\nabla F(a_0)[h_1] = 0$, and thus the final claim of the lemma follows. \qed 

%% file: Appendix/AppKoltCoup.tex

\section{Coupling via Koltchinskii (1994)} \label{sec:unifcoupling}

In this section we develop uniform coupling results for empirical processes that help verify Assumption \ref{ass:coupling}(i) in specific applications.
Our analysis is based on the Hungarian construction of \cite{massart:1989} and \cite{koltchinskii:1994}, and we state the results in a notation that abstracts from the rest of the paper due to their potential independent interest.
Thus, in what follows we consider $V \in \mathbf R^{d}$ to be a generic random variable distributed according to $P \in \mathbf P$, denote its support under $P$ by $\Omega(P) \subset \mathbf R^{d}$, and let $\lambda$ denote the Lebesgue measure on $\mathbf R^{d}$. For any function $f$ we further set
$$\mathbb G_n(f) \equiv \frac{1}{\sqrt n}\sum_{i=1}^n (f(V_i) - E_P[f(V)]).$$
The rates obtained through a Hungarian construction crucially depend on the ability of the functions indexing the empirical process to be approximated by a suitable Haar basis.
Here, we follow \cite{koltchinskii:1994} and control the relevant approximation errors through primitive conditions stated in terms of the integral modulus of continuity.
For a measure $P$ and a function $f:\mathbf R^{d} \rightarrow \mathbf R$, the integral modulus of continuity of $f$ is the function $\varpi(f, \cdot, P):\mathbf R_+\rightarrow \mathbf R_+$ defined for every $h\in \mathbf R_+$ as
\begin{equation}\label{appcoup:eq1}
\varpi(f, h, P) \equiv \sup_{\|s\|\leq h} (\int_{\Omega(P)} (f(v+s) - f(v))^21\{v+s\in \Omega(P)\} dP(v))^{\frac{1}{2}} .
\end{equation}
Intuitively, the integral modulus of continuity quantifies the ``smoothness" of a function $f$ by examining the difference between $f$ and its own translation.
For example, it is straightforward to verify that $\varpi(f,h,P) \lesssim h$ whenever $f$ is Lipschitz.
In contrast indicator functions such as $f(v) = 1\{v\leq t\}$ typically satisfy $\varpi(f,h,P) \lesssim h^{1/2}$.

We impose the following assumptions to establish the uniform coupling results.

\begin{assumption}\label{ass:coupreg}
(i) For all $P\in \mathbf P$, $P\ll \lambda$ and $\Omega(P)\subset \mathbf R^d$ is compact; (ii) The densities $dP/d\lambda$ satisfy $\sup_{P\in \mathbf P} \sup_{v\in \Omega(P)} \frac{dP}{d\lambda}(v) < \infty$ and $\inf_{P\in \mathbf P} \inf_{v\in \Omega(P)} \frac{dP}{d\lambda}(v) > 0$.
\end{assumption}

\begin{assumption}\label{ass:supportreg}
(i) For each $P\in \mathbf P$ there is a continuously differentiable bijection $T_P : [0,1]^{d} \rightarrow \Omega(P)$; (ii) The Jacobian $JT_P$ and its determinant $|JT_P|$ satisfy $\inf_{P\in \mathbf P} \inf_{v \in [0,1]^{d}} |JT_P(v)| > 0$ and $\sup_{P\in \mathbf P} \sup_{v\in [0,1]^{d}} \|JT_P(v)\|_{o} < \infty$.
\end{assumption}

\begin{assumption}\label{ass:coupprocessreg}
The classes $\mathcal F_n$ satisfy: (i) $\sup_{P\in \mathbf P}\sup_{f\in \mathcal F_n}\varpi(f, h, P) \leq \varphi_n(h)$ for some $\varphi_n : \mathbf R_+\rightarrow \mathbf R_+$ satisfying $\varphi_n(C h) \leq C^\kappa\varphi_n(h)$ for all $n$, $C \geq 1$, and some $\kappa > 0$; and (ii) $ \sup_{f\in \mathcal F_n} \|f\|_{\infty} \leq K_n$ for some $K_n > 0$.
\end{assumption}

In Assumption \ref{ass:coupreg} we impose that $V\sim P$ be continuously distributed for all $P\in \mathbf P$, with uniformly (in $P$) bounded supports and densities bounded from above and away from zero.
Assumption \ref{ass:supportreg} requires that the support of $V$ under each $P$ be ``smooth" in the sense that it be a differentiable transformation of the unit square.
Together, Assumptions \ref{ass:coupreg} and \ref{ass:supportreg} enable us to construct partitions of $\Omega(P)$ such that the diameter of each set in the partition is controlled uniformly in $P$; see Lemma \ref{aux:coup}. As a result, the approximation error by the Haar basis implied by each partition can be controlled uniformly by the integral modulus of continuity; see Lemma \ref{aux:haar}.
Together with Assumption \ref{ass:coupprocessreg}, which imposes conditions on the integral modulus of continuity of $\mathcal F_n$ uniformly in $P$, we can obtain a uniform coupling result through the analysis in \cite{koltchinskii:1994}.
We note that the homogeneity condition on $\varphi_n$ in Assumption \ref{ass:coupprocessreg}(i) is not necessary, but we impose it to simplify the bound.

The next theorem establishes a coupling for the empirical process $\mathbb G_n$.

\begin{theorem}\label{th:coup}
Let Assumptions \ref{ass:coupreg}-\ref{ass:coupprocessreg} hold, $\{V_i\}_{i=1}^n$ be i.i.d.\ with $V_i \sim P \in \mathbf P$ and for any $\delta_n \downarrow 0$ let $N_n \equiv  \sup_{P\in \mathbf P} N_{[\hspace{0.02 in}]}(\delta_n,\mathcal F_n,\|\cdot\|_{P,2})$, $J_n \equiv \sup_{P\in \mathbf P} J_{[\hspace{0.02 in}]}(\delta_n,\mathcal F_n,\|\cdot\|_{P,2})$,
\begin{equation}\label{th:coupdisp1}
S_n \equiv ( \sum_{i=0}^{ \lceil\log_2 n\rceil} 2^{i}\varphi_n^2(2^{-\frac{i}{d}}) )^{\frac{1}{2}} .
\end{equation}
If $N_n \uparrow \infty$, there is a Gaussian $\Iso$ (possible depending on $n$) so that uniformly in $P\in \mathbf P$
\begin{equation}\label{th:coupdisp2}
\|\mathbb G_{n} - \Iso\|_{\mathcal F_n} = O_P(\frac{K_n\log(nN_n)}{\sqrt n} + \frac{ K_n\sqrt{\log(nN_n)\log(n)}S_n}{\sqrt{n}} + J_n(1 + \frac{J_n K_n}{\delta_n^2\sqrt n})) .
\end{equation}
\end{theorem}

Theorem \ref{th:coup} is a mild modification of the results in \cite{koltchinskii:1994}.
The proof of Theorem \ref{th:coup} relies on a coupling of the empirical process on a sequence of grids of cardinality $N_n$, and employs the equicontinuity of $\mathbb G_{n}$ and $\Iso$ to obtain a coupling on the entire class $\mathcal F_n$.
The conclusion of Theorem \ref{th:coup} applies to any choice of grid accuracy $\delta_n$.
In order to obtain the best rate, $\delta_n$ must be chosen to balance the terms in \eqref{th:coupdisp2} and thus depends on the metric entropy of  $\mathcal F_n$ through the terms $N_n$ and $J_n$.

Below, we include the proof of Theorem \ref{th:coup} and auxiliary results.

\noindent {\sc Proof of Theorem \ref{th:coup}:} Let $\{\Delta_i(P)\}$ be the partitions of $\Omega(P)$ in Lemma \ref{aux:coup} and $\mathcal B_{P,i}$ the $\sigma$-algebra generated by $\Delta_i(P)$.
By Lemma \ref{aux:haar} and Assumption \ref{ass:coupprocessreg},
\begin{multline}\label{th:coup1}
\sup_{P\in \mathbf P} \sup_{f\in \mathcal F_n} (\sum_{i=0}^{ \lceil\log_2 n\rceil} 2^i E_P[(f(V) - E_P[f(V)|\mathcal B_{P,i}])^2])^{\frac{1}{2}} \\  \leq C_1( \sum_{i=0}^{ \lceil\log_2 n\rceil} 2^{i}\varphi_n^2(2^{-\frac{i}{d}}) )^{\frac{1}{2}} \equiv C_1S_n
\end{multline}
for some constant $C_1 > 0$ and for $S_n$ as defined in \eqref{th:coupdisp1}. Next, let $\mathcal F_{P, n,\delta_n} \subseteq \mathcal F_n$ denote a finite $\delta_n$-net of $\mathcal F_n$ with respect to $\|\cdot\|_{P,2}$. Since $N(\epsilon,\mathcal F_n,\|\cdot\|_{P,2})\leq N_{[\hspace{0.02 in }]}(\epsilon,\mathcal F_n,\|\cdot\|_{P,2})$, it follows from the definition of $N_n$ that we may choose $\mathcal F_{P,n,\delta_n}$ so that
\begin{equation}\label{th:coup2}
\sup_{P\in \mathbf P}\text{card}(\mathcal F_{P,n,\delta_n}) \leq \sup_{P\in \mathbf P} N_{[\hspace{0.02 in}]}(\delta_n,\mathcal F_n,\|\cdot\|_{P,2}) \equiv N_n.
\end{equation}
By Theorem 3.5 in \cite{koltchinskii:1994}, \eqref{th:coup1} and \eqref{th:coup2}, it follows that for each $n \geq 1$ there exists an isonormal process $\Iso$, such that for all $\eta_1>0$, $\eta_2 >0$
\begin{multline}\label{th:coup3}
\sup_{P\in \mathbf P} P(\frac{\sqrt{n}}{K_n}\|\mathbb G_{n} - \Iso\|_{\mathcal F_{P,n,\delta_n}} \geq \eta_1 + \sqrt{\eta_1}\sqrt{\eta_2}(C_1S_n  + 1))  \\
\lesssim N_n \exp\{-C_2\eta_1\} + n \exp\{-C_2\eta_2\} ,
\end{multline}
for some $C_2 > 0$.
Since $N_n \uparrow \infty$, \eqref{th:coup3} implies for any $\varepsilon > 0$ there are $C_3 > 0$, $C_4 > 0$ sufficiently large, such that setting $\eta_1 \equiv C_3\log(N_n)$ and $\eta_2 \equiv C_3\log(n)$ yields
\begin{equation}\label{th:coup4}
\sup_{P\in \mathbf P} P(\|\mathbb G_{n} - \Iso\|_{\mathcal F_{P,n,\delta_n}} \geq C_4 K_n\times\frac{\log(nN_n) + \sqrt{\log(N_n)\log(n)}S_n}{\sqrt{n}}) < \varepsilon .
\end{equation}
Next, note that by definition of $\mathcal F_{P,n,\delta_n}$, there exists a $\Gamma_{n,P} : \mathcal F_{n}\rightarrow \mathcal F_{P,n,\delta_n}$ such that $\sup_{P\in \mathbf P}\sup_{f\in \mathcal F_n}\|f - \Gamma_{n,P} (f)\|_{P,2}\leq \delta_n$.
Let $D(\epsilon,\mathcal F_n,\|\cdot\|_{P,2})$ denote the $\epsilon$-packing number for $\mathcal F_n$ under $\|\cdot\|_{P,2}$, and note $D(\epsilon,\mathcal F_n,\|\cdot\|_{P,2}) \leq N_{[\hspace{0.02 in}]}(\epsilon,\mathcal F_n,\|\cdot\|_{P,2})$.
Therefore, by Corollary 2.2.8 in \cite{vandervaart:wellner:1996} we can conclude that
\begin{multline}\label{th:coup5}
\sup_{P\in \mathbf P}E_P[\|\Iso- \Iso \circ \Gamma_{n,P}\|_{\mathcal F_n}] \\
\lesssim \sup_{P\in \mathbf P}\int_0^{\delta_n} \sqrt{\log D(\epsilon,\mathcal F_n, \|\cdot\|_{P,2})}d\epsilon \leq \sup_{P\in \mathbf P}J_{[\hspace{0.02 in}]}(\delta_n,\mathcal F_n,\|\cdot\|_{P,2}) \equiv J_n .
\end{multline}
Similarly, employing Lemma 3.4.2 in \cite{vandervaart:wellner:1996} yields that
\begin{multline}\label{th:coup6}
\sup_{P\in \mathbf P} E_P[\|\mathbb G_{n} - \mathbb G_{n} \circ \Gamma_{n,P}\|_{\mathcal F_n}] \\ \lesssim \sup_{P\in \mathbf P}J_{[\hspace{0.02 in}]}(\delta_n,\mathcal F_n,\|\cdot\|_{P,2})(1 + \sup_{P\in \mathbf P} \frac{J_{[\hspace{0.02 in}]}(\delta_n,\mathcal F_n,\|\cdot\|_{P,2})K_n}{\delta_n^2\sqrt{n}}) \equiv J_n(1 + \frac{J_nK_n}{\delta_n^2\sqrt{n}}).
\end{multline}
Therefore, combining \eqref{th:coup4}, \eqref{th:coup5}, and \eqref{th:coup6} together with the decomposition
\begin{equation*}
\|\mathbb G_{n} - \Iso\|_{\mathcal F_n} \\ \leq \|\mathbb G_{n} - \Iso\|_{\mathcal F_{P,n,\delta_n}} + \|\mathbb G_{n} - \mathbb G_{n}\circ \Gamma_{n,P}\|_{\mathcal F_n} + \|\Iso - \Iso \circ \Gamma_{n,P}\|_{\mathcal F_n} ,
\end{equation*}
establishes the claim of the theorem by Markov's inequality. \qed

\begin{lemma}\label{aux:coup}
Let $\mathcal B_P$ denote the completion of the Borel $\sigma-$algebra on $\Omega(P)$ with respect to $P$. If Assumptions \ref{ass:coupreg} and \ref{ass:supportreg} hold, then for each $P\in \mathbf P$ there exists a sequence $\{\Delta_i(P)\}$ of partitions of $(\Omega(P), \mathcal B_P, P)$ such that:

\vspace{-0.2in}
\begin{packed_enum}
    \item[(i)] $\Delta_i(P) = \{\Delta_{i,k}(P) : k = 0,\ldots, 2^i-1\}$, $\Delta_{i,k}(P) \in \mathcal B_P$, and $\Delta_{0,0}(P) = \Omega(P)$.
    \item[(ii)] $\Delta_{i,k}(P) = \Delta_{i+1,2k}(P) \cup \Delta_{i+1,2k+1}(P)$ and $\Delta_{i+1,2k}(P)\cap \Delta_{i+1,2k+1}(P) = \emptyset$ for any integers $k = 0,\ldots 2^i-1$ and $i\geq 0$.
    \item[(iii)] $P(\Delta_{i+1,2k}(P)) = P(\Delta_{i+1,2k+1}(P)) = 2^{-i-1}$ for $k = 0,\ldots 2^i-1$, $i\geq 0$.
    \item[(iv)] $\sup_{P\in \mathbf P}\max_{0\leq k \leq 2^i-1}\sup_{v,v^\prime \in \Delta_{i,k}(P)} \|v - v^\prime\|_2 = O(2^{-\frac{i}{d}})$.
    \item[(v)] $\mathcal B_P$ equals the completion with respect to $P$ of the $\sigma$-algebra generated by $\bigcup_{i} \Delta_i(P)$.
\end{packed_enum}
\end{lemma}

\noindent {\sc Proof:} Let $\mathcal A$ denote the Borel $\sigma$-algebra on $[0,1]^{d}$, and for any $A\in \mathcal A$ define
\begin{equation}\label{auxcoup1}
Q_P(A) \equiv P(T_P(A)) ,
\end{equation}
where $T_P(A) \in \mathcal B_P$ due to $T_P^{-1}$ being measurable.
Moreover, $Q_P([0,1]^{d}) = 1$ due to $T_P$ being surjective, and $Q_P$ is $\sigma$-additive due to $T_P$ being injective.
Hence, we conclude $Q_P$ defined by \eqref{auxcoup1} is a probability measure on $([0,1]^{d}, \mathcal A$).
In addition, for $\lambda$ the Lebesgue measure, we obtain from Theorem 3.7.1 in \cite{bogachev1:2007} that
\begin{equation}\label{auxcoup2}
Q_P(A) = P(T_P(A)) = \int_{T_P(A)} \frac{dP}{d\lambda}(v)d\lambda(v) = \int_A \frac{dP}{d\lambda}(T_P(a))|JT_P(a)|d\lambda(a) ,
\end{equation}
where $|JT_P(a)|$ denotes the Jacobian of $T_P$ at any point $a\in[0,1]^{d}$.
Hence, $Q_P$ has density with respect to Lebesgue measure given by $g_P(a) \equiv \frac{dP}{d\lambda}(T_P(a))|JT_P(a)|$ for any $a \in [0,1]^{d}$.
Next, let $a = (a_1,\ldots, a_{d})^\prime \in [0,1]^{d}$ and define for any $t\in [0,1]$
\begin{equation}\label{auxcoup3}
G_{l,P}(t|A) \equiv \frac{Q_P(a \in A : a_l \leq t)}{Q_P(A)} ,
\end{equation}
for any $A \in \mathcal A$ with $Q_P(A) > 0$ and $1\leq l \leq d$.
Also let $\text{m}(i) \equiv i - \lfloor \frac{i-1}{d}\rfloor \times d$ (i.e.\ $\text{m}(i)$ equals $i$ modulo $d$), set $\tilde \Delta_{0,0}(P) = [0,1]^{d}$, and inductively define the partitions of $[0,1]^{d}$
\begin{align}\label{auxcoup4}
\tilde \Delta_{i+1,2k}(P) & \equiv \{a \in \tilde \Delta_{i,k}(P) : G_{\text{m}(i+1),P}(a_{\text{m}(i+1)}| \tilde \Delta_{i,k}(P)) \leq \frac{1}{2}\} \nonumber \\ \tilde \Delta_{i+1,2k+1}(P)  & \equiv \tilde \Delta_{i,k}(P)\setminus \tilde \Delta_{i+1,2k}(P)
\end{align}
for $0\leq k \leq 2^i - 1$.
For $\text{cl}(\tilde \Delta_{i,k}(P))$ the closure of $\tilde \Delta_{i,k}(P)$, we then note that by construction each $\tilde \Delta_{i,k}(P)$ is a hyper-rectangle in $[0,1]^{d}$ -- i.e.\ it is of the general form
\begin{equation*}
\text{cl}(\tilde \Delta_{i,k}(P)) = \prod_{j=1}^{d} [l_{i,k,j}(P), u_{i,k,j}(P)] .
\end{equation*}
Moreover, since $g_P$ is positive on $[0,1]^{d}$ by Assumptions \ref{ass:coupreg}(ii) and \ref{ass:supportreg}(ii), it follows that for any $i \geq 0$, $0 \leq k \leq 2^i -1$ and $1\leq j \leq d$, we have
\begin{align}\label{auxcoup6}
l_{i+1,2k,j}(P) & = l_{i,k,j}(P) \nonumber \\ u_{i+1,2k,j}(P) & = \left\{
\begin{array}{l}
u_{i,k,j}(P) \text{ if }  j \neq \text{m}(i+1) \\
\text{ solves } G_{\text{m}(i+1),P}(u_{i+1,2k,j}(P)|\tilde \Delta_{i,k}(P)) = \frac{1}{2} \text{ if } j = \text{m}(i+1)
\end{array} \right. 
\end{align}
Similarly, since $\tilde \Delta_{i+1,2k+1}(P) = \tilde \Delta_{i,k}(P)\setminus \tilde \Delta_{i+1,2k}(P)$, it additionally follows that
\begin{equation}\label{auxcoup7}
u_{i+1,2k+1,j}(P) = u_{i,k,j}(P) \hspace{0.3 in} l_{i+1,2k+1,j}(P) = \left\{
\begin{array}{l}
l_{i,k,j}(P) \text{ if }  j \neq \text{m}(i+1) \\
u_{i+1,2k,j}(P) \text{ if } j = \text{m}(i+1)
\end{array} \right. 
\end{equation}
Since $Q_P(\text{cl}(\tilde \Delta_{i+1,2k}(P))) = Q_P(\tilde \Delta_{i+1,2k}(P))$ by $Q_P \ll \lambda$, \eqref{auxcoup3} and \eqref{auxcoup6} yield
\begin{align*}
Q_P(\tilde \Delta_{i+1,2k}(P)) & = Q_P(a \in \tilde \Delta_{i,k}(P) : a_{\text{m}(i+1)} \leq u_{i+1,2k,\text{m}(i+1)}(P)) \nonumber \\ & = G_{\text{m}(i+1),P}(u_{i+1,2k,\text{m}(i+1)}(P)|\tilde \Delta_{i,k}(P))Q_P(\tilde \Delta_{i,k}(P)) \nonumber \\ &  = \frac{1}{2}Q_P(\tilde \Delta_{i,k}(P)) .
\end{align*}
Therefore, since $\tilde \Delta_{i,k}(P) = \tilde \Delta_{i+1,2k}(P) \cup \tilde \Delta_{i+1,2k+1}(P)$, it follows  $Q_P(\tilde \Delta_{i+1,2k+1}(P)) = \frac{1}{2}Q_P(\tilde \Delta_{i,k}(P))$ for $0\leq k \leq 2^i -1$ as well.
In particular, $Q_P(\tilde \Delta_{0,0}(P)) = 1$ implies that
\begin{equation}\label{auxcoup9}
Q_P(\tilde \Delta_{i,k}(P)) = \frac{1}{2^{i}}
\end{equation}
for any integers $i\geq 1$ and $0\leq k \leq 2^i-1$.
Moreover, we note that result \eqref{auxcoup2} and Assumptions \ref{ass:coupreg}(ii) and \ref{ass:supportreg}(ii) together imply that the density $g_P$ of $Q_P$ satisfies
\begin{equation}\label{auxcoup10}
0 < \inf_{P\in \mathbf P} \inf_{a \in [0,1]^{d}} g_P(a) < \sup_{P\in \mathbf P}\sup_{a\in [0,1]^{d}} g_P(a) < \infty ,
\end{equation}
and therefore $Q_P(A) \asymp \lambda(A)$ uniformly in $A\in \mathcal A$ and $P\in \mathbf P$.
Hence, since by \eqref{auxcoup6} $u_{i+1,2k,j}(P) = u_{i,k,j}(P)$ and $l_{i+1,2k,j}(P) = l_{i,k,j}(P)$ for all $j\neq \text{m}(i+1)$, we obtain
\begin{multline}\label{auxcoup11}
\frac{(u_{i+1,2k,\text{m}(i+1)}(P) - l_{i+1,2k,\text{m}(i+1)}(P))}{(u_{i,k,\text{m}(i+1)}(P) - l_{i,k,\text{m}(i+1)}(P))} = \frac{\prod_{j=1}^{d}(u_{i+1,2k,j}(P) - l_{i+1,2k,j}(P))}{\prod_{j=1}^{d}(u_{i,k,j}(P) - l_{i,k,j}(P))} \\ = \frac{\lambda(\tilde \Delta_{i+1,2k}(P))}{\lambda(\tilde \Delta_{i,k}(P))} \asymp \frac{Q_P(\tilde \Delta_{i+1,2k}(P))}{Q_P(\tilde \Delta_{i,k}(P))} = \frac{1}{2}
\end{multline}
uniformly in $P \in \mathbf P$, $i\geq0$, and $0\leq k \leq 2^i-1$ by results \eqref{auxcoup9} and \eqref{auxcoup10}.
Moreover, by identical arguments but using \eqref{auxcoup7} instead of \eqref{auxcoup6} we conclude
\begin{equation}\label{auxcoup12}
\frac{(u_{i+1,2k+1,\text{m}(i+1)}(P) - l_{i+1,2k+1,\text{m}(i+1)}(P))}{(u_{i,k,\text{m}(i+1)}(P) - l_{i,k,\text{m}(i+1)}(P))} \asymp \frac{1}{2}
\end{equation}
also uniformly in $P\in \mathbf P$, $i\geq 0$ and $0\leq k \leq 2^i - 1$.
Thus, since $(u_{i+1,2k,j}(P) - l_{i+1,2k,j}(P)) =(u_{i+1,2k+1,j}(P) - l_{i+1,2k+1,j}(P)) = (u_{i,k,j}(P) - l_{i,k,j}(P))$ for all $j \neq \text{m}(i+1)$, and $u_{0,0,j}(P) - l_{0,0,j}(P) = 1$ for all $1\leq j \leq d$ we obtain from $\text{m}(i) = i -\lfloor \frac{i-1}{d}\rfloor\times d$, results \eqref{auxcoup11} and \eqref{auxcoup12}, and proceeding inductively that
\begin{equation}\label{auxcoup13}
(u_{i,k,j}(P) - l_{i,k,j}(P))\asymp 2^{-\frac{i}{d}} ,
\end{equation}
uniformly in $P\in \mathbf P$, $i\geq 0$, $0\leq k \leq 2^i - 1$, and $1\leq j \leq d$.
Thus, result \eqref{auxcoup13} yields
\begin{multline}\label{auxcoup14}
\sup_{P\in \mathbf P}\max_{0\leq k\leq 2^i -1}\sup_{a,a^\prime \in \tilde \Delta_{i,k}(P)} \|a - a^\prime\| \\ \leq \sup_{P\in \mathbf P} \max_{0\leq k \leq 2^i -1}\max_{1\leq j \leq d} \sqrt{d}\times (u_{i,j,k}(P) - l_{i,j,k}(P)) = O(2^{-\frac{i}{d}}) .
\end{multline}

We next obtain the desired sequence of partitions $\{\Delta_i(P)\}$ of $(\Omega(P),\mathcal B_P, P)$ by constructing them from the partition $\{\tilde \Delta_{i,k}(P)\}$ of $[0,1]^{d}$. To this end, set
\begin{equation*}
\Delta_{i,k}(P) \equiv T_P(\tilde \Delta_{i,k}(P))
\end{equation*}
for all $i\geq 0$ and $0\leq k \leq 2^i-1$. Note that $\{\Delta_i(P)\}$ satisfies conditions (i) and (ii) due to $T_P^{-1}$ being a measurable map, $T_P$ being bijective, and result \eqref{auxcoup4}.
In addition, $\{\Delta_i(P)\}$ satisfies condition (iii) since by definition \eqref{auxcoup1} and result \eqref{auxcoup9} we have
\begin{equation*}
P(\Delta_{i,k}(P)) = P(T_P(\tilde \Delta_{i,k}(P))) = Q_P(\tilde \Delta_{i,k}(P)) = 2^{-i} ,
\end{equation*}
for all $0\leq k \leq 2^i-1$.
Moreover, by Assumption \ref{ass:supportreg}(ii), $\sup_{P\in \mathbf P}\sup_{a\in [0,1]^{d}} \|JT_P(a)\|_{o} < \infty$, and hence by the mean value theorem we can conclude that
\begin{multline*}
\sup_{P\in \mathbf P}\max_{0\leq k \leq 2^i-1} \sup_{v,v^\prime \in \Delta_{i,k}(P)} \|v -v^\prime\|  = \sup_{P\in \mathbf P}\max_{0\leq k \leq 2^i-1} \sup_{a,a^\prime \in \tilde \Delta_{i,k}(P)} \|T_P(a) -T_P(a^\prime)\| \\ \lesssim \sup_{P\in \mathbf P} \max_{0\leq k \leq 2^i-1} \sup_{a,a^\prime \in \tilde \Delta_{i,k}(P)} \|a -a^\prime\| = O(2^{-\frac{i}{d}}) ,
\end{multline*}
by result \eqref{auxcoup14}, which verifies that $\{\Delta_i(P)\}$ satisfies condition (iv).
Also note that to verify $\{\Delta_i(P)\}$ satisfies condition (v) it suffices to show that $\bigcup_{i\geq 0} \Delta_i(P)$ generates the Borel $\sigma$-algebra on $\Omega(P)$.
To this end, we first aim to show that
\begin{equation}\label{auxcoup18}
\mathcal A = \sigma(\bigcup_{i\geq 0} \tilde \Delta_i(P)) ,
\end{equation}
where for a collection of sets $\mathcal C$, $\sigma(\mathcal C)$ denotes the $\sigma$-algebra generated by $\mathcal C$.
For any closed set $A \in \mathcal A$, then define $D_i(P)$ to be given by
\begin{equation*}
D_i(P) \equiv \bigcup_{k : \tilde \Delta_{i,k}(P)\cap A \neq \emptyset } \tilde \Delta_{i,k}(P) .
\end{equation*}
Notice that since $\{\tilde \Delta_i(P)\}$ is a partition of $[0,1]^{d}$, $A \subseteq D_i(P)$ for all $i \geq 0$ and hence $A \subseteq \bigcap_{i\geq 0} D_i(P)$.
Moreover, if $a_0\in A^c$, then $A^c$ being open and \eqref{auxcoup14} imply $a_0 \notin D_i(P)$ for $i$ sufficiently large.
Hence, $A^c \cap (\bigcap_{i\geq 0} D_i(P)) = \emptyset$ and therefore $A = \bigcap_{i\geq 0} D_i(P)$.
It follows that if $A$ is closed, then $A \in \sigma (\bigcup_{i\geq 0}\tilde \Delta_i(P))$, which implies $\mathcal A \subseteq \sigma (\bigcup_{i\geq 0}\tilde \Delta_i(P))$.
On the other hand, since $\tilde \Delta_{i,k}(P)$ is Borel for all $i\geq 0$ and $0\leq k \leq 2^i -1$, we also have $\sigma(\bigcup_{i\geq 0} \tilde \Delta_i(P))\subseteq \mathcal A$, and hence \eqref{auxcoup18} follows.
To conclude, we then note that
\begin{equation}\label{auxcoup20}
\sigma(\bigcup_{i\geq 0} \Delta_i(P)) = \sigma(\bigcup_{i\geq 0} T_P(\tilde \Delta_i(P))) = T_P(\sigma(\bigcup_{i\geq 0} \tilde \Delta_i(P))) = T_P(\mathcal A) ,
\end{equation}
by Corollary 1.2.9 in \cite{bogachev1:2007}. However, $T_P$ and $T_P^{-1}$ being continuous implies $T_P(\mathcal A)$ equals the Borel $\sigma$-algebra in $\Omega(P)$, and therefore \eqref{auxcoup20} implies $\{\Delta_i(P)\}$ satisfies condition (v) establishing the lemma. \qed

\begin{lemma}\label{aux:haar}
Let $\{\Delta_i(P)\}$ be as in Lemma \ref{aux:coup}, and $\mathcal B_{P,i}$ denote the $\sigma$-algebra generated by $\Delta_i(P)$. If Assumptions \ref{ass:coupreg} and \ref{ass:supportreg} hold, then there are $K_0 > 0$, $K_1 \geq 1$ such that for all $P\in \mathbf P$ and any $f$ satisfying $f\in L^2_P$ for all $P\in \mathbf P$:
$$E_P[(f(V) - E_P[f(V)|\mathcal B_{P,i}])^2] \leq K_0\times \varpi^2(f, K_1\times 2^{-\frac{i}{d}},P) .$$
\end{lemma}

\noindent {\sc Proof:} Since $\Delta_i(P)$ is a partition of $\Omega(P)$ and $P(\Delta_{i,k}(P)) = 2^{-i}$ for all $i\geq 0$ and $0\leq k\leq 2^i-1$, we may express $E_P[f(V)|\mathcal B_{P,i}]$ as an element of $L^2_P$ by
\begin{equation*}
E_P[f(V)|\mathcal B_{P,i}] = 2^i \sum_{k=0}^{2^i-1} 1\{V\in \Delta_{i,k}(P)\} \int_{\Delta_{i,k}(P)} f(v)dP(v) .
\end{equation*}
Hence, employing that $P(\Delta_{i,k}(P)) = 2^{-i}$ for all $i\geq 0$ and $0\leq k\leq 2^i-1$ together with $\Delta_i(P)$ being a partition of $\Omega(P)$, and applying Holder's inequality to the term $(f(v)-f(\tilde v))1\{v\in \Omega(P)\} 1\{\tilde v\in \Delta_{i,k}(P)\}$ we can conclude that
\begin{align*}
E_P[(f(V) & - E_P[f(V)|\mathcal B_{P,i}])^2] \nonumber \\ & = \sum_{k=0}^{2^i-1} \int_{\Delta_{i,k}(P)} (f(v) - 2^i \int_{\Delta_{i,k}(P)} f(\tilde v)dP(\tilde v))^2 dP(v) \nonumber \\
& = \sum_{k=0}^{2^i-1}2^{2i} \int_{\Delta_{i,k}(P)} (\int_{\Delta_{i,k}(P)}(f(v) - f(\tilde v))1\{ v\in \Omega(P)\} dP(\tilde v))^2dP(v) \nonumber \\
& \leq \sum_{k=0}^{2^i-1}2^{2i}P(\Delta_{i,k}(P)) \int_{\Delta_{i,k}(P)}  \int_{\Delta_{i,k}(P)}(f(v) - f(\tilde v))^21\{ v\in \Omega(P)\} dP(\tilde v)dP(v) \nonumber \\
& = \sum_{k=0}^{2^i-1}2^{i}\int_{\Delta_{i,k}(P)}  \int_{\Delta_{i,k}(P)}(f(v) - f(\tilde v))^21\{ v\in \Omega(P)\} dP(\tilde v)dP(v) .
\end{align*}
Let $D_i \equiv \sup_{P\in \mathbf P}\max_{0\leq k\leq 2^i-1} \text{diam}\{\Delta_{i,k}(P)\}$, where $\text{diam}\{\Delta_{i,k}(P)\}$ is the diameter of $\Delta_{i,k}(P)$.
Further note that by Lemma \ref{aux:coup}(iv), $D_i = O(2^{-\frac{i}{d}})$ and hence we have $\lambda(\{s \in \mathbf R^{d} : \|s\| \leq D_i\}) \leq M_12^{-i}$ for some $M_1 > 0$ and $\lambda$ the Lebesgue measure.
Noting that $\sup_{P\in \mathbf P} \sup_{v\in \Omega(P)} \frac{dP}{d\lambda}(v) < \infty$ by Assumption \ref{ass:coupreg}(ii), and doing the change of variables $s = v - \tilde v$ we then obtain for some constant $M_0 > 0$ that
\begin{align}\label{auxhaar3}
E_P[& (f(V) - E_P[f(V)|\mathcal B_{P,i}])^2] \nonumber \\ & \leq M_0 \sum_{k=0}^{2^i-1}2^{i}\int_{\Delta_{i,k}(P)}  \int_{\|s\| \leq D_i}(f(\tilde v+s) - f(\tilde v))^21\{ s+\tilde v \in \Omega(P)\} d\lambda(s)dP(\tilde v) \nonumber \\  & \leq M_0M_1\sup_{\|s\| \leq D_i}\sum_{k=0}^{2^i-1} \int_{\Delta_{i,k}(P)}(f(\tilde v+ s) - f(\tilde v))^21\{ \tilde v+s\in \Omega(P)\} dP(\tilde v) .
\end{align}
Hence, since $\{\Delta_{i,k}(P) : k = 0\ldots 2^i-1\}$ is a partition of $\Omega(P)$, $\varpi(f,h,P)$ is decreasing in $h$, and $D_i \leq K_1 2^{-\frac{i}{d}}$ for some $K_1 \geq 1$ by Lemma \ref{aux:coup}(iv), we obtain
\begin{equation}\label{auxhaar4}
E_P[(f(V) - E_P[f(V)|\mathcal B_{P,i}])^2] \leq M_0M_1 \times \varpi^2(f, K_1\times 2^{-\frac{i}{d}},P)
\end{equation}
by \eqref{auxhaar3}. Setting $K_0 \equiv M_0\times M_1$ in \eqref{auxhaar4} then establishes the lemma. \qed

%% file: Appendix/AppBootCoup.tex

\section{Uniform Bootstrap Coupling}\label{sec:bootcoup}

We next provide uniform coupling results for the multiplier bootstrap that allow us to verify Assumption \ref{ass:bootcoupling} in a variety of problems.
The results in this appendix may be of independent interest, as they extend the validity of the multiplier bootstrap to suitable non-Donsker classes $\mathcal F_n$.
For this reason, as in Section \ref{sec:unifcoupling}, we state the results in a notation that abstracts from the rest of the paper.
Hence, here $V \in \mathbf R^{d}$ should be interpreted as a generic random variable whose distribution is given by $P\in \mathbf P$. 
For $\{\omega_i\}_{i=1}^n$ i.i.d.\ standard normal random variables independent of $\{V_i\}_{i=1}^n$ we also set
$$\Wboot(f) \equiv \frac{1}{\sqrt n}\sum_{i=1}^n \omega_i\{f(V_i) - \frac{1}{n}\sum_{j=1}^n f(V_j)\}.$$

Our coupling results rely on a series approximation to the elements of $\mathcal F_n$.
To this end, we will assume that for each $P \in \mathbf P$ there is a basis $\{f_{d,n,P}\}_{d=1}^{d_n}$, with $d_n$ possibly diverging to infinity, that provides a suitable approximation to every $f\in \mathcal F_n$.
Formally, for $f_{n,P}^{d_n}(v) \equiv (f_{1,n,P}(v),\ldots, f_{d_n,n,P}(v))^\prime$, we impose the following:

\begin{assumption}\label{ass:4seriescoup}
For each $P\in \mathbf P$ there is an array of functions $\{f_{d,n,P}\}_{d=1}^{d_n} \subset L^2_P$ such that: (i) The eigenvalues of $E_P[f_{n,P}^{d_n}(V)f_{n,P}^{d_n}(V)^\prime]$ are bounded by $1 \leq C_n$ uniformly in $P\in \mathbf P$; (ii) $\sup_{P\in \mathbf P} \max_{1\leq d\leq d_n} \|f_{d,n,P}\|_{\infty} \leq K_n$ with $1\leq K_n$ finite.
\end{assumption}

\begin{assumption}\label{ass:4seriesreg}
For every $f\in \mathcal F_n$ and $P\in \mathbf P$ there is a $\beta_{n,P}(f)\in \mathbf R^{d_n}$ such that: (i)  The class $\mathcal G_{n,P} \equiv \{(f - \int fdP)  - f_{n,P}^{d_n\prime}\beta_{n,P}(f) : f\in \mathcal F_n\}$ has envelope $G_{n,P}$ which  satisfies $\|g\|_{P,2} \leq \delta_n \|G_{n,P}\|_{P,2}$ for all $P\in \mathbf P$, $g\in \mathcal G_{n,P}$, and some $\delta_n > 0$ with
$$J_{1n} \equiv \sup_{P\in \mathbf P} \{J_{[\hspace{0.02 in}]}(\delta_n\|G_{n,P}\|_{P,2},\mathcal G_{n,P},\|\cdot\|_{P,2}) + \sqrt{n}E_P[G_{n,P}(V)\exp\{-\frac{ n\delta_n^2 \|G_{n,P}\|_{P,2}^2}{G_{n,P}^2(V)\eta_{n,P}}\}]\}$$
finite and $\eta_{n,P} \equiv 1+\log N_{[\hspace{0.02 in}]}(\delta_n\|G_{n,P}\|_{P,2},\mathcal G_{n,P},\|\cdot\|_{P,2})$;
(ii) The set $\mathcal B_{n} \equiv  \{\beta_{n,P}(f) : f\in \mathcal F_n, P \in \mathbf P\}\cup\{0\}$ satisfies $J_{2n} \equiv \int_0^\infty \sqrt{\log(N(\epsilon,\mathcal B_{n},\|\cdot\|_2))} d\epsilon < \infty$.
\end{assumption}

Assumption \ref{ass:4seriescoup} imposes our regularity conditions on the approximating functions $\{f_{d,n,P}\}_{d=1}^{d_n}$.
We emphasize that the functions $\{f_{d,n,P}\}_{d=1}^{d_n}$ need not be known: They are only employed in the theoretical construction of the coupling.
In certain applications, such as when $\mathcal F_n$ is finite dimensional, a basis $\{f_{d,n,P}\}_{d=1}^{d_n}$ may be naturally available.
The approximating requirements on $\{f_{d,n,P}\}_{d=1}^{d_n}$ are imposed in Assumption \ref{ass:4seriesreg}.
In particular, Assumption \ref{ass:4seriesreg}(i) requires that the remainder of the approximation of $\mathcal F_n$ by $\{f_{d,n,P}\}_{d=1}^{d_n}$ not be ``too large."
Intuitively, Assumption \ref{ass:4seriesreg}(i) controls the ``bias" in a series approximation of $\mathcal F_n$ by linear combinations of $\{f_{d,n,P}\}_{d=1}^{d_n}$.
Assumption \ref{ass:4seriesreg}(ii) in turn controls the ``variance" of the series approximation by demanding that the class of approximating functions have a finite entropy.

We next show Assumptions \ref{ass:4seriescoup} and \ref{ass:4seriesreg} suffice for coupling $\Wboot$.

\begin{theorem}\label{th:mainbootcoup}
Let Assumptions \ref{ass:4seriescoup}, \ref{ass:4seriesreg} hold, $\{(\omega_i,V_i)\}_{i=1}^n$ be i.i.d.\ with $V_i \sim P\in \mathbf P$, $\omega_i \sim N(0,1)$, $\omega_i$ and $V_i$ independent, and  $d_n\log(1+d_n)K_n^2C_n = o(n)$. Then: (i) There is a linear Gaussian $\IsoW$ (possibly depending on $n$) independent of $\{V_i\}_{i=1}^n$ with
$$\|\Wboot - \IsoW\|_{\mathcal F_n} = O_P(J_{2n}\{\frac{K_n^2C_n d_n\log(1+d_n)}{n}\}^{1/4} + J_{1n})$$
uniformly in $P\in \mathbf P$ with $E[\IsoW(f)] = 0$ and $E[(\IsoW(f)\IsoW(g)] = {\rm Cov}_P\{f(V),g(V)\}$ for any $f,g\in \mathcal F_n$. 
(ii) If in addition $\sup_{P\in \mathbf P}\| (\text{\rm Var}_P\{f_{n,P}^{d_n}(V)\})^{-1}\|_{o,2}\leq \xi_n < \infty$ and $\xi_n\sqrt{d_n\log(1+d_n)C_n}K_n/\sqrt n = o(1)$, then uniformly in $P\in \mathbf P$
$$\|\Wboot - \IsoW\|_{\mathcal F_n} = O_P(\frac{J_{2n}K_n \sqrt{\xi_nC_nd_n\log(1+d_n)}}{\sqrt n} + J_{1n}).$$
\end{theorem}

Theorem \ref{th:mainbootcoup}(i) derives a rate of convergence for the coupled process, while Theorem \ref{th:mainbootcoup}(ii) improves on the rate under the additional requirement that $\text{Var}_P\{f_{n,P}^{d_n}(V)\}$ be bounded away from singularity.
The rates of both Theorems \ref{th:mainbootcoup}(i) and \ref{th:mainbootcoup}(ii) depend on the selected sequence $d_n$, which should be chosen optimally.
Heuristically, the proof of Theorem \ref{th:mainbootcoup} proceeds in two steps.
First, we construct a multivariate normal random variable $\IsoW(f_{n,P}^{d_n}) \in \mathbf R^{d_n}$ that is coupled with $\Wboot(f_{n,P}^{d_n}) \in \mathbf R^{d_n}$, and then employ the linearity of $\Wboot$ to obtain a suitable coupling on the subspace $\mathbb S_{n,P} \equiv \overline{\text{span}}\{f_{n,P}^{d_n}\}$.
Second, we employ Assumption \ref{ass:4seriesreg}(i) to show that a successful coupling on $\mathbb S_{n,P}$ leads to the desired construction since $\mathcal F_n$ is well approximated by $\{f_{d,n,P}\}_{d=1}^{d_n}$.

Below, we include the proof of Theorem \ref{th:mainbootcoup} and auxiliary results.

\noindent {\sc Proof of Theorem \ref{th:mainbootcoup}:}
We first couple $\Wboot$ on a finite dimensional subspace and then show that such a result suffices for coupling $\Wboot$ and $\IsoW$ on $\mathcal F_n$.
To this end, let $\mathbb S_{n,P}\equiv \overline{\text{span}}\{f_{n,P}^{d_n}\}$ and note that Assumption \ref{ass:4seriesreg}(ii) and Lemma \ref{lm:finbootcoup} imply that there exists a linear Gaussian process $\mathbb G_{P}^{(1)}$ on $\mathbb S_{n,P}$ and a sequence $R_n = o(1)$ such that
\begin{equation}\label{th:mainbootcoup1}
\sup_{\beta \in \mathcal B_n} |\Wboot (f_{n,P}^{d_n\prime} \beta) - \mathbb G^{(1)}_{P}(f_{n,P}^{d_n\prime} \beta)| = O_P(J_{2n}R_n)
\end{equation}
uniformly in $P\in \mathbf P$, $E[\mathbb G_P^{(1)}(f_{n,P}^{d_n\prime} \beta)]=0$, and also $E[(\mathbb G_P^{(1)}(f_{n,P}^{d_n\prime} \beta_1))(\mathbb G_P^{(1)}(f_{n,P}^{d_n\prime} \beta_2))] = \text{Cov}_P(f_{n,P}^{d_n}(V)^\prime \beta_1,f_{n,P}^{d_n}(V)^\prime \beta_2\}$.
To establish part (i) of the theorem we will set $R_n = (d_n\log(1+d_n)C_nK_n^2/n)^{1/4}$ and employ Lemma \ref{lm:finbootcoup}(i), while to establish part (ii) we will set $R_n = (\xi_n d_n \log(1+d_n)C_nK_n^2/n)^{1/2}$ and employ Lemma \ref{lm:finbootcoup}(ii) instead.

Next note that since $\Wboot(f - c) = \Wboot(f)$ for any $c\in \mathbf R$ and $f\in L^2_P$, we may assume without loss of generality that $E_P[f(V)] = 0$ for all $f\in \mathcal F_n$. 
For any closed linear subspace $\mathbb A_{n,P}\subseteq L^2_P$ let $\text{Proj}\{f |\mathbb A_{n,P}\}$ denote the $\|\cdot\|_{P,2}$ projection of $f$ onto $\mathbb A_{n,P}$ and $\mathbb A^\perp_{n,P} \equiv \{f \in L^2_P : f = g - \text{Proj}\{g |\mathbb A_{n,P}\} \text{ for some } g \in L^2_P\}$.
Assuming the underlying probability space is suitably enlarged to carry a linear isonormal Gaussian process $\mathbb G_{P}^{(2)}$ on $\{\text{Proj}\{f|\mathbb S_{n,P}^\perp\} : f\in \mathcal F_n\cup \mathcal G_{n,P}\}$ independent of $\mathbb G_{P}^{(1)}$ and $\{V_i\}_{i=1}^n$, we then set
\begin{equation*}
\IsoW(f) \equiv \mathbb G_{P}^{(1)} (\text{Proj}\{f |\mathbb S_{n,P}\}) + \mathbb G_{P}^{(2)} (\text{Proj}\{f |\mathbb S_{n,P}^\perp\}) ,
\end{equation*}
which is linear in $f$ by linearity of $f\mapsto \text{Proj}\{f | \mathbb S_{n,P}\}$, $\mathbb G_{P}^{(1)}$, and $\mathbb G_{P}^{(2)}$, and satisfies $E[\IsoW(f)] = 0$ and $E[\IsoW(f)\IsoW(g)] = \text{Cov}_P\{f(V),g(V)\}$.
Moreover, since $\IsoW$ is sub-Gaussian with respect to $\|\cdot\|_{P,2}$, it follows from Corollary 2.2.8 in \cite{vandervaart:wellner:1996}, $N(\delta_n\|G_{n,P}\|_{P,2},\mathcal G_{n,P},\|\cdot\|_{P,2}) = 1$ due to $\|g\|_{P,2}\leq \delta_n\|G_{n,P}\|_{P,2}$ for all $g\in \mathcal G_{n,P}$ and $P\in \mathbf P$, bracketing numbers being larger than covering numbers, Jensen's inequality, and the definition of $J_{1n}$ in Assumption \ref{ass:4seriesreg}(i) that
\begin{multline}\label{th:mainbootcoup3}
E_P[\sup_{g\in \mathcal G_{n,P}}|\IsoW(g)|] \lesssim \delta_n \|G_{n,P}\|_{P,2}  + \int_0^\infty \sqrt{\log(N(\epsilon,\mathcal G_{n,P},\|\cdot\|_{P,2}))}d\epsilon \\
\leq \delta_n \|G_{n,P}\|_{P,2} + \int_0^{\delta_n \|G_{n,P}\|_{P,2}} \sqrt{1+\log(N_{[\hspace{0.02 in}]}(\epsilon,\mathcal G_{n,P},\|\cdot\|_{P,2}))}d\epsilon \lesssim J_{1n} .
\end{multline}
To obtain an analogous bound for $\Wboot$, note $\sup_{g\in \mathcal G_{n,P}} \|g\|_{P,2} \leq \delta_n \|G_{n,P}\|_{P,2}$ by Assumption \ref{ass:4seriesreg}(i) and $|E_P[g(V)]| \leq \|g\|_{P,2}$ by Jensen's inequality imply that
\begin{multline}\label{th:mainbootcoup4}
\sup_{g\in \mathcal G_{n,P}} |\Wboot(g)| \leq \sup_{g\in \mathcal G_{n,P}} |\frac{1}{\sqrt n} \sum_{i=1}^n \omega_i g(V_i)|  \\ + |\frac{1}{\sqrt n} \sum_{i=1}^n \omega_i| \times \{ \sup_{g\in \mathcal G_{n,P}} |\frac{1}{n}\sum_{i=1}^n g(V_i) - E_P[g(V)]| + \delta_n\|G_{n,P}\|_{P,2}\}.
\end{multline}
Next, define the class $\tilde {\mathcal G}_{n,P} \equiv \{(\omega,v)\mapsto \omega g(v) : g \in \mathcal G_{n,P}\}$, and with some abuse of notation let $P$ index the joint distribution of $(V,\omega)$.
Further note that if $\{[g_{i,l,P}, g_{i,u,P}]\}_i$ is a bracket for $\mathcal G_{n,P}$, then the functions $\{[\tilde g_{i,l,P},\tilde g_{i,u,P}]\}$ given by
\begin{align*}
\tilde g_{i,l,P}(\omega,v) & \equiv  \max\{\omega,0\}g_{i,l,P}(v) + \min\{\omega,0\} g_{i,u,P}(v) \notag \\
\tilde g_{i,u,P}(\omega,v) & \equiv  \min\{\omega,0\}g_{i,l,P}(v) + \max\{\omega,0\} g_{i,u,P}(v)
\end{align*}
form a bracket for $\tilde {\mathcal G}_{n,P}$.
Moreover, since $E[\omega^2] = 1$ and $\omega$ and $V$ are independent, it follows that $\|\tilde g_{i,u,P} - \tilde g_{i,l,P}\|_{P,2} = \|g_{i,u,P} - g_{i,l,P}\|_{P,2}$.
Setting $\tilde G_{n,P}(\omega,v) \equiv |\omega| G_{n,P}(v)$, then note that $\tilde G_{n,P}$ is an envelope for $\tilde {\mathcal G}_{n,P}$, which satisfies $\|\tilde G_{n,P}\|_{P,2} = \|G_{n,P}\|_{P,2}$.
For $\eta_{n,P} \equiv 1+\log N_{[\hspace{0.02 in}]}(\delta_n \|G_{n,P}\|_{P,2},\mathcal G_{n,P},\|\cdot\|_{P,2})$ (as in Assumption \ref{ass:4seriesreg}(i)) we then obtain by Theorem 2.14.2 in \cite{vandervaart:wellner:1996} that
\begin{multline}\label{th:mainbootcoup6p1}
E_P[\sup_{g\in \mathcal G_{n,P}}|\frac{1}{\sqrt n}\sum_{i=1}^n\omega_i g(V_i)|] \lesssim J_{[\hspace{0.02 in}]}(\delta_n\|G_{n,P}\|_{P,2},{\mathcal G}_{n,P},\|\cdot\|_{P,2}) \\ + \sqrt{n}E_P[|\omega|G_{n,P}(V)1\{|\omega|\frac{G_{n,P}(V)}{\|G_{n,P}\|_{P,2}} > \frac{\sqrt n \delta_n}{\sqrt{\eta_{n,P}}}\}].
\end{multline}
Moreover, since $\omega$ follows a standard normal distribution, we have $E[|\omega|1\{|\omega| > a\}] \lesssim \exp\{-a^2/2\}$ for any $a \geq 0$.
Therefore, the independence of $\omega$ and $V$ implies
$$E_P[|\omega|G_{n,P}(V)1\{|\omega|  \frac{G_{n,P}(V)}{\|G_{n,P}\|_{P,2}}> \frac{\sqrt{n} \delta_n}{\sqrt {\eta_{n,P}}}\}]  \lesssim E_P[G_{n,P}(V) \exp\{-\frac{n\delta_n^2\|G_{n,P}\|_{P,2}^2}{2G_{n,P}^2(V)\eta_{n,P}}\}]$$
which together with result \eqref{th:mainbootcoup6p1} and the definition of $J_{1n}$ in Assumption \ref{ass:4seriesreg}(i) yields
\begin{equation}\label{th:mainbootcoup6}
E_P[\sup_{g\in \mathcal G_{n,P}}|\frac{1}{\sqrt n}\sum_{i=1}^n\omega_i g(V_i)|] \lesssim J_{1n}.
\end{equation}
Moreover, by Lemmas 2.3.1 and 2.9.1 in \cite{vandervaart:wellner:1996} we have
\begin{multline}\label{th:mainbootcoup7}
E_P[\sup_{g\in \mathcal G_{n,P}} |\frac{1}{\sqrt n}\sum_{i=1}^n g(V_i) - E_P[g(V)]|] + \delta_n \|G_{n,P}\|_{P,2}\\
\lesssim E_P[\sup_{g\in \mathcal G_{n,P}}|\frac{1}{\sqrt n}\sum_{i=1}^n \omega_i g(V_i)|] + \delta_n \|G_{n,P}\|_{P,2} \lesssim J_{1n},
\end{multline}
where the final inequality follows from \eqref{th:mainbootcoup6} and the definition of $J_{1n}$.
Thus, \eqref{th:mainbootcoup4}, \eqref{th:mainbootcoup6}, and \eqref{th:mainbootcoup7} together with Markov's inequality imply that uniformly in $P\in \mathbf P$
\begin{equation}\label{th:mainbootcoup8}
 \|\Wboot\|_{\mathcal G_{n,P}} = O_P(J_{1n}).
\end{equation}
Next, we use the linearity of the processes $f\mapsto \Wboot(f)$ and $f\mapsto \IsoW(f)$ to obtain that
\begin{multline*}
\|\Wboot - \IsoW\|_{\mathcal F_n} \leq \sup_{f\in \mathcal F_n} |\Wboot(f_{n,P}^{d_n\prime}\beta_{n,P}(f)) - \IsoW(f_{n,P}^{d_n\prime}(\beta_{n,P}(f)))| + \|\Wboot - \IsoW\|_{\mathcal G_{n,P}} \\
\leq \sup_{\beta \in \mathcal B_n}|\Wboot(f_{n,P}^{d_n\prime}\beta) - \IsoW(f_{n,P}^{d_n\prime}\beta)| + O_P(J_{1n}) = O_P(J_{2n}R_n + J_{1n}),
\end{multline*}
where the second inequality holds uniformly in $P\in \mathbf P$ by \eqref{th:mainbootcoup3} and Markov's inequality, result \eqref{th:mainbootcoup8}, and set inclusion, while the equality holds uniformly in $P\in \mathbf P$ by result \eqref{th:mainbootcoup1}.
The first claim of the theorem then follows by using Lemma \ref{lm:finbootcoup}(i) to set $R_n = (d_n\log(1+d_n)C_nK_n^2/n)^{1/4}$ in \eqref{th:mainbootcoup1}, and the second part of the theorem follows from using Lemma \ref{lm:finbootcoup}(ii) to set $R_n = (\xi_nd_n \log(1+d_n)C_nK_n^2/n)^{1/2}$ instead. \qed

\begin{lemma}\label{lm:finbootcoup}
Let $\{(\omega_i,V_i)\}_{i=1}^n$ be i.i.d.\ with $V_i \sim P\in \mathbf P$, $\omega_i \sim N(0,1)$, and $\omega_i$ and $V_i$ independent.
Suppose Assumption \ref{ass:4seriescoup} holds, $d_n\log(1+d_n) K_n^2 C_n = o(n)$, and $\mathcal B_n \subset \mathbf R^{d_n}$ satisfies $0\in \mathcal B_n$ and $J_{2n} \equiv \int_0^\infty \sqrt{\log(N(\epsilon,\mathcal B_n,\|\cdot\|_2))}d\epsilon < \infty$.
Then: (i) There is a linear Gaussian process $\IsoW$ on $\mathbb S_{n,P} \equiv \overline{\text{\rm span}}\{f^{d_n}_{n,P}\}$ independent of $\{V_i\}_{i=1}^n$ with
\begin{equation*}
\sup_{\beta \in \mathcal B_n} |\Wboot (f_{n,P}^{d_n\prime} \beta) - \IsoW (f_{n,P}^{d_n\prime} \beta)| = O_P(J_{2n} \{\frac{d_n\log(1+d_n) C_n K_n^2}{n}\}^{1/4})
\end{equation*}
uniformly in $P\in \mathbf P$ and satisfying $E[\IsoW(f_{n,P}^{d_n\prime}\beta)] = 0$ and $E[\IsoW(f_{n,P}^{d_n\prime}\beta_1)\IsoW(f_{n,P}^{d_n\prime}\beta_2)] = {\rm Cov}_P\{f_{n,P}^{d_n}(V)^\prime \beta_1,f_{n,P}^{d_n}(V)^\prime \beta_2\}$.
(ii) If in addition $\sup_{P\in \mathbf P}\|\text{\rm Var}_P^{-1}\{f_{n,P}^{d_n}(V)\}\|_{o,2} \leq \xi_n < \infty$ and $\xi_n\sqrt{d_n\log(1+d_n)C_n}K_n/\sqrt n = o(1)$, then uniformly in $P\in \mathbf P$
\begin{equation*}
\sup_{\beta \in \mathcal B_n} |\Wboot (f_{n,P}^{d_n\prime} \beta) - \IsoW(f_{n,P}^{d_n\prime} \beta)| = O_P(\frac{J_{2n} \sqrt {\xi_n d_n\log(1+d_n)C_n}K_n}{\sqrt n}).
\end{equation*}
\end{lemma}

\noindent {\sc Proof:} First note that $\Wboot(f - c) = \Wboot(f)$ for any $c\in \mathbf R$ and $f\in L^2_P$.
We may therefore assume without loss of generality that $E_P[f_{n,P}^{d_n}(V)] = 0$, and for every $P\in \mathbf P$ we let $\Sigma_n(P) \equiv \text{Var}_P\{f_{n,P}^{d_n}(V)\} = E_P[f_{n,P}^{d_n}(V)f_{n,P}^{d_n}(V)^\prime]$ and define
\begin{equation*}
\hat \Sigma_n(P) \equiv \frac{1}{n} \sum_{i=1}^n (f_{n,P}^{d_n}(V_i)-\frac{1}{n}\sum_{j=1}^n f_{n,P}^{d_n}(V_j))(f_{n,P}^{d_n}(V_i)-\frac{1}{n}\sum_{j=1}^n f_{n,P}^{d_n}(V_j))^\prime .
\end{equation*}
For a sequence $R_n$ with $R_n = o(1)$, and any constant $M > 0$ and $P\in \mathbf P$ define the event
\begin{equation}\label{lm:finbootcoup2}
A_{n,P}(M) \equiv \{\|\hat \Sigma_n^{1/2}(P) - \Sigma_n^{1/2}(P)\|_{o,2} \leq M R_n \} .
\end{equation}
Further note that by Lemma \ref{lm:oprate} it follows we may select $R_n = o(1)$ such that we have
\begin{equation}\label{lm:finbootcoup3}
\liminf_{M\uparrow \infty} \liminf_{n\rightarrow \infty} \inf_{P\in \mathbf P} P(\{V_i\}_{i=1}^n \in  A_{n,P}(M) ) = 1 .
\end{equation}
In particular, to establish part (i) we will set $R_n = (d_n\log(1+d_n)C_nK_n^2/n)^{1/4}$ and employ Lemma \ref{lm:oprate}(i), while to establish part (ii) we will set $R_n = (\xi_n d_n \log(1+d_n)C_nK_n^2/n)^{1/2}$ and employ Lemma \ref{lm:oprate}(ii) instead.

Next, let $\mathcal N_{d_n} \in \mathbf R^{d_n}$ follow a standard normal distribution and be independent of $\{(\omega_i,V_i)\}_{i=1}^n$ (defined on the same suitably enlarged probability space).
Further let $\{\hat \nu_d\}_{d=1}^{d_n}$ denote eigenvectors of $\hat \Sigma_n(P)$, $\{\hat \lambda_d\}_{d=1}^{d_n}$ represent the corresponding (possibly zero) eigenvalues and define the random variable $\mathbb Z_{n,P} \in \mathbf R^{d_n}$ to be given by
\begin{equation}\label{lm:finbootcoup4}
\mathbb Z_{n,P} \equiv \sum_{d : \hat \lambda_d \neq 0} \hat \nu_d  \frac{(\hat \nu_d^\prime \Wboot(f_{n,P}^{d_n}))}{\hat \lambda_d^{1/2}} + \sum_{d : \hat \lambda_d = 0} \hat \nu_d (\hat \nu_d^\prime \mathcal N_{d_n}) .
\end{equation}
Then note that since $\Wboot(f_{n,P}^{d_n}) \sim N(0,\hat \Sigma_n(P))$ conditional on $\{V_i\}_{i=1}^n$, and $\mathcal N_{d_n}$ is independent of $\{(\omega_i,V_i)\}_{i=1}^n$, $\mathbb Z_{n,P}$ is Gaussian conditional on $\{V_i\}_{i=1}^n$.
Furthermore,
\begin{equation*}
E[\mathbb Z_{n,P}\mathbb Z_{n,P}^\prime|\{V_i\}_{i=1}^n] = \sum_{d = 1}^{d_n} \hat \nu_d \hat \nu_d^\prime = I_{d_n}
\end{equation*}
by direct calculation for $I_{d_n}$ the $d_n\times d_n$ identity matrix.
Hence, $\mathbb Z_{n,P} \sim N(0,I_{d_n})$ conditional on $\{V_i\}_{i=1}^n$ almost surely in $\{V_i\}_{i=1}^n$ and is thus independent of $\{V_i\}_{i=1}^n$.
Moreover, we also note that by Theorem 3.6.1 in \cite{bogachevgauss} and $\Wboot(f_{n,P}^{d_n}) \sim N(0,\hat \Sigma_n(P))$ conditional on $\{V_i\}_{i=1}^n$, it follows that $\Wboot(f_{n,P}^{d_n})$ belongs to the range of $\hat \Sigma_n(P) : \mathbf R^{d_n}\rightarrow \mathbf R^{d_n}$ almost surely in $\{(\omega_i,V_i)\}_{i=1}^n$.
Therefore, since $\{\hat \nu_d : \hat \lambda_d \neq 0 \}_{d=1}^{d_n}$ spans the range of $\hat \Sigma_n(P)$, we conclude from \eqref{lm:finbootcoup4} that for any $\beta \in \mathbf R^{d_n}$
\begin{equation*}
\beta^\prime \hat \Sigma_n^{1/2}(P) \mathbb Z_{n,P} = \beta^\prime \sum_{d:\hat{\lambda}_d \neq 0} \hat \nu_d (\hat \nu_d^\prime \Wboot(f_{n,P}^{d_n})) = \Wboot(\beta^\prime f_{n,P}^{d_n})  .
\end{equation*}
Analogously, we define for any $\beta \in \mathbf R^{d_n}$ the linear Gaussian process $\IsoW$ on $\mathbb S_{n,P}$ by
\begin{equation*}
\IsoW(\beta^\prime f_{n,P}^{d_n}) \equiv \beta^\prime \Sigma_n^{1/2}(P) \mathbb Z_{n,P} ,
\end{equation*}
which by construction is  independent of $\{V_i\}_{i=1}^n$  and satisfies $E[\IsoW(f_{n,P}^{d_n\prime}\beta)] = 0$ and $E[\IsoW(f_{n,P}^{d_n\prime}\beta_1)\IsoW(f_{n,P}^{d_n\prime}\beta_2)] = {\rm Cov}_P\{f_{n,P}^{d_n}(V)^\prime \beta_1,f_{n,P}^{d_n}(V)^\prime \beta_2\}$.
Setting
\begin{equation}\label{lm:finbootcoup8}
\bar {\mathbb G}_{P}(\beta) \equiv (\beta^\prime(\hat \Sigma_n^{1/2}(P) - \Sigma_n^{1/2}(P))\mathbb Z_{n,P}) 1\{ A_{n,P}(M)\} ,
\end{equation}
where  $1\{A_{n,P}(M)\}$ denotes the indicator for the event $\{V_i\}_{i=1}^n \in A_{n,P}(M)$, then note 
\begin{equation}\label{lm:finbootcoup9}
\sup_{\beta \in \mathcal B_n}|\Wboot(f_{n,P}^{d_n \prime}\beta) - \IsoW(f_{n,P}^{d_n \prime}\beta)| 1\{A_{n,P}(M)\}  = \sup_{\beta \in \mathcal B_n} |\bar{\mathbb G}_{P}(\beta)|.
\end{equation}
Moreover, we note that conditional on $\{V_i\}_{i=1}^n$, $\bar {\mathbb G}_{P}$ is sub-Gaussian under the semi-metric $\rho_n(\tilde \beta,\beta) \equiv \|(\hat \Sigma_n^{1/2}(P) - \Sigma_n^{1/2}(P))(\tilde \beta -  \beta)\|_{2}$.
Since $\|\hat \Sigma_n^{1/2}(P) - \Sigma_n^{1/2}(P)\|_{o,2} \leq MR_n$ whenever $1\{A_{n,P}(M)\} = 1$ we obtain, whenever $\{V_i\}_{i=1}^n \in A_{n,P}(M)$, that
\begin{align}\label{lm:finbootcoup10}
\int_0^\infty \sqrt{\log(N(\epsilon,\mathcal B_n,\rho_n))}d\epsilon & \leq \int_0^\infty \sqrt{\log(N(\epsilon/MR_n,\mathcal B_n,\|\cdot\|_2))}d\epsilon \notag \\
& = MR_n \int_0^\infty \sqrt{\log(N(u,\mathcal B_n,\|\cdot\|_2))}du,
\end{align}
where the equality follows from the change of variables $\epsilon = MR_n u$.
Therefore, since $0\in \mathcal B_n$, Corollary 2.2.8 in \cite{vandervaart:wellner:1996} and \eqref{lm:finbootcoup10} imply
\begin{equation}\label{lm:finbootcoup11}
E[\sup_{\beta \in \mathcal B_n} |\bar {\mathbb G}_{P}(\beta)| |\{V_i\}_{i=1}^n] \lesssim  MR_n \int_0^\infty \sqrt{\log(N(u,\mathcal B_n,\|\cdot\|_2))}du \equiv MR_n J_{2n}.
\end{equation}
Next, note \eqref{lm:finbootcoup8}, \eqref{lm:finbootcoup9}, and \eqref{lm:finbootcoup11} together with Markov's inequality imply that
\begin{multline}\label{lm:finbootcoup12}
P( \sup_{\beta \in \mathcal B_n} |\Wboot(f_{n,P}^{d_n\prime}\beta) - \IsoW(f_{n,P}^{d_n\prime}\beta) | > M^2 R_nJ_{2n}; ~ A_{n,P}(M)) \\ \leq P(  \sup_{\beta \in \mathcal B_n} |\bar{\mathbb G}_{P}(\beta)| > M^2 R_nJ_{2n})  \lesssim \frac{1}{M}
\end{multline}
for all $P\in \mathbf P$. Therefore, combining results \eqref{lm:finbootcoup3} and \eqref{lm:finbootcoup12}, we can finally conclude
\begin{multline*}
\limsup_{M\uparrow \infty} \limsup_{n\rightarrow \infty} \sup_{P\in \mathbf P} P(\sup_{\beta \in \mathcal B_n} |\Wboot(f_{n,P}^{d_n\prime}\beta) - \IsoW(f_{n,P}^{d_n\prime}\beta)| > M^2  R_nJ_{2n}) \\ \lesssim \limsup_{M\uparrow \infty} \limsup_{n\rightarrow \infty} \sup_{P\in \mathbf P} \{\frac{1}{M} + P(\{V_i\}_{i=1}^n \notin A_{n,P}(M)) \} = 0 .
\end{multline*}
The first claim of the lemma then follows by employing Lemma \ref{lm:oprate}(i) to set $R_n = (d_n\log(1+d_n)C_nK_n^2/n)^{1/4}$ in \eqref{lm:finbootcoup2}, while the second claim follows by employing Lemma \ref{lm:oprate}(ii) to set $R_n = (\xi_n d_n\log(1+d_n)C_nK_n^2/n)^{1/2}$. \qed

\begin{lemma}\label{lm:oprate}
Let $\{V_i\}_{i=1}^n$ be i.i.d.\ with $V \sim P \in \mathbf P$, suppose Assumption \ref{ass:4seriescoup} holds, define $\Sigma_n(P) \equiv \text{\rm Var}_P\{f_{n,P}^{d_n}(V)\}$ and its sample analogue $\hat \Sigma_n(P)$ to equal
$$\hat \Sigma_{n}(P)  \equiv \frac{1}{n} \sum_{i=1}^n (f_{n,P}^{d_n}(V_i) - \frac{1}{n}\sum_{j=1}^n f_{n,P}^{d_n}(V_j)) (f_{n,P}^{d_n}(V_i) - \frac{1}{n}\sum_{j=1}^n f_{n,P}^{d_n}(V_j))^\prime ,$$
and assume $d_n \log(1+d_n)K_n^2C_n = o(n)$.
(i) Then, it follows that uniformly in $P\in \mathbf P$
\begin{equation*}
\|\hat \Sigma_n^{1/2}(P) - \Sigma_n^{1/2}(P)\|_{o,2} = O_P(\{\frac{d_n\log(1+d_n)C_n K_n^2}{n}\}^{1/4}).
\end{equation*}
(ii) If in addition $\sup_{P\in \mathbf P} \|\Sigma_n^{-1}(P)\|_{o,2} \leq \xi_n < \infty$ and $\xi_n\sqrt{d_n\log(1+d_n)C_n}K_n/\sqrt n = o(1)$, then we can also conclude uniformly in $P\in \mathbf P$ that
\begin{equation*}
\|\hat \Sigma_n^{1/2}(P) - \Sigma_n^{1/2}(P)\|_{o,2} = O_P(\frac{\sqrt{\xi_n d_n\log(1+d_n)C_n} K_n}{\sqrt n}) .
\end{equation*}
\end{lemma}

\noindent {\sc Proof:}
First note that we may without loss of generality assume that $E_P[f_{n,P}^{d_n}(V)] = 0$. 
Next note that Assumption \ref{ass:4seriescoup}(ii) implies that for all $P\in \mathbf P$ we must have
\begin{equation}\label{lm:oprate1}
\|\frac{1}{n}\{f^{d_n}_{n,P}(V_i)f^{d_n}_{n,P}(V_i)^\prime - E_P[f^{d_n}_{n,P}(V)f_{n,P}^{d_n}(V)^\prime]\}\|_{o,2} \leq \frac{2d_nK_n^2}{n}
\end{equation}
almost surely for all $P\in \mathbf P$ since each entry of the matrix $f^{d_n}_{n,P}(V_i)f^{d_n}_{n,P}(V_i)^\prime$ is bounded by $K_n^2$.
Similarly, employing $\|f_{n,P}^{d_n}(V_i)f_{n,P}^{d_n}(V_i)^\prime\|_{o,2} \leq d_nK_n^2$ almost surely we obtain
\begin{equation}\label{lm:oprate2}
\|\frac{1}{n} E_P[\{f_{n,P}^{d_n}(V)f_{n,P}^{d_n}(V)^\prime - E_P[f_{n,P}^{d_n}(V) f_{n,P}^{d_n}(V)^\prime]\}^2]\|_{o,2}
\leq \frac{2d_n K_n^2 C_n}{n}.
\end{equation}
Thus, employing results \eqref{lm:oprate1} and \eqref{lm:oprate2}, together with $d_n\log(1+d_n)K_n^2C_n = o(n)$, we obtain by Theorem 6.1(ii) in \cite{tropp:2012} that for all $P\in \mathbf P$
\begin{multline}\label{lm:oprate3}
P(\|\frac{1}{n}\sum_{i=1}^n f_{n,P}^{d_n}(V_i)f_{n,P}^{d_n}(V_i)^\prime - E_P[f_{n,P}^{d_n}(V) f_{n,P}^{d_n}(V)^\prime]\|_{o,2} >  \frac{M\sqrt{d_n\log(1+d_n)C_n}K_n}{\sqrt n})  \\
\leq  d_n \exp\{ - \frac{M^2(d_n\log(1+d_n)K_n^2C_n)}{2n} \frac{n}{MBd_nK_n^2C_n}\}
\end{multline}
for some $B < \infty$.
Hence, we can conclude from \eqref{lm:oprate3} that uniformly in $P\in \mathbf P$
\begin{equation}\label{lm:oprate4}
\|\frac{1}{n}\sum_{i=1}^n f_{n,P}^{d_n}(V_i)f_{n,P}^{d_n}(V_i)^\prime - E_P[f_{n,P}^{d_n}(V) f_{n,P}^{d_n}(V)^\prime]\|_{o,2} = O_P(\frac{\sqrt{d_n\log(1+d_n)C_n}K_n}{\sqrt n}).
\end{equation}
Recalling that we had without loss of generality set $E_P[f_{n,P}^{d_n}(V)] = 0$, next note that $E_P[f_{d,n,P}^2(V)] \leq \|E_P[f_{n,P}^{d_n}(V)f_{n,P}^{d_n}(V)^\prime]\|_o \leq C_n$, Markov's inequality, and Lemmas 2.2.9 and 2.2.10 in \cite{vandervaart:wellner:1996} imply, uniformly in $P\in \mathbf P$, that
\begin{multline}\label{lm:oprate5}
\|\frac{1}{n}\sum_{i=1}^n f_{n,P}^{d_n}(V_i)\|_2 \leq \sqrt{d_n} \max_{1\leq d \leq d_n} |\frac{1}{n}\sum_{i=1}^n f_{d,n,P}(V_i) | \\ = O_P(\frac{K_n\log(1+d_n)\sqrt{d_n}}{n} + \frac{\sqrt{C_n d_n \log(1+d_n)}}{\sqrt n}).
\end{multline}
Therefore, since for any $a,b\in \mathbf R^{d_n}$ we have $\|ab^\prime\|_{o,2}\leq \|a\|_2 \|b\|_2$, results \eqref{lm:oprate4} and \eqref{lm:oprate5} together with the triangle inequality yield, uniformly in $P\in \mathbf P$, that
\begin{equation}\label{lm:oprate8}
\|\hat \Sigma_n(P) - \Sigma_n(P)\|_{o,2} = O_P(\frac{\sqrt{d_n\log(1+d_n)C_n}K_n}{\sqrt n}).
\end{equation}
Finally, since $\hat \Sigma_n(P) \geq 0 $ and $\Sigma_n(P)\geq 0$, Theorem X.1.1 in \cite{bhatia:1997} implies that
\begin{equation}\label{lm:oprate9}
\|\hat \Sigma_n^{1/2}(P) - \Sigma_n^{1/2}(P)\|_{o,2} \leq \|\hat \Sigma_n(P) - \Sigma_n(P)\|_{o,2}^{1/2}
\end{equation}
almost surely, and hence the first claim the lemma follows from \eqref{lm:oprate8} and \eqref{lm:oprate9}.

For the second claim, let \underline{eig}$\{A\}$ denote the smallest eigenvalue of any Hermitian matrix $A$.
Since $\|\Sigma_n^{-1}(P)\|_{o,2} = 1/\underline{\text{eig}}\{\Sigma_n(P)\}$, $\sup_{P\in \mathbf P}\|\Sigma_n^{-1}(P)\|_{o,2} \leq \xi_n$, result \eqref{lm:oprate8}, Corollary III.2.6 in \cite{bhatia:1997}, and $\xi_n\sqrt{d_n\log(1+d_n)C_n}K_n/\sqrt n = o(1)$ imply
\begin{multline*}
\liminf_{n\rightarrow \infty} \inf_{P\in \mathbf P} P( \underline{\text{eig}}\{\hat \Sigma_n(P)\} > \frac{1}{2\xi_n}) \\
\geq \liminf_{n\rightarrow \infty} \inf_{P\in \mathbf P} P( \underline{\text{eig}}\{\Sigma_n(P)\} > \frac{1}{2\xi_n} + \|\hat \Sigma_n(P) - \Sigma_n(P)\|_{o,2})  = 1.
\end{multline*}
Therefore, Applying Theorem X.3.8 in \cite{bhatia:1997} we can then conclude that
\begin{multline*}
\limsup_{M\uparrow \infty} \limsup_{n\rightarrow \infty} \sup_{P\in \mathbf P} P(\|\hat \Sigma_n^{1/2}(P) - \Sigma_n^{1/2}(P)\|_{o,2} > M\frac{\sqrt{\xi_n d_n\log(1+d_n)C_n} K_n}{\sqrt n})
\\ \leq \limsup_{M\uparrow \infty}  \limsup_{n\rightarrow \infty} \sup_{P\in \mathbf P} P(\|\hat \Sigma_n(P) - \Sigma_n(P)\|_{o,2} >  M\frac{\sqrt{ d_n\log(1+d_n)C_n} K_n}{\sqrt n}) = 0,
\end{multline*}
where the final equality follows from result \eqref{lm:oprate8}. \qed

\begin{lemma}\label{lm:auxbootcalc}
For any positive random variable $U$ with $E[U^2] < \infty$ and finite constant $A > 0$ it follows that $E[U\exp\{-A/U^2\}] \leq E[U]\exp\{-A/E[U^2]\} + E[U^2]/\sqrt{2A}.$
\end{lemma}

\noindent {\sc Proof:} First note $u\mapsto u \exp\{-A/u^2\}$ is convex on $u\in (0,\sqrt{2A}]$.
Therefore Jensen's inequality, $u \mapsto u\exp\{- A^2/u^2\}$ being increasing in $u \in (0,\infty)$, $E[1\{0<U <\sqrt{2A}\}U] \leq E[U]$ due to $U$ being positive a.s., and $\exp\{-A/U^2\} \leq 1$ due to $A > 0$, imply
\begin{align*}
E[U\exp\{-\frac{A}{U^2}\}] & = E[1\{0 < U \leq \sqrt{2A}\}U \exp\{-\frac{A}{U^2}\}] + E[1\{U > \sqrt{2A}\} U\exp\{-\frac{A}{U^2}\}] \\
& \leq E[U] \exp\{-\frac{A}{E[U^2]}\} + E[1\{U > \sqrt{2A}\} U].
\end{align*}
The claim of the lemma therefore follows from $E[1\{U > \sqrt{2A}\}U] \leq E[U^2]/\sqrt{2A}$ by the Cauchy Schwarz inequality and Markov's inequality. \qed